\documentclass[12pt,letter]{amsart}
\usepackage{amscd,amssymb,graphics,color,hyperref}
\usepackage[latin1]{inputenc}
\usepackage{srcltx}

\usepackage{mathrsfs}

\usepackage[all]{xy}
\newtheorem{theorem}{Theorem}[section]
\newtheorem{lemma}[theorem]{Lemma}
\newtheorem{proposition}[theorem]{Proposition}
\newtheorem{corollary}[theorem]{Corollary}
\newtheorem{conjecture}[theorem]{Conjecture}

\theoremstyle{definition}
\newtheorem{definition}[theorem]{Definition}
\newtheorem{construction}[theorem]{Construction}
\newtheorem{example}[theorem]{Example}

\theoremstyle{remark}
\newtheorem{remark}[theorem]{Remark}

\newcommand{\HH} {\mathbb{H}}
\newcommand{\NN} {\mathbb{N}}
\newcommand{\ZZ} {\mathbb{Z}}
\newcommand{\QQ} {\mathbb{Q}}
\newcommand{\RR} {\mathbb{R}}
\newcommand{\CC} {\mathbb{C}}

\newcommand{\PP} {\mathbb{P}}
\renewcommand{\AA} {\mathbb{A}}
\newcommand{\GG} {\mathbb{G}}

\newcommand {\shAff} {\mathcal{A}\text{\textit{ff}}}

\newcommand {\shExt} {\mathcal{E} \!\text{\textit{xt}}}

\newcommand {\shF} {\mathcal{F}}

\newcommand {\shH} {\mathcal{H}}
\newcommand {\shHom} {\mathcal{H}\!\text{\textit{om}}}
\newcommand {\shI} {\mathcal{I}}

\newcommand {\shL} {\mathcal{L}}
\newcommand {\shM} {\mathcal{M}}
\newcommand {\shMPL} {\mathcal{MPL}}

\newcommand {\shN} {\mathcal{N}}
\newcommand {\shPL} {\mathcal{PL}}
\newcommand {\shQ} {\mathcal{Q}}
\newcommand {\shR} {\mathcal{R}}
\newcommand {\shS} {\mathcal{S}}
\newcommand {\shT} {\mathcal{T}}

\newcommand {\shP} {\mathcal{P}}

\newcommand {\shX} {\mathcal{X}}
\newcommand {\shY} {\mathcal{Y}}

\newcommand {\foR} {\mathfrak{R}}

\newcommand {\A} {\boldsymbol{A}}

\newcommand {\Aff} {\operatorname{Aff}}

\newcommand {\Ann} {\operatorname{Ann}}

\renewcommand {\Bar} {\operatorname{Bar}}
\newcommand {\bct} {\mathrm{bct}}

\newcommand {\C} {\mathscr{C}}
\newcommand {\Cat} {\underline{\mathrm{Cat}}}

\newcommand {\cl} {\operatorname{cl}}
\newcommand {\codim} {\operatorname{codim}}

\newcommand {\coker} {\operatorname{coker}}

\newcommand {\dlog} {\operatorname{dlog}}

\newcommand {\Ext} {\operatorname{Ext}}

\newcommand {\GL} {\operatorname{GL}}
\newcommand {\gp} {{\operatorname{gp}}}

\newcommand {\Hom} {\operatorname{Hom}}

\newcommand {\id} {\operatorname{id}}
\newcommand {\im} {\operatorname{im}}

\newcommand {\Int} {\operatorname{Int}}

\newcommand {\Lin} {\operatorname{Lin}}

\newcommand {\LPoly} {\underline{\mathrm{LPoly}}}
\newcommand {\lra} {\longrightarrow}
\newcommand {\ls} {\dagger}

\newcommand {\M} {\mathcal{M}}

\newcommand {\maxid} {\mathfrak{m}}

\renewcommand{\O} {\mathcal{O}}

\newcommand {\ord} {\operatorname{ord}}

\renewcommand{\P} {\mathscr{P}}
\newcommand {\PA} {\operatorname{PA}}

\newcommand {\Pic} {\operatorname{Pic}}

\newcommand {\phd} {\mathrm{phd}}

\newcommand {\pre} {\mathrm{pre}}

\newcommand {\Proj} {\operatorname{Proj}}

\newcommand {\rank} {\operatorname{rank}}

\newcommand {\sgn} {\mathrm{sgn}}
\newcommand {\sing} {\mathrm{sing}}
\newcommand {\Sing} {\operatorname{Sing}}
\newcommand {\SL} {\operatorname{SL}}
\newcommand {\Spec} {\operatorname{Spec}}

\newcommand {\Strata} {\operatorname{Strata}}
\newcommand {\supp} {\operatorname{supp}}

\newcommand {\Trans} {\operatorname{Trans}}

\newcommand{\bbfamily}{\fontencoding{U}\fontfamily{bbold}\selectfont}
\newcommand{\textbb}[1]{{\bbfamily#1}}
\newcommand {\lfor} {\mbox{\textbb{[}}}
\newcommand {\rfor} {\mbox{\textbb{]}}}
\newcommand {\T} {\shT}
\newcommand {\X} {\shX}
\newcommand {\Y} {\shY}

\newcommand {\W} {\mathscr{W}}
\newcommand {\Gm} {\GG_m}
\newcommand {\shLS} {\mathcal{LS}}
\newcommand {\D} {\operatorname{D}}

\def\mydate{\ifcase\month \or January\or February\or March\or
April\or May\or June\or July\or August\or September\or October\or 
November\or December\fi \space\number\day,\space\number\year}

\begin{document}
\def\mapright#1{\smash{
 \mathop{\longrightarrow}\limits^{#1}}}
\def\mapleft#1{\smash{
 \mathop{\longleftarrow}\limits^{#1}}}
\def\exact#1#2#3{0\lra#1\lra#2\lra#3\lra0}
\def\mapup#1{\Big\uparrow
  \rlap{$\vcenter{\hbox{$\scriptstyle#1$}}$}}
\def\mapdown#1{\Big\downarrow
  \rlap{$\vcenter{\hbox{$\scriptstyle#1$}}$}}
\def\dual#1{#1^{\hspace{0.5pt}\scriptscriptstyle \vee}}
\def\dualvs#1{#1^{\hspace{0.5pt}*}}
\def\invlim{\mathop{\rm lim}\limits_{\longleftarrow}}

\title[Mirror Symmetry via Logarithmic Degeneration Data I]{
Mirror Symmetry via Logarithmic Degeneration Data I}

\author{Mark Gross} 
\address{UCSD Mathematics, 9500 Gilman Drive, La Jolla, CA 92093-0112, USA}
\email{mgross@math.ucsd.edu}
\thanks{This work was partially supported by NSF grant
0204326, EPSRC, the Heisenberg program of the DFG, and SPP 1094}

\author{Bernd Siebert} \address{Mathematisches Institut,
Albert-Ludwigs-Universit\"at, Eckerstra\ss e~1, 79104 Freiburg,
Germany}
\email{bernd.siebert@math.uni-freiburg.de}

\maketitle
\tableofcontents
\bigskip

\section*{Introduction.}

This paper is the first arising from our project announced in
\cite{Announce} aiming at establishing a new paradigm for mirror
symmetry. At the center of this approach is residual data associated
to certain maximally unipotent degenerations $f:\X\to\mathcal{S}$ of
Calabi-Yau varieties. The residual degeneration data consists of the
central fibre $\X_0$ of the degeneration together with the
log-structure induced from the inclusion $\X_0\subset\X$, a
polarization, and an element of a ``log K\"ahler moduli space''.  We
claim that mirror symmetry comes down to an involution acting on
residual degeneration data. In particular, degenerating families
should be mirror-dual if and only if they have dual expressions in
terms of residual data. The deepest result that we prove here is a
basic duality between logarithmic complex and K\"ahler moduli of
central fibres of degenerations as log spaces.  Finer consequences for
mirror symmetry will be addressed in a sequel to this paper.

The central idea is to restrict attention to what we call {\it toric
degenerations} of Calabi-Yau varieties (Definition \ref{toric
degen}). Roughly put, these are degenerations in which the singular
fibre is a union of toric varieties and the map to the base is log
smooth off of some bad set $Z$. In this context, we define a dual
intersection complex capturing key data about the degeneration. The
dual intersection complex we construct is an affine manifold $B$ with
singularities together with a polyhedral decomposition $\P$. The
affine structure depends on both the irreducible components of $\X_0$
and on information about the structure of $f$ at the most singular
points of $\X_0$.

If in addition, $\X$ is polarized with a choice of relatively ample
line bundle $\shL$, then the dual intersection complex $B$ comes 
equipped with a convex multi-valued piecewise linear function
$\varphi$. This function can be used to define a discrete Legendre
transform, which gives a new affine manifold $\check B$ which is, in
a suitable sense, dual to $B$. In addition $\check B$ carries a new
convex multi-valued piecewise linear function $\check\varphi$. Then
$\check B$ should be the dual intersection complex of a mirror
degeneration.

This part of the construction is relatively simple and conceptual, but
it covers only the discrete part of mirror symmetry. It is the
treatment of moduli that makes this paper so long. One apparent source
of moduli is a change of gluing of the components of $\X_0$. However,
an essential insight of this paper is that the correct limiting
version of complex moduli also involves a choice of logarithmic
structure. More justification for doing this will come from the study
of deformation theory in the sequel of this paper. On the other hand,
the source of the moduli on the K\"ahler side is perhaps less clear.
In this paper we have fairly much reverse-engineered from the complex
moduli side. The fact that the K\"ahler moduli space is then related
(in an expected fashion) to the logarithmic Picard group should be
seen as the first striking verification of our approach. In the sequel
to this paper, we will connect the K\"ahler moduli more directly to a
logarithmic version of $H^{1,1}$ and couple this with a base-change
theorem,  justifying our definition further \cite{sequel}.

Let us consider a simple example of the phenomena we have just outlined.
Consider $\X\subseteq\PP^3\times\AA^1$ defined by the equation
\[
f_4+tx_0x_1x_2x_3=0
\]
for $f_4$ a general choice of homogeneous polynomial of degree four,
with $t$ the coordinate on $\AA^1$. Let $f:\X\rightarrow\AA^1$ be the
projection. Then $\X_0$ is a union of toric varieties ($\PP^2$'s) meeting 
along toric strata, and the map is generically normal crossings, except
at $24$ points in $\X_0$ where the total space $\X$ is singular; call
this set of points $Z$.

We can build the dual intersection complex $B$ of $\X$ as follows.
Because $\X_0$ is normal crossings, $B$ will be a simplicial complex
and coincide with the traditional dual intersection complex: we have 
a vertex for each irreducible component of $\X_0$, and if 
$v_0,\ldots,v_k$ are vertices corresponding to components $X_{v_0},
\ldots,X_{v_k}$, then $\langle v_0,\ldots,v_k\rangle$ is a simplex of
the dual intersection complex if and only if $X_{v_0}\cap\cdots\cap
X_{v_k}\not=\emptyset$. Thus in this case, $B$ is the boundary of a
tetrahedron. $B$ carries a polyhedral decomposition $\P$, namely the
collection of simplices of this simplicial complex. To make $B$ an
affine manifold with singularities, we remove the midpoints of each
edge; these will be singularities of the affine structure. Identify
each face with the standard simplex in $\RR^2$; this gives an affine
structure in the interior of each face. Standard simplices appear here
because $f:\X\setminus Z\rightarrow\AA^1$ is normal crossings; in the
more general case we will replace simplices by lattice polytopes.

To define affine charts in a neighbourhood of each vertex $v$, we 
specify a \emph{fan structure} at $v$. This means we choose a 
complete rational polyhedral fan $\Sigma_v$ in $\RR^2$, and a
homeomorphism between an open neighbourhood of $v$ and an open
neighbourhood of $0$ in $\RR^2$, giving a one-to-one correspondence
between cells of $\P$ containing $v$ and cones in $\Sigma_v$, so that 
in the interior of each maximal cell of $\P$ containing $v$, the
homeomorphism is an element of $\Aff(\ZZ^2)$. To construct the dual
intersection complex, we take the fan $\Sigma_v$ to be the fan
defining the irreducible component $X_v$. This gives the following
picture:\\[1ex]\mbox{}
\begin{center}
\includegraphics{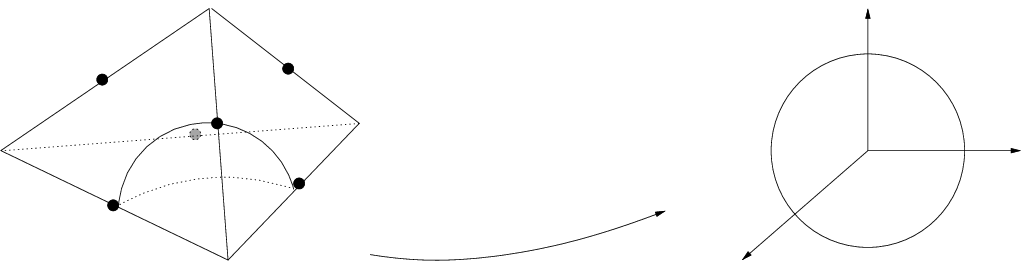}
\end{center}

To construct the intersection complex $\check B$ of $\X$, we need to
choose a polarization on $\X_0$, say given by the restriction of
$\O_{\PP^3}(n)$ for some $n$ to $\X_0$. Then restricting this line
bundle to each irreducible component of $\X_0$, we take the Newton
polytope of this line bundle, giving in our example the standard
simplex of $\RR^2$ rescaled by a factor of $n$. These are glued
together in the obvious way to reflect the intersections of the
components of $\X_0$, again yielding a tetrahedron in this example.
This gives a polyhedral decomposition $\check\P$ on $\check B$ which
is dual to $\P$. Finally we specify an affine structure with
singularities by specifying a fan structure at each vertex. This time
we take as fan the normal fan to the corresponding (maximal) cell of
the dual intersection complex. The example of the quartic is perhaps
misleading though because $\check B$ is the same thing as $B$,
rescaled by a factor of $n$. We will see in \S 4, however, that it is
always true that $B$ and  $\check B$ are related by a discrete
Legendre transform.

Now going from a toric degeneration to its intersection complex or dual
intersection complex is easy; going backwards is more difficult, and
we do not complete this task in this paper. Rather, we only show how to
go from $(B,\P)$ to a family of log Calabi-Yau spaces. To illustrate
this, consider the case that $B$ is two dimensional, as above, and the
polyhedral decomposition $\P$ subdivides $B$ into standard simplices. 
Then to each vertex $v\in\P$ is associated a fan $\Sigma_v$,
obtained by looking at $\P$ in a small neighbourhood of $v$ in $B$. 
We then obtain a corresponding toric variety $X_v$. These varieties
can then be glued together in a way whose combinatorics is specified
by $B$ and $\P$: we want to construct a scheme or algebraic space 
$X_0$ whose dual intersection complex is $(B,\P)$. There is
some moduli to this gluing, say specified by some data $s$.
We then obtain spaces $X_0(B,\P,s)$ given by this gluing. In this
simple case, $X_0(B,\P,s)$ is normal crossings, and as studied by
Friedman in \cite{Friedman}, the deformation theory of $X_0(B,\P,s)$
is controlled partly by a sheaf $\shN_D=
\shExt^1(\Omega^1_{X_0(B,\P,s)},\O_{X_0(B,\P,s)})$, a line bundle
on the singular locus $D$ of $X_0(B,\P,s)$. It turns out this sheaf
can be calculated in terms of $B$ and $s$, and in particular, the restriction
of $\shN_D$ to any irreducible component of $D$, corresponding to
an edge $\tau$ of $\P$, is $\O_{\PP^1}(n)$, where suitably defined monodromy
of the affine structure is $\begin{pmatrix}1&n\\0&1\end{pmatrix}$
around the singular point in the interior of $\tau$. (See 
Example~\ref{normal crossings}. In the quartic example above, $n=4$
for each singular point.) In order for $X_0(B,\P,s)$ to be smoothable,
$\shN_D$ must have a section not generically zero along any component
(and the zeroes will be related to the singular set $Z$ mentioned above).
This places restrictions on both the numbers $n$ (obviously we need
$n\ge 0$, which is the positivity condition of \S 1.5) and $s$.
It turns out in this case specifying a section of $\shN_D$ with no zeroes at 
singular points of $D$
is the same thing as choosing a log structure on $X_0(B,\P,s)$
and a morphism to a log point. Thus the full moduli space we are interested
in is, in this example, the moduli space of pairs $X_0(B,\P,s)$ and sections
of $\shN_D$ not vanishing at any singular point of 
$D$. This can be viewed as a rough approximation to the 
moduli space of a smoothing of $X_0(B,\P,s)$: they in fact have the same
dimension.

For more general $B,\P$ or more general types of toric degenerations
(not just normal crossings) the situation is much more complicated.
This is the difficulty alluded to above which makes this paper so long.
It is necessary to consider more general degenerations, for if we restricted
to normal crossings degenerations, the mirror degenerations could only involve
irreducible components isomorphic to $\PP^n$.

As much more motivation was already given in \cite{Announce}, we will
end this discussion here and rather give a summary and reader's guide
to the paper. A lot of the content is quite technical, so a good deal
can be skipped on a first reading, and we will try to explain what
is or is not essential for absorbing the basic ideas of the paper.

Section 1 is devoted entirely to affine manifolds, and the structures
on them which will play a role in the remainder of the paper. Section
1.1 is intended as a warm-up, reminding the reader of basic concepts
of affine manifolds, simple examples, and a review of the Legendre 
transform in this context. This provides important motivation, as the
Legendre transform is key to Hitchin's elaboration of the SYZ 
(Strominger-Yau-Zaslow)
approach to mirror symmetry \cite{SYZ}, \cite{Hit}. As we replace the Legendre
transform with the discrete Legendre transform, we are able to argue
that our approach is an algebro-geometric version of SYZ.

Section 1.2 introduces affine manifolds with singularities, and gives 
some basic examples, such as Example \ref{polytope}. These first two 
sections are straightforward and are essential reading.

In \S1.3, we introduce the basic combinatorial objects: polyhedral
decompositions of integral affine manifolds with singularities. This
is essentially just a decomposition of an integral affine manifold
with  singularities into lattice polytopes, though there is some
subtlety in how such polytopes are allowed to interact with the
singular set of the affine manifold. The definition, Definition
\ref{PD2}, is essential for this paper, and Construction
\ref{basicconstruct} is the prime example: we will use this to
construct the discrete Legendre transform and the dual intersection
complex. We note that both this definition and construction have
already been given in \cite{Announce}. The remainder of the section
explores properties of polyhedral decompositions and auxiliary
constructions. These are all quite important for the paper, but
proofs may be skipped without great harm.

Section 1.4 describes the discrete Legendre transform and its
properties: Proposition \ref{dlt} and Proposition \ref{dlt2}
summarize these. We make connections with the standard discrete
Legendre transform on $\RR^n$ in Example \ref{torusdlt}, and with
Batyrev duality in Example \ref{Batyrevdual}.

Section 1.5 defines two properties of affine manifolds with
singularities. The first is a positivity property, a generalisation
of the phenomenon which occurs in elliptic fibrations which forces
suitably normalized monodromy around Kodaira-type $I_n$ fibres to be
positive in a certain sense. The second is simplicity, an analogue of
the notion introduced for torus fibrations in \cite{SlagI}. These
definitions, especially that of simplicity, are a bit technical, and
the whole section could be skipped until these definitions are used
in \S\S4 and 5.

In \S2, we begin the process of constructing toric log Calabi-Yau
spaces from pairs $(B,\P)$, where $B$ is an integral affine manifold
with singularities and $\P$ is a polyhedral decomposition on $B$. We
actually give two dual constructions in \S\S2.1 and 2.2, depending
on whether we view $(B,\P)$ as the intersection complex or dual
intersection complex of the degeneration. The former case, which we
refer to as the cone picture, is technically easier and in fact
produces a projective scheme (if $B$ is compact), 
whereas the latter case, which we
refer to as the fan picture, is technically harder and produces only
an algebraic space. In fact, most of \S2.2 is devoted to
constructing an explicit \' etale open cover of this algebraic space,
which we will need for later aspects of the construction. This
construction has a number of subtleties having to do with the role of
the singularities in $B$. We recommend strongly reading \S2 through
Example \ref{twopyramids}, as this will illustrate these subtleties.
Most of the rest of the section is devoted to the details of the
construction of an algebraic space $X_0(B,\P,s)$ starting with data
$(B,\P)$ and so-called open gluing data $s$ (Definition
\ref{open_gluing_data}). This can be skimmed, but open gluing data
will play an important role. At the end of \S2.2, we compare the fan
and cone pictures, and compute the dualizing sheaf and basic
invariants of $X_0(B,\P,s)$. These statements should be read.

Section 3 introduces log structures. The main new idea in this paper
is the significance of log structures in mirror symmetry. While the
definition of a log structure is very simple, learning to work with
them can take some time. So \S3.1 provides an introduction to the
sorts of log structures we use. Given a toric degeneration
$f:\X\rightarrow\shS$, the special fibre $\X_0$ inherits a canonical
log structure from the inclusion $\X_0\subseteq\X$. This log
structure should be thought of as providing key information about the
smoothing. Furthermore, while our construction of the dual
intersection complex as an affine manifold depends not just on $\X_0$
but the map $f:\X\rightarrow\shS$, in fact all the information needed
about $f$ is contained in the log structure on $\X_0$. So our dual
intersection complex construction makes sense not just for toric
degenerations, but also for what we call toric log Calabi-Yau spaces,
log spaces which look like degenerate fibres of toric degenerations,
defined over the standard log point.  Thus, if one wants to
understand mirror symmetry by studying  singular fibres of
degenerations, one cannot work just with the  singular fibres, but
must also involve the log structure. In fact, the discrete Legendre
transform will interchange information about the irreducible
components of a toric log Calabi-Yau space with information about the
log structure of its mirror. Hence, log structures lie at the heart
of our construction. Section 3.1 is thus essential reading.

In \S\S3.2 and 3.3 we begin to address the question: given an 
affine manifold with singularities $B$ with polyhedral decomposition
$\P$ and open gluing data $s$, how do we put a log structure on
$X_0(B,\P,s)$ in order to make it a log Calabi-Yau space? We are able
to construct a sheaf of sets on $X_0(B,\P,s)$ whose sections define
suitable log structures. This construction is explained in \S3.2,
and in \S3.3, we actually compute this sheaf explicitly. This latter
section is the technical heart of the paper, and can  probably be
skipped on a first reading. The eventual goal is to identify global
sections of this sheaf; this is finally accomplished in the case that
$B$ is positive and simple in Theorem \ref{lifted}, with one of
the main theorems of the paper being Theorem \ref{moduli_theorem}. 
This allows us, in the simple case, to calculate the moduli space of
log Calabi-Yau spaces with a given dual intersection complex
$(B,\P)$. This moduli space will play the role of complex moduli in
mirror symmetry for toric log Calabi-Yau spaces. 

Section 4 introduces the notions of toric degenerations and log
Calabi-Yau spaces, and in \S4.1 we give these definitions and the
construction of the dual intersection complex, reversing the
constructions of the previous two chapters. The construction of the
dual intersection  complex was already explained in less detail in
\cite{Announce}. A brief \S4.2 discusses the polarized case and the
intersection complex. Section 4.3 makes the connection with
positivity: we show the dual intersection complex of a toric
degeneration of Calabi-Yau varieties is always positive, hence
justifying the definition of positive. Finally, in  \S4.4, we continue
some of the calculations leading up to the proof of
Theorems~\ref{lifted} and \ref{moduli_theorem}. This is again
technical, but Examples~\ref{obstructed} and \ref{degenexamples}
should prove informative.

Section 5 ties together all the strands so far. We complete the 
calculation of the moduli of log Calabi-Yau spaces with a given dual
intersection complex in the simple case. 
The proof relies on most of the technical
aspects developed so far in the paper. However, the answer
is elegant: the moduli space in fact coincides with a
cohomology group of a natural sheaf on $B$, determined canonically by
the affine structure. Furthermore, this is precisely the group
expected from previous experience with the Strominger-Yau-Zaslow
conjecture. In \S5.2, we compute the log Picard group;
this can probably be skipped altogether as it merely provides
motivation for the definition of the log K\"ahler moduli space in
\S5.3. The log K\"ahler moduli space will be canonically isomorphic
to the moduli space of log Calabi-Yau spaces for the mirror. There is
no technical content in \S5.3, and in fact the reader may wish to
read this section early on! There we also describe the many aspects
of this program left undone.

\begin{remark}
\label{convention}
Contrary to \cite{Oda} or \cite{Fulton}, we use
the following convention. Let $\Sigma$ be a fan defining a toric
variety $X$, with toric Weil divisors $D_1,\ldots,D_n$ corresponding
to rays $R_1,\ldots,R_n$ of $\Sigma$ with primitive generators
$v_1,\ldots,v_n$. If $D=\sum a_i D_i$ is a Cartier divisor, we take
the piecewise linear function $\psi$ on the fan $\Sigma$
corresponding to $D$ to take the values $\psi(v_i)=a_i$, rather than
$\psi(v_i)=-a_i$, as in \cite{Oda} or \cite{Fulton}. This will affect
various signs throughout.
\qed
\end{remark}

We would like to thank 
Klaus Altmann,
Robert Friedman,
Maxim Kontsevich,
Marco K\"uhnel,
Arthur Ogus,
Martin Olsson,
Simone Pavanelli,
Stefan Schr\"oer,
Balazs Szendr\"oi,
Richard Thomas,
and Ilia Zharkov for useful conversations during the progress
of this work.

Part of this work has been done while the first author visited RIMS in
Kyoto and the second author stayed at the Institut de Math\'ematiques
de Jussieu in Paris. They thank these named institutions for
hospitality, and in particular their respective hosts, Kyoji Saito and
Claire Voisin. We would also especially like to thank the referee for
his invaluable comments.

\newpage

\section{Affine Manifolds}
\label{section1}

\subsection{Affine Manifolds and Invariants}

We will start by reviewing some basic notions concerning affine
manifolds and their relation to mirror symmetry. For basic
information on affine manifolds, we follow Goldman and Hirsch's paper
\cite{GH}.

We fix $M=\ZZ^n$ an abelian group, $N=\Hom_{\ZZ}(M,\ZZ)$,
$M_{\RR}=M\otimes_{\ZZ} \RR$, $N_{\RR}=N\otimes_{\ZZ} {\RR}$. Then
$N=\dualvs{M}$ and $N_\RR=\dualvs{M_\RR}$ if
$\dualvs{\Lambda}=\Hom_\ZZ(\Lambda ,\ZZ)$ and
$\dualvs{V}=\Hom_\RR(V,\RR)$ for an abelian group $\Lambda$ and an
$\RR$-vector space $V$ respectively, and $\otimes_\ZZ \RR$ is viewed
as functor taking abelian groups to $\RR$-vector spaces.

We set
\[
\Aff(M_{\RR})=M_{\RR}\rtimes \GL_n(\RR)
\]
to be the group of affine transformations of $M_{\RR}$, with subgroup
\[
\Aff(M)=M\rtimes \GL_n(\ZZ).
\]

If $M$ and $M'$ are two different lattices, then we denote by
$\Aff(M_{\RR},M'_{\RR})$ the $\RR$-vector space of affine maps
between $M_{\RR}$ and $M'_{\RR}$. Here
\[
\Aff(M_{\RR},M'_{\RR})=M_{\RR}'\times\Hom(M_{\RR},M'_{\RR}).
\]
Similarly 
\[
\Aff(M,M')=M'\times\Hom(M,M')
\]
is the $\ZZ$-module of affine maps between $M$ and $M'$.

\begin{definition}
Let $B$ be an $n$-dimensional manifold. An \emph{affine structure} on
$B$ is given by an open cover $\{U_i\}$ along with coordinate charts
$\psi_i:U_i\to M_{\RR}$, whose transition functions
$\psi_i\circ\psi_j^{-1}$ lie in $\Aff(M_{\RR})$. The affine structure
is \emph{integral} if the transition functions lie in $\Aff(M)$.  If
$B$ and $B'$ are (integral) affine manifolds of dimensions $n$ and
$n'$ respectively, then a continuous map $f:B\to B'$ is (integral)
affine if locally $f$ is given by elements of $\Aff(\RR^n,\RR^{n'})$
($\Aff(\ZZ^n, \ZZ^{n'})$). If in addition $f$ is a local
diffeomorphism, we say $f$ is \' etale (integral) affine.
\end{definition}

\begin{remark} In other papers on the role of affine manifolds in
mirror symmetry \cite{KS}, \cite{HZ}, affine manifolds with
transition maps in $M_{\RR}\rtimes\GL_n(\ZZ)$ are considered. We restrict
to the integral case here because these are the affine manifolds
which arise when studying degenerations; this is roughly equivalent
on the mirror side to focusing on symplectic manifolds with \emph{integral}
symplectic forms. Dual intersection complexes arising from toric degenerations
are integral affine manifolds with singularities, as are intersection
complexes arising from a choice of ample line bundle on a toric
degeneration. In particular, the examples given in the introduction
associated to the quartic are integral.
\end{remark}

\begin{proposition}
 Let $\pi:\tilde B\to B$ be the universal covering of an (integral)
affine manifold $B$, inducing an (integral) affine structure on
$\tilde B$. Then there is an \'etale (integral) affine map
$\delta:\tilde B\to M_{\RR}$, called the \emph{developing map}, and
any two such maps differ only by an (integral) affine transformation.
\end{proposition}

\proof
This is standard, see \cite{GH}, pg. 641 for a proof.
\qed\medskip

Note that there is no need for the developing map to be injective or
a covering space; it is only a local isomorphism in general.

\begin{definition}
The fundamental group $\pi_1(B)$ acts on $\tilde B$ by deck
transformations; for $\gamma\in\pi_1(B)$, let $\Psi_{\gamma}:\tilde
B\to\tilde B$ be the corresponding deck transformation with
$\Psi_{\gamma_2}\circ \Psi_{\gamma_1}=\Psi_{\gamma_1\gamma_2}$. Then
by the uniqueness of the developing map, there exists a
$\rho(\gamma)\in \Aff(M_{\RR})$ such that $\rho(\gamma)\circ
\delta\circ \Psi_{\gamma}=\delta$. The map $\rho:\pi_1(B)\to
\Aff(M_{\RR})$ is called the \emph{holonomy representation}. If the
affine structure is integral, then $\im\rho \subseteq \Aff(M)$.
\end{definition}

Recall that we compose loops $\gamma_1$ and $\gamma_2$ in $\pi_1(B)$
so that $\gamma_1\gamma_2$ is obtained by first following $\gamma_1$
and then $\gamma_2$. Hence $\Psi_{\gamma_1\gamma_2}=
\Psi_{\gamma_2}\circ \Psi_{\gamma_1}$, from which it follows we have
defined $\rho$ to be a group homomorphism. 
A different way to view this is by observing that $B$ is naturally
endowed with an \emph{affine} connection (as opposed to the 
ubiquitous linear connections) by pulling back the standard
affine connection on $M_{\RR}\cong \RR^n$ via charts. Then
$\rho(\gamma)$ is given by parallel transport along 
$\gamma^{-1}$ with respect to the affine connection. (See
\cite{GH}.)

Note conversely that given an immersion $\delta:\tilde B\to M_{\RR}$
and a representation $\rho:\pi_1(B)\to \Aff(M_{\RR})$ such that
$\rho(\gamma)\circ \delta\circ \Psi_{\gamma}=\delta$, these data
induce an affine structure on $B$.

\begin{example}
\label{torusexamples}
(1) If $B$ is an (integral) affine manifold and $G$ is a group acting
properly and discontinuously on  $B$ via (integral) affine
transformations, then $B/G$ inherits an (integral) affine structure
from $B$.

(2) As $M_{\RR}$ is naturally an affine manifold, with the developing
map being the identity, if $\Gamma\subseteq M_{\RR}$ is any lattice
acting by translations on $M_{\RR}$, we obtain an affine structure on
$M_{\RR}/\Gamma$. Here
for $\lambda\in \pi_1(M_{\RR}/\Gamma) =\Gamma$, $\rho(\lambda)$ is
translation by $-\lambda$. Thus the affine  structure is integral if
and only if $\Gamma\subseteq M$.

(3) This example is from \cite{Baues}. Take $M=\ZZ^2$, and consider
the subgroup $G\subseteq \Aff(M_{\RR})$ defined by
\[
G=\{A\in \Aff(M_{\RR})| \hbox{$A(m_1,m_2)=(m_1+vm_2+u+v(v-1)/2,m_2+v)$ for
$u,v\in\RR$}\}.
\]
(This is not quite the form given in \cite{Baues}, but rather $G$ has
been conjugated by translation by $(0,-1/2)$ to obtain better
integrality properties). $G$ is isomorphic to $\RR^2$, and if we
choose any lattice $\Gamma\subseteq G$, then $\Gamma$ acts properly
and discontinuously, so that $M_{\RR}/\Gamma$ is an affine manifold,
topologically a two-torus. This is the only other affine structure on
the two-torus obtained from $M_{\RR}$ by dividing out by a properly
discontinuous group action. (See \cite{Baues}, Theorem 4.5). Note the
affine structure is integral with respect to the integral structure
$M\subseteq M_{\RR}$ if and only if 
\[
\Gamma\subseteq
\{A\in \Aff(M_{\RR})| \hbox{$A(m_1,m_2)= (m_1+vm_2+u+v(v-1)/2,m_2+v)$
for $u\in\ZZ$, $v\in \ZZ$}\}.
\]
\qed
\end{example}

We recall the notion of \emph{radiance obstruction} from \cite{GH}.

\begin{definition}
Let $\Lin:\Aff(M_{\RR})\to \GL_n(\RR)$ and $\Trans:\Aff(M_{\RR})\to
M_{\RR}$ be the projections. Here $\Lin$ is a homomorphism, but
$\Trans$ is a crossed homomorphism with respect to the regular
representation of $\GL_n({\RR})$ on $M_{\RR}$, i.e.
\[
\Trans(A_1A_2)=\Lin(A_1)(\Trans(A_2))+\Trans(A_1).
\]
Given an affine representation $\rho:G\to \Aff(M_{\RR})$ of a group
$G$, set $\tilde\rho=\Lin\circ\rho$. Then $\Trans\circ\rho$ can be
interpreted as an element $c_{\rho}\in H^1(G, M_{\RR}^{\tilde\rho})$,
where $M_{\RR}^{\tilde\rho}$ denotes the  $G$-module $M_{\RR}$
defined by the representation $\tilde\rho$. The class $c_{\rho}$ is
called the \emph{radiance obstruction} of $\rho$. If $B$ is an affine
manifold and $\rho:\pi_1(B)\rightarrow\Aff(M_{\RR})$ is the holonomy
representation of $B$, then the \emph{radiance obstruction} of $B$
is $c_{B_0}:=c_{\rho}$.
\qed
\end{definition}

The radiance obstruction is an important invariant of an affine
manifold. It can be viewed as being analogous to the cohomology class
of the symplectic form on a symplectic manifold.

\begin{theorem}
Two affine representations $\rho_1,\rho_2:G\to \Aff(M_{\RR})$ with
$\tilde\rho_1=\tilde\rho_2$ are conjugate by a translation if and
only if $c_{\rho_1}=c_{\rho_2}$. 
\end{theorem}

\proof \cite{GH}, page 631. \qed\medskip

This is important for the following reason.  If we compose an affine
structure with a translation, i.e.\ replace $\delta$ with $\tau_a\circ
\delta$ where $\tau_a$ denotes translation by $a$, then the holonomy
representation $\rho$ is replaced by $\rho'$ with
$\rho'(\gamma)=\tau_a\circ \rho(\gamma)\circ  \tau_a^{-1}$, and then
$\Trans\circ\rho$ and $\Trans\circ\rho'$ are cohomologous:
\[
\Trans(\rho'(\gamma))=a-\tilde\rho(\gamma)(a)+\Trans(\rho(\gamma)).
\]
Thus holonomy representations from affine structures related by
translation are conjugate by a translation, so the radiance
obstruction helps classify  holonomy representations. In particular,
this allows us to identify integral affine structures:

\begin{proposition}
An affine representation $\rho:G\to M_{\RR}\rtimes
\GL_n(\ZZ)\subseteq \Aff(M_{\RR})$ is conjugate by a translation to a
representation $\rho':G\to \Aff(M)$ if and only if the radiance
obstruction $c_{\rho}\in H^1(G,M_{\RR}^ {\tilde\rho})$ is in the
image of $H^1(G,M^{\tilde\rho})\to H^1(G,M_{\RR}^{\tilde\rho})$.
\end{proposition}

\proof
If $\rho$ is conjugate to a representation $\rho':G\to  \Aff(M)$ then
$u=\Trans\circ \rho':G\to M$ is a crossed homomorphism representing an
element of $H^1(G,M^{\tilde\rho})$ whose image in
$H^1(G,M^{\tilde\rho}_{\RR})$ is $c_{\rho}$. Conversely, if
$c_{\rho}$ is in the image of $H^1(G,M^{\tilde\rho})$, then $u$ is
cohomologous to a $v:G\to M^{\tilde\rho}$, i.e. there exists $a\in
M_{\RR}$ such that $u(g)-v(g)=a-\tilde\rho(g)a$. If $\tau_a$ denotes
translation by $a$, then this says
$\tau_a^{-1}\circ\rho(g)\circ\tau_a$ is in $\Aff(M)$ for all $g\in
G$.
\qed
\medskip

The radiance obstruction in fact also plays a role in understanding
the Legendre transform, as we shall see, and this is the primary
reason for introducing it here.

\begin{definition}
If $B$ is an affine manifold, there is a flat linear connection $\nabla$ on
$\shT_B$, where if $y_1,\ldots,y_n$ are local affine coordinates,
$\partial/\partial {y_1},\ldots,\partial/\partial{y_n}$ are a frame
of flat sections of $\shT_B$. Denote by $\Lambda_{\RR}$
the local system of flat sections and $\check\Lambda_{\RR}$
the dual local system of flat sections of the
dual connection on $\T_B$. If furthermore the holonomy of $B$ is contained in
$M_{\RR}\rtimes \GL_n(\ZZ)$, rather than $M_{\RR}\rtimes\GL_n(\RR)$,
then there exist integral subsystems
$\Lambda\subseteq\Lambda_{\RR}$ and $\check\Lambda \subseteq
\check\Lambda_{\RR}$ coming from the inclusions $M\subseteq M_{\RR}$
and $N\subseteq N_{\RR}$. 
\end{definition}

We note that the monodromy of the local system $\Lambda_{\RR}$
is precisely the linear part of the holonomy representation. Thus
the radiance obstruction can be viewed as measuring the difference
between the monodromy of $\Lambda_{\RR}$ and the holonomy representation.

\begin{remark}
\label{realising}
\cite{GH} gives a number of ways of realising the radiance obstruction.
The \v Cech realisation will also be of use to us.

Choose an open covering $\{U_i\}$ of $B$ along with affine charts
$\psi_i:U_i\to M_{\RR}$. Such a chart allows us to identify
$\T_{U_i}$ canonically with $U_i\times M_{\RR}$ and the graph of
$\psi_i$ can be viewed as a section $s_i\in\Gamma(U_i,\T_{U_i})$,
which is parallel for the affine connection, and hence is 
independent of $\psi_i$ up to addition by flat
sections of $\T_{U_i}$. Then $(s_j-s_i)_{ij}$ form a \v Cech
1-cocycle for $\Lambda_{\RR}$, hence represents an element of
$H^1(B,\Lambda_{\RR})$. This group is naturally isomorphic to
$H^1(\pi_1(B),M_{\RR}^{\tilde\rho})$, and this \v Cech 1-cocycle
represents the radiance obstruction under this isomorphism. If the
charts $\psi_i$ are integral, then $s_j-s_i\in \Gamma(U_i\cap
U_j,\Lambda)$, yielding the radiance obstruction in
$H^1(B,\Lambda)$.
\qed
\end{remark}

\begin{definition}
\label{fundseq1}
Let $\shAff_{\RR}(B,{\RR})$ denote the sheaf of affine maps from $B$ to
$\RR$, i.e. functions locally of the form $f\circ \delta$, where
$f\in \shAff(M_{\RR},{\RR})$ and $\delta$ is the developing map. This
sheaf fits into a natural exact sequence
\[
\exact{{\RR}}{\shAff_{\RR}(B,{\RR})}{\check\Lambda_{\RR}},
\]
analogous to the exact sequence of $\RR$-vector spaces \emph{affine}
\[
\exact{{\RR}}{\shAff(M_{\RR},{\RR})}{N_{\RR}}.
\]
Similarly, if $B$ is an integral affine manifold, define
$\shAff(B,\ZZ)$ to be the sheaf of affine functions on $B$ locally of
the form $f\circ \delta$, where $f\in \shAff(M,\ZZ)$. There is an 
exact sequence
\[
\exact{\ZZ}{\shAff(B,\ZZ)}{\check\Lambda}.
\]
\qed
\end{definition}

Another description of the radiance obstruction:

\begin{proposition}
\label{extclass}
The extension class of
\[
0\lra\RR\lra\shAff_{\RR}(B,\RR)\lra
\check\Lambda_{\RR}\lra 0
\]
in $\Ext^1(\check\Lambda_{\RR},\RR)=H^1(B,\Lambda_{\RR})$
coincides with the radiance obstruction of $B$.
\end{proposition}

\proof Let $\{U_i\}$ be a cover of $B$ of contractible open
sets, and on each $U_i$ choose a splitting $\alpha_i:
\check\Lambda_{\RR}\rightarrow\shAff_{\RR}(B,\RR)$. Then 
the extension class is determined by the \v Cech cocycle
$(\alpha_j-\alpha_i)_{ij}$, $\alpha_j-\alpha_i:\check\Lambda_{\RR}
\rightarrow\RR$. To compare this with the radiance obstruction, note
$\alpha_i$ determines an affine chart $\psi_i:U_i\rightarrow
\Gamma(U_i,\Lambda_{\RR})=\Hom(\Gamma(U_i,\check\Lambda_{\RR}),\RR)$
by $\psi_i(b)(n)=\alpha_i(n)(b)$, for $b\in U_i$, $n\in \Gamma(U_i,
\check\Lambda_{\RR})$ so that $\alpha_i(n)$ is an affine linear
function on $U_i$. As in Remark \ref{realising}, each chart $\psi_i$
determines a well-defined section $s_i\in\Gamma(U_i,\T_{U_i})$, the
graph of $\psi_i$, and the radiance obstruction
is represented by $(s_j-s_i)_{ij}$. But for $b\in U_i\cap U_j$, 
$(s_j-s_i)(b)$ is the functional on $\check\Lambda_{\RR,b}
=\T_{B,b}^*$ given by 
$(s_j-s_i)(b)(n)=\alpha_j(n)(b)-\alpha_i(n)(b)$. But this is
precisely $\alpha_j-\alpha_i:\check\Lambda_{\RR,b}\rightarrow\RR$,
which doesn't depend on $b$, so $s_j-s_i$ and $\alpha_j-\alpha_i$
coincide.
\qed
 \medskip

To make contact with the metric form of SYZ, we briefly discuss
metrics on affine manifolds.

\begin{definition} 
\label{hessian}
Let $B$ be an affine manifold. A Hessian metric $g$ on $B$ is a
Riemannian metric on $B$ such that locally, for affine coordinates
$(y_1,\ldots,y_n)$, there is a potential function $K$ such that
$g_{ij}=\partial^2 K/\partial y_i\partial y_j$. The pair $(B,g)$ is
called an \emph{affine K\"ahler} manifold or a \emph{Hessian} 
manifold. This can be defined in a coordinate independent way
as $g=\nabla dK$, a section of the bundle 
$S^2\T^*\subseteq \T^*\otimes \T^*$. 
\end{definition}

Affine K\"ahler metrics were first studied by Cheng and Yau in \cite{ChengYau},
where metrics whose potential in addition
satisfies the Monge-Amp\`ere equation were studied. Such metrics should
be especially important in the study of mirror symmetry, but we do not pursue
this further here. Affine K\"ahler manifolds
were called Hessian manifolds in \cite{SY}.

The potential function $K$ is only defined locally, up to an affine
function. Thus $K$ can actually be defined as a multi-valued function
on $B$ well defined up to affine functions, and we consider such
functions more generally.

Let $B$ be an affine manifold, $\pi:\tilde B\to B$ the universal
cover and $\delta:\tilde B\to M_{\RR}$ the developing map. We
consider continuous maps $\varphi:\tilde B\to \RR$ which satisfy the
condition
\[
\varphi-\varphi\circ \Psi_{\gamma}=\alpha(\gamma),
\]
where $\alpha$ is a map $\alpha:\pi_1(B)\to \shAff(M_{\RR},\RR)$. Here
an element of $\shAff(M_{\RR},\RR)$ induces a map on $\tilde B$ via
composition with $\delta$. Of course $\shAff(M_{\RR},\RR)$ is a left
$\pi_1(B)$-module, with $\gamma$ acting by composition on the right
with $\rho(\gamma^{-1})$ (or thinking of an element of $\shAff(M_{\RR},
\RR)$ as a map on $\tilde B\to {\RR}$, $\gamma$ acts by composition
with $\Psi_{\gamma}$). Then $\alpha$ is a crossed homomorphism:
\begin{eqnarray*}
\alpha(\gamma_1\gamma_2)&=&\varphi-\varphi\circ
\Psi_{\gamma_2}\circ \Psi_{\gamma_1}\\
&=&\varphi-(\varphi-\alpha(\gamma_2))\circ \Psi_{\gamma_1}\\
&=&\gamma_1(\alpha(\gamma_2))+\alpha(\gamma_1).
\end{eqnarray*}
Thus $\alpha$ defines an element $[\alpha]\in
H^1(\pi_1(B),\shAff(M_{\RR},\RR))$. Furthermore, adding an element of
$\shAff(M_{\RR},\RR)$ to $\varphi$ replaces $\alpha$ with a
cohomologous $\alpha'$, and any representative of $[\alpha]$ can be
obtained in this way.

Alternatively, we can define $[\alpha]\in H^1(B,\shAff_{\RR}(B,{\RR}))
\cong H^1(\pi_1(B),\shAff(M_{\RR},\RR))$. Choosing a covering of $B$ by
simply connected sets $U_i$ and choosing a representative $\varphi_i$
for $\varphi$ on $U_i$, $\varphi_i-\varphi_j\in \shAff_{\RR} (U_i\cap
U_j, \RR)$ and hence we obtain a \v Cech cocycle in
$H^1(B,\shAff_{\RR}(B,{\RR}))$.

Suppose now that $B$ carries a metric of Hessian form. Any two local
potential functions $K$ for the metric differ by an affine function,
so we can patch to get $K:\tilde B\to \RR$. If $K-K\circ
\Psi_{\gamma}=\alpha(\gamma)$, the class $[\alpha]\in
H^1(\pi_1(B),\shAff(M_{\RR},\RR))\cong H^1(B,\shAff_{\RR}(B,\RR))$ is
called the \emph{class} of the metric $g$, as defined by Kontsevich
and Soibelman \cite{KS}.

The developing map yields an isomorphism
\[
\delta^*:\delta^*\T^*_{M_{\RR}}\to\T_{\tilde B}^*.
\]
Since $\T_{M_{\RR}}^*=M_{\RR}\times N_{\RR}$, $\delta^*$ gives
an isomorphism of $\tilde B\times N_{\RR}$ with $\T_{\tilde B}^*$.
Let $q:\tilde B\times N_{\RR}\to N_{\RR}$ be the projection.
Then $(\delta^*)^{-1}(dK)$ is a section of $\delta^*\T_{M_{\RR}}^*$,
so we can view
\[
\check\delta:=q\circ (\delta^*)^{-1}(dK)
\]
as a function $\check\delta:\tilde B\to N_{\RR}$. This is just the
differential $dK$ under these identifications. Because the Hessian of
$K$ is positive definite, $\check\delta$ is an immersion, hence defining
a new affine structure on $B$. Note that
$d\alpha(\gamma)$ is naturally identified with an element of
$N_{\RR}$. Under this identification, one can check that
\[
{}^t \tilde\rho(\gamma^{-1})\circ \check\delta\circ \Psi_{\gamma}
+d\alpha(\gamma)=\check\delta,
\]
and thus $\pi_1(B)$ acts on the new affine structure on $\tilde B$ by
affine transformations. Dividing $\tilde B$ by this action, we obtain a new
affine structure on $B$, which we denote by $\check B$. So the
holonomy representation  $\check\rho:\pi_1(\check B)\to
\Aff(N_{\RR})$ of the affine structure given by $\check\delta$ has
linear part $\dualvs{\tilde\rho}$ dual to the representation
$\tilde\rho$. Also $d\alpha(\gamma)\in N_{\RR}$ is just the
projection of $\alpha(\gamma)\in \shAff(M_{\RR},\RR)$ onto $N_{\RR}$.
Thus the radiance obstruction of the affine structure $\check \delta$
is just the projection of $[\alpha]\in H^1(\pi_1(\check B),
\shAff(M_{\RR},\RR))$ to $H^1(\pi_1(B),N_{\RR})$.

We can also define the Legendre transform of $K$ as, for $x\in \tilde B$,
\[
\check K(x)=\langle \check\delta(x),\delta(x)\rangle -K(x).
\]
Then the Hessian of $\check K$ defines the same metric as $K$ on $B$,
and $d\check K=\delta$. For a proof of this, see \cite{Hit}, \S 5.
We call $(\check B,\check K)$ the
\emph{Legendre transform} of $(B,K)$. We note we have shown the
radiance obstruction of $\check B$ is determined by the class of the
metric. Conversely, the class of the metric on $\check B$ is at least
partially determined by the radiance obstruction of $B$, since
$d\check K=\delta$.

The role of the Legendre transform in the SYZ picture of mirror
symmetry is well-understood, and was first explained by Hitchin
\cite{Hit}. See also \cite{Leung}, which develops this point of view
further. It is precisely the presence of the potential function $K$
which allows us to pass between $B$ and $\check B$, thus dualizing
affine manifolds. One of the key points of this paper is that we can
in fact formulate a discrete version of this, and replace $K$ with a
piecewise linear function. This will enable us to obtain an
algebro-geometric analog of SYZ.

\begin{example}
In Example~\ref{torusexamples}, (2), we can take any convex quadratic
function $K:M_{\RR}\to\RR$ to serve as a potential. Now $K$ satisfies
the periodicity condition
\[
K(x+\gamma)=K(x)+\alpha(\gamma)(x)
\]
for $\gamma\in\Gamma$, $x\in M_{\RR}$, and some $\alpha(\gamma) \in
\shAff(M_{\RR},\RR)$. Taking differentials of this equation, we get
$\check\delta(x+\gamma)=\check\delta(x)+d\alpha(\gamma)\in N_{\RR}$.
Thus if we set
\[
\check\Gamma=\{d\alpha(\gamma)\in N_{\RR}|\gamma\in\Gamma\},
\]
then the dual torus can be identified as $N_{\RR}/\check\Gamma$.

In Example~\ref{torusexamples}, (3), in fact there is no convex
function $K:M_{\RR} \to \RR$ defining a Hessian metric on $B$. (The
easiest way to see this is to note that the torus bundle
$X(B):=\T_B/\Lambda$ has a natural complex structure on it making 
$X(B)$ isomorphic to a primary Kodaira surface, and the pull-back of
$K$ to $X(B)$ is the K\"ahler potential of a K\"ahler metric. However
a primary Kodaira surface is not K\"ahler). \qed
\end{example}

\subsection{Affine manifolds with singularities}

\begin{definition}
An affine manifold with singularities is a topological manifold $B$
along with a closed set $\Delta\subseteq B$ which is locally a finite union of
locally closed submanifolds of codimension at least 2, and an affine
structure on $B_0=B\setminus \Delta$. An affine manifold with
singularities is \emph{integral} if the affine structure on $B_0$ is
integral.  We always denote by $i:B_0\to B$ the inclusion map. A
continuous map $f:B\to B'$ of (integral) affine manifolds with
singularities is (integral) affine if $f^{-1}(B_0')\cap B_0$ is dense
in $B$ and
\[
f|_{f^{-1}(B_0')\cap B_0}:f^{-1}(B_0')\cap B_0\to B_0'
\]
is (integral) affine.
\qed
\end{definition}

\begin{example}
\label{twodimaffine}
We give a variety of two-dimensional examples.

(1) Let $B=\{z\in \CC| |z|<1\}$, $B_0=B\setminus \{0\}$. Let
$\shH\rightarrow B_0$ be the universal cover, where $\shH$ is the upper
half-plane with coordinate $w$ and covering map given by $w\mapsto e^{2\pi 
i w}$. Then
\[
\delta(w)=\left({\rm Re}\left( e^{2\pi i w}\right), {\rm
Re}\left(ne^{2\pi i w} \left(w-{1\over 2\pi i}\right)\right)\right)
\]
defines an affine structure on $B_0$ (this is the developing map
of $\shH$ into ${\RR }^2$). If $\gamma$ is a
counterclockwise simple loop around the origin, then $\rho(\gamma)$
is linear, given by the matrix $\begin{pmatrix}1&0\cr
-n&1\cr\end{pmatrix}$.

To see how this affine structure arises, consider the family $f_0:
B_0\times\CC/\langle 1,\tau(z)\rangle\rightarrow B_0$ of elliptic
curves over $B_0$ with period $\tau(z)={n\over 2\pi i}\log z$. 
This elliptic curve can be compactified over $B$ by adding a Kodaira-type
$I_n$ fibre, i.e. a cycle of $n$ projective lines. Let $y$ be the fibre
coordinate in this family: then $dz\wedge dy$ is a holomorphic $2$-form,
and ${\rm Re}(dz\wedge dy)$ is a symplectic form, making the
elliptic fibration into a Lagrangian fibration. The Arnold-Liouville
Theorem yields an affine structure on $B_0$. A choice of affine coordinates
$x_1,x_2$ in a neighbourhood $U\subseteq B_0$ is given by a choice
of continuously varying local basis $\gamma_1,\gamma_2$ for $H_1(f^{-1}(b),
\ZZ)$ for $b\in U$. Then at $b\in U$,
\[
dx_j=dz_1\left(\int_{\gamma_j}(\iota(\partial/\partial z_1){\rm Re}(dz\wedge dy)
)|_{f_0^{-1}(b)}
\right)+dz_2\left(\int_{\gamma_j}(\iota(\partial/\partial z_2){\rm Re}
(dz\wedge dy))|_{f_0^{-1}(b)}\right),
\]
for $z=z_1+iz_2$. (See \cite{SlagII}, \S 2 for more details concerning
the Arnold-Liouville theorem in the special Lagrangian situation). 
Taking $\gamma_1$ to be given by the period 1 and $\gamma_2$ to be given
by the period $\tau(z)$, then $dx_1={\rm Re}(dz)$ and $dx_2={\rm Re}(\tau(z)
dz)$, so we can take $x_1={\rm Re} z$ and $x_2={\rm Re}\left( {n\over
2\pi i}(z\log z-z)\right)$ as given.

We note that ${\rm Re}z=0$ gives a well-defined line through the origin. 
This allows one to identify this affine structure in a neighbourhood
of zero for $n=1$ with a neighbourhood
of the singularity defined in (2).

(2) Next we will give an affine manifold (with boundary)
with singularities by gluing together polyhedra. The affine
manifold will be a  union of two triangles, as depicted in either the
left or right hand pictures below, but we define the affine structure
by drawing the affine embedding of two open sets covering $B$,
obtained by making cuts as shown. Here the solid lines denote cuts,
the ``$\times$'' being the singular point, and affine coordinates are
given at the vertices of the two triangles. Thus the two pictures give
two systems of affine coordinates, linear on each triangle. The
intersection of the two coordinate charts is just $B$ minus the common
edge of the two triangles. On the left-hand component the change of
coordinates is the identity, but on the right-hand component the
change of coordinates is given by $(x_1,x_2)\mapsto (x_1,x_1+x_2)$.
Alternatively, the right-hand triangles in the left and right-hand pictures
are the same, identified by the linear transformation
$\begin{pmatrix} 1&0\\1&1\end{pmatrix}$. 
\begin{center}
\input{sect4fig1.pstex_t}
\end{center}

(3) The dual intersection complex of the quartic degeneration given in
the introduction is an integral affine manifold with singularities.

(4) A cone with a cone angle $0<\alpha\le 2\pi$ is an example of an
affine manifold with singularities: take the wedge  $\{z\in\CC| 0\le
\arg z\le \alpha\}$ in $\CC$, and identify the two edges via
rotation. The holonomy is given by rotation by the angle $-\alpha$,
i.e.
\[
\begin{pmatrix} \cos\alpha & \sin\alpha\\ -\sin\alpha&\cos\alpha\end{pmatrix}
\]
Such a rotation is integral with respect to any lattice in $\CC$ if 
$\alpha=\pi$, is integral with respect to the lattice $\ZZ[i]\subseteq
\CC$ if $\alpha= n\pi/2$, $1\le n<4$, and is integral with respect to
the lattice $\ZZ[(1+i\sqrt{3})/2]$ if $\alpha=n\pi/3$, $1\le n<6$.
In no other cases is the cone integral.

In the case these affine structures are integral, these are all
interesting singularities. However, in this paper, the only
two-dimensional singularities we will deal with in our general
construction will be those occuring in (1)-(3).  In particular, the
examples of (4) do not have polyhedral decompositions in the sense defined
in the next section. A more general approach to building
degenerations from affine manifolds with singularities should also
deal with singularities of the second type. This seems possible in
two dimensions, and there is some hope of generalising this to higher
dimensions.
\qed
\end{example}

\begin{example}
\label{ellipticK3}
Let $f:X\to\PP^1$ be an elliptically fibreed K3 surface with
holomorphic 2-form $\Omega=\Omega_1+i\Omega_2$. The fibres of $f$ are
Lagrangian with respect to both $\Omega_1$ and $\Omega_2$. As in Example
\ref{twodimaffine}, (1), by the 
Arnold-Liouville theorem, one then obtains from each of these
symplectic forms an affine structure on $B_0=\PP^1\setminus\Delta$,
where $\Delta$ is the discriminant locus of $f$. These affine
structures will in fact be related by a Legendre transform. 
Furthermore, if $[\Omega_i] \in H^2(X,\ZZ)$, then one can show the
affine structure induced by $\Omega_i$ will be integral.
\end{example}

\begin{example}
\label{polytope}
Let $\Xi\subseteq M_{\RR}\cong \RR^n$ be a polytope containing $0$ in
its interior, and let $B=\partial\Xi$. Denote by $\Bar(B)$ the first
barycentric subdivision of the boundary of $\Xi$ (see Definition
\ref{bary}). Let $\Delta$ be the union of the simplices in $\Bar(B)$
not containing any vertex of $\Xi$ nor containing the barycenter of
any $n-1$-dimensional face of $\Xi$. Set $B_0=B\setminus\Delta$. 

We can define an affine structure on $B_0$ as follows. For every 
$n-1$-dimensional face $\sigma$ of $\Xi$, let
$W_{\sigma}=\Int(\sigma)$. For each vertex $v$ of $\Xi$, let $\check
v\subseteq B$ be the union of all simplices of $\Bar(B)$ containing
$v$, and set $W_v=\Int(\check v)$. Then $\{W_{\sigma}\}\cup \{W_v\}$
form an open cover for $B_0$. We can then define affine charts for
each $n-1$-dimensional face $\sigma$, 
\[
\psi_{\sigma}:W_{\sigma}\to \AA_{\sigma}\subseteq M_{\RR}
\]
to be the inclusion of $\sigma$ inside the affine hyperplane
$\AA_{\sigma}$ spanned by $\sigma$. For a vertex $v$, define a chart
\[
\varphi_v:W_v\to M_{\RR}/\RR v
\]
by projection. It is a simple exercise to show this defines an affine
structure on $B_0$, which is integral when $\Xi$ is a reflexive
lattice polytope \cite{Bat}. See \cite{HZ} for a generalisation of this
construction.
\end{example}

\begin{example}
\label{aspmorr}
Let $M=\ZZ^4+{1\over 5}(1,2,3,4)$, and let $\Xi\subseteq M_{\RR}$ be
the polytope with vertices 
\[
(-1,-1,-1,-1), (1,0,0,0), (0,1,0,0), (0,0,1,0), (0,0,0,1).
\]
(Adding the fractional lattice point is not necessary for this
example, but will be used in Example~\ref{degenexamples}.) By
Example~\ref{polytope}, we obtain an integral affine manifold with
singularities structure on $\partial\Xi$. Furthermore, the group
$\ZZ/5\ZZ$ acts on $M$ with generator
\[
\begin{pmatrix}
0&0&0&-1\\
1&0&0&-1\\
0&1&0&-1\\
0&0&1&-1
\end{pmatrix}
\]
This cyclically permutes the vertices of $\Xi$, and one can check it
induces integral affine automorphisms of the boundary. Thus if
$B=\partial\Xi/(\ZZ/5\ZZ)$,  $B$ is also an integral affine manifold
with singularities. In this case $B$ is in fact a lens space.
\end{example}

\begin{example}
\label{enriques}
Let $\Xi\subseteq M_{\RR}=\RR^3$ be the octahedron with vertices
\[
(\pm 1,0,0),(0,\pm 1,0), (0,0,\pm 1).
\]
Then $\ZZ/2\ZZ$ acts on $\Xi$ by negation, and thus acts on
$\partial\Xi$ by affine transformations. Take $B=
\partial\Xi/(\ZZ/2\ZZ)$. This is a real projective plane, with six
singular points. 
\end{example}

\subsection{Polyhedral decompositions}
We would like to define a polyhedral decomposition of an affine
manifold $B$ with singularities. Intuitively, we want this to be a
cell decomposition of $B$ into polyhedra in affine space.  There are
two subtleties which make the definition slightly complicated. The
first is that we would like our cells to be able to be
self-intersecting: for example, if $B=\RR/\ZZ$, we would like to take
$B$ to be a maximal cell (viewed as the interval $[0,1]$ with
endpoints identified) and $0\mod \ZZ$ a zero-dimensional cell. This is
the dual intersection complex of a nodal elliptic curve, and we do not
wish to rule out such a basic example. As a result, we first define a
polyhedral decomposition of a region in $M_{\RR}$, and then use this
definition locally on $B$. Secondly, there is some subtlety in how
cells are allowed to interact with the discriminant locus, and this
will show up in the additional restriction given below for {\it toric}
polyhedral decompositions.

\begin{definition} 
\label{PD1}
A \emph{polyhedral decomposition} of a closed set $R\subseteq
M_{\RR}$ is a locally finite covering $\P$ of $R$ by closed convex
polytopes (called \emph{cells}) with the property that
\item{(1)} if $\sigma\in\P$ and $\tau\subseteq\sigma$ is a face then
$\tau\in\P$;
\item{(2)} if $\sigma,\sigma'\in\P$, then $\sigma\cap\sigma'$ is a
common face of $\sigma$ and $\sigma'$.

We say the decomposition is \emph{integral} if all vertices
(0-dimensional elements of $\P$) are contained in $M$.
\qed
\end{definition}

For a polyhedral decomposition $\P$ and $\sigma\in\P$ we define the
\emph{(relative) interior} of $\sigma$ as
\[
\Int(\sigma)=\sigma\setminus
\bigcup_{\tau\in\P,\tau\subsetneq\sigma}\tau.
\]

\begin{definition}
\label{PD2}
Let $B$ be an integral affine manifold with singularities. A
\emph{polyhedral decomposition} of $B$ is a collection $\P$ of closed
subsets of $B$ (called \emph{cells}) covering $B$ which satisfies the
following properties. If $\{v\}\in\P$ for some point $v\in B$, then
$v\not\in\Delta$ and there exists an integral polyhedral decomposition
$\P_v$ of a closed neighbourhood of the origin $R_v\subseteq
\Lambda_{\RR,v}\cong\T_{B,v}$  (the stalk of the local system
$\Lambda_{\RR}$ at $v$ or equivalently the tangent space of $B$ at
$v$) which is the closure of an open neighbourhood of the origin, and
a continuous map ${\rm exp}_v:R_v\to B$, ${\rm exp}_v(0)=v$,
satisfying
\item{(1)} ${\rm exp}_v$ is locally a homeomorphism onto its image, is
injective on $\Int(\tau)$ for all $\tau\in\P_v$, and is an integral
affine map in some neighbourhood of the origin. 
\item{(2)} For every top-dimensional $\tilde\sigma\in\P_v$,
$\exp_v(\Int(\tilde\sigma))\cap\Delta=\emptyset$ and the restriction
of ${\rm exp}_v$ to $\Int(\tilde\sigma)$ is integral affine. 
Furthermore, $\exp_v(\tilde\tau)\in\P$ for all $\tilde\tau\in\P_v$.
\item{(3)} $\hbox{$\sigma\in\P$ and $v\in\sigma$} \Leftrightarrow
\hbox{$\sigma={\rm exp_v}(\tilde\sigma)$ for some
$\tilde\sigma\in\P_v$ with $0\in\tilde\sigma$.}$
\item{(4)} Every $\sigma\in\P$ contains a point $v\in\sigma$ with 
$\{v\}\in\P$.

In addition we say the polyhedral decomposition is \emph{toric} if it satisfies
the additional condition
\item{(5)}
For each $\sigma\in\P$, there is a neighbourhood $U_{\sigma}
\subseteq B$ of $\Int(\sigma)$ and an integral affine submersion
$S_{\sigma}:U_{\sigma}\to M'_{\RR}$ where $M'$ is a lattice of rank
equal to $\dim B-\dim \sigma$ and $S_{\sigma}(\sigma\cap
U_{\sigma})=\{0\}$.
\qed
\end{definition}

We will write $\P_{\max}=\{\sigma\in\P|\hbox{$\sigma$ is maximal,
i.e. $\dim\sigma=\dim B$}\}$. 

\begin{example} 
\label{decompexamp}
(1) If $B=M_{\RR}$, $\Delta=\emptyset$, a polyhedral decomposition of
$B$ is just an integral polyhedral decomposition of $M_{\RR}$ in the
sense of Definition \ref{PD1}.  In this case, when the discriminant
locus $\Delta$ is empty, the toric condition is vacuous. If
$B=M_{\RR}/\Gamma$ for some lattice $\Gamma$, then a polyhedral
decomposition of $B$ is induced by a polyhedral decomposition of
$M_{\RR}$ invariant under $\Gamma$. Similarly, in
Example~\ref{torusexamples} (3), the same holds.  For example,
$M_{\RR}/\Gamma$ has a polyhedral decomposition containing only one
maximal cell, coming from a fundamental domain for the action of
$\Gamma$ on $M_{\RR}$. 

(2) When $B$ has singularities, the definition of toric polyhedral
decomposition imposes some slightly subtle restrictions on how the
cells of $\P$ interact with $\Delta$. In particular, it places strong
conditions on the holonomy of $B$ locally near $\Delta$: see
Proposition~\ref{monodromy2}. However, it is a useful exercise to
verify in Example~\ref{polytope} that if one takes $\P$ to be the
collection of all proper subfaces of $\Xi$, then $\P$ is a toric
polyhedral decomposition of $B$. Furthermore, in
Examples~\ref{aspmorr} and~\ref{enriques}, this polyhedral
decomposition on $\partial\Xi$ descends to a polyhedral decomposition
of $B$. In particular, in the case of Example~\ref{aspmorr}, there is
only one cell of dimension three in the decomposition, and one
vertex. One reason for the complexity of the definition of polyhedral
decomposition is that we wish to allow cells to be  self-intersecting,
i.e. be polytopes with some sides identified. These examples show the
necessity of this. \qed
\end{example}

\begin{remark}
\label{fan1}
Given a polyhedral decomposition $\P$ on $B$, if $v$ is a vertex of
$\P$, we can look at the polyhedral decomposition of $R_v$ in a small
neighbourhood of the origin in $\Lambda_{\RR,v}$. This clearly
coincides with a small neighbourhood of the origin of a complete
rational polyhedral fan $\Sigma_v$ in $\Lambda_{\RR,v}$:
\begin{center}
\includegraphics{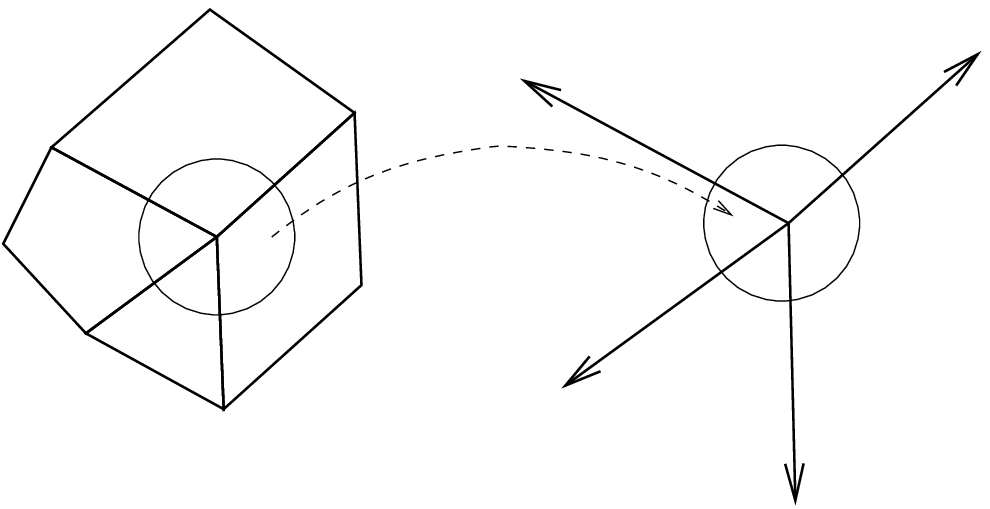}
\end{center}
In fact, we shall shortly see that the data of the affine structures on
maximal cells of $\P$ and a fan structure at each vertex
$v$ essentially determine the affine structure on $B$.
\qed
\end{remark}

\begin{definition} 
\label{bary}
Recall that if $\sigma\subset M_{\RR}$ is a polytope, then the
barycenter  $\Bar(\sigma)$ of $\sigma$ is the average of the vertices
of $\sigma$, and thus is invariant under affine transformations. The
first barycentric subdivision of $\sigma$ is then the triangulation
of $\sigma$ consisting of all simplices spanned by barycenters of
ascending chains of faces of $\sigma$. Thus given a polyhedral
decomposition $\P$ of an affine manifold with singularities $B$, we
can define the first barycentric subdivision $\Bar(\P)$ of $\P$ to be
the triangulation consisting of all images of simplices in the first
barycentric subdivisions of all $\tilde\sigma\in \P_v$ for all
vertices $v$. Because barycentric subdivisions are affine invariants,
this gives a well-defined triangulation of $B$. Note that $\Bar(\P)$
need not be a polyhedral decomposition of $B$ if
$\Delta\not=\emptyset$. For example, a vertex of $\Bar(\P)$
might be contained in $\Delta$.

For any $\tau\in\P$, let $W_{\tau}$ be the union of the interiors of
all simplices of $\Bar(\P)$ intersecting $\Int(\tau)$ (i.e. having
the barycenter of $\tau$ as a vertex). This can be thought of as the
(open) star of the barycenter of $\tau$ in the triangulation
$\Bar(\P)$ of $B$. Set $\W=\{W_{\tau}| \tau\in\P\}$; this is an open
covering of $B$.
\qed
\end{definition}

\begin{construction}
\label{basicconstruct} 
We will now describe a standard procedure for constructing affine
manifolds with singularities along with polyhedral decompositions.
Let $\P'$ be a collection of $n$-dimensional integral polytopes in
$M_{\RR}$. Suppose we are given integral affine identifications of
various proper faces of the polytopes in $\P'$ in such a way that once we
glue the polytopes using these identifications, we obtain a manifold
$B$, along with a decomposition $\P$ consisting of images of faces of
polytopes in $\P'$. In particular, we have the identification map
\[
\pi:\coprod_{\sigma'\in\P'} \sigma'\to B.
\]

Now $B$ is not yet an affine manifold with singularities. It only has
an affine structure defined in the interiors of maximal cells. When
two polytopes of $\P'$ are identified along faces, we have an affine
structure on that face, but no affine structure in the directions
``transversal'' to that face. We cannot, however, expect an affine
structure on all of $B$, and we need to choose a discriminant locus.
This still does not give enough data to specify an affine structure on
$B$ away from the discriminant locus. We will also need to choose a
\emph{fan structure} at each vertex of $\P$. We define these data 
below, but first define the discriminant locus.

Let $\Bar(\P)$ be the first barycentric subdivision of $\P$. Then we
define $\Delta'\subseteq B$ to be the union of all simplices  in
$\Bar(\P)$ not containing a vertex of $\P$ or the barycenter of a
maximal cell of $\P$. 

For a vertex $v$ of $\P$, let $W_v$ be the union of the interiors of
all simplices in $\Bar(\P)$ containing $v$. Then $W_v$ is an open
neighbourhood of $v$, and
\[
\{W_v|\hbox{$v$ a vertex of $\P$}\}\cup \{\Int(\sigma)|
\sigma\in\P_{\max}\}
\]
forms an open covering of $B\setminus\Delta'$. To define an affine
structure on $B\setminus\Delta'$, we need to choose affine charts on
$W_v$.

For a vertex $v$ of $\P$, let
\[
\P'_v=\{(v',\sigma')|\hbox{$v'\in\sigma'\in\P'$ a vertex, $\pi(v')=v$}\}.
\]
Let $R_v$ be the quotient of $\coprod_{(v',\sigma')\in\P'_v}\sigma'$
by the following equivalence relation. Let $(v_i',\sigma_i')\in\P_v'$
and $y_i\in\sigma_i'$ for $i=1,2$. Let $\omega_i'\subseteq\sigma_i'$
be the smallest subface of $\sigma_i'$ containing $y_i$. Then
$y_1\sim y_2$ if
\begin{enumerate}
\item
$\pi(y_1)=\pi(y_2)$.
\item
$v_i'\in\omega_i'$.
\item $\pi$ identifies germs of $\omega_1'$ at $v_1'$ and $\omega_2'$
at $v_2'$. By this we mean there exist arbitrarily small open
neighbourhoods $U_i$ of $v_i'$ in $\omega_i'$ such that
$\pi(U_1)=\pi(U_2)$.
\end{enumerate}
(This equivalence relation was not quite correct in \cite{Announce}.)
For example, if $\P'$ consists of the unit square in $\RR^2$, and $B$
is obtained by identifying opposite sides, we have a unique vertex
$v$ in $\P$ and the picture
\begin{center}
\input{squares.pstex_t}
\end{center}
Condition (3) is necessary in this example, as otherwise there would be
further identifications of the interior one-dimensional faces
of $R_v$.

There is a continuous map 
\[
\pi_v:R_v\to B
\]
defined by taking $b\in\sigma'\subseteq R_v$ to $\pi(b)$, and it is
easy to see that if $U_v$ is the connected component of
$\pi_v^{-1}(W_v)$ containing the equivalence class of points
$\{v'\in\sigma'| (v',\sigma')\in\P_v'\}$ in $R_v$, then $U_v\to W_v$ is a
homeomorphism. We write this equivalence class of points as $v'$.

$R_v$ has an abstract polyhedral decomposition $\P_v$, with $v'$ the
unique vertex in $U_v$ mapping to $v$. We will need to find an
embedding $i_v:R_v\to M_{\RR}$. If this is done appropriately, then a
coordinate chart $\psi_v:W_v\to M_{\RR}$ can be defined as
$i_v\circ\pi_v^{-1}|_{W_v}$, and ${\rm exp}_v:i_v(R_v)\to B$ can be
defined as $\pi_v\circ i_v^{-1}$, giving both an affine structure on
$B\setminus\Delta'$ and a proof that $\P$ is a polyhedral
decomposition of $B$.

To do this, we need to choose a \emph{fan structure} at each vertex
$v$ of $\P$. This means for each $v$ we choose a complete rational
polyhedral fan $\Sigma_v$ in $M_{\RR}$ (see \cite{Oda} for a definition)
and a one-to-one inclusion
preserving correspondence between elements of $\P_v$ containing $v'$
and elements of $\Sigma_v$ which we write as $\sigma\mapsto
\sigma_{v'}$.  Furthermore, this correspondence should have the
property that there exist integral affine isomorphisms $i_{\sigma}$
between the tangent wedge of $\sigma$ at $v'$ and $\sigma_{v'}$ which
preserves the correspondence. Such an isomorphism, if it exists, is
unique (it is determined by specifying what it does to primitive
integral generators of each ray, and by integrality, these must be
sent to primitive integral generators). By this uniqueness, the maps
$i_{\sigma}$ glue together to give a map
\[
i_v:R_v\to M_{\RR}
\]
which is a homeomorphism onto its image. Then it is easy to see that
using $\psi_v$ and ${\rm exp}_v$ as defined above one obtains an
integral affine structure on $B_0:=B\setminus\Delta'$ 
and one sees that $\P$ is a
polyhedral decomposition.
\qed
\end{construction}

Suppose we have carried out Construction~\ref{basicconstruct},  and
so obtained an affine structure on $B\setminus\Delta'$ and a
polyhedral decomposition $\P$. It often happens that our choice of
$\Delta'$ is too crude, and we can still extend the affine structure
to a larger open set of $B$. Our original discriminant  locus
$\Delta'$ is a union of all codimension two simplices of $\Bar(\P)$
not containing a vertex of $\P$ or the barycenter of a maximal cell
of $\P$. Let $\tau$ be such a simplex. Locally at a point $b\in
\Int(\tau)$, $B\setminus\Delta'$ takes the form $\RR^{n-2}\times
(\RR^2\setminus \{(0,0)\})$ topologically. We say the holonomy around
the simplex $\tau$ is trivial if the  holonomy of the affine
structure on $B\setminus\Delta'$ about a simple loop around the
origin in $\RR^2\setminus\{(0,0)\}\subset \RR^{n-2}\times
(\RR^2\setminus \{(0,0)\})$ is trivial.

Let $\Delta$ be the union of all codimension two simplices in
$\Delta'$ about which the holonomy is non-trivial.

\begin{proposition}
\label{extaff}
The affine structure on $B\setminus\Delta'$
extends uniquely to an affine structure on $B\setminus\Delta$.
\end{proposition}

\proof 
Let $\omega$ be a simplex of $\Bar(\P)$ contained in $\Delta'$ but
not $\Delta$. Let $\tau$ be the minimal cell of $\P$ containing
$\omega$. Then $0<\dim\tau<\dim B$. Choose a vertex $v$ of $\tau$
such that $\omega\subseteq \overline{W_v}$, and we obtain
$R_v\subseteq\Lambda_{\RR,v}$, $\P_v$ and $\exp_v:R_v\to B$. Let
$\tilde\tau\in\P_v$ with $0\in\tilde\tau$ and
$\exp_v(\tilde\tau)=\tau$. Let  $U_v$ be the connected component of
$\exp_v^{-1}(W_v)$ containing $0$, and let $\tilde \omega\subseteq
\tilde\tau\cap \overline{U_v}$ be such that
$\exp_v(\tilde\omega)=\omega$.  Then $\tilde\tau$ is contained in 
maximal cells $\tilde\sigma_1,\ldots,\tilde\sigma_m\in\P_v$. Let $U$
be a small open neighbourhood of $\Int(\tilde\omega)$ which is
contained in $\bigcup\tilde\sigma_i$ and such that $\exp_v(U)$ is
disjoint from any simplex in $\Delta'$ not containing $\omega$. Then
$U\subseteq \Lambda_{\RR,v}$ inherits an integral affine structure.

Let $U'=U\setminus\exp_v^{-1}(\Delta')$. If we can show
$\exp_v|_{U'}$ is an integral affine isomorphism onto its image in
$B\setminus\Delta'$, then we can patch $U$ and $B\setminus\Delta'$
and so extend the integral affine structure across $\Int(\omega)$.

It is clear from our choice of $U$ that $\exp_v|_{U'}$ is a
homeomorphism onto its image, and by construction it is integral
affine on
\[
U'\cap (U_v\cup \Int(\tilde\sigma_1)\cup\cdots \cup
\Int(\tilde\sigma_m)).
\]
If $y\in U'\setminus(U_v\cup\Int(\tilde\sigma_1)\cup\cdots\cup
\Int(\tilde\sigma_m))$,  we need to show $\exp_v$ is integral affine
at $y$. Let $V$ be a small open neighbourhood of $y$ in $U'$, mapping
$V$ onto an open set $\exp_v(V)$ which we can identify via its affine
structure as an open set in $\Lambda_{\RR,\exp_v(y)}$. Then $\exp_v$
is affine on $V\cap\tilde\sigma_i$, say given by affine
transformations
\[
\psi_i:\Lambda_{\RR,v}
\to \Lambda_{\RR,\exp_v(y)}.
\]
Let $\gamma$ be the image of a loop based at $y$ which passes from
$y$ into $\tilde\sigma_i$ to $0$, then into $\tilde\sigma_j$ and back
to $y$. Then up to conjugation by a translation,
$\rho(\gamma)=\psi_i\circ\psi_j^{-1}$.  Now $\gamma$ can be viewed as
a loop around a codimension 2 simplex of $\Delta'$ containing
$\omega$, so $\rho(\gamma)$ is assumed to be trivial. This says
$\psi_i=\psi_j$. This holds for all $1\le i,j\le m$ and so  $\exp_v$
is affine at $y$ coinciding with $\psi_i$ for any $i$.

Uniqueness follows since for a map between connected open subsets of
$\RR^n$ to be affine can be checked on any connected, dense open subset.
\qed

\begin{definition} 
Let $B$ be an integral affine manifold with singularities and a
polyhedral decomposition $\P$. We define a sheaf $\Lambda_{\P,\RR}$
on $B_0$ as a subsheaf of $\Lambda_{\RR}$ by setting
\[
\Gamma(U,\Lambda_{\P,\RR })=\left\{ v\in\Gamma(U,\Lambda_{\RR
})\bigg| \hbox{$\forall y\in U$, $\sigma\in\P$ with $y\in\sigma$, $v$
is tangent to $\sigma$ at $y$} \right\}.
\]
Similarly, we can define $\Lambda_{\P}$ as above as a subsheaf of
$\Lambda$. Note that if $y\in \Int(\sigma)$, the stalk
$(\Lambda_{\P,\RR})_y$ is the tangent space to $\sigma$ at $y$, and
$\Lambda_{\P,y}$ is a lattice in $(\Lambda_{\P,\RR})_y$.
\qed
\end{definition}

The existence of a polyhedral decomposition of an affine manifold
with singularities places restrictions on the nature of its
singularities, and whether the polyhedral decomposition is toric can
be detected via its holonomy.

\begin{proposition} 
\label{monodromy1}
Let $\P$ be a polyhedral decomposition of $B$. an integral affine
manifold with singularities. For any
$\tau\in\P$, there exists an open neighbourhood $U_{\tau}$ of
$\Int(\tau)$ such that if $y\in \Int(\tau)\setminus\Delta$ and
\[
\rho:\pi_1(U_{\tau}\setminus\Delta,y)\to \Aff(\Lambda_{\RR,y})
\]
is the local holonomy representation, 
\[
\tilde\rho:\pi_1(U_{\tau}\setminus\Delta,y)\to \GL(\Lambda_{\RR,y})
\]
the linear part of the holonomy representation (hence the monodromy
representation for the local system $\Lambda$),
then the radiance obstruction of
$\rho$ is trivial and 
\[
(\tilde\rho(\gamma)-I)((\Lambda_{\P,\RR})_y)=0
\]
for all $\gamma\in\pi_1(U_{\tau}\setminus\Delta,y)$.
\end{proposition}

\proof
Let $v$ be a vertex of $\tau$, $R_v,\P_v$ and $\exp_v$ as usual, and
let $\tilde\tau\in\P_v$ with $0\in\tilde\tau$ and
$\exp_v(\tilde\tau)=\tau$. Let $\tilde U_{\tau}$ be a sufficiently
small open neighbourhood of $\Int(\tilde\tau)$ such that $\exp_v$
maps  $\tilde U_{\tau}$ homeomorphically to an open neighbourhood
$U_{\tau}$ of $\Int(\tau)$.  Let $\tilde y\in R_v$ map to $y\in
\Int(\tau)\setminus\Delta$.  Consider a loop in
$U_{\tau}\setminus\Delta$ based at $y$, which we can lift to a loop
$\tilde\gamma$ in $\tilde U_{\tau}$ based at $\tilde y$. We can always
assume, by deforming $\gamma$, that $\tilde\gamma$ intersects
non-maximal cells of $\P_v$ only in a finite number of points.
Indexing these points by $\ZZ/m\ZZ$, we write them as $\tilde y=\tilde
z_0,\tilde z_1,\ldots,\tilde z_{m-1}$ in order through which
$\tilde\gamma$ passes through them, so that we arrive at $\tilde
z_m=\tilde z_0$ when we have finished traversing $\tilde\gamma$. Let
$\tilde\sigma_i\in\P_v$ be the maximal cell that $\tilde\gamma$ passes
through between $\tilde z_i$ and $\tilde z_{i+1}$. Finally, let
$z_i=\exp_v(\tilde z_i)$. Identifying a small open neighbourhood $V_i$
of $z_i$ with an open neighbourhood of the origin in
$\Lambda_{\RR,z_i}$, the map $\exp_v$ is given by affine maps 
\begin{eqnarray*}
\psi_i^{i-1}:\Lambda_{\RR,v}&\to&\Lambda_{\RR,z_i}\\
\psi_i^i:\Lambda_{\RR,v}&\to&\Lambda_{\RR,z_i}
\end{eqnarray*}
on $\exp_v^{-1}(V_i)\cap\tilde\sigma_{i-1}$ and $\exp_v^{-1}(V_i)\cap
\tilde\sigma_i$ respectively. If $\tilde z_i\in \Int(\tilde\omega_i)$
for some cell $\tilde\omega_i\in\P_v$, then $\psi_i^{i-1}$ and
$\psi_i^i$ must agree on the linear space $\RR\tilde\omega_i$ spanned
by $\tilde\omega_i$, and since we can assume $U_{\tau}$ is
sufficiently small, we can assume $\RR\tilde\omega_i$ contains
$\RR\tilde\tau$. Thus $(\psi_i^i)^{-1}\circ \psi_i^{i-1}$ is the
identity on $\RR\tilde\tau$, while $\psi_0^{m-1}$ and $\psi_0^0$
identify $\RR\tilde\tau$ with $(\Lambda_{\P,\RR})_{y}$. Now we see
that 
\[
\rho(\gamma)^{-1}=\psi_{0}^{m-1}\circ (\psi_{m-1}^{m-1})^{-1}
\circ\psi_{m-1}^{m-2} \circ\cdots\circ(\psi_1^1)^{-1}
\circ\psi_1^0\circ(\psi_0^0)^{-1},
\]
which is the identity on $(\Lambda_{\P,\RR})_y$. In particular,
$0\in(\Lambda_{\P,\RR})_y$ is fixed by $\rho(\gamma)$ so the
translational part of $\rho(\gamma)$ is zero and the radiance
obstruction vanishes. \qed

\begin{remark} This local triviality of the radiance obstruction
can be viewed as saying that the radiance obstruction $c_{B_0}\in
H^1(B_0,\Lambda)$ is in the image of the inclusion
$H^1(B,i_*\Lambda)\hookrightarrow H^1(B_0,\Lambda)$ induced by the
Leray spectral sequence for $i:B_0\hookrightarrow B$. Indeed, local
triviality says that the radiance obstruction is zero in 
$H^0(B,R^1i_*\Lambda)$. 
\end{remark}

\begin{definition}
Let $B$ be an integral affine manifold with singularities with a
polyhedral decomposition $\P$. If $\sigma\in\P$, then the subspaces
$(\Lambda_{\P,\RR})_y$ of $\Lambda_{\RR,y}$ for all $y\in
\Int(\sigma)\setminus\Delta$ are canonically identified via parallel
transport in a neighbourhood of $\Int(\sigma)$ by Proposition
\ref{monodromy1}. For any $y\in\Int(\sigma)\setminus \Delta$ we denote
this subspace of $\Lambda_{\RR,y}$ by $\Lambda_{\sigma,\RR}$.
Similarly, we denote by $\Lambda_{\sigma}$ the corresponding subspace
of $\Lambda_y$ for $y\in \Int(\sigma)\setminus\Delta$.
\end{definition}

\begin{proposition}
\label{monodromy2}
Let $\P$ be a polyhedral decomposition of $B$. Then $\P$ is toric if
and only if for each $\tau\in\P$, there exists an open neighbourhood
$U_{\tau}$ of $\Int(\tau)$ such that if $y\in
\Int(\tau)\setminus\Delta$ with
\[
\rho:\pi_1(U_{\tau}\setminus\Delta,y)\to \Aff(\Lambda_{\RR,y}),
\]
the local holonomy representation, then
$(\tilde\rho(\gamma)-I)(\Lambda_{\RR,y})\subseteq
\Lambda_{\tau,\RR}$.
\end{proposition}

\proof First suppose $(\tilde\rho(\gamma)-I)(\Lambda_{\RR,y})
\subseteq \Lambda_{\tau,\RR}$ in $U_{\tau}$,  and we will show,
possibly after shrinking $U_{\tau}$, that there is an affine
submersion $S_{\tau}:U_{\tau}\to\Lambda_{\RR,y}/ \Lambda_{\tau,\RR}$
with $S_{\tau}(\Int(\tau))=0$, thus showing $\P$ is toric. Let
$v\in\tau$ be a vertex, $R_v$, $\P_v$ as usual, with $\tilde\tau
\in\P_v$ with $0\in\tilde\tau$ and $\exp_v(\tilde\tau)=\tau$.
Possibly by shrinking $U_{\tau}$, we can assume there is an open
neighbourhood of $\Int(\tilde\tau)$, $\tilde U_{\tau}\subseteq R_v$,
with $\exp_v$ a homeomorphism of $\tilde U_{\tau}$ onto $U_{\tau}$.
If we choose $v$ to be in the same connected component of
$\tau\setminus\Delta$ as $y$, we can identify $\Lambda_{\RR,y}$ and
$\Lambda_{\RR,v}$ by parallel transport along a path in $\tau$, and
under this identification, the linear space $\RR\tilde\tau$ spanned
by $\tilde\tau$ is identified with $\Lambda_{\tau,\RR}$.

We can then define 
\[
S_{\tau}=p_{\tau}\circ \exp_v^{-1}
\]
where $p_{\tau}:\Lambda_{\RR,v}\to\Lambda_{\RR,v}/\RR\tilde\tau$ is
the projection. This is a continuous map, and clearly
$S_{\tau}(\Int(\tau)) =0$, so we just need to show $S_{\tau}$ is an
affine submersion.

To do so, it suffices to show it is an affine submersion at any
point $y'\in \Int(\tau)\setminus\Delta$. Identify a small
neighbourhood $V$ of such a point with a neighbourhood in
$\Lambda_{\RR,y'}$. If $\tilde\sigma_1,\ldots,\tilde\sigma_r$ are the
maximal cells of $\P_v$ containing $\tilde\tau$, then $\exp_v$ is
given by an affine map $\psi_i:\Lambda_{\RR,v}\to \Lambda_{\RR,y'}$
on $\tilde\sigma_i \cap V$. On the other hand, holonomy around a
simple loop based at $y$ into $\tilde\sigma_i$ through $y'$ into
$\tilde\sigma_j$ back to $y$ is given by $\psi_i^{-1}\circ\psi_j$, so
by assumption on the local holonomy, it is clear that
$p_{\tau}\circ\psi_i^{-1}=p_{\tau}\circ\psi_j^{-1}$. Thus $S_{\tau}$
is affine at $y'$, as $S_{\tau}$ coincides with the affine submersion
$p_{\tau}\circ\psi_i^{-1}$ for any $i$.

Conversely, suppose we have an affine submersion $S_{\tau}:U_{\tau}
\to M'_{\RR}$ as in Definition~\ref{PD2}, (5). Then for  $y\in
\Int(\tau)\setminus\Delta$, $S_{\tau*}:\Lambda_{\RR,y} \to M_{\RR}'$
clearly yields a factorization of $S_{\tau*}$ as $\Lambda_{\RR,y}\to
\Lambda_{\RR,y}/\Lambda_{\tau,\RR} \mapright{\cong} M_{\RR}'$. Thus
on $U_{\tau}\setminus\Delta$, the sheaf kernel of
$S_{\tau*}:\Lambda_{\RR}\to M_{\RR}'$ is just obtained from
$\Lambda_{\tau,\RR}$ by parallel transport. Thus for any
$\gamma\in\pi_1(U_{\tau}\setminus\Delta,y)$, as
$\tilde\rho(\gamma)(\Lambda_{\tau,\RR})=\Lambda_{\tau,\RR}$, the 
transformation $\tilde\rho(\gamma)$ induces the identity on 
$\Lambda_{\RR,y}/\Lambda_{\tau,\RR}\cong M'_{\RR}$. Thus
$(\tilde\rho(\gamma)-I)(\Lambda_{\RR,y})\subseteq \Lambda_{\tau,\RR}$
as desired.
\qed
\medskip

There will also be information in the normal directions of a cell
$\tau\in\P$, for example the fan $\Sigma_\tau$ introduced in
Definition~\ref{fan2} below. The normal spaces of the cells also glue
into a sheaf, but now special care has to be taken with the
discriminant locus because of monodromy. The following definition
turns out to work best.

\begin{definition}
Let $B$ be an integral affine manifold with singularities and a toric
polyhedral decomposition $\P$. We define a subsheaf $\shQ_{\P,\RR}$ of
$i_*(\Lambda_{\RR}/\Lambda_{\P,\RR})$ as follows. A section
$s\in\Gamma(U,i_*(\Lambda_{\RR}/\Lambda_{\P,\RR}))$ is a section of
$\shQ_{\P,\RR}$ if there is an open cover $\{U_i\}$ of $U$ such that
for any $i$ and for any path $\gamma:[0,1]\to U_i\setminus \Delta$,
$\gamma^*s\in\Gamma([0,1],\gamma^*(\Lambda_{\RR}/\Lambda_{\P,\RR}))$
is induced by a section $t\in\Gamma([0,1],\gamma^*\Lambda_{\RR})$ (so
local sections are constant under parallel transport). We similarly
define a subsheaf $\shQ_{\P}$ of $i_*(\Lambda/\Lambda_{\P})$.
\end{definition}

We note that if $y\in \Int(\sigma)$ for some $\sigma\in\P$, and $U$
is a connected open neighbourhood of $y$ disjoint from any
$\tau\in\P$ with $\sigma\not\subseteq\tau$, then 
\[
\Gamma(U,\shQ_{\P,\RR})=\Lambda_{\RR,y'}/\Lambda_{\sigma,\RR}
\]
and
\[
\Gamma(U,\shQ_{\P})=\Lambda_{y'}/\Lambda_{\sigma}
\]
where $y'\in\sigma\cap (U\setminus \Delta)$. By
Proposition~\ref{monodromy2}, this is independent of the choice of
$y'$. We will write
\begin{eqnarray*}
\shQ_{\sigma,\RR}&=&(\shQ_{\P,\RR})_y\\
\shQ_{\sigma}&=&(\shQ_{\P})_y
\end{eqnarray*}
for any $y\in \Int(\sigma)$, being independent of the choice of $y$.

\begin{example}
Suppose $\dim B=2$, $\sigma_1,\sigma_2\in\P_{\max}$ with
$\sigma_1\cap \sigma_2=\tau$ a cell of dimension one, and suppose
$\Delta\cap\tau=\{p\}$. Thus $\Int(\tau)\setminus\Delta$ has two
connected components. Then $\Lambda/\Lambda_{\P}$ is the constant
sheaf $\ZZ$ on $\Int(\tau)\setminus\Delta$, and
$i_*(\Lambda/\Lambda_{\P})$ then has stalk $\ZZ\oplus\ZZ$ at $p$,
which is bigger than we want. On the other hand, if $\Lambda$ has
non-trivial monodromy around $p$, then
$\Gamma(U\setminus\{p\},\Lambda)=\Gamma(U\setminus\{p\},\Lambda_{\P})$
for $U$ any sufficiently small neighbourhood of $p$. Thus
$i_*\Lambda/ i_*\Lambda_{\P}$ has stalk $0$ at $p$. The above
definition is designed to give an intermediate value for the stalk,
namely $\ZZ$.
\end{example}

\begin{definition}
\label{fan2}
\emph{The fan $\Sigma_{\tau}$.} Suppose $B$ is integral affine with
singularities and $\P$ is a toric polyhedral decomposition. Let $\tau
\in\P$. We will construct a rational polyhedral fan in
$\shQ_{\tau,\RR}$ (with integral structure $\shQ_{\tau}$) associated
to $\tau$.

Let $S_{\tau}:U_{\tau}\to M'_{\RR}$ be as in Definition \ref{PD2},
(5). For any point $y\in \Int(\tau)\setminus\Delta$, $S_{\tau,*}$
identifies $\Lambda_{\RR,y}/\Lambda_{\tau,\RR}= \shQ_{\tau,\RR}$ with
$M'_{\RR}$, so we can take $M'_{\RR}=\shQ_{\tau,\RR}$. The idea now
is that if  we take $U_{\tau}$ small enough, then inside $U_{\tau}$
the decomposition $\P$ looks like the pull back of a fan in
$M'_{\RR}$ under the map $S_{\tau}$. Indeed, let $\P_{\tau}$ denote
the collection of closures of connected components of the sets
$\Int(\sigma)\cap U_{\tau}$, where $\sigma$ ranges over all elements
of $\P$ containing $\tau$.  For any $\sigma\in\P_{\tau}$, set
\[
K_\sigma=
\{aS_{\tau}(x)|\hbox{$a\in\RR_{\ge0}$ and $x\in\sigma\cap U_{\tau}$}\}
=\RR_{\ge 0}\cdot S_\tau(\sigma).
\]
Because $S_{\tau}$ contracts $\Int(\tau)$ to a point, this can be
viewed as the tangent wedge to $\sigma$ at a point $y\in \Int(\tau)$
modulo the linear space $\Lambda_{\tau,\RR}$.  We set
\[
\Sigma_{\tau}=\{K_\sigma|\sigma\in\P_{\tau}\}.
\]
It is clear that $\Sigma_{\tau}$ is a complete rational polyhedral
fan in $\shQ_{\tau,\RR}$. Note that if $\tau$ is a vertex $v$, then
this coincides with $\Sigma_v$ of Remark~\ref{fan1}.
\qed
\end{definition}

Note we do not necessarily have a one-to-one correspondence between
cells of $\P$ containing $\tau$ and cones of $\Sigma_{\tau}$. This is
because we have allowed the possibility that cells can be
self-intersecting, as in Example~\ref{decompexamp}, (1), where a
fundamental domain for the action of $\Gamma$ on $M_{\RR}$ is used to
induce a polyhedral decomposition of $M_{\RR}/\Gamma$. It is because
of this possibility that we use the following notion to generalize
the concept of inclusion of cells:

\begin{definition}
Let $\Cat(\P)$ be the category whose  set of objects is the set $\P$,
and morphisms are defined as follows: if $\tau\not\subseteq\sigma$
then $\Hom(\tau,\sigma)=\emptyset$; if $\tau=\sigma$ then
$\Hom(\tau,\sigma)=\{id\}$, and if $\tau\subset\sigma$, then
$\Hom(\tau,\sigma)$ is the set of 1-simplices in the triangulation 
$\Bar(\P)$ with endpoints $\Bar(\tau)$ and $\Bar(\sigma)$. If $e_1$
is an edge joining $\Bar(\tau)$ and $\Bar(\sigma)$ and $e_2$ an edge
joining $\Bar(\sigma)$ and  $\Bar(\omega)$, then $e_2\circ e_1$ is
defined to be the third edge of the unique $2$-simplex containing
$e_1$ and $e_2$.
\end{definition}

The key point here is that if the neighbourhood $U_{\tau}$ of
$\Int(\tau)$ is taken to be sufficiently small, then there is a
natural one-to-one correspondence between the set of connected
components of $\Int(\sigma)\cap U_{\tau}$ and the set
$\Hom(\tau,\sigma)$, with each component of $\Int(\sigma)\cap
U_{\tau}$ intersecting precisely one element of $\Hom(\tau,\sigma)$. 
Thus there is one cone in the fan $\Sigma_{\tau}$ coming from
$\sigma$ for each element of $\Hom(\tau,\sigma)$, and we obtain a
one-to-one correspondence between cones of $\Sigma_{\tau}$ and the
set
\[
\coprod_{\sigma\in\P} \Hom(\tau,\sigma).
\]

Now we can obtain the relationship between $\Sigma_{\tau}$ and
$\Sigma_{\sigma}$ when $\tau\subset\sigma$; this depends on an
$e\in\Hom(\tau,\sigma)$ corresponding to a connected component $T_e$
of $\Int(\sigma)\cap U_{\tau}$. Choose a path $\gamma$ in $U_{\tau}$
from some $x\in \Int(\tau)\setminus\Delta$ to some $y\in
T_e\setminus\Delta$. By parallel transport this defines an
isomorphism $\Lambda_{\RR,x}\to \Lambda_{\RR,y}$, and hence a
surjection $p_e:\shQ_{\tau,\RR}\to \shQ_{\sigma,\RR}$, which is
defined independently of the path $\gamma$ by
Proposition~\ref{monodromy2}.  If $K_e$ denotes the cone of
$\Sigma_{\tau}$ corresponding to the connected component $T_e$, then
the map $p_e$ is easily seen to identify  $\Sigma_{\sigma}$ with the
quotient fan $\Sigma_{\tau}(K_e)$, which is defined as follows. This
identification is completely canonical, depending only on $e$.

\begin{definition}
\label{quotfan}
Let $\Sigma$ be a fan on a vector space $M_{\RR }$, and let
$\tau\in\Sigma$ be a cone. We denote by $\Sigma(\tau)$ the fan in
$M_{\RR }/{\RR }\tau$ consisting of the cones
\[
\{(\sigma+{\RR }\tau)/{\RR }\tau|\sigma\supseteq\tau, \sigma\in\Sigma\}.
\]
We say $\Sigma(\tau)$ is the \emph{quotient fan} of $\Sigma$ by the
cone $\tau$.

We also denote by $\tau^{-1}\Sigma$ the fan of (not strictly convex!)
cones
\[
\{\sigma+{\RR }\tau|\sigma\supseteq\tau, \sigma\in\Sigma\}.
\]
We call this the \emph{localisation} of the fan $\Sigma$ at $\tau$.
\qed
\end{definition}

\begin{definition}
\label{normalfan}
\emph{The normal fan $\check\Sigma_{\tau}$.}  Let $\tau\in\P$, and
let $y\in \Int(\tau)\setminus\Delta$ be any point. Let $v$ be a
vertex of $\tau$ contained in the same connected component of
$\tau\setminus \Delta$ as $y$, and let $R_v$, $\P_v$, and $\exp_v$ be
as usual, and $\tilde\tau \in\P_v$ with $0\in\tilde\tau$ and
$\exp_v(\tilde\tau)=\tau$. Then via parallel translation along
$\tau$, $\Lambda_{\RR,v}$ can be canonically identified with
$\Lambda_{\RR,y}$ so that $\RR\tilde\tau$ and $\Lambda_{\tau,\RR}$
are identified. We can then view $\tilde\tau\subseteq
\Lambda_{\tau,\RR}$ as a polytope under this identification, with
$\dim\tilde\tau=\dim\Lambda_{\tau,\RR}$. In
particular, $\tilde\tau$ is well-defined up to translation in
$\Lambda_{\tau,\RR}$ independently of the choice of $y$. This
associates to each cell a canonical (up to translation) polytope.

We define $\check\Sigma_{\tau}$ to be the normal fan in
$\dualvs{\Lambda_{\tau,\RR}}$ to $\tilde\tau$ in the usual way,
as follows. For each face  $\tilde\omega\subseteq\tilde\tau$, let 
\[
\check K_\omega=\{f\in\dualvs{\Lambda_{\tau,\RR}}|\hbox{
$f|_{\tilde\omega}$  is constant and $\langle f,z\rangle\ge \langle
f,x\rangle$  for $z\in\tilde\tau$, $x\in\tilde\omega$}\}.
\]
Then the collection of cones $\{\check K_{\omega}\}$ form the normal
fan $\check\Sigma_{\tau}$ in $\dualvs{\Lambda_{\tau,\RR}}$,  which is
completely well-defined. It is then easy to see there is a one-to-one
correspondence between cones of $\check\Sigma_{\tau}$ and the set
\[
\coprod_{\omega\in\P}\Hom(\omega,\tau).
\]
\qed
\end{definition}

\begin{definition} 
Let $\P$ be a polyhedral decomposition of $B$. Let
$\shAff_{\RR}(B,\RR)$ denote the sheaf of continuous functions on $B$
to $\RR$ which are affine when restricted to $B_0$. We define $\shAff
(B,\ZZ)$ to be the subsheaf of $\shAff_{\RR}(B,\RR)$ of functions which
when restricted to $B_0$ are in $\shAff(B_0,\ZZ)$
(Definition~\ref{fundseq1}).
\end{definition}

Global affine functions are either unbounded or constant.

\begin{proposition}
\label{no global affine functions}
Let $\P$ be a polyhedral decomposition of $B$ and $\varphi:B\to \RR$ a
bounded affine function. Then $\varphi$ is constant.
\end{proposition}
\proof
Assume $\varphi$ is non-constant. Then its linear part is non-zero at
every $x\in B_0$. We construct a piecewise integral affine map
$\gamma:\RR_{\ge 0}\to B$ with intervals of linearity mapping to edges
of $\P$ and such that $\varphi\circ\gamma$ is strictly monotone.
Then $\varphi\circ \gamma$ is unbounded because $\gamma$ traces
through infinitely many edges and on each edge $\varphi\circ\gamma$
increases at least by $1$.

First note that for each edge $\tau\in\P$ the restriction
$\varphi|_{\Int(\tau)}$ is affine, regardless of $\tau$ intersecting
$\Delta$ or not. In fact, there exists a top-dimensional cell
$\sigma\in\P$ containing $\tau$ and $\varphi|_{\Int(\sigma)}$ is
affine linear. By continuity of $\varphi$ the pull-back to some
polytope $\tilde\sigma$ covering $\sigma$ is then also affine, and
there is an edge $\tilde\tau\subset\tilde\sigma$ mapping to $\tau$.

Now we construct $\gamma$ inductively; in the $n$-th step $\gamma_n$
passes through $n$ edges. For $n=0$ just choose a vertex $v\in B$ and
let $\gamma_0:\{0\}\to B$ have image $v$. To construct $\gamma_{n+1}$
let $\tau$ be the last edge passed through by $\gamma_n$. Thus if
$\gamma_n$ is defined on $[0,a]$ and the integral length of $\tau$
is $l$ then $\gamma_n|_{[a-l,a]}$ is given by
$[a-l,a]\simeq\tilde\tau\to B$. The identification of $[a-l,a]$ with
$\tilde\tau$ is done in such a way that $a-l$ maps to
$\gamma_n(a-l)=\gamma_{n-1}(a-l)$. At the endpoint $w=\gamma_n(a)$ the
level set of $\varphi$ locally separates $B$ into two connected
components with $\gamma_n([a-\epsilon,a))$ being contained in only one
of them for $\epsilon$ small. Now because $\P$ defines a fan structure
at $w$ there exists an edge $\tau'\in\P$ emanating from $w$ into the
other connected component. Then $\varphi|_{\tau'}$ is strictly
increasing when starting from $w$. Hence we can attach an interval
equal to the integral length of $\tau'$ to the domain of definition
of $\gamma_n$ and map it in an integral affine way onto $\tau'$
to obtain $\gamma_{n+1}$.
\qed

\begin{proposition}
\label{extendafffun}
If $i:B_0\hookrightarrow B$ denotes the inclusion and $B$ has a
polyhedral decomposition $\P$, then $\shAff_{\RR }(B,{\RR })
=i_*\shAff_{\RR }(B_0,{\RR })$ and $\shAff (B,\ZZ)=i_*\shAff(B_0,\ZZ)$.
\end{proposition}

\proof Let $U\subseteq B$ be open, $U_0=U\setminus\Delta$. Then 
there is certainly an inclusion given by restriction $\shAff_{\RR }(U,{\RR })
\to \shAff_{\RR }(U_0,{\RR })$, because $U_0$ is dense in $U$ and thus
a function on $U_0$ has at most one continuous extension to $U$. Conversely,
given an affine function on $U_0$, to extend it to $U$ we need only
to find extensions on an open covering of $U$; by uniqueness these
extensions will glue. Thus we can assume $U$ to be as small as we like.
In particular, we can assume $U$ is a small neighbourhood of
a point $y\in \Int(\tau)\cap\Delta$, $\tau\in\P$. Let $v$ be a vertex
of $\tau$, $R_v$, $\P_v$, $\exp_v$ as usual, $\tilde\tau\in\P_v$
with $0\in\tilde\tau$ and $\exp_v(\tilde\tau)=\tau$. Let $\tilde y
\in \Int(\tilde\tau)$ satisfy $\exp_v(\tilde y)=y$, and let
$\tilde U$ be a lift of $U$ to an open neighbourhood of $\tilde y$.
If $\tilde\sigma_1,\ldots,\tilde\sigma_m$ are the maximal cells
of $\P_v$ containing $\tilde\tau$, then 
\[
f_i:=f\circ\exp_v|_{(\tilde U\cap\tilde\sigma_i)
\setminus\exp_v^{-1}(\Delta)}
\]
extends to an affine map $\psi_i:\tilde U\cap\tilde\sigma_i\subseteq
\Lambda_{\RR,v}\to\RR$. Furthermore, because $f_i$ and $f_j$ agree on
$(\tilde U\cap \tilde \sigma_i\cap\tilde\sigma_j)
\setminus\exp_v^{-1}(\Delta)$, by continuity, $\psi_i$ and $\psi_j$
agree on $\tilde\sigma_i\cap \tilde\sigma_j$. Thus the $\psi_i$'s
glue to give a function $\psi: \tilde U\to\RR$. It follows that $f$ can be
extended across $\Delta$ as $\psi\circ \exp_v^{-1}|_U$. The same
argument works in the integral case.
\qed

\begin{proposition}
\label{affexact}
If $B$ has a polyhedral decomposition, then there are exact sequences
\[
\exact{\RR}{\shAff_{\RR}(B,\RR)}{i_*\check\Lambda_{\RR}}
\]
and
\[
\exact{\ZZ}{\shAff(B,\ZZ)}{i_*\check\Lambda}.
\]
\end{proposition}

\proof
We have from Definition~\ref{fundseq1} an exact sequence
\[
\exact{{\RR }}{\shAff_{\RR }(B_0,{\RR })}{\check\Lambda_{\RR }}
\]
on $B_0$,
hence an exact sequence
\[
0\to {\RR }\to
{\shAff_{\RR }(B,{\RR })}\to {i_*\check\Lambda_{\RR }}
\]
by Proposition~\ref{extendafffun}, so we just need to show
surjectivity on the right. By Proposition~\ref{extclass} and Proposition~\ref{monodromy1}, 
there is an open cover $\{U_{\tau}|\tau\in\P\}$
of $B$ such that the exact sequence for $\shAff(B_0,\RR)$ splits on
$U_{\tau}\cap B_0$ for each $\tau\in\P$. But then on $U_{\tau}$,
$\shAff(B,\RR)=\RR\oplus i_*\check\Lambda_{\RR}$, so we have surjectivity.
\qed\medskip

Next we wish to define piecewise linear functions. One natural
definition would be to define a piecewise linear function on $B$ with
respect to a polyhedral decomposition $\P$ as a continuous function
on $B$ which is affine on the interior of every maximal simplex.
However, this only depends on the affine structure of the polytopes
and ignores singularities, and this proposed definition can behave a
bit perversely near singularities.  Instead, we refine this by

\begin{definition} 
\label{PLfunc}
Let $\P$ be a polyhedral decomposition of $B$. If $U\subseteq B$ is an
open set, then a \emph{piecewise linear} function on $U$ is a
continuous function $\varphi:U\to \RR$ which is affine on $U \cap
\Int(\sigma)$ for each $\sigma\in\P_{\max}$, and satisfies the
following property: for any $y\in U$, $y\in \Int(\sigma)$ for some
$\sigma\in\P$, there exists a neighbourhood $V$ of $y$ and a $f\in
\Gamma(V,\shAff_{\RR}(B,\RR))$ such that $\varphi-f$ is zero on $V\cap
\Int(\sigma)$. We denote by $\shPL_{\P,\RR}(B,\RR)$ the sheaf on $B$ of
such functions. Similarly, we define the subsheaf of integral
piecewise linear functions $\shPL_{\P}(B,\ZZ) \subseteq
\shPL_{\P,\RR}(B,\RR)$ so that $\varphi\in \shPL_{\P}(B,\RR)$ if
$\varphi|_{\Int(\sigma)}\in \shAff(\Int(\sigma),\ZZ)$ for each
$\sigma\in\P_{\max}$, and $f$ can be taken to be integral in a
neighbourhood of each $y$.
\end{definition}

The restriction on the behaviour near the singularities 
of the affine structure makes the following definition possible.

\begin{definition}
\label{PLfan}
Let $\P$ be a toric polyhedral decomposition of $B$ and $\varphi$ be
a  piecewise linear function on $B$. For each $\sigma\in\P$, we
obtain by Definition~\ref{fan2} a fan $\Sigma_{\sigma}$ in
$\shQ_{\sigma,\RR}$. Choose any point  $y\in \Int(\sigma)$, and a
small open neighbourhood $U$ of $y$. By Definition~\ref{PLfunc} there
exists an affine function $f:U\to\RR$ such that $f|_{U\cap
\Int(\sigma)}=\varphi|_{U\cap \Int(\sigma)}$. Then $\varphi-f$ is the
pullback by $S_{\sigma}:U\to \shQ_{\sigma,\RR}$ of a function on
$\shQ_{\sigma,\RR}$ which is piecewise linear with respect to the fan
$\Sigma_{\sigma}$. We write this function as $\varphi_{\sigma}$,
well-defined up to a linear function.
\end{definition}

\begin{definition}
If $\P$ is a polyhedral decomposition of $B$, then
there is a commutative diagram of exact sequences
\[
\begin{matrix}&&0&&0&&&&\cr
&&\mapdown{}&&\mapdown{}&&&&\cr
&&\ZZ&\mapright{=}&\ZZ&&&&\cr
&&\mapdown{}&&\mapdown{}&&&&\cr
0&\mapright{}&\shAff(B,\ZZ)&\mapright{}&\shPL_{\P}
(B,\ZZ)&\mapright{}
&\shMPL_{\P}&\mapright{}&0\cr
&&\mapdown{}&&\mapdown{}&&\mapdown{=}&&\cr
0&\mapright{}&i_*\check\Lambda&\mapright{}&\shPL_{\P}
(B,\ZZ)/\ZZ&\mapright{}
&\shMPL_{\P}&\mapright{}&0\cr
&&\mapdown{}&&\mapdown{}&&&&\cr
&&0&&0&&&&\cr\end{matrix}
\]
defining a sheaf $\shMPL_{\P}$. A similar diagram with $\RR$ subscripts
defines a sheaf $\shMPL_{\P,\RR}$.
\qed
\end{definition}

\begin{definition}
A multi-valued piecewise linear function on an open subset
$U\subseteq B$ with respect to $\P$ is an element $\varphi\in
\Gamma(U,\shMPL_{\P,\RR})$; it is integral if it is in
$\Gamma(U,\shMPL_{\P})$. Such a section should be thought of as a
collection of piecewise linear functions on an open covering $U_i$ of
$U$ which differ only by affine linear functions on overlaps. The
\emph{first Chern class} of $\varphi\in H^0(B,\shMPL_{\P,\RR})$ is the
image of $\varphi$ under the boundary map 
\[
c_1:H^0(B,\shMPL_{\P,\RR})\to H^1(B,i_*\check\Lambda_{\RR})
\]
(or
\[
c_1:H^0(B,\shMPL_{\P})\to H^1(B,i_*\check\Lambda)
\]
in the integral case). Using Definition~\ref{PLfan}, a multi-valued
piecewise linear function on $B$ determines a piecewise linear function
$\varphi_{\sigma}:\shQ_{\sigma,\RR}\to\RR$ with respect to
$\Sigma_{\sigma}$, well-defined up to a linear function.
\qed
\end{definition}

\begin{definition} 
A multi-valued piecewise linear function $\varphi$ on $B$ with respect
to $\P$ is \emph{(strictly) convex} if $\varphi_{\tau}$ is a
(strictly) convex piecewise linear function on the fan $\Sigma_{\tau}$
for all $\tau\in\P$. More precisely, if $\sigma$ is a maximal cone of
$\Sigma_{\tau}$ and $n_{\sigma}\in  \dualvs{\shQ_{\tau,\RR}}$
satisfies $\langle n_{\sigma},y\rangle=\varphi_{\tau}(y)$ for all
$y\in\sigma$, then $\langle n_{\sigma},y\rangle\le\varphi_{\tau}(y)$
($\langle n_{\sigma},y\rangle < \varphi_{\tau}(y)$) for all
$y\not\in\sigma$. Note this is independent of the choice of
$\varphi_{\tau}$.
\end{definition}

While we will not use it in this paper, the following concept
is a natural generalisation of regular triangulations:

\begin{definition}
A polyhedral decomposition $\P$ is \emph{regular} if there is a
strictly convex multi-valued piecewise linear function on $B$ with
respect to $\P$.
\end{definition}

One of the key points of this paper is that strictly convex piecewise
linear functions can play the same role as the affine K\"ahler
potentials discussed in \S1: in other words, they allow us to dualize
the affine structure using a discrete Legendre transform which we will
explain in the next section.

\begin{remark}
\label{straight}
\emph{Straightening the discriminant locus.} 
When we construct an affine manifold with singularities
using Construction~\ref{basicconstruct}, the discriminant
locus has been chosen to be a union of simplices in the barycentric
subdivision of $\P$. The intersection and dual intersection complexes
we construct in \S 4 are special cases of Construction~\ref{basicconstruct}. 
However, in real life, one
might encounter affine manifolds with singularities arising from
other sources, such as elliptically fibred K3 surfaces
(Example~\ref{ellipticK3}). Often these examples carry polyhedral
decompositions, but there is no reason to expect the discriminant locus
to be a union of barycentric simplices. It is convenient to
note these affine manifolds can then be deformed into examples
of the sort constructed in Construction~\ref{basicconstruct}.

Thus if $B,\P$ does not arise through Construction~\ref{basicconstruct},
it will often make sense to replace $B$ and
$\P$ with a closely related $B_{str}$ and $\P_{str}$ as follows. Each
maximal cell  $\sigma$ of $\P$ determines a lattice polytope
$\tilde\sigma \subseteq \Lambda_{\RR,y}$ for any $y\in \Int(\sigma)$
using Definition~\ref{normalfan}. Furthermore, $B$ is clearly
obtained as a topological space via integral affine identifications
between various faces of the polytopes of
\[
\P'=\{\tilde\sigma|\sigma\in\P_{\max}\}.
\]
Also, we have a fan structure $\Sigma_v$ at each vertex of
$\P$, and by Construction~\ref{basicconstruct} and
Proposition~\ref{extaff}, we obtain a new affine manifold with
singularities structure on $B$, which we call $B_{str}$, with
discriminant locus $\Delta_{str}$ contained in $\Delta'_{str}$, the
union of all codimension two simplices of $\Bar(\P)$ not containing
vertices or intersecting interiors of maximal cells of $\P$. Then as
a manifold $B_{str}=B$, and $\P_{str}=\P$, but the structure of
affine manifolds with singularities is different because
$\Delta_{str}\not=\Delta$. In the following figure, we see the effects of
straightening. The dotted lines show the first barycentric
subdivision.
\begin{center}
\includegraphics{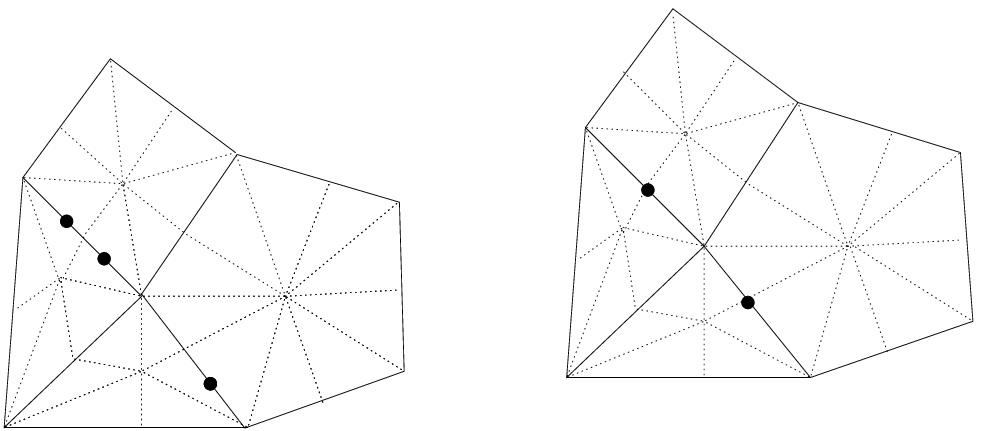}\\[-3ex]
\end{center}

All of our examples of affine manifolds with singularities will
essentially arise through Construction~\ref{basicconstruct}, so this
straightening will not be necessary in this paper.
\qed
\end{remark}

\subsection{The discrete Legendre transform}

We now introduce a key concept for relating an algebro-geometric form
of mirror symmetry to the Strominger-Yau-Zaslow approach. Let $B$ be
an integral affine manifold with singularities with a toric
polyhedral decomposition $\P$, and assume the discriminant locus of
$B$ has been straightened, or equivalently, $(B,\P)$ arises
from Construction~\ref{basicconstruct}. Let $\varphi$ be a  strictly convex
multi-valued integral piecewise linear  function on $B$. We will
construct a new integral affine manifold with singularities $\check
B$ with discriminant locus $\check\Delta$. As manifolds, $B=\check B$
and $\Delta=\check\Delta$, but the affine structure is dual. In
addition, we obtain a toric polyhedral decomposition $\check\P$ and a
strictly convex multi-valued integral piecewise linear function
$\check\varphi$ on $\check B$. We will say $(\check B,\check\P,
\check\varphi)$ is the \emph{discrete Legendre transform} of
$(B,\P,\varphi)$.

First we define $\check\P$. For any $\sigma\in\P$, define
$\check\sigma$ to be the union of all simplices in $\Bar(\P)$
intersecting $\sigma$ but  disjoint from any proper subcell of
$\sigma$. Put
\[
\check\P=\{\check\sigma|\sigma\in\P\}.
\]
This is the usual dual cell to $\sigma$, with $\dim\check\sigma=n-
\dim\sigma$. The figure below shows what $\check\P$ looks like for the
example of Remark~\ref{straight}
\begin{center}
\includegraphics{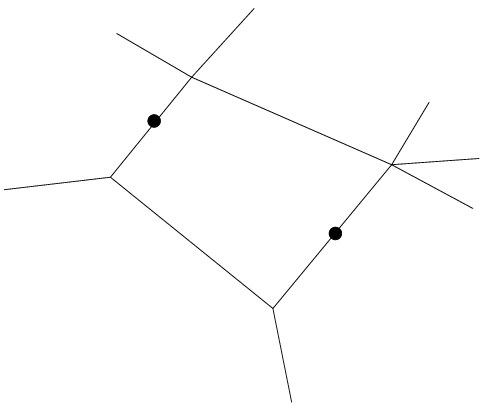}\\[-3ex]
\end{center}
Of course, $\check\sigma$ is usually not a polyhedron with respect to the
affine  structure on $B$, and we will build a new affine structure
using Construction
\ref{basicconstruct}.

For any $\sigma\in\P$, we obtain a fan $\Sigma_{\sigma}$ living in
$\shQ_{\sigma,\RR}$ and a piecewise linear function
$\varphi_{\sigma}$ on $\Sigma_{\sigma}$ by Definition~\ref{PLfan},
defined up to choice of a linear function. This function is strictly
convex by assumption, and we can consider the corresponding Newton
polytope, i.e. set
\[
\widetilde{\check\sigma}=
\{x\in\dualvs{\shQ_{\sigma,\RR}}|\langle x,y\rangle
\ge -\varphi_{\sigma}(y)\quad \forall y\in\shQ_{\sigma,\RR}\}.
\]
Note that because $\varphi_{\sigma}$ is strictly convex there is a 
one-to-one inclusion  reversing correspondence between the faces of
$\widetilde{\check\sigma}$ and cones in $\Sigma_{\sigma}$. Given
$e\in \Hom(\sigma,\tau)$ for some $\tau\in\P$, there is a corresponding
cone $\tau_e\in\Sigma_{\sigma}$, with corresponding cell
$\check\tau_e\subseteq\widetilde{\check\sigma}$ given by
\[
\check\tau_e
=\{x\in\widetilde{\check\sigma}|\langle x,y\rangle=-\varphi_{\sigma}(y)\quad
\forall y\in\tau_e\}.
\]
In addition, $\widetilde{\check\sigma}$ is an integral polytope
because $\varphi_{\sigma}$ has integral slopes.

We can glue these polyhedra together as follows. In the category
$\Cat(\P)$, let $e\in\Hom(\sigma,\tau)$, determining a face $\check\tau_{e}$ of
$\widetilde{\check\sigma}$. In addition, because
$\Sigma_{\tau}=\Sigma_{\sigma} (\tau_{e})$, it is easy to see that
$\widetilde{\check\tau}\subseteq  \dualvs{\shQ_{\tau,\RR}}\subseteq
\dualvs{\shQ_{\sigma,\RR}}$ is just a translate of the face
$\check\tau_{e}$ of $\widetilde{\check\sigma}$. This yields a
contravariant functor $F$ from $\Cat(\P)$ to the category of
topological spaces with $F(\tau)=\widetilde{\check\tau}$. For
$e\in\Hom(\sigma,\tau)$, $F(e):\widetilde{\check\tau}
\to\widetilde{\check\sigma}$ is the unique identification via
translation of $\widetilde{\check\tau}$ with the face
$\check\tau_{e}$ of $\widetilde{\check\sigma}$. We can then
construct the limit of $F$, i.e. the quotient space of
$\coprod_{\sigma\in\P}\widetilde{\check\sigma}$ by these
identifications. We will write this quotient map as
\[
\check\pi:\coprod_{\sigma\in\P}\widetilde{\check\sigma}\to \check B.
\]

It is easy to see this quotient space is homeomorphic to  $B$ itself,
with the image of $\widetilde{\check\sigma}$ in $B$ being
$\check\sigma$. This is just a combinatorial identification, but we
can achieve this as a $C^0$ map by mapping
$\widetilde{\check\sigma}$ to $\check\sigma\subseteq B$ in a piecewise linear way on
the barycentric subdivisions. In particular, we have now constructed
$\check B$ by identifying via affine transformations the faces of the
lattice polytopes $\widetilde{\check\sigma}$ where $\sigma$ runs
over minimal cells of $\P$. This is the first step in constructing an
affine manifold via the method of Construction \ref{basicconstruct}.
Note that under this piecewise linear identification of $B$ and
$\check B$, we have 
\[
\Bar(\check\P)=\Bar(\P)
\]
and
\[
\Cat(\check\P)=\Cat(\P)^{\rm op}.
\]

We also take the discriminant locus $\check\Delta'$ as usual as the
union of codimension two simplices not intersecting vertices of
$\check\P$ or interiors of maximal cells of $\check\P$. Of course
$\Delta\subseteq \check\Delta'$. The latter discriminant locus is
only provisional, and we will shrink it in Proposition~\ref{dlt} below.

To finish specifying an integral affine structure on $\check
B\setminus \check\Delta'$, we just need to specify a fan structure at
each vertex $\check\sigma$ of $\check\P$ (for $\sigma$ a maximal cell
of $\P$).

Let $y$ be the barycenter of a maximal cell $\sigma$, so
$\check\sigma =\{y\}$. The cones of the desired fan
$\Sigma_{\check\sigma}$ should be in one-to-one correspondence with
\[
\coprod_{\check\tau\in\check\P} \Hom(\check\sigma,\check\tau)
=
\coprod_{\tau\in\P} \Hom(\tau,\sigma),
\]
and the latter set is in one-to-one correspondence with cones in the
normal fan $\check\Sigma_{\sigma}$ (see Definition~\ref{normalfan}).
So it is enough to give integral affine identifications as follows.
Let $e\in\Hom(\tau,\sigma)$, so $y$ is an endpoint of $e$ and the
other endpoint of $e$ is the barycenter of $\tau$. Then by
Definition~\ref{normalfan}, we obtain $\tilde\tau\subseteq
\tilde\sigma\subseteq\Lambda_{\RR,y}$, well-defined up to
translation, and thus we can assume $0\in\tilde\tau$, getting the 
corresponding cone of $\check\Sigma_{\sigma}$ being
\begin{eqnarray*}
\check K_{e}&=&\{f\in\dualvs{\Lambda_{\RR,y}}|
\hbox{$f|_{\tilde\tau}$ is constant and $\langle f,z\rangle\ge
\langle f,x\rangle$ for all $z\in\tilde\sigma,x\in\tilde\tau$}\}\\
&=&\{f\in(\RR\tilde\tau)^{\perp}|\langle f,z\rangle\ge 0\quad
\forall z\in\tilde\sigma\}.
\end{eqnarray*}
On the other hand, the choice of $e$ selects a maximal cone 
$\sigma_{e}\in\Sigma_{\tau}$ and hence a vertex $\check\sigma_{e}$ of
$\widetilde{\check\tau}$, and $\widetilde{\check\tau}$ has tangent
wedge at this vertex given by the dual cone to $\sigma_e$, namely
\[
\{f\in\dualvs{\shQ_{\tau,\RR}}|\langle f,z\rangle\ge 0
\quad\forall z\in\sigma_{e}\subseteq \shQ_{\tau,\RR}\}.
\]
We must compare these two cones as follows. We can deform $e$
slightly to a path $\gamma$ connecting $y\in \Int(\sigma)$ to $y'\in
\Int(\tau)\setminus\Delta$; parallel transport along $\gamma$
identifies $\Lambda_{\RR,y}$ and $\Lambda_{\RR,y'}$ in such a  way that
$\RR\tilde\tau$ and $(\Lambda_{\P,\RR})_{y'}$ are identified. Under
this identification, $(\RR\tilde\tau)^{\perp}$ is identified with
$\dualvs{\shQ_{\tau,\RR}}$, and it is then clear the two cones defined
above coincide. This gives the integral affine transformation
identifying $\check K_e$ with the corresponding tangent wedge of
$\widetilde{\check\tau}$. This is sufficient to define the fan
structure of $\check B$ at $y$, and thus we obtain an integral affine
structure on $\check B \setminus\check\Delta'$.

\begin{proposition}
\label{dlt} 
Let $B$ be an integral affine manifold with singularities with a
toric polyhedral decomposition $\P$, and a strictly convex
multi-valued integral piecewise linear function $\varphi$ on $B$. Let
$\check B$ be as constructed above. Then
\item{(1)} Under the canonical identification $\check B=B$ the linear part
of the holonomy representation $\rho^{\check B}:\pi_1(\check
B\setminus \check\Delta') \to \Aff(N_{\RR})$ is dual to the linear
part of the holonomy representation $\rho^B:
\pi_1(B\setminus\Delta')\to \Aff(M_{\RR})$ (where
$\Delta':=\check\Delta'$). In other words,
\[
\tilde\rho^{\check B}(\gamma)={}^t\tilde\rho^B(\gamma^{-1})
\]
for all $\gamma\in\pi_1(B\setminus\Delta')$. Thus in particular,  by
Proposition~\ref{extaff}, the affine structure on $\check B\setminus
\check\Delta'$ extends to $\check B_0=\check B\setminus\check\Delta$,
where $\check\Delta=\Delta$. In addition, 
\[
\Lambda^{B_0}\cong \check\Lambda^{\check B_0}
\]
and
\[
\check\Lambda^{B_0}\cong \Lambda^{\check B_0},
\]
where the superscripts denote the affine structure with respect to
which these local systems are defined. (These are isomorphisms of
abstract sheaves, not as subsheaves of the tangent or cotangent
bundles, as the  identification $B_0=\check B_0$ is not differentiable.)
\item{(2)} $\check\P$ is a toric polyhedral decomposition of $\check B$.
\item{(3)} The radiance obstruction of $\check B_0$ coincides with
$c_1(\varphi)\in H^1(B,i_*\check\Lambda^{B_0})$ under a natural inclusion
\[
H^1(B,i_*\check\Lambda^{B_0})
\hookrightarrow H^1(\check B_0,\Lambda^{\check B_0}).
\]
\item{(4)} For $\sigma\in\P$, $\check\Sigma_{\sigma}$ and 
$\Sigma_{\check\sigma}$ coincide, and $\check\Sigma_{\check\sigma}$ and 
$\Sigma_{\sigma}$ coincide.
\end{proposition}

\proof (1) A path $\gamma$ in $\check B\setminus\check\Delta'$ 
passes successively through open sets $\Int(\check v_1)$ to
$W_{\check\sigma_1}$ to $\Int(\check v_2)$ to $W_{\check\sigma_2}$
etc., where $v_1,\ldots,v_n$ are vertices of $\P$ and
$\sigma_1,\ldots,\sigma_n$ are maximal cells of $\P$; here
$W_{\check\sigma_i}$ is the open neighbourhood of the vertex
$\check\sigma_i$ given in Definition~\ref{bary}, and in fact
$W_{\check\sigma_i}=\Int(\sigma_i)$ and $\Int(\check v_i)=W_{v_i}$.
Consider parallel transport from $\Int(\check v_i)$ to $W_{\check
\sigma_i}$ in the local system $\Lambda^{\check B_0}_{\RR}$ given by
the  affine structure on $\check B_0$. Now $v_i$ is the barycenter of
$\Int(\check v_i)$ and $\check\sigma_i=\{y_i\}$, with $y_i$ the
barycenter of $\sigma_i$. Then by construction, $\Lambda^{\check
B_0}_{\RR,v_i}$ is identified with $\dualvs{\shQ_{v_i,\RR}}
=\check\Lambda_{\RR,v_i}$ and $\Lambda^{\check B_0}_{\RR,y_i}$ is
identified with $\check\Lambda_{\RR,y_i}$. Furthermore, according to
the description of the fan structure above, $\check\Lambda_{\RR,v_i}$
and $\check\Lambda_{\RR,y_i}$ are identified using parallel
translation in $\Lambda_{\RR}$ between $v_i$ and $y_i$. Thus parallel
translation in $\Lambda^{\check B_0}_{\RR}$ from $v_i$ to $y_i$ is
the inverse transpose of parallel translation in $\Lambda_{\RR}$ from
$v_i$ to $y_i$. This makes it clear that $\tilde\rho^{\check
B}(\gamma)={}^t\tilde\rho^B(\gamma^{-1})$ for $\gamma$ a loop.

In particular, the linear part of the holonomy about a codimension
two simplex of $\check\Delta'$ is trivial if and only if the linear
part of the holonomy of $B$ is trivial around the simplex. In
addition, by Proposition~\ref{monodromy1}, the local radiance
obstruction near the simplex vanishes, so if the linear part of the
holonomy is trivial, so is the holonomy itself. Hence by
Proposition~\ref{extaff} we can take $\check\Delta=\Delta$. In
addition, the identifications of local systems follows from the
statement about $\tilde\rho^B$ and $\tilde\rho^{\check B}$.

(2) Since we constructed $\check B$ and $\check\P$ by means of
Construction~\ref{basicconstruct}, $\check\P$ is a polyhedral
decomposition. That it is toric follows from (1) and
Proposition~\ref{monodromy2}. Indeed, if $y\in
\Int(\check\tau)\setminus\check\Delta$, it is enough to show that for
any simple loop based at $y$ around a codimension two simplex
$\omega$ of $\check\Delta$ intersecting $\Int(\check\tau)$,
\[
(\tilde\rho^{\check B}(\gamma)-I)(\Lambda^{\check B_0}_{\RR,y})
\subseteq (\Lambda^{\check B_0}_{\check\P,\RR})_y=\Lambda^{\check
B_0}_{\check\tau, \RR}.
\]
But $\Lambda_{\RR,y}^{\check B_0}$ can be identified with 
$\check\Lambda_{\RR,y},$ and $\Lambda^{\check B_0}_{\check\tau,\RR}$
can be identified with $\dualvs{\shQ_{\tau,\RR}}$. Also
$\tilde\rho^{\check B}(\gamma) ={}^t\tilde\rho^B(\gamma^{-1})$. Now
$\omega$ contains the barycenter of $\check\tau$, which is also the
barycenter of $\tau$, so by Proposition \ref{monodromy1},
\[
(\tilde\rho^B(\gamma^{-1})-I)((\Lambda_{\P,\RR})_{y})=0,
\]
from which it immediately follows that
\[
(\tilde\rho^{\check B}(\gamma)-I)(\check\Lambda_{\RR,y})
\subseteq\dualvs{\shQ_{\tau,\RR}},
\]
which is what is desired. 

(3) We will compare the radiance obstruction of the affine structure
on $\check B\setminus\check\Delta'$ in  $H^1(\check
B\setminus\check\Delta',\Lambda^{\check B_0})=H^1(B\setminus\Delta',
\check\Lambda^{B_0})$ and $c_1(\varphi)\in H^1(B\setminus\Delta',
\check\Lambda^{B_0})$. Write $\check\Lambda$ for
$\check\Lambda^{B_0}$. This computation will be sufficient since if
$i':B\setminus\Delta'\hookrightarrow B_0$ and $i:B_0\hookrightarrow
B$ are the inclusions, then it is easy to see that
$i_*'(\check\Lambda|_{B\setminus \Delta'})=\check\Lambda$ on $B_0$,
and then by the Leray spectral sequences for $i$ and $i'$ we get a
series of inclusions
\[
H^1(B,i_*\check\Lambda^{B_0})\hookrightarrow H^1(B_0,\check\Lambda^{B_0})
\hookrightarrow H^1(B\setminus\Delta',\check\Lambda^{B_0}).
\]
Thus if $c_1(\varphi)$ and the radiance obstruction coincide in
$H^1(B\setminus\Delta',\check\Lambda^{B_0})$, then they in fact coincide
in $H^1(B,i_*\check\Lambda^{B_0})$.

To calculate the radiance obstruction and $c_1(\varphi)$ on $B\setminus\Delta'$,
we can use \v Cech cohomology with respect to the open covering
of $B\setminus\Delta'$ given by
\[
\{\Int(\check v)=W_v| \hbox{$v$ a vertex in $\P$}\}\cup
\{\Int(\sigma)=W_{\check\sigma}|\sigma\in \P_{\max}\}.
\]
We can choose the zero function as a representative for $\varphi$ on
each $\Int(\sigma)$, as $\varphi$ is linear there. On each $W_v$, we
choose a representative $\varphi_v$ for $v$ with $\varphi_v(v)=0$; we
can then view $\varphi_v$ as a piecewise linear function on the fan
$\Sigma_v$ in $\Lambda_{\RR,v}$. As such, $\varphi_v$ is given by an
element $n_e\in \check\Lambda_{\RR,v}$ for each $e\in 
\coprod_{v\in\sigma} \Hom(v,\sigma)$ corresponding to maximal cones
of $\Sigma_v$. The \v Cech 1-cocycle representing $c_1(\varphi)$ then
takes the value $n_e$ on the connected component of $W_v\cap
\Int(\sigma)$ corresponding to $e\in\Hom(v,\sigma)$.

On the other hand, to compute the radiance obstruction of $\check B
\setminus\Delta'$, we can use the \v Cech realisation given in
Remark~\ref{realising}. There is a canonical coordinate chart
$\psi_{\check\sigma}:W_{\check\sigma}\to\check
\Lambda^{B_0}_{\RR,\check\sigma} $ taking $\check\sigma$ to zero, the
graph of which gives a section of $\T_{W_{\check\sigma}}$. On the
other hand, the choice of function $\varphi_v$ determines the
polyhedron $\widetilde{\check v}\subseteq \check
\Lambda_{\RR,v}^{B_0}$, with vertices $-n_e$ for $e\in
\coprod_{v\in\sigma}\Hom(v,\sigma)$. This gives a chart $\psi_{\check
v}:\Int(\check v)\to\check\Lambda^{B_0}_{\RR,v}$ identifying
$\Int(\check v)$ with $\Int\big(\hspace{1pt}\widetilde{\check v}
\hspace{1pt}\big)$. Again the graph of this chart gives a section of
the tangent bundle $\T_{W_{\check v}}$. It is then clear, using
parallel translation along $e$ to identify
$\check\Lambda^{B_0}_{\RR,v}$ and  $\check
\Lambda^{B_0}_{\RR,\check\sigma}$, that the difference of these two
sections is $n_e$ on the connected component of $W_v\cap
W_{\check\sigma}$ corresponding to $e$. Thus $c_1(\varphi)$ and the
radiance obstruction coincide.

(4) follows easily from the definitions.
\qed
\medskip

Finally, we wish to define $\check\varphi$, the Legendre transform of
$\varphi$. For each $\sigma\in\P$, choose $\tilde\sigma\subseteq
\Lambda_{\sigma,\RR}$ as in Definition~\ref{normalfan}, well-defined
up to translation. Let $\check\varphi_{\check\sigma}$ be the function
defined on the normal fan $\check\Sigma_{\sigma}$ given by
\[
\check\varphi_{\check\sigma}(y)=-\inf\{\langle y,x\rangle|x\in\tilde\sigma\}
\]
for $y\in \dualvs{\Lambda_{\sigma,\RR}}$. This is a piecewise linear function
on the fan $\check\Sigma_{\sigma}=\Sigma_{\check\sigma}$, and it is a
standard easy fact that it is strictly convex, (see \cite{Oda}, A.3, with
opposite sign conventions) with the Newton polytope
of $\check\varphi_{\check\sigma}$ being $\tilde\sigma$. If $\tilde\sigma$
is changed by a translation, $\check\varphi_{\check\sigma}$ is changed
by a linear function, so it is well-defined modulo linear functions.

Now let $S_{\check\sigma}:U_{\check\sigma}\to \dualvs{\Lambda_{\sigma,
\RR}}$ be the affine submersion given by the fact that $\check\P$ is
toric, with $U_{\check\sigma}$ sufficiently small. Then we can denote
also by $\check\varphi_{\check\sigma}$ the composition
$\check\varphi_{\check\sigma}\circ S_{\check\sigma}$. It is then easy
to see

\begin{proposition}
\label{dlt2}
\item{(1)} The $(U_{\check\sigma},\check\varphi_{\check\sigma})$
determine an integral multi-valued piecewise linear function on $\check B$
which is strictly convex. We call this $\check\varphi$, and say the discrete 
Legendre transform of $(B,\P,\varphi)$ is $(\check B,\check\P,\check\varphi)$.
\item{(2)} The discrete  Legendre transform of $(\check B,\check
\P,\check\varphi)$ is $(B,\P,\varphi)$.
\end{proposition}

\begin{example}
\label{torusdlt}
Let us relate this version of the discrete Legendre transform to the
more traditional form for polyhedral decompositions of $M_{\RR}$. We
do this by taking $B=M_{\RR}/\Gamma$ for some integral lattice
$\Gamma\subseteq M$. Then a polyhedral decomposition $\P$ of $B$ is
induced by a periodic polyhedral decomposition of $M_{\RR}$ (which we
will also call $\P$). A multi-valued strictly convex piecewise linear
function $\varphi$ on $B$ can then be viewed, by patching together
choices of values for $\varphi$, as a single-valued strictly convex
piecewise linear function $\varphi:M_{\RR} \to \RR$ satisfying the
periodicity condition
\[
\varphi(x+\gamma)=\varphi(x)+\alpha(\gamma)(x),
\]
for some map $\alpha:\Gamma\to \shAff(M_{\RR},\RR)$. We then define the
dual polyhedral decomposition $\check\P$ in  $N_{\RR}$ as follows.
For each $\sigma\in\P_{\max}$, let  $n_{\sigma}\in N_{\RR}$ be the
slope of $\varphi|_{\sigma}$. For each $\tau\in\P$, let $\check\tau$
be the convex hull in $N_{\RR}$ of $\{n_{\sigma}|\hbox{$\tau
\subseteq \sigma$ maximal}\}$. It is easy to see that
$\check\P=\{\check\tau| \tau\in\P\}$ forms a polyhedral decomposition
of $N_{\RR}$. Furthermore, $\Gamma$ embeds in $N_{\RR}$ by mapping
$\gamma\in \Gamma$ to the linear part of $\alpha(\gamma)$. It is easy
to see that $\check\P$ is  periodic with respect to $\Gamma\subseteq
N_{\RR}$ and hence descends to a decomposition of $\check B=
N_{\RR}/\Gamma$. If $\varphi$ has integral slopes, then
$\Gamma\subseteq N$, and $\check\P$ is an integral polyhedral
decomposition. Finally, we can define the Legendre transform of
$\varphi$ by, for $n\in \check v\subseteq N_{\RR}$ a maximal cell of
$\check\P$,
\[
\check \varphi(n)=\langle n,v\rangle -\varphi(v).
\]
This in turn defines a multi-valued function $\check\varphi$ on
$\check B$. It is easy to check that $(\check
B,\check\P,\check\varphi)$ agrees with the discrete Legendre
transform defined earlier of $(B,\P,\varphi)$ up to sign conventions.

This discrete Legendre transform appears to be useful in applied
mathematics: see for example \cite{chyn}.
\qed
\end{example}

\begin{example}
\label{Batyrevdual}
The discrete Legendre transform is also a generalization of Batyrev
duality \cite{Bat}. 
Let $\Xi\subseteq M_{\RR}$ be a reflexive polytope, with $0$
the unique interior lattice point. Then $B=\partial\Xi$ has an affine
manifold with singularities structure as in Example~\ref{polytope}.
Let $\Sigma$ be the fan in $M_{\RR}$ obtained from $\Xi$ by taking
$\Sigma$ to consist of cones over proper faces of $\Xi$. Then in fact
any piecewise linear function $\varphi$ on $\Sigma$ gives rise to a
multi-valued piecewise linear function $\varphi_B$ on $B$ as
follows. 

For any proper face $\tau$ of $\Xi$, or equivalently $\tau\in\P$, let
$W_{\tau}$ be as in Definition~\ref{bary},
so that $\{W_{\tau}\}$ form an open cover of $B$. For each
$\tau$, choose $n_{\tau}\in N_{\RR}$ such that $\varphi(m)=\langle
n_{\tau},m\rangle$ for all $m\in\tau$, and let
$\varphi_{\tau}:W_{\tau} \to\RR$ be the restriction of
$\varphi-n_{\tau}$ to $W_{\tau}$. It is then easy to see that
$\varphi_B=\{(W_{\tau},\varphi_{\tau}) |\tau\in\P\}$ represents a
multi-valued piecewise linear function on  $B$, strictly convex if
$\varphi$ is. If $\varphi$ is integral, then $n_{\tau}$ can be chosen
in $N$, and $\varphi_B$ is also integral. Note that $\varphi_B$ is
not the same thing as the restriction of $\varphi$ to $\partial\Xi$.

If $\varphi$ is taken to represent the anti-canonical class of the
toric variety defined by $\Sigma$, i.e. $\varphi(v)=1$ for all
vertices $v$ of $\Xi$, then $\varphi_B$ is strictly convex. As an
exercise, check in this case that the discrete Legendre transform of
$(B,\P,\varphi_B)$ is obtained as the boundary of the Batyrev dual
$\Xi^*\subseteq N_{\RR}$, defined by 
\[
\Xi^*=\{x\in N_{\RR}| \hbox{$\langle x,y\rangle \ge -1$ for all 
$y\in\Xi$}\},
\]
with affine structure as in Example~\ref{polytope}.
\end{example}

\subsection{Positivity and simplicity}

We will now introduce several possible restrictions on integral
affine manifolds with singularities and polyhedral decompositions.
These restrictions will be used in \S\S4 and 5, though we will give
some motivation here for introducing them.

Let $B$ be an integral affine $n$-dimensional manifold with
singularities and $\P$ a toric polyhedral decomposition. We will,
once and for all, choose for each $\omega\in\P$ with $\dim\omega=1$
an integral generator $d_{\omega}$ of $\Lambda_{\omega}$. In
addition, for each $\rho\in\P$ with $\dim\rho=n-1$, we choose an
integral generator $\check d_{\rho}$ of $\dualvs{\shQ_\rho}$. This is a
rank 1 lattice contained naturally in $\check\Lambda_y$ for any $y\in
\Int(\rho)\setminus\Delta$. 

These choices give us a canonical ordering of the vertices of
$\omega$ or the maximal cells containing $\rho$. Indeed, for
$\dim\omega=1$, the morphisms $e:v\to\omega$ are in one-to-one
correspondence with the maximal cones of the dual fan
$\check\Sigma_{\omega}$ living in the one-dimensional space
$\dualvs{\Lambda_{\omega,\RR}}$. Then as $d_{\omega}\in
\Lambda_{\omega}$, there is a well-defined $e_{\omega}^+:
v^+_{\omega}\to\omega$ corresponding to the cone on which
$d_{\omega}$ is positive, and an $e_{\omega}^-: v_{\omega}^-
\to\omega$ corresponding to the cone on which  $d_{\omega}$ is
negative. 

The choice of $\check d_{\rho}$ distinguishes the morphisms from
$\rho$ to maximal cells: $g_{\rho}^+:\rho\to\sigma^+_{\rho}$ and
$g^-_{\rho}:\rho\to\sigma^-_{\rho}$ as follows. The fan
$\Sigma_{\rho}$ in $\shQ_{\rho,\RR}$ consists of two cones
corresponding to the two different morphisms to maximal cells. We
take $g^+_{\rho}$ to correspond to the cone on which $\check
d_{\rho}$ is positive and $g^-_{\rho}$  to that cone on which $\check
d_{\rho}$ is negative.

Now in general, if $\tau\in\P$, $f_i:v_i\to\tau$, $e_i:\tau\to
\sigma_i$ for $i=1,2$, $v_i$ vertices of $\P$ and $\sigma_i\in
\P_{\max}$, we obtain a loop $\gamma^{e_1e_2}_{f_1f_2}$ in $B_0$
based at $v_1$, by traversing in order the edges $e_1\circ
f_1:v_1\to\sigma_1$,  $e_1\circ f_2:v_2\to\sigma_1$ (backwards), $e_2\circ
f_2:v_2\to\sigma_2$, $e_2\circ f_1:v_1\to\sigma_2$  (backwards) of $\Bar(\P)$:
\begin{center}
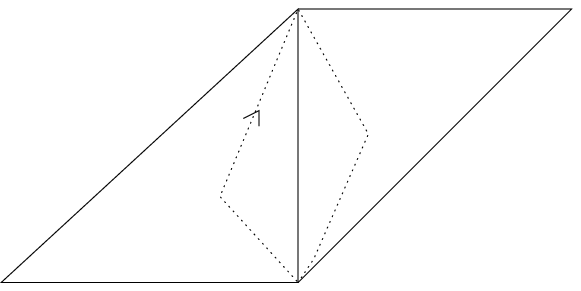
\end{center}

This induces a monodromy transformation
\[
T_{f_1f_2}^{e_1e_2}:\Lambda_{v_1}\to\Lambda_{v_1}
\]
via parallel transport around $\gamma^{e_1e_2}_{f_1f_2}$.  It is
often convenient to parallel transport to a point $y\in
\Int(\tau)\setminus\Delta$ near $v_1$, to get
$T_{f_1f_2}^{e_1e_2}:\Lambda_y\to\Lambda_y$. By Propositions
\ref{monodromy1} and~\ref{monodromy2},
$T^{e_1e_2}_{f_1f_2}|_{\Lambda_{\tau}}$ is the identity, and
$(T_{f_1f_2}^{e_1e_2}-I)(\Lambda_y) \subseteq \Lambda_{\tau}$.

There are two special cases where this gives useful information. If
$\tau= \omega\in\P$ is dimension one, we have
$e_{\omega}^{\pm}:v_{\omega}^{\pm} \to\omega$. So given
$e_i:\omega\to\sigma_i$, $i=1,2$, we get
\[
T_{\omega}^{e_1e_2}:=T^{e_1e_2}_{e^+_{\omega}e^-_{\omega}}.
\]
It follows that there exists an $n_{\omega}^{e_1e_2}\in
\Lambda_{\omega}^{\perp}=\dualvs{\shQ_{\omega}}\subseteq\check\Lambda_y$
such that $T_{\omega}^{e_1e_2}:\Lambda_y\to\Lambda_y$ is given by 
\[
T^{e_1e_2}_{\omega}(m)=m+\langle n_{\omega}^{e_1e_2},m\rangle d_{\omega}.
\]
On the other hand, if $\tau=\rho\in\P$ is dimension $n-1$, 
$f_i:v_i\to\rho$, we get
\[
T_{f_1f_2}^{\rho}:=T^{g^+_{\rho}g^-_{\rho}}_{f_1f_2}.
\]
Then there exists an $m^{\rho}_{f_1f_2} \in\Lambda_{\rho}
\subseteq\Lambda_y$ such that
\[
T^{\rho}_{f_1f_2}(m)=m+\langle\check d_{\rho},m\rangle 
m^{\rho}_{f_1f_2}.
\]
Combining these two cases, if we are given a map $e:\omega\to\rho$
with $\dim\omega=1$ and $\dim\rho=n-1$, we can define
\[
T_e:=T_{\omega}^{g^+_{\rho}\circ e, g^-_{\rho}\circ e}
=T^{\rho}_{e\circ e_{\omega}^+, e\circ e_{\omega}^-}
\]
and then we have
\[
n_e:=n_{\omega}^{g^+_{\rho}\circ e,g_{\rho}^-\circ e}, \quad
m_e:=m^{\rho}_{e\circ e^+_{\omega}, e\circ e^-_{\omega}}.
\]
Note $n_e$ is proportional to $\check d_{\rho}$, and $m_e$ is
proportional to $d_{\omega}$, and the constant of proportionality is
the same.

It is also not hard to see that all the monodromy data can be encoded
in these numbers. In particular, the $n_e$ (or $m_e$) for all
$e:\omega\to\rho$ determine the minimal discriminant locus $\Delta$,
if $B$ is obtained from Construction~\ref{basicconstruct}. Indeed,
for $n>2$ if $\tau$ is an $n-2$ dimensional simplex in $\Bar(\P)$
contained in $\Delta'$ and containing the edge $e:\omega\to\rho$,
then a loop about $\tau$ in $B_0$ is homotopic (up to orientation) to
$\gamma^{g^+_{\rho}, g^-_{\rho}}_{e\circ e^+_{\omega}, e\circ
e^-_{\omega}}$. Thus if $n_e\not=0$, all $n-2$ dimensional simplices
of $\Delta'$ containing $e$ must be included in $\Delta$; if $n_e=0$,
they are not included in $\Delta$ (see Proposition~\ref{extaff}). Of
course, the same is true for $m_e$. For $n=2$ analogous statements
hold but now $\omega=\tau=\rho$ is an edge.

\begin{definition}
\label{positive def}
Let $B$ be an integral affine manifold with singularities and a toric 
polyhedral decomposition $\P$. We say $(B,\P)$ is \emph{positive} if
for all $e:\omega\to\rho$ with $\dim\omega=1$, $\dim\rho=n-1$,
$m_e$ is a non-negative multiple of $d_{\omega}$, or equivalently
$n_e$ is a non-negative multiple of $\check d_{\rho}$.
\end{definition}

This may seem mysterious at first, though the importance of
positivity will become clearer in \S4. Note that it is independent
of the choices of $d_{\omega}$ or $\check d_{\rho}$: indeed, changing
the sign of one or the other also changes the orientation of the
loop.

This condition can be explained very roughly in dimension two as
follows. Given a 2-torus fibration $f:X\to B$ from an oriented
four-manifold to a two-dimensional base, with degenerate fibres which
are pinched tori, there are actually two sorts of such degenerate
fibres, distinguished by the sign of the intersection number of the
two sheets of the singular fibre at the pinch point. Of course, if
$f:X\to B$ is an elliptic fibration (or special Lagrangian
fibration), only the positive case occurs. This places a restriction
on the monodromy of the fibration, and it is precisely this
restriction that the notion of positivity captures.

\begin{proposition}
If $\varphi$ is a strictly convex multi-valued piecewise linear
function on $(B,\P)$, $(\check B,\check\P,\check\varphi)$ the
discrete Legendre transform of $(B,\P,\varphi)$, then $(B,\P)$ is positive
if and only if $(\check B,\check\P)$ is positive.
\end{proposition}

\proof Under the identification of Proposition~\ref{dlt}, (1), of
$\Lambda^{B_0}$ with $\check\Lambda^{\check B_0}$, we can identify,
for $y\in \Int(\omega)\setminus\Delta$, $\Lambda^{B_0}_{\omega}$ with
$\dualvs{(\shQ^{\check B}_{\check\omega})}\subseteq
\check\Lambda_y^{\check B_0}$, and thus we can equate a choice of
$d_{\omega}$ with a choice of $\check d_{\check\omega}$. Similarly,
we can identify $\dualvs{(\shQ^B_{\rho})}$ with $\Lambda^{\check
B_0}_{\check \rho}$ for $\check y\in
\Int(\check\rho)\setminus\check\Delta$. Thus we can equate a choice
of $\check d_{\rho}$ with a choice of $d_{\check\rho}$. 

Now suppose we are given $e:\omega\to\rho$, and hence $\check
e:\check\rho\to\check\omega$. We are then given
$e^{\pm}_{\omega}:v^{\pm}_{\omega}\to\omega$ and
$g_{\rho}^{\pm}:\rho\to\sigma_{\rho}^{\pm}$. It is easy to check that
$\check g_{\rho}^{\pm}:\check\sigma_{\rho}^{\pm}\to\check\rho$
coincides with $e^{\pm}_{\check\rho}:v^{\pm}_{\check\rho}\to\check
\rho$, preserving the sign, and similarly, $\check e^{\pm}_{\omega}:
\check\omega\to \check v_{\omega}^{\pm}$ coincides with
$g^{\pm}_{\check \omega}: \check\omega\to
\sigma^{\pm}_{\check\omega}$.  Thus in $B$, the monodromy
transformation $T_e$ is given by transport along the loop $\gamma$
which goes from $v^+_{\omega}$ to the barycenter of $\sigma_{\rho}^+$
to $v^-_{\omega}$ to the barycenter of  $\sigma_{\rho}^-$ to
$v_{\omega}^+$, while in $\check B$, the monodromy transformation
$T_{\check e}$ is given by transport along $\check\gamma$ which goes
from $\check\sigma_{\rho}^+$ (the barycenter of $\sigma_{\rho}^+$) to
the barycenter of $\check v_{\omega}^+$ (i.e. $v^+_{\omega}$) to
$\check\sigma_{\rho}^-$ to the barycenter of $\check v^-_{\omega}$
and then to $\check\sigma_{\rho}^+$. But this is precisely the
opposite direction of $\gamma$. It then follows from
Proposition~\ref{dlt}, (1), that $T_{\check e}=T_e^t$. Thus if
$T_e(m)=m+\langle n_e,m\rangle d_{\omega}$ then
$T_{\check e}(n)=n+\langle n,d_{\omega}\rangle n_e$.
Thus $m_{\check e}=n_e$, and one is a positive multiple of $\check
d_{\rho}$ if and only if the other is.
\qed

\begin{remark}
\label{psiomega}
Positivity can also be interpreted as follows. For a given
$\omega\in\P$ with $\dim\omega=1$, the data $n_{\omega}^{e_1e_2}$,
$e_i:\omega \to\sigma_i\in\P_{\max}$ determine a piecewise linear
function $\psi_{\omega}$ on the fan $\Sigma_{\omega}$, well-defined
up to a linear function, as follows. For any $e:\omega\to\tau\in\P$,
there is a cone in $\Sigma_{\omega}$ corresponding to $e$, which we
write as $K_e$. Now fix some  $e:\omega\to\sigma\in\P_{\max}$, and
take $\psi_{\omega}|_{K_e}=0$. On any other maximal cone
$K_{e'}\in\Sigma_{\omega}$, take  $\psi_{\omega}|_{K_{e'}}=
-n_{\omega}^{ee'}$. We just need to check this is continuous. If two
maximal cones $K_{e_1}, K_{e_2}$ of $\Sigma_{\omega}$ intersect in a
cone $K_f$ of $\Sigma_{\omega}$, giving a diagram
\[
\xymatrix@C=30pt
{&&\sigma_1\\
\omega\ar@/^/[rru]^{e_1}\ar[r]^(.7){f}\ar@/_/[drr]_{e_2}&\tau
\ar[ru]\ar[rd]&\\
&&\sigma_2
}
\]
then we need to check $n_{\omega}^{ee_1}=n_{\omega}^{ee_2}$ on $K_f$.
But it is easy to see that
\[
T_{\omega}^{e_1e_2}\circ T_{\omega}^{ee_1}=T_{\omega}^{ee_2},
\]
and thus 
\[
n_{\omega}^{e_1e_2}=n_{\omega}^{ee_2}-n_{\omega}^{ee_1}.
\]
But $n_{\omega}^{e_1e_2}\in\Lambda^{\perp}_{\tau}$ by 
Proposition~\ref{monodromy1}, so is zero on
$K_f$. This gives continuity.

Note that changing the initial choice of $e$ just changes
$\psi_{\omega}$ by a linear function.

The condition of positivity is now easily seen to be equivalent to 
$\psi_{\omega}$ being convex for all $\omega$ (though not necessarily
strictly convex).

Similarly the data $m_{f_1f_2}^{\rho}$ determine a piecewise linear
function $\check\psi_{\rho}$ on the normal fan $\check\Sigma_{\rho}$;
again, positivity is seen to be equivalent to $\check\psi_{\rho}$ convex
for all $\rho$.
\qed
\end{remark}

\begin{example}
One can check that the construction of Example~\ref{polytope} only
produces positive $(B,\P)$. To do this, given $\Xi\subseteq M_{\RR}$, 
$\omega\subseteq\rho$ faces of dimension one and of codimension one
respectively, and a choice of $d_{\omega}$ and $\check d_{\rho}$, we
obtain the endpoints $v^+_{\omega}$ and $v^-_{\omega}$ of $\omega$ and
the maximal faces $\sigma_{\rho}^+$ and $\sigma_{\rho}^-$ containing
$\rho$. There exists $n^+,n^-\in N$ such that $\langle
n^+,\sigma_{\rho}^+\rangle=\langle n^-,\sigma_{\rho}^- \rangle=1$, by
reflexivity of $\Xi$. On the other hand, $\Lambda_{v^+_{\omega}}$ can
be identified with $M/v^+_{\omega}$, and with $e:\omega\to\rho$ the
unique morphism, 
\[
T_e(m)=m+\langle m,n_+-n_-\rangle (v_{\omega}^- -v_{\omega}^+),
\]
as can be computed as in \cite{Ruan}, Theorem 8.2, or \cite{HZ},
Lemma 2.4. See also \cite{GBB}, \S 2.
But since  $\Xi$ is convex, $n_+-n_-$ is positive on
$\sigma_{\rho}^+$. Also, $v_{\omega}^- -v_{\omega}^+$ is a positive
multiple of $d_{\omega}$. Hence $B$ is positive.
\end{example}
\bigskip

We will define a certain class of $(B,\P)$ which we will call
\emph{simple}. At this point in time, this definition is a bit
difficult to motivate. However, it should be viewed as analogous to
the notion of simplicity defined in \cite{SlagI} and used extensively
in  \cite{SlagII} and \cite{TMS}. Given a torus fibration $f:X\to B$,
with discriminant locus $\Delta$, $B_0=B\setminus\Delta$, $X_0
=f^{-1}(B_0)$, $f_0=f|_{X_0}$ and $i:B_0\hookrightarrow B$ the
inclusion,  we said $f$ was $G$-simple for a group $G$ if
$i_*R^pf_{0*}G= R^pf_*G$ for all $p$. This essentially says that the
cohomology of any singular fibre is isomorphic to the monodromy
invariant part of the cohomology of a nearby non-singular fibre. It
was only with this assumption that the Strominger-Yau-Zaslow
conjecture implies the usual connection between $h^{1,1}$ and
$h^{1,2}$ of mirror pairs of Calabi-Yau threefolds. If the fibrations
are not simple, the connection between these numbers is more tenuous,
and in particular we may not expect a clear isomorphism between
complexified K\"ahler moduli and complex moduli of the mirror.

One would like some analogue of simplicity for affine manifolds with
singularities. We cannot use the definition directly, because we don't
have torus fibrations. Instead we will use a definition which places
restrictions on the monodromy of the local system $\Lambda$ on
standard open sets $W_{\tau}\cap B_0$. In dimensions two and three,
this will essentially coincide with the restrictions on monodromy
implied by the old definition of simplicity. However, our new notion
of simplicity also depends on $\P$: the polyhedral decomposition
determines the ``scale'' of detail we see in the discriminant locus.
The definition of simplicity only depends on the monodromy
transformations $T_e$ for various $e$. For example, if $\dim B=2$ and
there are several points in $\Delta$ on a single edge, simplicity
places a restriction on the monodromy of a loop enclosing all these
points, and not on loops around each point on the edge. So even in the
two-dimensional case, simplicity in the old sense only implies
simplicity in the new sense if  $\P$ is sufficiently fine. Again, it
is useful to assume the discriminant locus has already been
straightened to avoid this problem. In any event, with these caveats,
our definition clearly gives the right thing in dimensions two and
three, and we shall see more clearly the naturality of this definition
in \S 5, as well as in \cite{GBB} and \cite{sequel}.

Before giving our definition of simplicity, we define some auxiliary
polytopes associated with $(B,\P)$ which encode key data about the
monodromy.

\begin{definition}
\label{NDelta}
For a given $\omega\in\P$ with $\dim\omega=1$, fix
$e:\omega\to\sigma\in\P_{\max}$, and let $\check\Delta(\omega)$ be the
convex hull in $\Lambda_{\omega,\RR}^{\perp}
\subseteq\check\Lambda_{\RR,y}$ (for $y\in
\Int(\omega)\setminus\Delta$ near $v^+_{\omega}$) of the set
\[
\{n_{\omega}^{ee'}|e':\omega\to\sigma'\in\P_{\max}\}.
\]
Similarly, for a given $\rho$, fix $f:v\to\rho$ and let 
$\Delta(\rho)$ be the convex hull in 
$\Lambda_{\rho,\RR}\subseteq\Lambda_{\RR,y}$ (for $y\in
\Int(\rho)\setminus \Delta$ near $v$) of
$\{m^{\rho}_{ff'}|f':v'\to\rho\}$.

These polytopes are well-defined up to translation, independent of
the choice of $\sigma$ or $v$. 

If $\tau\in\P$ with $1\le\dim\tau\le n-1$, we can also define
auxiliary  polytopes related to $\tau$. Given $e:\omega\to\tau$ with 
$\dim\omega=1$, choose some $f:\tau\to\sigma\in\P_{\max}$, and set
$\check\Delta_e(\tau)$ to be the convex hull of
\[
\{n_{\omega}^{f\circ e,f'\circ e}| f':\tau\to\sigma'\in\P_{\max}\}
\]
in $\Lambda^{\perp}_{\tau,\RR}
=\dualvs{\shQ_{\tau,\RR}}\subseteq\check\Lambda_{\RR,y}$, for $y\in
\Int(\tau)\setminus\Delta$ near $v^+_{\omega}$. Thus $\check\Delta_e(\tau)$
is just a face of $\check\Delta(\omega)$.

Similarly, given $e:\tau\to\rho$, $\dim\rho=n-1$, choose some 
$f:v\to\tau$ and set $\Delta_e(\tau)$ to be the convex hull of
$\{m^{\rho}_{e\circ f,e\circ f'}| f':v'\to\tau\}$. Again,
$\Delta_e(\tau)$ is a face of $\Delta(\rho)$.
\qed
\end{definition}

\begin{remark}
\label{polytopefacts}
1)\ \ 
These polytopes are only of interest in the positive case, in which
case $\check\Delta(\omega)$ and $\Delta(\rho)$ can be viewed as the
Newton polytopes determined by $\psi_{\omega}$ and
$\check\psi_{\rho}$ (Remark \ref{psiomega}) respectively. In this
case, $\psi_{\omega}$ can be recovered, (up to a linear function) as
\[
\psi_{\omega}(y)=-\inf\{\langle y,x\rangle | x\in\check\Delta(\omega)\}
\]
and similarly
\[
\varphi_{\rho}(x)=-\inf\{\langle y,x\rangle | y\in\Delta(\rho)\}.
\]
Furthermore, by construction of $\psi_{\omega}$, if
$e_i:\omega\rightarrow \sigma_i\in\P_{\max}$, $i=1,2$ define cones
$K_{e_i}\in\Sigma_{\omega}$, and $\psi_{\omega}$ plus any linear
function is given by $-p_{e_i}$ on $K_{e_i}$, then
$n_{\omega}^{e_1e_2}=p_{e_2}-p_{e_1}$. Thus all the 
$n_{\omega}^{e_1e_2}$'s can be recovered from knowing
$\check\Delta(\omega)$ in the positive case, and similarly the
$m^{\rho}_{e_1e_2}$'s can be recovered from $\Delta(\rho)$ for
$e_i:v_i\rightarrow\rho$.

The same comments hold for the $\check\Delta_e(\tau)$ and
$\Delta_e(\tau)$: given $e:\omega\rightarrow\tau$ with
$\dim\omega=1$, the $n_{\omega}^{f_1\circ e,f_2\circ e}$ can be
recovered from $\check\Delta_e(\tau)$ for $f_i:\tau\rightarrow
\sigma_i\in\P_{\max}$; similarly given $f:\tau\rightarrow\rho$ with
$\codim\rho=1$, the $m^{\rho}_{f\circ e_1,f\circ e_2}$ can be
recovered from $\Delta_e(\tau)$ for $e_i:v_i\rightarrow\tau$.\\[1ex]
2)\ \ From these remarks it also follows readily that $\Sigma_\omega$
is a refinement of the normal fan of $\check \Delta(\omega)$. The
normal vector $\check d_\rho$ of a codimension-one cell $\tilde\rho$
containing $\tilde\omega$ is parallel to an edge of
$\check\Delta(\omega)$ iff $n_e\neq
0$ for $e:\omega\to\rho$ the corresponding morphism. This is the case
iff the corresponding edge ($n=2:$ vertex) of the barycentric
decomposition is contained in $\Delta$.

Similarly, $\check \Sigma_\rho$ is a refinement of the normal fan of
$\Delta(\rho)$. Hence $\Delta(\rho)$ is a polytope with faces parallel
to faces of $\tilde\rho\subset\Lambda_{\rho,\RR}$ itself. Again, an
edge $\tilde\omega\subset\tilde\rho$ occurs in this list iff the
corresponding edge of the barycentric subdivision is part of $\Delta$.

Note that if all $n_e$ (or, equivalently, all $m_e$) are primitive
it follows that all these polytopes are determined uniquely just by
$\P$ and the combinatorics of the discriminant locus relative
to $\P$. 
\qed
\end{remark}

\begin{definition}
\label{simplicity}
Let $B$ be an integral affine manifold with singularities and $\P$ a
toric polyhedral decomposition of $B$. Suppose $(B,\P)$ is positive.
Then we say $(B,\P)$ is \emph{simple} if for every $\tau\in\P$ with
$1\le\dim\tau\le n-1$, the following condition holds:

Let 
\[
\P_1(\tau)=\{e:\omega\to\tau|\dim\omega=1\}
\]
and
\[
\P_{n-1}(\tau)=\{f:\tau\to\rho|\dim\rho=n-1\}.
\]
\pagebreak

\noindent
Then there exist disjoint subsets 
\[
\Omega_1,\ldots,\Omega_p\subseteq 
\P_1(\tau)
\]
and
\[
R_1,\ldots,R_p\subseteq \P_{n-1}(\tau)
\]
for some $p\ge 0$ such that
\item{(1)} For $e\in\P_1(\tau)$ and $f\in\P_{n-1}(\tau)$, 
$n_{f\circ e}=0$ unless $e\in\Omega_i$, $f\in R_i$ for some $i$.
\item{(2)} For each $1\le i\le p$, $\check\Delta_e(\tau)$ is independent
(up to translation) of $e\in\Omega_i$, and we denote $\check\Delta_e(\tau)$
by $\check\Delta_i$ for $e\in\Omega_i$; similarly, $\Delta_f(\tau)$
is independent (up to translation) of $f\in R_i$, and we denote
$\Delta_f(\tau)$ by $\Delta_i$.
\item{(3)} Let $e_1,\ldots,e_p$ be the standard basis of $\ZZ^p$, and
let $\check\Delta(\tau)$ be the convex hull of
\[
\bigcup_{i=1}^p \check\Delta_i\times\{e_i\}
\]
in $(\dualvs{\shQ_{\tau}}\oplus \ZZ^p)\otimes\RR$. Then $\check\Delta(\tau)$
is an elementary simplex of some dimension (i.e. a simplex whose only
integral points are its vertices). Similarly, if $\Delta(\tau)$
is the convex hull of
\[
\bigcup_{i=1}^p\Delta_i\times\{e_i\}
\]
in $(\Lambda_{\tau}\oplus\ZZ^p)\otimes\RR$, then $\Delta(\tau)$
is an elementary simplex.

These last two conditions are vacuous if $p=0$.
\qed
\end{definition}

In case $\tau$ is one-dimensional the polytope $\check\Delta(\tau)$
agrees with the polytope defined in Definition~\ref{NDelta}, and
similarly in the one-codimensional case with $\Delta(\tau)$. Thus
the notation is consistent with previous usage.

\begin{remark}
\label{simplicityfacts}
We give a number of elementary observations related to this
definition, which will hopefully build some feeling for a rather
complex situation. We always suppose that $(B,\P)$ is positive and
simple.

(1) For any $\tau\in\P$ with $0<\dim\tau<\dim B$,
$e_i:\tau\rightarrow \sigma_i\in\P_{\max}$, $f:\omega\rightarrow\tau$
with $\dim\omega=1$, $f\in \Omega_i$, then $n_{\omega}^{e_1\circ
f,e_2\circ f}$ is  independent of $f\in\Omega_i$. Indeed, by
Remark~\ref{polytopefacts}, $n_{\omega}^{e_1\circ f,e_2\circ f}$ is
determined by $\check\Delta_f(\tau)= \check\Delta_i$, independently of
$f\in\Omega_i$. Similarly for $e:\tau\rightarrow \rho$ with
$\codim\rho=1$, $e\in R_i$, $f_i:v_i\rightarrow\tau$,
$m^{\rho}_{e\circ f_1,e\circ f_2}$ is independent of $e$ in $R_i$.

(2) If $f:\omega\to\tau$ is not in $\bigcup_{i=1}^p\Omega_i$ then
$n_{\omega}^{e\circ f,e'\circ f} =0$. In fact, by (1) in the
definition, $n_g =0$ for all $g:\omega\to\rho$ factoring through $f$
for $\dim\rho=n-1$. In general there exists a sequence
$g_1,\ldots,g_m$, $g_i:\omega\to \rho_i$ factoring through $f$, 
$\dim\rho_i=n-1$, such that $T_{\omega}^{e\circ f,e'\circ f}
=T_{g_1}^{\pm 1}\circ\cdots\circ T_{g_m}^{\pm 1}$, hence
$n_{\omega}^{e\circ f,e'\circ f}=\pm n_{g_1}\pm\cdots\pm n_{g_m}=0$.
\pagebreak

(3) By taking $p$ to be as small as possible, we can assume that
$n_{g\circ f}\not=0$, $m_{g\circ f}\not=0$ for any $g\in
R_i,f\in\Omega_i$. Indeed, by (1), if $n_{g\circ f}=0$ for some
$f\in\Omega_i$, in fact $n_{g\circ f}=0$ for all $f\in\Omega_i$. But
if $n_{g\circ f}=0$, we also have $m_{g\circ f}=0$, so $m_{g\circ
f}=0$ for all $g\in R_i$. Thus $n_{g\circ f}=m_{g\circ f}=0$ for all
$g\in R_i,f\in\Omega_i$. So we could just as well leave off $R_i$ and
$\Omega_i$, and condition (1) of Definition~\ref{simplicity} still
holds.

Note also we can always assume $\dim\check\Delta_i>0$ and
$\dim\Delta_i>0$, for if one of these is a point, then
$n_{g\circ f}=0$ for all $g\in R_i, f\in\Omega_i$.

(4) Let $\Xi\subseteq M_{\RR}$ be an elementary simplex, and
$\check\Delta_1,\ldots,\check\Delta_p\subseteq\Xi$ disjoint faces. Let
$T_{\check\Delta_i} \subseteq M_{\RR}$ be the tangent space to
$\check\Delta_i$. It is easy to see that $T_{\check\Delta_1}
+\cdots+T_{\check\Delta_p}$ forms an internal direct sum in $M_{\RR}$
(check this for a standard simplex!) Thus, in the case of
Definition~\ref{simplicity}, if  $T_{\check\Delta_i}
\subseteq\dualvs{\shQ_{\tau,\RR}}$ is the tangent space to
$\check\Delta_i$, then since $\check\Delta_1 \times\{e_1\}, \ldots,
\check\Delta_p\times\{e_p\}$ are disjoint faces of
$\check\Delta(\tau)$, $T_{\check\Delta_1}+ \cdots+T_{\check\Delta_p}$
is an internal direct sum in $\dualvs{\shQ_{\tau,\RR}}$. The dual
statement holds for the $\Delta_i$.

(5) If $p$ is taken to be as small as possible as in (2), then 
$p\le\min(\dim\tau,\codim\tau)$. Indeed, $\dim\dualvs
{\shQ_{\tau,\RR}}=\codim \tau$, and by the assumption $p$ is as small
as possible, $\dim\check\Delta_i\ge 1$ for all $i$ by (2). Thus by (3),
$p\le\codim\tau$. Dually, $p\le\dim\tau$.
\qed
\end{remark}

\begin{example} 
Let us consider what simplicity means for $\dim B\le 3$. It is vacuous if
$\dim B=1$. For $\dim B=2$, let $\tau\in\P$ with $\dim\tau=1$. Then 
\[
\P_1(\tau)=\P_{n-1}(\tau)=\{id:\tau\to\tau\}.
\]
If $n_{id}=0$, i.e. if $\Int(\tau)\cap\Delta=\emptyset$, then we take
$p=0$, and conditions (1)-(3) are vacuous.

If $n_{id}\not=0$, then $n_{id}=n\cdot \check d_{\omega}$ for some
$n\in\ZZ$, $n>0$ by positivity, and we take $\Omega_1=R_1=\{id\}$,
and $\check\Delta(\tau)$ and  $\Delta(\tau)$ are both line segments
of length $n$. (The monodromy operator $T_{id}$ is $\begin{pmatrix}
1&n\\ 0&1\end{pmatrix}$ in an appropriate basis). Then the condition
that $\check\Delta(\tau)\times\{e_1\}$ and
$\Delta(\tau)\times\{e_1\}$ are elementary simplices is just
stating that $n=1$. Thus, if the discriminant locus of $B$ is
straightened, simplicity is equivalent to monodromy about each point
in $\Delta$ being $\begin{pmatrix} 1& 1\\ 0&1\end{pmatrix}$.

In dimension 3, consider first $\dim\tau=1$, so $\P_1(\tau)=\{id\}$.
Then either $n_e=0$ for all $e\in\P_2(\tau)$, in which case we take
$p=0$, and then $\tau\cap\Delta=\emptyset$, or else we have to take
$\Omega_1=\P_1(\tau)$, $R_1=\{e\in\P_2(\tau)|n_e\not=0\}$. In the
latter case, condition (1) is satisfied, and (2) and (3) together
imply $\Delta_e(\tau)$ is a line segment of length 1 for each $e\in
R_1$, while $\check\Delta(\tau)$ is a standard simplex of dimension 1
or 2. If $\dim\check\Delta(\tau)=1$, it is then easy to see that near
$\tau$, $\Delta$ is just a line segment. By the definition of
polyhedral decomposition, $\Delta$ must be contained in a union of
2-cells of $\P$. However, each 2-cell containing a one-dimensional
part of this piece of the discriminant locus must lie in the unique
(local) affine plane containing $\Delta$ which is invariant under
holonomy, and this plane contains $\tau$, as depicted below. If
$\dim\check\Delta(\tau)=2$, $\Delta$ is just a trivalent graph, as
depicted.
\[
\hbox{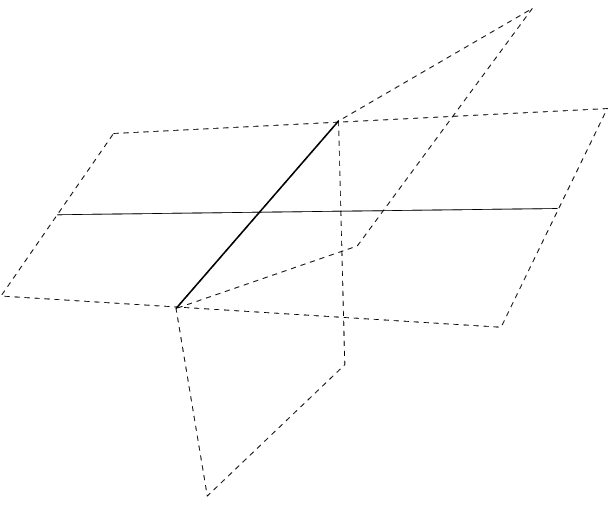
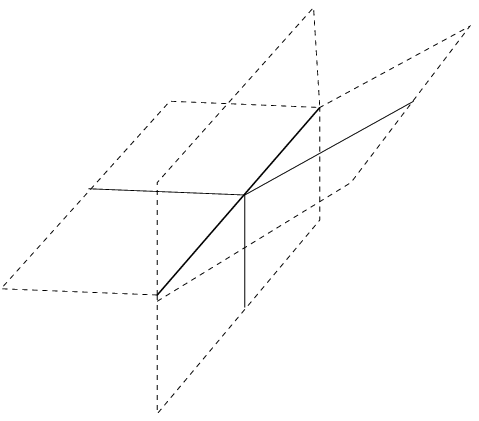}
\]
Furthermore, in a suitable basis with $d_{\tau}=(1,0,0)$, the
monodromy about $\Delta$ in the first case is 
\[
\begin{pmatrix} 1&1&0\\ 0&1&0\\0&0&1\end{pmatrix},
\]
whereas in the second case, the monodromy about the three legs is
\[
\hbox{
$\begin{pmatrix} 1&1&0\\ 0&1&0\\ 0&0&1\end{pmatrix}$,
$\begin{pmatrix} 1&0&1\\ 0&1&0\\ 0&0&1\end{pmatrix}$, and
$\begin{pmatrix} 1&-1&-1\\ 0&1&0\\ 0&0&1\end{pmatrix}$}
\]
respectively.

For $\dim\tau=2$, the picture is very similar, with two pictures
occuring depending on the dimension of $\Delta(\tau)$, namely
\begin{center}
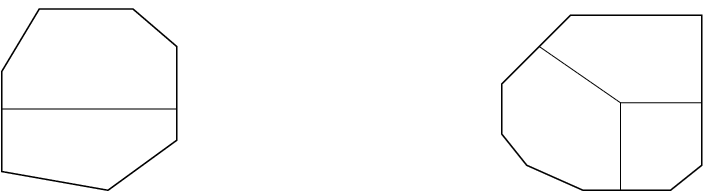
\end{center}
and the monodromy is the transpose of the ones appearing in the
$\dim\tau=1$ case. This agrees with the behaviour of semi-stable,
simple torus fibrations discussed in \cite{TMS}. However, in that
paper, another type of behaviour occured. It was possible to have a
semi-stable, simple torus fibration with a quadrivalent vertex, over
which occured a fibre of type $(1,1)$. Notice that we do not have
quadrivalent vertices in our situation; it appears that the
straightening process would remove such a vertex, replacing it by two
trivalent vertices.

The role of the partition really only becomes important in dimension
4 and higher. Rather than give a complete description in the
four-dimensional case, we just give an example of the behaviour we
expect. We may have, for example, $\tau\in\P$ with $\dim\tau=2$,
$\tau$ a parallelogram, say with vertices $(0,0)$, $(1,0)$, $(a,b)$,
and $(a+1,b)$, where $\gcd(a,b)=1$, $b\not=0$. On the other hand, we
may have $\Sigma_{\tau}$ the two-dimensional fan with $1$-dimensional
cones generated by $(\check b,-\check a)$, $(-\check b,\check a)$,
$(0,1)$ and $(0,-1)$, where again $\gcd(\check a,\check b)=1$ and
$\check b\not=0$. Then it is possible to arrange monodromy so that
the definition of simplicity is satisfied, using the following choice
of $\Omega_i, R_i$.  We take  $\Omega_1$ to be the two parallel edges
of $\tau$ joining $(0,0)$ with $(1,0)$, and $(a,b)$ with $(a+1,b)$,
while $\Omega_2$ consists of the two other edges of $\tau$.
Meanwhile, $R_1$ consists of the two 3-cells corresponding to the
rays of $\Sigma_{\tau}$ generated by $(0,\pm 1)$, while $R_2$
consists of the two 3-cells corresponding to the rays generated by
$\pm (\check b,-\check a)$. Furthermore, we may then  have $\check\Delta_1$
being the line segment with endpoints $(0,0)$ and $(1,0)$; $\check\Delta_2$
the line segment with endpoints $(0,0)$ and $(\check a,\check b)$;
while $\Delta_1$ is the line segment with endpoints $(0,0)$ and
$(1,0)$; $\Delta_2$ is the line segment with endpoints $(0,0)$
and $(a,b)$. If this is the case, then the definition of simplicity
is in fact satisfied. Furthermore, in this case, near $\tau$ the
discriminant locus consists of two two-dimensional manifolds
intersecting at a point. This is in fact the sort of behaviour we
would expect for the discriminant locus of an abelian surface
fibration over a complex surface, so it is quite possible this sort
of picture would appear in the study of degenerations of
hyperk\"ahler manifolds.

Note in this case the monodromy about the two branches of the
discriminant locus takes the form, in an appropriate basis,
\[
\hbox{
$\begin{pmatrix}
1&0&1&0\\
0&1&0&0\\
0&0&1&0\\
0&0&0&1
\end{pmatrix}$
and
$\begin{pmatrix}
1&0&\check a a&\check b a\\
0&1&\check a b&\check b b\\
0&0&1&0\\
0&0&0&1
\end{pmatrix}$
}
\]
\qed
\end{example}

\begin{example}
The $(B,\P)$ constructed in Example~\ref{decompexamp} (2) are
\emph{not} in general simple. To get simple $(B,\P)$ from reflexive
polytopes, one must use a more sophisticated construction. This
construction has in fact been given by Haase and Zharkov in
\cite{HZ}. In the notation of that paper, simplicity is implied by a
sufficiently generic choice of the vectors $\lambda$ and $\nu$, which
will correspond to choosing ample polarizations on maximally
projective crepant resolutions of the toric varieties $\PP_{\Xi}$ and
$\PP_{\Xi^*}$. This is shown in \cite{GBB}, where the more general
case of complete intersections in toric varieties is considered.
So one could think of simplicity as
a sort of ampleness condition.
\end{example}

\section{From polyhedral decompositions to algebraic spaces}
\label{section2}

In this section, given an integral affine manifold with singularities
$B$ with toric polyhedral decomposition $\P$, we will construct
algebraic spaces obtained by gluing toric varieties.

We fix an integral affine manifold with singularities $B$, and a
toric polyhedral decomposition $\P$ once and for all in this
section.  There are two constructions we can make, depending on
whether we view $B$ as an intersection complex or a dual intersection
complex of the central fibre of a degeneration. We shall first consider
the case where we view $B$ as an intersection complex, in which case
the  maximal cells of $\P$ correspond to irreducible components of a
scheme, with a polarisation. This is conceptually simpler
than the dual picture, as we shall see, but should be viewed as
somewhat less important for our point of view. We refer to this
construction as the \emph{cone picture}, and will refer to the dual
construction, in which vertices of $\P$ correspond to irreducible
components, as the \emph{fan picture}.

Whenever $P$ is a monoid, $A$ a ring, we denote by $A[P]$ the monoid
algebra over $A$ defined by
\[
A[P]=\bigoplus_{m\in P} A\cdot z^m
\]
with
\[
z^m\cdot z^{m'}=z^{m+m'}.
\]

\subsection{The cone picture}

\begin{definition}
\label{cones}
If $\sigma\in\P$,  let $\tilde\sigma\subseteq\Lambda_{\sigma,\RR}$ be
as in Definition \ref{normalfan}, well-defined up to translation. Set
\[
C(\sigma)=\{(rp,r)| r\in\RR_{\ge 0}, p\in\tilde\sigma\}
\subseteq \Lambda_{\sigma,\RR}\oplus\RR
\]
be the cone over $\tilde\sigma$. Let $\check P_{\sigma}$ be the
monoid $C(\sigma)\cap (\Lambda_{\sigma}\oplus\ZZ)$. The monoid
$\check P_{\sigma}$ is well-defined up to unique isomorphism:  if
$\tilde\sigma$ is translated, there is a well-defined isomorphism
identifying the two choices of $\check P_{\sigma}$. The projection
$\check P_{\sigma} \to \ZZ$ defines a $\ZZ$-grading on the monoid
ring $R_{\sigma}:=\ZZ[\check P_{\sigma}]$. Define
\[
\check X_{\sigma}:=\Proj R_{\sigma}.
\]
This also defines a line bundle $\O_{\check X_{\sigma}}(1)$ on
$\check X_{\sigma}$.
\qed
\end{definition}

Note that $R_{\sigma}$ and $\check X_{\sigma}$ are well-defined
up to isomorphism. 

We can define a contravariant functor $\check F$ from $\Cat(\P)$ to
the category of graded rings as follows. We take $\check F(\tau)
=R_{\tau}$. If $\tau\subseteq\sigma$, then $e\in\Hom(\tau,\sigma)$
determines a face $\tilde\tau$ of $\tilde\sigma
\subseteq\Lambda_{\sigma,\RR}$, so that the edge of the barycentric
subdivision of $\tilde\sigma$ joining the barycenters of
$\tilde\sigma$ and $\tilde\tau$ maps to $e$. This face $\tilde\tau$
can be identified, by translation, with $\tilde\tau\subseteq 
\Lambda_{\tau,\RR}$. This then gives a well-defined inclusion of
$\check P_{\tau}$ in $\check P_{\sigma}$.  We then define $\check
F(e)$ to be the graded ring homomorphism  $\ZZ[\check
P_{\sigma}]\to\ZZ[\check P_{\tau}]$ defined by 
\[
z^m\longmapsto\begin{cases}0&m\not\in \check P_{\tau}\cr
z^m&m\in \check P_{\tau}.\cr\end{cases}
\]

\begin{example}
\label{squarecone}
Let $B={\RR }^2/\ZZ^2$, with $M=\ZZ^2$, and take as decomposition
$\P$ of $B$ the decomposition with only one maximal cell, coming from
the unit square in $\RR^2$. The set $\P$ consists of one square
$\sigma$, two edges $\tau_1,\tau_2$, and one point $v$.  We see
\[
R_{\sigma}=\ZZ[x,y,z,w]/(xy-zw)
\]
and $\check X_{\sigma}=\PP^1\times\PP^1$, while
$R_{\tau_1}=\ZZ[u,v]$, $R_{\tau_2}=\ZZ[s,t]$ and $R_{v}=\ZZ[q]$. The
two different maps $R_{\sigma}\to R_{\tau_1}$ are given by
\[
x\mapsto u, z\mapsto v, y\mapsto 0, w\mapsto 0
\]
and
\[
x\mapsto 0, z\mapsto 0, y\mapsto v, w\mapsto u
\]
respectively, and the two different maps $R_{\sigma}\to R_{\tau_2}$
are given by
\[
x\mapsto s, w\mapsto t, y\mapsto 0, z\mapsto 0
\]
and
\[
x\mapsto 0, w\mapsto 0, y\mapsto t, z\mapsto s
\]
respectively.
\end{example}

The inverse limit of $\check F$ exists in the category of graded
rings. (See \cite{McLane}, III.3 and III.4 for the notions of inverse and
direct limits of functors.) Calling this limit the graded ring
$R=:\invlim \check F$, we get a scheme $\check
X_0(B,\P)=\Proj R$, along with a line bundle $\O_{\check
X_0(B,\P)}(1)$. If $\P$ is finite, then $R$ is a finitely generated
$\ZZ$-algebra and $\Proj R$ is a projective scheme with $\O_{\check
X_0(B,\P)}(1)$ ample.
The normalization of  $\check X_0(B,\P)$ consists of
$\coprod_{\sigma\in\P_{\max}} \check X_{\sigma}$. However, we would
also like a  twisted version of this construction, so that we get an
entire family of gluings. We will do this over a ring $A$.

\begin{definition}
\label{conegluingdata}
Let $A$ be a ring. Then \emph{gluing data (for the cone picture) for
$\P$ over $A$} are data $\check s=(\check
s_e)_{e\in\coprod_{\tau,\sigma\in\P}\Hom(\tau,\sigma)}$, where if
$e:\tau\to\sigma$ then $\check s_e\in \Hom(\check P_{\tau},\Gm(A))$,
satisfying the conditions
\begin{enumerate}
\item $\check s_{id}(p)=1$ for all $p\in\check P_{\sigma}$, $id:\sigma
\to\sigma$ the identity.
\item
If $e_1:\sigma_1\to\sigma_2$ and $e_2:\sigma_2\to \sigma_3$,
$e_3=e_2\circ e_1$, then $\check s_{e_1}\cdot \check 
s_{e_2}|_{\check P_{\sigma_1}} =\check s_{e_3}$ in $\Hom(\check
P_{\sigma_1},\Gm(A))$.
\end{enumerate}

Given such gluing data, we can then define a contravariant functor
$\check F_{A,\check s}$ from $\Cat(\P)$ to the category of graded
$A$-algebras via
\[
\check F_{A,\check s}(\sigma)=R_{\sigma}\otimes_{\ZZ} A
\]
and if $e:\sigma\to\tau$ then
\[
\check F_{A,\check s}(e)(z^p\otimes 1)=\check F(e)(z^p)\otimes \check
s_e(p).
\]
\qed
\end{definition}

\begin{definition}
Given gluing data $\check s$, we let 
\[
\check X_0(B,\P,\check s)=\Proj \left( \lim_{\longleftarrow} \check
F_{A,\check s}\right).
\]
If $\P$ is finite, this is a 
projective scheme over $\Spec A$, with ample line bundle
$\O_{\check X_0(B,\P,\check s)}(1)$.
\qed
\end{definition}

\begin{example}
\label{squarecone contd}
Continuing with Example~\ref{squarecone}, we can take $A=k$ an
algebraically closed field. Using arbitrary gluing data, we get four
maps 
\[
R_{\sigma}\otimes k\lra R_v\otimes k
\]
given by
\begin{eqnarray*}
x\mapsto \lambda_x q,& y,z,w\mapsto 0\\
y\mapsto \lambda_y q,& x,z,w\mapsto 0\\
z\mapsto \lambda_z q,& x,y,w\mapsto 0\\
w\mapsto \lambda_w q,& x,y,z\mapsto 0
\end{eqnarray*}
The two maps
\[
R_{\tau_1}\otimes k \lra R_v\otimes k
\]
are given by
\begin{eqnarray*}
u\mapsto \lambda_1 q, & v\mapsto 0\\
u\mapsto 0, & v\mapsto \lambda_2 q
\end{eqnarray*}
and the two maps 
\[
R_{\tau_2}\otimes k\lra R_v\otimes k
\]
are given by
\begin{eqnarray*}
s\mapsto \mu_1 q, & t\mapsto 0\\
s\mapsto 0, & t\mapsto \mu_2 q.
\end{eqnarray*}
The cocycle condition of Definition \ref{conegluingdata} then
dictates that the maps $R_{\sigma}\otimes k\to R_{\tau_1}\otimes k$
are
\[
x\mapsto \lambda_x\lambda^{-1}_1 u, z\mapsto \lambda_z\lambda_2^{-1}
v, y\mapsto 0, w\mapsto 0
\]
and
\[
x\mapsto 0, z\mapsto 0, y\mapsto\lambda_y\lambda_2^{-1} v, 
w\mapsto \lambda_w\lambda_1^{-1} u
\]
respectively, and the two different maps $R_{\sigma}\otimes k\to
R_{\tau_2}\otimes k$ are given by
\[
x\mapsto \lambda_x\mu_1^{-1} s, w\mapsto \lambda_w\mu_2^{-1} t,
y\mapsto 0, z\mapsto 0
\]
and
\[
x\mapsto 0, w\mapsto 0, y\mapsto \lambda_y
\mu_2^{-1} t, z\mapsto \lambda_z\mu_1^{-1} s
\]
respectively. Then, as a scheme,  $\check X_0(B,\P,s)$ can be thought
of as $\PP^1\times\PP^1$ with $\{0\}\times\PP^1$ and
$\{\infty\}\times\PP^1$ glued using $(u,v)\mapsto
(\lambda_w\lambda_x^{-1}u,\lambda_y\lambda_z^{-1}v)$, and
$\PP^1\times\{0\}$ and  $\PP^1\times\{\infty\}$ glued using
$(s,t)\mapsto (\lambda_z^{-1}\lambda_xs, \lambda_y^{-1}\lambda_wt)$,
or equivalently, if $\lambda= \lambda_x^{-1}\lambda_y^{-1}
\lambda_z\lambda_w$, we glue using $(u,v) \mapsto (\lambda u,v)$ and
$(s,t)\mapsto (\lambda^{-1} s,t)$. Thus the gluing is not arbitrary,
and as far as the underlying scheme is concerned, we only have a
one-parameter family of gluings. For a more general gluing the ample
line bundle of bidegree $(1,1)$ would not descend. As we will see in
Example~\ref{not a scheme} many gluings do not even give a scheme but
an algebraic space. Also, different choices of gluing may give the
same underlying scheme, but a different choice of line bundle.
\end{example}

\begin{example}
\label{polytope2}
If $B$ is obtained from a reflexive polytope $\Xi\subseteq M_{\RR}$
as in Example~\ref{polytope} and $\P$ is the polyhedral decomposition
whose elements are proper faces of $\Xi$, then we can describe
$\check X_0(B,\P,1)$ (where $1$ denotes the trivial gluing data) as
follows. The polytope $\Xi$ defines a projective toric variety
$(\PP_{\Xi},\O_{\PP_{\Xi}}(1))$, where $\Xi$ is the Newton polytope
of $\O_{\PP_{\Xi}}(1)$. The interior point $0\in\Xi$ corresponds to a
section $s_0\in\Gamma(\PP_{\Xi},\O_{\PP_{\Xi}}(1))$ which vanishes
precisely once on each toric divisor of $\PP_{\Xi}$. The zero locus
$s_0=0$ coincides with $\check X_0(B,\P,1)$, with  polarization
induced by $\O_{\PP_{\Xi}}(1)$. Different choices of gluing data will
give possibly non-isomorphic schemes.
\end{example}

In fact in general this may give too many choices of gluings. 
An arbitrary gluing will give a scheme which cannot appear as the
central fibre of a toric degeneration (see Definition~\ref{toric degen}),
even locally. The
only gluings we are interested in are those for which there is  a
special sort of open covering, which we shall describe in detail in
the fan picture. In theory, one should develop the cone and fan
pictures on completely equal footing, but this would greatly increase
the size of this paper. Furthermore, the two pictures are related by
a discrete Legendre transform, so for us it is not really worth the
effort to develop both pictures in equal generality. So we move on to
the fan picture.

\subsection{The fan picture}

For any rational polyhedral fan $\Sigma$, we will write $X(\Sigma)$
for the toric variety over $\Spec \ZZ$ defined by $\Sigma$.

\begin{definition}
If $\sigma\in\P$, let $X_{\sigma}:=X(\Sigma_{\sigma})$, where
$\Sigma_{\sigma}$ is the fan in $\shQ_{\sigma,\RR}$ defined in 
Definition~\ref{fan2}.
\end{definition}

Note that if $\Sigma$ is a fan in $M_{\RR}$, and $\tau\in\Sigma$ with
$\Sigma(\tau)$ the quotient fan (Definition~\ref{quotfan}), then
there is a canonical closed embedding $X(\Sigma(\tau))\to X(\Sigma)$.
This can be defined as follows: given $\sigma\in\Sigma$, $\sigma$
determines an open subset $X(\sigma)\subseteq X(\Sigma)$ given by
$X(\sigma)= \Spec \ZZ[\dual{\sigma}\cap N]$, where 
\[
\dual{\sigma}=\{n\in N_{\RR}|\langle n,m\rangle\ge 0\quad\forall m\in\sigma\}
\]
as usual. If
$\sigma\supseteq\tau$, then
\[
X((\sigma+\RR\tau)/\RR\tau)=\Spec \ZZ[\dual{
((\sigma+\RR\tau)/\RR\tau)}\cap N].
\]
Here $(\sigma+\RR\tau)/\RR\tau$ is a cone in $M_{\RR}/\RR\tau$, and
we identify the dual of this space with $(\RR\tau)^{\perp}\subseteq
N_{\RR}$. The inclusion 
\[
X((\sigma+\RR\tau)/\RR\tau)\to X(\sigma)
\]
is determined by the ring map given by
\[
z^m\longmapsto \begin{cases} z^m& \hbox{if $m\in \dual{\sigma}\cap 
(\RR\tau)^{\perp}=\dual{((\sigma+\RR\tau)/\RR\tau)}$;}
\cr
0&\hbox{otherwise}.\cr\end{cases}
\]
These inclusions patch on open sets and give a well-defined inclusion
$X(\Sigma(\tau))$ in $X(\Sigma)$.

We can then define a contravariant functor $F:\Cat(\P)\to
\underline{\rm Sch}/\ZZ$ as follows. We set
\[
F(\tau):=X_{\tau}.
\]
For $e\in\Hom(\tau,\sigma)$, we obtain as before
Definition~\ref{quotfan} a surjection $p_e:
\shQ_{\tau,\RR}\to\shQ_{\sigma,\RR}$ which identifies
$\Sigma_{\sigma}$ with $\Sigma_{\tau}(K_e)$, where $K_e$ is the cone
of $\Sigma_{\tau}$ corresponding to $e$. This gives a well-defined
inclusion map
\[
F(e):X_{\sigma}\to X_{\tau}.
\]

\begin{example}
Let $B=\RR^2/(\ZZ(1,2)+\ZZ(2,1))$, with polyhedral decomposition coming
from the quotient of the following periodic decomposition of $\RR^2$:
\begin{center}
\includegraphics{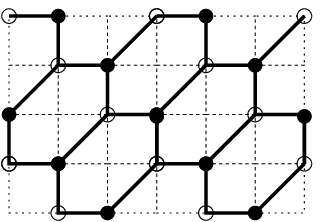}
\end{center}
There are two vertices, $v_1$ and $v_2$, indicated by the two types
of vertices in the figure. It is clear that $X_{v_1}$ and $X_{v_2}$
are both isomorphic to $\PP^2$. There are three edges $\tau_1,\tau_2$
and $\tau_3$ with $X_{\tau_i}=\PP^1$. The maps $X_{\tau_i}\to
X_{v_j}$ identify $X_{\tau_i}$ with one of the various coordinate
lines of $X_{v_j}$, as dictated by the combinatorics of the picture.
There is one maximal cell $\sigma$, $X_{\sigma}=pt$, and three
different maps $X_{\sigma}\to X_{v_1}$ and three different maps 
$X_{\sigma}\to X_{v_2}$, each identifying $X_{\sigma}$ with one of
the three zero-dimensional strata of $X_{v_i}$.
\end{example}

We now proceed to twisting and gluing. Let $S$ be any scheme. We
denote by $\Gm(S)$ the group of invertible regular functions on $S$.

Recalling the open covering $\W$ from Definition~\ref{bary}, we
denote by $H^p(\W,\cdot)$ the \v Cech cohomology groups with respect
to the open covering $\W$. We write 
\[
W_{\tau_1\cdots\tau_p}:=W_{\tau_1}\cap\cdots\cap W_{\tau_p}.
\]

\begin{lemma}
\label{opencovfacts}
$W_{\tau_1\cdots\tau_p}\not=\emptyset$ if and only if
$\tau_1\subseteq\cdots\subseteq\tau_p$ after reordering. Also, if
$\tau\subsetneq\sigma$, then for each $e\in\Hom(\tau,\sigma)$,
$W_{\tau\sigma}$ contains a unique connected component intersecting
$e$, which we write as $W_{e}$. This gives a one-to-one
correspondence between components of $W_{\tau\sigma}$ and
$\Hom(\tau,\sigma)$. More generally, the connected components of
$W_{\tau_1\cdots\tau_p}$ are in one-to-one correspondence with $p-1$
dimensional simplices of $\Bar(\P)$ with vertices at the barycenters
of $\tau_1,\ldots,\tau_p$.
\end{lemma}

\proof The first statement follows immediately from the fact that
$W_{\tau_1\cdots\tau_p}$ is the union of the interiors of all
simplices of $\Bar(\P)$ containing the barycenters of
$\tau_1,\ldots,\tau_p$. To prove the second statement, let $W_{e}$ be
the union of the interiors of all simplices of $\Bar(\P)$ containing
$e$. It is easy to see each $W_{e}$ is connected, and an open subset
of $B$. Furthermore, any two simplices of $\Bar(\P)$ have disjoint
interior, so $W_{\tau\sigma}=\coprod_{e\in\Hom(\tau,\sigma)} W_e$.
Thus $\{W_{e}|e\in\Hom(\tau,\sigma)\}$ is the set of connected
components of $W_{\tau\sigma}$. The general case is similar. \qed

\begin{definition}
Let $S$ be a scheme. We call a \v Cech 1-cocycle
$s=(s_{e})_{e\in\coprod \Hom(\tau,\sigma)}$ for the sheaf
$\shQ_{\P}\otimes \Gm(S)$ on $B$ with respect to the open covering
$\W$ \emph{closed gluing data (for the fan picture) for $\P$ over
$S$}. Here $s_e\in \Gamma(W_e, \shQ_{\P}\otimes
\Gm(S))=\shQ_{\sigma}\otimes \Gm(S)$ for $e:\tau\rightarrow\sigma$.
\end{definition}

We sometimes write $s_{\tau\sigma}$ instead of $s_e$ when $e:\tau\to
\sigma$ is clear from the context.

Since $\Sigma_{\sigma}$ is a fan living in $\shQ_{\sigma,\RR}$, 
$\shQ_{\sigma}\otimes\Gm(S)$ acts on $X_{\sigma}\times S$, and hence
$s_{e}$ gives an automorphism
\[
s_{e}:X_{\sigma}\times S\to X_{\sigma}\times S.
\]
Thus we have

\begin{definition}
\label{functor}
For closed gluing data $s$ for $\P$, we define a functor
\[
F_{S,s}:\Cat(\P)\to \underline{\rm Sch}/S
\]
by 
\[
F_{S,s}(\tau):=X_{\tau}\times S
\]
and for $e\in\Hom(\tau,\sigma)$,
\[
F_{S,s}(e):=(F(e)\times id)\circ s_{e}: X_{\sigma}\times S
\to X_{\tau}\times S.
\]
\qed
\end{definition}

Let $X$ be an algebraic space, $A$ a closed subspace and $A\to \bar
A$ a finite morphism. Then the fibred coproduct $X\amalg_A \bar A$
in the category of algebraic spaces is an algebraic space. This is a
corollary of Artin's theory of the existence of formal modifications,
see \cite{artin}, Theorem~6.1. In particular, the identification of
two algebraic spaces along a closed subspace is again an algebraic
space. Applying this fact inductively it is possible to construct the
direct limit
\[
X_0 (B,\P,s)=\lim_{\longrightarrow} F_{S,s}
\]
as an algebraic space locally of finite type over $S$. In other
words, for every $\tau\in \P$ there is a morphism $q_\tau:
X_\tau\times S\to X_0 (B,\P,s)$ and for every $e\in\Hom (\tau,
\sigma)$ we have $q_\sigma= q_\tau \circ F_{S,s}(e)$; moreover,
$X_0(B,\P,s)$ is a universal object in the category of algebraic
spaces with this property. The constructed space is of finite type
over $S$ whenever $\P$ is finite. 

Instead of providing full details for the use of Artin's result, we
are going to construct $X_0(B,\P,s)$ directly as a quotient of a
scheme by an \'etale equivalence relation. In other words, instead of
constructing $X_0(B,\P,s)$ by gluing proper toric varieties together
along closed substrata, we glue together a collection of open sets to
obtain $X_0(B,\P,s)$. However, this cannot be done for all choices of
$s$ (see Remark~\ref{openclosedgluing}), and in fact those choices of
$s$ for which this cannot be done are irrelevant for our theory. The
existence of an open covering of $X_0(B,\P,s)$ by some standard open
sets is crucial. Thus we will only discuss the open gluings, and this
topic occupies the rest of this section. 
\medskip

First let us explain what these standard open sets are.

\begin{definition}
\label{fanopen}
Let $\sigma\in\P_{\max}$ be a maximal cell, $y\in \Int(\sigma)$ and
$\tilde\sigma\subseteq \Lambda_{\RR,y}$ as in Definition~\ref{cones}.
We denote by $\dual{C(\sigma)}$ the dual cone to $C(\sigma)$ in 
$\check\Lambda_{\RR,y}\oplus\RR$, and set
\[
P_{\sigma}=\dual{C(\sigma)}\cap (\check\Lambda_y\oplus\ZZ),
\]
so that the affine toric variety defined by the fan of faces of
$C(\sigma)$ is 
\[
U(\sigma):=\Spec \ZZ[P_{\sigma}].
\]
In addition, the projection $\Lambda_v\oplus \ZZ \to \ZZ$ defines an
element $\rho_\sigma\in P_\sigma\subseteq
\check\Lambda_v\oplus\ZZ$, hence a monomial  regular function
$z^{\rho_{\sigma}}$ on $U(\sigma)$. Define $V(\sigma)$ to be the zero
scheme of this regular function.  Note that by construction
$z^{\rho_{\sigma}}$ vanishes precisely once on each toric divisor of
$U(\sigma)$, so $V(\sigma)$ is the reduced union of all toric
divisors in $U(\sigma)$. By the uniqueness property of $\tilde\sigma$
as in Definition~\ref{cones},  $U(\sigma)$ and $V(\sigma)$ are unique
up to unique toric isomorphism.
\qed
\end{definition}

$V(\sigma)$ can be described in terms of a monoid ring also:

\begin{definition}
Let $\sigma$ be a (not necessarily strictly convex) cone in $N_{\RR}$,
for $N$ a lattice, and let $P$ be the monoid $\sigma\cap N$. We
write
\[
\partial P:=(\partial\sigma\cap N)\cup\{\infty\},
\]
with addition
\[
p+q=\begin{cases}
p+q&\hbox{if $p,q,p+q\in\partial\sigma$}\\
\infty&\hbox{otherwise.}
\end{cases}
\]
We define the monoid ring $\ZZ[\partial P]$ by using the convention that
$z^{\infty}=0$.
\qed
\end{definition}

It is easy to see that in the notation of Definition~\ref{fanopen}, 
$V(\sigma)=\Spec\ZZ[\partial P_{\sigma}]$.

The collection $\{V(\sigma)\times S|\sigma\in \P_{\max}\}$ will 
serve as our \'etale cover of $X_0(B,\P,s)$. The next task is to
write down the \'etale equivalence relation on $\coprod_{\sigma\in
\P_{\max}} V(\sigma)\times S$. This will be a closed subspace
\[
\foR\subset
\Big(\coprod_{\sigma_1\in \P_{\max}} \!\!\!V(\sigma_1)\times S\Big)
\times_S
\Big(\coprod_{\sigma_2\in \P_{\max}} \!\!\!V(\sigma_2)\times S\Big)
=\coprod_{\sigma_1,\sigma_2\in \P_{\max}} \!\!\!
V(\sigma_1) \times V(\sigma_2)\times S.
\]
\begin{remark}
The guiding principle for the construction of this equivalence
relation comes from the relation of the $V(\sigma)$ with the proper
schemes $X_v$ from before. A vertex $v\in\tilde\sigma$ gives the
$n$-dimensional face $(v,1)^\perp\cap \dual{C(\sigma)}$ of
$\dual{C(\sigma)}$. The integral points of this face form a  monoid
$P_{\sigma,v}$, and $\Spec \ZZ[P_{\sigma,v}]$ is one of the
irreducible components of $V(\sigma)$. An alternative way to see this
component is as follows: for $y\in \Int(\sigma)$, identify
$\Lambda_y$ with $\Lambda_v$ by parallel transport in $\sigma$ and
project
\[
(v,1)^\perp\cap \dual{C(\sigma)}\subset \dualvs{(\Lambda_{\RR,v}\oplus
\RR)} = \check\Lambda_{\RR,v}\oplus\RR
\]
onto the first factor. This identifies  $(v,1)^\perp\cap
\dual{C(\sigma)}$ with the dual of the tangent wedge to $\tilde\sigma$
at $v$. Thus the irreducible component $\Spec \ZZ[P_{\sigma,v}]$ is
naturally identified with the affine open set of the toric variety
$X_v=X(\Sigma_v)$ corresponding to the cone of $\Sigma_v$ determined
by $v\in\tilde\sigma$.

The problem in writing this down is that the gluing of open subsets
for maximal cells $\sigma_1,\sigma_2$ intersecting in more than one
point does not generally come from gluing of open sets in the ambient
spaces $U(\sigma_1), U(\sigma_2)$. Rather each vertex in the
intersection prescribes one gluing of an open subset of
$U(\sigma_i)$; the restrictive form of the monodromy assures
compatibility of the gluing only after restriction to the divisor
$V(\sigma_i)\subset U(\sigma_i)$.
\qed
\end{remark}

\begin{construction}
\label{basicgluing}
Let us first discuss the case of only two maximal cells $\sigma_1,
\sigma_2$ without self-intersections and intersecting in one cell
$\tau=\sigma_1\cap \sigma_2$. Also let us assume the gluing parameter
$s$ to be trivial for the time being. Choose vertices $v_i\in
\sigma_i$ and embeddings $\sigma_i= \tilde\sigma_i \subset
\Lambda_{\RR,v_i}$. We are going to use the notations
\[
M_i=\Lambda_{v_i},\quad
N_i=M_i^*,
\quad P_i:=P_{\sigma_i}=\dual{C(\sigma_i)}\cap (N_i\oplus \ZZ),
\quad \rho_i:=\rho_{\sigma_i}=(0,1)\in N_i\oplus\ZZ.
\]
For each vertex $w\in \tau$, there are linear identifications
\[
\psi_w: M_2\lra M_1, \quad \psi_w^t: N_1\lra N_2
\]
given by parallel transport from $v_2$ through $w$ to $v_1$. We will
distinguish $\tau$ as a face of $\sigma_1$ or a face of $\sigma_2$ by
using the notation $\tau_i\subset \sigma_i$. 

We wish to identify open subsets $V(\tau_i)\subseteq V(\sigma_i)$
and a natural isomorphism between $V(\tau_1)$ and $V(\tau_2)$.

To begin, it is worth understanding $V(\sigma_i)$ a bit better.
The monoid $\partial P_i$ has an interpretation in terms of
fans. Let $\check \Sigma_i$ be the normal fan to $\sigma_i$, living in
$N_i\otimes\RR$. Then by projecting $N_i\oplus\ZZ\rightarrow N_i$,
we can in fact identify $\partial P_i$ with the monoid $N_i\cup\{\infty\}$
with addition
\[
p+q=\left\{ \begin{array}{ll}
p+q&\text{if $p,q$ are in a common cone of 
$\check \Sigma_i$}\\
\infty& \text{otherwise.}
\end{array}\right.
\]
Note that for any face $\omega\subseteq\sigma_i$, there
is a corresponding face $\check\omega$ of $\dual{C(\sigma_i)}$
dual to $C(\omega)$, i.e.
\[
\check\omega=\{n\in \dual{C(\sigma_i)}| \hbox{$\langle n,(m,1)\rangle =0$
for all $m\in\omega$}\}.
\]
We will also denote by $\check\omega$ the projection of $\check\omega$
to $N_i$; this is the cone of $\check\Sigma_i$ denoted $\check
K_\omega$ in Definition~\ref{normalfan}. In particular, $\Spec
\ZZ[\check\omega\cap N_i]$ is then a closed subscheme of
$V(\sigma_i)$, a toric stratum corresponding to $\omega$.

In future $\check\omega$ may refer both to the face of
$C(\sigma_i)^\vee\subset N_i\oplus \ZZ$ and to its projection to
$N_i$. It will always be clear from the context which alternative
applies.

Now define
\[
U(\tau_i)= \Spec \ZZ[Q_i] 
\]
with
\[
Q_i:= \dual{C(\tau_i)} \cap (N_i\oplus\ZZ)
= P_i+ \big(\RR\check\tau_i \cap (N_i\oplus \ZZ)\big).
\]
As an abstract toric variety $U(\tau_i)$ is isomorphic to
\[
\Spec \ZZ[
\Hom(C(\tau_i)\cap (M_i\oplus\ZZ), \NN)]\times \GG_m^{\codim
\tau_i};
\]
in particular, the isomorphism class as an abstract toric variety
depends only on $\tau$.
It is easy to see
that $\ZZ[\partial Q_i]\cong\ZZ[Q_i]/(z^{\rho_i})$. Thus in analogy
with the definition for maximal cones we now define
\[
V(\tau_i)=\Spec \ZZ[\partial Q_i].
\]
The inclusion $P_i\subseteq Q_i$ induces the inclusion $U(\tau_i)
\subseteq U(\sigma_i)$, and this restricts to an inclusion
$V(\tau_i)\subseteq V(\sigma_i)$.

It is worth underlining a very important point here. We were able to
define $V(\sigma)$ for $\sigma$ maximal in a completely canonical way.
On the other hand, for $\tau$ non-maximal, we have not defined
$V(\tau)$ but only $V(\tau_i)$ as a subset of the canonically defined
$V(\sigma_i)$. The point is that $V(\tau)$ is not well-defined up to a
unique isomorphism because of holonomy in a neighbourhood of
$\Int(\tau)$, i.e. there is no canonical affine structure on
$\sigma_1\cup\sigma_2$. It is only when this holonomy is trivial and
the affine structure of $B$ is defined in a neighbourhood of
$Int(\tau)$ that $V(\tau)$ becomes well-defined. It is this
indeterminacy of $V(\tau)$ when there is holonomy  which causes the
gluing between $V(\tau_1)$ and $V(\tau_2)$ below to be non-trivial.
Believing in the independent existence of $V(\tau)$ is an easy trap to
fall into, one we have made countless times.
In particular, there is no canonical
isomorphism between $U(\tau_1)\subset U(\sigma_1)$ and
$U(\tau_2)\subset U(\sigma_2)$, and we cannot glue these schemes,
even though $V(\tau_1)$ and $V(\tau_2)$ will glue.

The monoid $\partial Q_i$ also has an interpretation in terms of
fans as before. By projecting $\partial \dual{C(\tau_i)}$
to $N_i$, we obtain an identification
$\partial Q_i= N_i\cup \{\infty\}$ with addition
\[
p+q=\left\{ \begin{array}{ll}
p+q&\text{if $p,q$ are in a common cone of $\check \tau_i^{-1}
\check \Sigma_i$}\\
\infty& \text{otherwise.}
\end{array}\right.
\]
Here $\check\tau_i^{-1}\check\Sigma_i$ is the localisation of the fan
$\check\Sigma_i$ at $\check\tau_i$ as defined in Definition~\ref{quotfan}.

Here is the key point of the construction. We want to identify the
sets $V(\tau_1)$ and $V(\tau_2)$ in order to glue together
$V(\sigma_1)$ and $V(\sigma_2)$. This is done by specifying an
isomorphism of monoids
\[
\phi:\partial Q_1\lra \partial Q_2.
\]
Via the embedding $\partial Q_i\setminus \{\infty\} \subset N_i$ such
a morphism is equivalent to a piecewise linear map $\phi: N_1\to N_2$
linear on cones of $\check\tau_1^{-1}\check\Sigma_1$, mapping
$\check\tau_1^{-1} \check\Sigma_1$ to $\check\tau_2^{-1}
\check\Sigma_2$. To define $\phi$ on the maximal cone in $\partial
Q_1$ corresponding to a vertex $w\in\tau$ use the linear map
${}^t\psi_w: N_1\to N_2$ defined by parallel transport through $w$ as
above. Parallel transport identifies the $\sigma_i$ with maximal
cells in the induced polyhedral decomposition of
$R_w\subseteq\Lambda_{\RR,w}$, whose intersection maps to $\tau$.
Since $\check\tau_i ^{-1}\check\Sigma_i$ depends only on
$\tau_i\subset M_i$ this provides the desired identification of fans.

We need to check that our definition of $\phi$ is well-defined on
intersections of cones. Let $w_1,w_2$ be vertices of $\tau$ and
$\check w_i\subset N_1$ the corresponding maximal cones of
$\check\tau_1^{-1} \check\Sigma_1$. Then $\check\omega:= \check
w_1\cap \check w_2$ is dual to a subface $\omega\subset\tau$
containing $w_1,w_2$. For $n\in \check w_1\cap \check w_2$ and $m\in
M_1$ we compute
\[
\langle {}^t\psi_{w_2}^{-1}\,{}^t\psi_{w_1}(n),m\rangle
=\langle n,\psi_{w_1}\psi_{w_2}^{-1}(m)\rangle
=\langle n,T_{\gamma}(m)\rangle,
\]
where $T_\gamma$ is parallel transport along the loop passing from
$v_1$ to $w_2$ inside $\sigma_1$, then to $w_1$ inside $\sigma_2$ and
finally back to $v_1$. By Proposition~\ref{monodromy1}, $T_\gamma(m)$
is congruent to $m$ modulo the tangent space $\check \omega^\perp$ of
$\omega$. Hence $\langle n, T_\gamma(m)\rangle = \langle n,m\rangle$
for all $m\in M_1$ and
\[
{}^t\psi_{w_2}(n)={}^t\psi_{w_1}(n).
\]
Thus $\phi$ is well-defined. This gives our gluing isomorphism
\[
\Phi_{\sigma_1\sigma_2}: V(\tau_2)\lra V(\tau_1).
\]
\qed
\end{construction}

\begin{example}
\label{twotriangles}
We illustrate the definition of $\phi$ using Example \ref{twodimaffine},
(2).
Let $\sigma_1$ and $\sigma_2$ denote the left and right triangles
in Example \ref{twodimaffine}, (2) respectively, 
in either the left or right-hand pictures, and 
let $\tau=\sigma_1\cap\sigma_2$. Then using the
left-hand picture, $\dual{C(\sigma_1)}$ is the cone in $N_\RR\oplus
\RR$ generated by 
\[
(-1,0,0),(0,1,0),(1,-1,1)
\]
and $\dual{C(\sigma_2)}$ is the cone generated by 
\[
(1,0,0), (0,1,0),(-1,-1,1).
\]
In both instances $\check\tau_i$ is the cone generated by the first
vector. Then $\check\Sigma_1$ and $\check\Sigma_2$ are 
\begin{center}
\includegraphics{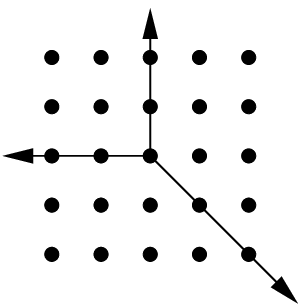}$\quad\quad\quad$
\includegraphics{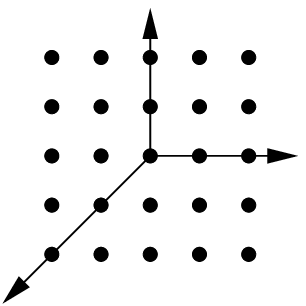}
\end{center}
obtained by forgetting the last coordinate. Indeed, this is just the
normal fan of $\sigma_i$. Finally, localizing at $\check\tau$ gives
$\partial Q_1$ and $\partial Q_2$ determined by the picture
\begin{center}
\includegraphics{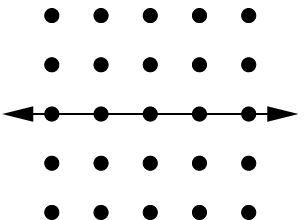}
\end{center}
The upper half plane corresponds to the vertex $(0,0)$ and the lower
half-plane corresponds to the vertex $(0,1)$. Now given that we have
been using the left-hand chart in Example~\ref{twodimaffine}, (2),
the map $\psi_{(0,0)}: M_2\to M_1$
should be viewed as the identity, while $\psi_{(0,1)}:M_2\to M_1$ is
given by the matrix $\begin{pmatrix}1&0\cr 1&1\end{pmatrix}$.  This
is the transformation which takes $\sigma_2$ in the left-hand chart
to $\sigma_2$ in the right-hand chart. Thus if we denote by $x$ the
element $(1,0)\in\partial Q_i$, $y$ the element $(0,1)$, $z$ the
element $(0,-1)$, then
\[
\phi:\partial Q_1\to\partial Q_2
\]
takes $x\mapsto x$, $y\mapsto y$ and $z\mapsto x^{-1}z$ since
\[
\psi^t_{(0,1)}\begin{pmatrix} 0\\ -1\end{pmatrix}
=\begin{pmatrix}-1\\-1\end{pmatrix}=-\begin{pmatrix}1\\0
\end{pmatrix}+\begin{pmatrix}0\\-1\end{pmatrix}.
\]
The map $\phi$ induces the gluing map
\[
\Phi_{\sigma_1\sigma_2}: 
\Spec \ZZ[x,x^{-1},y,z]/(yz)\to \Spec \ZZ[x,x^{-1},y,z]/(yz).
\]
\qed
\end{example}

\begin{example}
\label{twopyramids}
Next we do a three-dimensional example. Here again we use two
different cuts, and the shaded area shows where we make cuts.
\begin{center}
\input{sect4fig5.pstex_t}
\end{center}
This time we just draw the fan $\check\Sigma_i(\check\tau_i)$ (the
fan $\check\tau_i^{-1}\check\Sigma_i$ is the pullback of this to
$M_i$) living in $M_i/(1,0,0)\ZZ$:
\begin{center}
\input{sect4fig6.pstex_t}
\end{center}
with the first through fourth quadrants corresponding to the vertices
$(0,0,0)$, $(0,1,0)$, $(0,1,1)$ and $(0,0,1)$ of $\tau$ respectively.
Then we can write 
\[
\ZZ[\partial Q_i]=\ZZ[u,u^{-1},x,y,z,w]/(yw,xz),
\]
where $x,y,z,w$ are primitive generators of the rays as labelled in
the diagram and $u$ corresponds to $(1,0,0)\in M_i$. We take the maps
$\psi_{(0,0,0)}$ and $\psi_{(0,0,1)}$ to be the identity, and then
$\psi_{(0,1,0)}=\psi_{(0,1,1)}$ is given by the matrix
\[
\begin{pmatrix}1&0&0\cr 1&1&0\cr 0&0&1\end{pmatrix}.
\]
The gluing map $\phi: \ZZ[\partial Q_1]\to \ZZ[\partial Q_2]$ is
\[
u\mapsto u, x\mapsto x, y\mapsto y, z\mapsto z, w\mapsto u^{-1}w.
\]
\qed
\end{example}
\medskip

\noindent
Let us next discuss twisted gluings. Let $S$ be an arbitrary scheme,
and we wish to glue the $V(\sigma_i)\times S$ along $V(\tau_i)\times
S$ using a twisting of $\Phi_{\sigma_1\sigma_2}\times\id_S$.
Twistings are given as follows. Let 
\[
s_i:\partial Q_i\setminus\{\infty\}\to \Gm(S)
\]
be a map satisfying
\[
\hbox{$s_i(p+q)=s_i(p)\cdot s_i(q)$ if $p+q\not=\infty$.}
\]
We can view this as a piecewise multiplicative function on $N_i$.
Then $s_i$ defines a map $s_i:V(\tau_i)\times S\to V(\tau_i)\times S$
induced by 
\[
z^p\longmapsto s_i(p) z^p.
\]
We can then glue $V(\tau_1)\times S$ and $V(\tau_2)\times S$ via
$s_1^{-1}\circ(\Phi_{\sigma_1\sigma_2}\times\id_S)\circ s_2$.
\vspace{3ex}

So far the discussion was essentially local, focusing on only two
maximal cells intersecting in one common face. In general the part
$\foR_{\sigma_1\sigma_2}$ of the equivalence relation in
$V(\sigma_2)\times V(\sigma_1)$ may have more than one component. 
For any $\sigma\in\P_{\max}$, and a vertex $v\in\sigma$, we have
$R_v\subseteq\Lambda_{\RR,v}$, $\P_v$ and $\exp_v$ as usual, and
$\tilde\sigma\in\P_v$ mapping to $\sigma$ via $\exp_v$. Denoting this
map $\pi:\tilde\sigma \to\sigma$, it is an integral affine
isomorphism in the interior of $\tilde\sigma$ and is finite-to-one on
the  boundary. This map is unique up to integral affine isomorphisms
between different choices of $\tilde\sigma$.

Now given $\sigma_1,\sigma_2\in\P_{\max}$, let
$\tilde\sigma_i\subseteq M_i\otimes\RR$,
$M_i=\Lambda_{\sigma_i}$, $\pi_i: \tilde\sigma_i\to\sigma_i$ be
as above. Consider the fibred product $\tilde\sigma_1\times_B
\tilde\sigma_2$. Each connected component of this space will give one
component of $\foR_{\sigma_1\sigma_2}$.

\begin{lemma}\label{sigma_1 times_B sigma_2}
Let $\sigma_1,\sigma_2\in \P_{\max}$. Then the connected components
of $\tilde\sigma_1\times_B \tilde\sigma_2$ are graphs of integral
affine isomorphisms  $\tilde\tau_1\to \tilde\tau_2$ over $B$ for
faces $\tilde\tau_i\subset \tilde\sigma_i$. In particular, the
$\tilde\tau_i$ cover the same $\tau\subset\sigma_1 \cap \sigma_2$.
Exactly those isomorphisms occur that do not extend to larger subsets
of $\tilde\sigma_i$. 
\end{lemma}
\proof
We have disjoint decompositions $\tilde\sigma_i= \coprod_{\tilde\tau
\subset\tilde\sigma_i} \Int(\tilde\tau)$. By the definition of
polyhedral decompositions $\pi_i|_{\Int(\tilde\tau)}$ is a
homeomorphism onto its image. For any pair $\tilde\tau_1,
\tilde\tau_2$ covering the same $\tau \subset \sigma_1 \cap \sigma_2$
the composition $(\pi_2|_{\Int(\tilde\tau_2)})^{-1}\circ
\pi_1|_{\Int(\tilde\tau_1)}$ is affine. Hence it extends to an affine
linear map $\lambda$ from $\RR \tilde\tau_1\subset M_1\otimes\RR$ to
$\RR\tilde\tau_2 \subset M_2\otimes\RR$. This shows already that
$\tilde\sigma_1 \times_B \tilde\sigma_2$ has a decomposition into
polyhedral subsets that are graphs of affine transformations
$\lambda:\tilde\tau_1\to \tilde\tau_2$ commuting with $\pi_i$.

To see how the strata fit together choose a vertex $w_1\in
\tilde\tau_1$ and put $w_2= \lambda(w_1)$. Let $R_w\subseteq
\Lambda_{\RR,w}$, $\P_w$, $\exp_w$ be a chart for $(B,\P)$ at
$w=\pi_1(w_1) =\pi_2(w_2)$. By the definition of polyhedral
decomposition there exist unique affine embeddings
\[
\tilde\sigma_i\lra R_w
\]
respecting the polyhedral decompositions and mapping $w_i$ to the
origin. This induces integral affine maps $M_i\otimes\RR\to
\Lambda_{\RR,w}$ for $i=1,2$, and the composition
\[
M_1\otimes\RR\lra \Lambda_{\RR,w}\lra M_2\otimes\RR
\]
restricts to $\lambda$ on $\tilde\tau_1$. This shows that
$\tilde\sigma_1\times_B \tilde\sigma_2$ has a decomposition as a union
of polyhedra with each maximal cell the graph of an integral affine
transformation $\lambda:\tilde\tau_1\to\tilde\tau_2$ that commutes
with $\pi_i$ and does not extend to a map between larger faces of
$\tilde\sigma_i$. Moreover, we see that locally
$\tilde\sigma_1\times_B \tilde\sigma_2$ is isomorphic to the
intersection of two (possibly equal) cells of the polyhedral
decomposition of $R_w$. Hence the maximal cells in the polyhedral
decomposition of $\tilde\sigma_1 \times_B\tilde\sigma_2$ are disjoint.
Therefore we have the claimed decomposition of $\tilde\sigma_1
\times_B \tilde\sigma_2$ into pairwise disjoint polyhedra.
\qed

\begin{remark}
\label{shyalgebraic}
The reader shy of algebraic spaces may read the remainder of this
section with the restriction on $\P$ that the maps
$\pi:\tilde\tau\rightarrow\tau \subseteq B$ are injective for all
$\tau\in\P$. It is then the case that the algebraic space constructed
below only involves gluings of affine schemes along Zariski open
subsets rather than \'etale open subsets, and hence produces a
scheme.
\qed
\end{remark}

Let $\tilde\tau_1\to \tilde\tau_2$ be an integral affine
transformation whose graph is a connected component of
$\tilde\sigma_1\times_B \tilde\sigma_2$, as established by the lemma.
Except for the trivial case $\tilde\tau_1=\tilde\sigma_1
=\tilde\sigma_2= \tilde \tau_2$ there is a unique edge of the
barycentric subdivision of $\tilde\sigma_i$ joining the barycenter of
$\tilde\tau_i$ with the barycenter of $\tilde\sigma_i$. This edge
determines $\tilde\tau_i$ uniquely.  It also corresponds uniquely to
an edge of the barycentric subdivision of $\sigma_i$. We may thus
index the components of $\tilde\sigma_1\times_B \tilde\sigma_2$ by
$(e_1,e_2)$ with $e_i\in \Hom(\tau,\sigma_i)$ and $\tau\subset
\sigma_1\cap\sigma_2$. The mentioned trivial case is covered by
$(\id,\id)$. Not all pairs occur, but only those that are
maximal for the partial order defined by
\[
(e_1,e_2)<(f_1,f_2)\quad\Longleftrightarrow\quad
\exists\, h\in\Hom(\tau,\omega): e_i=f_i\circ h.
\]
\begin{example}
Consider $B=\RR^2/\ZZ^2$ with the affine structure induced from
$\RR^2$, and the polyhedral decomposition with one maximal cell
$\sigma$ covered by $\tilde\sigma= [0,1]\times[0,1] \subset \RR^2$. The
following figure shows the fibred product of $\tilde\sigma$ with itself
over $B$. 
\medskip
\begin{center}
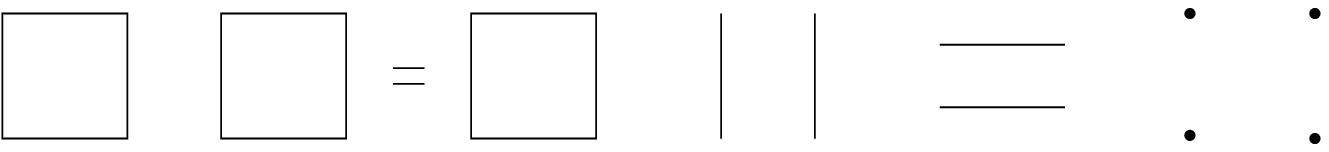
\end{center}
\end{example}
There are $4\cdot 4=16$ vertices in $\tilde\sigma\times_B
\tilde\sigma$; a pair of vertices $(v_i,v_j)$ is isolated in
$\tilde\sigma\times_B \tilde\sigma$ iff $v_i$ and $v_j$ do not belong
to a proper face of $\tilde\sigma$. This is the case iff $i-j\equiv2$
modulo $4$. Similarly, the pairs $(\tau_i,\tau_i)$ and
$(\omega_i,\omega_i)$ are not maximal: the identity isomorphisms
$\id_{\tau_i}$ and $\id_{\omega_i}$ extend to the identity on
$\tilde\sigma$.
\qed

\begin{definition}
\label{generalgluing1}
We have thus indexed the connected components of
$\tilde\sigma_1\times_B \tilde\sigma_2$ by a subset of pairs
$(e_1,e_2)\in \Hom(\tau,\sigma_1) \times \Hom(\tau,\sigma_2)$ for
$\tau\subset \sigma_1\cap\sigma_2$. By Lemma~\ref{sigma_1 times_B
sigma_2} it makes sense to call pairs corresponding to connected
components \emph{maximal}.
\end{definition}

We will need the following:

\begin{definition}
\label{Ve}
Let $e:\tau\rightarrow\sigma$. Then $e$ determines a cone $K_e$ in
$\Sigma_{\tau}$, and hence determines an open subset  $V_e\subseteq
X_{\tau}$ with $V_e=\Spec \ZZ[\dual{K_e}\cap \dualvs{\shQ_{\tau}}]$.
The embedding $V_e\subseteq X_{\tau}$ is canonical.
\end{definition}

For any maximal pair $(e_1,e_2)$ let $\tilde\tau_i \subseteq
\tilde\sigma_i$ be the face given by $e_i:\tau\to\sigma_i$. For any
common choice of vertex $f:w\to \tau$ of $\tilde\tau_i$ the
compositions $e_i\circ f$ define a path joining the interior of
$\sigma_2$ to the interior of $\sigma_1$ passing through $w$. Parallel
transport along this path defines an affine map
$\psi_w:\Lambda_{\sigma_2} \rightarrow\Lambda_{\sigma_1}$; note this
depends on $(e_1,e_2)$ and not just $\sigma_1$, $\sigma_2$ and $w$. As
in the case covered in Construction~\ref{basicgluing} without
self-intersections this induces an isomorphism between the irreducible
component $V_{e_1\circ f}\subset V(\tau_1)$ and $V_{e_2\circ f}\subset
V(\tau_2)$. The monodromy argument as in
Construction~\ref{basicgluing} shows that these isomorphisms are
compatible on intersections. We obtain an isomorphism
\[
\Phi_{e_1 e_2}: V(\tau_2)\lra V(\tau_1),
\]
of the open subschemes $V(\tau_i)\subset V(\sigma_i)$, defined above,
which is toric on each irreducible component.

\begin{construction}
\label{generalgluing2}
We next want to twist the isomorphisms $\Phi_{e_1e_2}$ over an
arbitrary scheme $S$. We need to describe automorphisms of the
$V(\tau)\times S$ slightly more abstractly. Let $\tau\in\P$ be any
cell, and let $\tilde\tau \subseteq \Lambda_{\tau,\RR}$ be as in
Definition~\ref{normalfan}, with $\pi:\tilde\tau\to\tau$ the map to
$B$.

Consider the sheaf $\pi^*(\shQ_{\P}\otimes\Gm(S))$. This can be
viewed as follows. If we choose any $\sigma\in\P_{\max}$ with
$e:\tau\to \sigma$ a morphism, then we can identify $\tilde\tau$ with
a face of $\tilde\sigma\subseteq\Lambda_{\sigma,\RR}=\Lambda_{\RR,z}$
for any $z\in \Int(\sigma)$. Under this identification, at a point
$y\in \Int(\tilde\omega)$, for $\tilde\omega\subseteq\tilde\tau$ a
face, the stalk of  $\pi^*(\shQ_{\P}\otimes\Gm(S))$ at $y$ is
\[
(\Lambda_z/\Lambda_\omega)\otimes\Gm(S)
=\Hom(\Lambda_\omega^{\perp},\Gm(S)).
\]
Thus an element $s\in\Gamma(\tilde\tau,\pi^*(\shQ_{\P}\otimes
\Gm(S)))$ can be viewed as given by a collection $(s_v)$, $v$ ranging
over the vertices of $\tilde\tau$, with $s_v\in
\Hom(\check\Lambda_z,\Gm(S))$, such that for any face
$\tilde\omega\subseteq\tilde\tau$ containing vertices $v$ and $w$, 
\[
s_v|_{\Lambda_\omega^{\perp}}
=s_w|_{\Lambda_\omega^{\perp}}.
\]
Now consider the fan $\check\Sigma_{\tau}$ in
$\dualvs{\Lambda_{\tau,\RR}}$; by parallel transport to $z\in
\Int(\sigma)$, we can consider $\dualvs{\Lambda_{\tau,\RR}}$ as a
quotient of $\check\Lambda_{\RR,z}$, and  we can pull back the fan
$\check\Sigma_{\tau}$ to $\check\Lambda_{\RR,z}$ to get a fan of (not
strictly) convex cones. We call this fan  $\check\tau^{-1}
\check\Sigma_{\sigma}$ as in Construction~\ref{basicgluing}; it of
course depends on $e$. The maximal cones of this fan are in
one-to-one correspondence with vertices $v$ of $\tilde\tau$. Then $s$
defines a map $s:\check\Lambda_z\to \Gm(S)$ by defining $s$ to be
$s_v$ on the cone of  $\check\tau^{-1}\check\Sigma_{\sigma}$
corresponding to $v\in\tilde\tau$.  This is a piecewise
multiplicative function, and it is clear the set of piecewise
multiplicative functions with respect to the fan
$\check\tau^{-1}\check\Sigma_{\sigma}$ coincides with
$\Gamma(\tilde\tau, \pi^*(\shQ_{\P}\otimes\Gm(S)))$. We then have
$V(\tau)\subseteq V(\sigma)$ as usual, and $s$ then induces an
automorphism of $V(\tau)\times S$, which we also denote by $s$.

Note that the group $\Gamma(\tilde\tau,
\pi^*(\shQ_{\P}\otimes\Gm(S)))$ does not depend on $e:\tau\to\sigma$,
so we only need to understand the dependence of the action on
$V(\tau)\times S$ described above. Let $e_i:\tau\to\sigma_i$, $i=1,2$
with $\sigma_i\in\P_{\max}$. Then as usual, we obtain
$V(\tau_1)\subseteq V(\sigma_1)$, $V(\tau_2)\subseteq V(\sigma_2)$
and a map
\[
\Phi_{e_1e_2}:V(\tau_2)\to V(\tau_1),
\]
(even if the pair $(e_1,e_2)$ is not maximal). Using the notation of
Construction~\ref{basicgluing}, the map $\Phi_{e_1e_2}$ is induced by
a map of monoids 
\[
\phi:\partial Q_1\to \partial Q_2
\]
obtained on each cone corresponding to a vertex of $\tau_i$ by
parallel translation through that vertex. Now
$s\in\Gamma(\tilde\tau,\pi^*(\shQ_{\P}\otimes\Gm(S)))$ defines both
$s_1:\partial Q_1\to \Gm(S)$ and $s_2:\partial Q_2\to \Gm(S)$; as
this has been done by parallel transport into the interiors of
$\sigma_1$ and $\sigma_2$, it is easy to see that
\[
s_1(p)=s_2(\phi(p)),
\]
and hence we have
\[
s_1\circ(\Phi_{e_1e_2}\times\id_S)=(\Phi_{e_1e_2}\times\id_S)\circ s_2.
\]
Thus, even though we have not defined $V(\tau)$ abstractly, the
action of $\Gamma(\tilde\tau,\pi^*(\shQ_{\P}\otimes\Gm(S)))$ on
$V(\tau)$ is well-defined with respect to the gluing maps
$\Phi_{e_1e_2}\times\id_S$.
\end{construction}

\begin{definition}
\label{PMdef}
For $\tau\in\P$, we will write $\check{PM}(\tau)$ for the group 
$\Gamma(\tilde\tau,\pi^*(\shQ_{\P}\otimes\Gm(S)))$ described above.
Here, $PM$ stands for piecewise multiplicative.
\end{definition}

We can now define open gluing data.

\begin{definition}
\label{open_gluing_data}
\emph{Open gluing data} for $\P$ over $S$ are data
$s=(s_e)_{e\in\coprod\Hom(\tau,\sigma)}$ with
\[
s_e\in \check{PM}(\tau)
\]
for $e:\tau\to \sigma$. They must satisfy
\begin{enumerate}
\item $s_{id}=1$ for $id:\sigma\to\sigma$ the identity.
\item For $e_1:\sigma_1\to\sigma_2$, $e_2:\sigma_2\to\sigma_3$ and
$e_3=e_2\circ e_1$, yielding $\tilde\sigma_1\subseteq
\tilde\sigma_2\subseteq \tilde\sigma_3$, we have
\[
s_{e_2}|_{\tilde\sigma_1}\cdot s_{e_1}=s_{e_3}.
\]
\end{enumerate}
Let $Z^1(\P,\shQ_{\P}\otimes\Gm(S))$ denote the set of all open
gluing data for $\P$ over $S$. This forms a group under
multiplication.

We define trivial open gluing data to be open gluing data $s=(s_e)$
such that there exists $t=(t_{\sigma})_{\sigma\in\P}$ with
$t_{\sigma}\in \check{PM}(\sigma) =\Gamma(\tilde\sigma,
\pi^*(\shQ_{\P}\otimes\Gm(S)))$ such that for any $e:\tau\to\sigma$, 
\[
s_e=t_{\tau}^{-1}t_{\sigma}|_{\tau}.
\]
Here restriction denotes the restriction of  $t_{\sigma}$ to the face
$\tilde\tau$ of $\tilde\sigma$ corresponding to $e$. 

We call the set of trivial open gluing data $B^1(\P,\shQ_{\P}\otimes
\Gm(S))$, and write
\[
H^1(\P,\shQ_{\P}\otimes\Gm(S))={Z^1(\P,\shQ_{\P}\otimes\Gm(S))
\over B^1(\P,\shQ_{\P}\otimes\Gm(S))}.
\]

Given open gluing data $s=(s_e)$, we can associate to it closed
gluing data $\bar s=(\bar s_e)$ as follows. By the construction of
$\shQ_{\P}$, we have $(\shQ_{\P})_y\cong\shQ_{\tau}\cong
\Gamma(W_{\tau},\shQ_{\P})$ for $y\in \Int(\tau)$, $\tau\in\P$. Thus
for any $e:\tau\to\sigma$ there is a map given by the composition of
natural maps
\[
\Gamma(\tilde\tau,\pi^*(\shQ_{\P}\otimes \Gm(S))) \to
(\shQ_{\P})_y\otimes\Gm(S)\mapright{\cong}
\Gamma(W_{\tau},\shQ_{\P}\otimes\Gm(S)) \mapright{}
\Gamma(W_{e},\shQ_{\P}\otimes\Gm(S))
\]
where the last map is restriction. We set $\bar s_e$ to be the image
of $s_e$ under this map. This gives a \v Cech 1-cocycle, because of
the 1-cocycle condition on open gluing data.
\qed
\end{definition}

For a discussion of the difference between open and closed gluing data
see Remark~\ref{openclosedgluing}. The notation as a cohomology group
is motivated by the fact that the definition fits into the framework
of barycentric complexes developed in the appendix, however for the
dual cell complex.

Given open gluing data $s$, we can define, for each maximal pair
$(e_1,e_2)$,
\[
\Phi_{e_1 e_2}(s):= s_{e_1}^{-1} \circ (\Phi_{e_1e_2}\times\id_S)
\circ s_{e_2}.
\]
Define $\foR_{e_1 e_2}\subset V(\sigma_2)\times V(\sigma_1)\times S$
as the graph of $\Phi_{e_1 e_2}(s)$ relative $S$.

\begin{example}
The simplest example is the case $B=\RR/\ZZ$, and $\P$ consists of
one maximal cell $\sigma$ and one vertex $v$. Identifying
$\tilde\sigma$ with the interval $[0,1]$, $\tilde\sigma\times_B
\tilde\sigma$ has three connected components: the diagonal and the
points $(0,1)$ and $(1,0)$.  In this case $V(\sigma)\cong \Spec
\ZZ[x,y]/(xy)$, and if there is no twisting, the three corresponding
components of $\foR$ in  $V(\sigma)\times V(\sigma)=\Spec
\ZZ[x,y,x',y']/(xy,x'y')$ are the diagonal, the subvariety given by
the equation $x'y=1$, and the subvariety given by the equation
$xy'=1$. The quotient of $V(\sigma)$ by this \'etale equivalence
relation is a nodal elliptic curve. One can also take
$B=\RR/n\ZZ$ with $n>1$ and subdivide $B$ into $n$ intervals of 
length one. In this case one obtains $n$ copies of $\Spec \ZZ[x,y]/(xy)$,
and after quotienting one obtains a cycle of $n$ rational curves,
i.e. a Kodaira type $I_n$ fibre.
\end{example}

\begin{lemma}
\label{equivrel}
The $\foR_{e_1 e_2}$ are pairwise disjoint; their union
\[
\foR=\bigcup_{(e_1, e_2)} \foR_{e_1 e_2}
\]
defines an \'etale, quasi-compact equivalence relation on
$\coprod_{\sigma\in \P_{\max}} V(\sigma)\times S$. Furthermore,
$\foR$ is closed in $\coprod_{\sigma_1,\sigma_2\in\P_{\max}} 
V(\sigma_1)\times V(\sigma_2)\times S$.
\end{lemma}
\proof
By construction $\foR_{e_1 e_2}$ is the graph of an isomorphism
between open subsets of $V(\sigma_1)\times S$ and $V(\sigma_2)\times
S$, and thus the projections of $\foR_{e_1e_2}$ to $V(\sigma_i)\times
S$ are local isomorphisms, hence \'etale. Closedness of
$\foR_{e_1e_2} \subseteq V(\sigma_2)\times V(\sigma_1)\times S$
follows from maximality of the pair $(e_1,e_2)$, so that the
isomorphism cannot be extended. Quasicompactness holds by local
finiteness of $\P$.

To check that the $\foR_{e_1 e_2}$ are pairwise disjoint
it suffices to show this after pull-back to the irreducible components
\[
\Big( V(\sigma_2)\cap X_{w_2} \Big) \times \Big( V(\sigma_1)\cap
X_{w_1}\Big) \times S,
\]
where the $w_i\in \tilde\sigma_i$ run over the set of vertices.  (The
notation $V(\sigma_i)\cap X_{w_i}$ is a slight abuse of notation: by
this we mean the irreducible component of $V(\sigma_i)$ corresponding
to the vertex $w_i$.) But by the local description of
$\tilde\sigma_1\times_B \tilde\sigma_2$ a pair of vertices
$(w_1,w_2)$ belongs to at most one maximal cell, indexed by
$(e_1,e_2)$. This translates into the statement that only $\foR_{e_1
e_2}$ intersects the irreducible component in question non-trivially.

By construction $\foR$ is symmetric and contains the diagonal. It
remains to show transitivity. Let $\sigma_1,\sigma_2,\sigma_3 \in
\P_{\max}$ with $\pi_i:\tilde\sigma_i \to \sigma_i$ as before. It
suffices to verify the transitivity relation on each irreducible
component of $V(\sigma_1)$ separately. Let $w_1\in \tilde\sigma_1$ be
a vertex and $\foR_{e_2 e_1}$ a component of $\foR$ that is
non-trivial on $V(\sigma_1) \cap X_{w_1}$. Let $R_w\to B$ be a chart
at $w=\pi_1(w_1)$ and $\tilde\sigma_i\to R_w$ for $i=1,2$ the affine
embeddings providing the local identification $\Phi_{e_2 e_1}(s)$ 
whose graph is $\foR_{e_2 e_1}$. Let $w_2\in \tilde\sigma_2$ be the
vertex mapping to $0\in R_w$. In verifying transitivity we may now
restrict to a component $\foR_{f_3 f_2}$ of $\foR$ that is non-trivial
on $V(\sigma_2)\cap X_{w_2}$. For these there is another affine
embedding $\tilde\sigma_3\to R_w$ that together with the already given
embedding of $\tilde\sigma_2$ provides the local isomorphism
$\Phi_{f_3 f_2}(s)$ with graph $\foR_{f_3 f_2}$. Let $w_3\in
\tilde\sigma_3$ be the vertex mapping to $0\in R_w$. Define $\tau$ to
be the image of $\tilde\sigma_1\cap \tilde\sigma_3 \subset R_w$ in
$B$, and let $g_i\in \Hom(\tau, \sigma_i)$, $i=1,3$ be the edges of
the barycentric decomposition belonging to the faces
$\tilde\sigma_1\cap \tilde\sigma_3\subset \tilde\sigma_i$. The pair
$(g_3,g_1)$ parametrizes a unique component $\foR_{g_3 g_1}$ of
$\foR$, which is the graph of $\Phi_{g_3 g_1}(s)$. 

\begin{center}
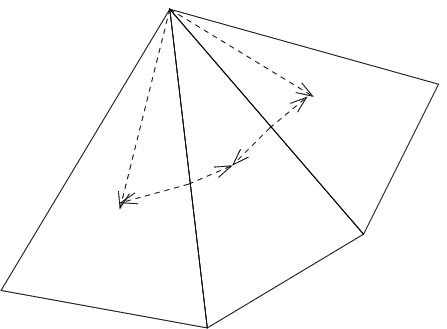
\end{center}

It is clear from toric geometry that the composition of $\Phi_{f_3
f_2}|_{X_{w_2}}$ with $\Phi_{e_2 e_1}|_{X_{w_1}}$ equals $\Phi_{g_3
g_1}|_{X_{w_3}}$ on the common domain of definition. As for the
twisting by $s$, we observe that the cocycle condition implies
that when restricted to $w_i$,
\[
s_{e_1}= s_{0 \tilde\sigma_1} s_{0\omega}^{-1}, \quad
s_{e_2}= s_{0 \tilde\sigma_2} s_{0\omega}^{-1},
\]
where $\omega= \tilde\sigma_1\cap\tilde\sigma_2$ and with indices
referring to the barycentric subdivision of $R_w$. In the definition
of $\Phi_{e_2 e_1}(s)$ on $X_{w_1}$ the term $s_{0\omega}^{-1}$
therefore cancels:
\[
\Phi_{e_2 e_1}(s)|_{X_{w_1}} = s_{0 \tilde\sigma_2}^{-1} \circ
\Phi_{e_2 e_1} \circ s_{0 \tilde\sigma_1}.
\]
Similar formulae hold for $\Phi_{f_3 f_2}$ and $\Phi_{g_3 g_1}$. Hence
\[
\Phi_{f_3 f_2}(s)\circ \Phi_{e_2 e_1}(s)|_{X_{w_1}}
= s_{0 \tilde\sigma_3}^{-1} \circ
\Phi_{g_3 g_1} \circ s_{0 \tilde\sigma_1}|_{X_{w_1}}
= \Phi_{g_3 g_1}(s)|_{X_{w_1}}.
\]
This establishes transitivity.
\qed
\medskip

\begin{definition}
If $s$ are open gluing data for $\P$ over $S$, we set
\[
X_0(B,\P,s):= \Big(\coprod_{\sigma\in \P_{\max}} 
V(\sigma)\times S\Big) \Big/\foR.
\]
By Lemma \ref{equivrel} this 
is an algebraic space locally of finite type over $S$. 
\qed
\end{definition}

It remains to show that this space is a direct limit for the functor
$F_{S,s}$. We denote by
\[
p: \coprod_{\sigma\in \P_{\max}}V(\sigma)\times S \lra X_0(B,\P,s)
\]
the quotient morphism, and $p_{\sigma}:V(\sigma)\times S\to
X_0(B,\P,s)$ the restriction of $p$ to $V(\sigma)\times S$.

\begin{lemma}\label{q_tau}
Let $s$ be open gluing data for $\P$ over $S$, and $\bar s$ the
associated closed gluing data. There exists a system of maps $q_\tau:
X_\tau\times S \to X_0(B, \P,s)$ that is compatible with the functor
$F_{S,\bar s}$; that is, for every $\tau_1,\tau_2\in \P$ and $e\in
\Hom (\tau_1,\tau_2)$ we have
\[
q_{\tau_2}= q_{\tau_1} \circ F_{S,\bar s}(e).
\]
\end{lemma}

\proof
We first define $q_w$ for $w\in \P$ a vertex. The set of maximal
cones of the fan $\Sigma_w$ defining $X_w$ is in one-to-one
correspondence with $\coprod_{\sigma\in \P_{\max}} \Hom(\{w\},
\sigma)$, so $X_w$ has an open covering indexed by this set. The
affine open set $V_e\subset X_w$ corresponding to $e\in \Hom(\{w\},
\sigma)$ is what we denoted by $V(\sigma)\cap X_{\tilde w}$ in the
proof of Lemma \ref{equivrel} for the unique lift $\tilde w\in
\tilde\sigma$ distinguished by $e$. This gives a closed embedding
\[
\iota_e: V_e\lra V(\sigma)
\]
identifying $V_e$ with the irreducible component of $V(\sigma)$
corresponding to $\tilde w$. On $V_e$ define $q_w$ as composition
\[
q_w: V_e\times S \stackrel{\iota_e\times \id_S}{\lra}
V(\sigma)\times S\stackrel{s_e^{-1}}{\lra} V(\sigma)\times S
\stackrel{p}{\lra} X_0(B,\P,s).
\]
(Technically, $s_e^{-1}$ only acts on the image of
$\iota_e\times\id_S$, but it can be extended in any way we like to
$V(\sigma)\times S$ for convenience of notation in this proof.) To
check well-definedness let $\sigma_1,\sigma_2\in \P_{\max}$ and
$e_i\in \Hom(\{w\}, \sigma_i)$. Let $\tilde\sigma_i \subset R_w$ be
the embeddings of the coverings of $\sigma_i$ into the chart at $w$
corresponding to $e_i$ as in the definition of $\foR$. Let
$\tilde\tau_i=\tilde\sigma_1\cap \tilde\sigma_2$, viewed as subcell
of $\tilde \sigma_i$, $\tau\in\P$ the image of $\tilde\tau_i$,
$f_i:\tau\to\sigma_i$, $g:w\to\tau$ the induced morphisms, with
$f_i\circ g=e_i$.  By the definition of $\foR$
\[
p|_{V(\tilde\tau_2)\times S} = p|_{V(\tilde\tau_1)\times S} \circ \Phi_{f_1
f_2}(s).
\]
On the other hand, $\iota_{e_i}$ identifies $V_{e_1}\cap V_{e_2}$ with
an irreducible component of $V(\tilde\tau_i)\subset V(\sigma_i)$, and
these two embeddings are related by $\Phi_{f_1 f_2}$:
\[
\iota_{e_1}|_{V_{e_1}\cap V_{e_2}} = \Phi_{f_1 f_2}\circ
\iota_{e_2}|_{V_{e_1}\cap V_{e_2}}.
\]
On $V_{e_1}\cap V_{e_2}$ it follows
\begin{eqnarray*}
(p|_{V(\tilde\tau_1)\times S})\circ s_{e_1}^{-1}\circ
(\iota_{e_1}\times\id_S)
&=&(p|_{V(\tilde\tau_1)\times S})\circ s_{e_1}^{-1}\circ 
(\Phi_{f_1 f_2}\times \id_S)\circ (\iota_{e_2}\times\id_S)\\
&=&(p|_{V(\tilde\tau_1)\times S})\circ
s_{e_1}^{-1}s_{f_1}\circ\Phi_{f_1 f_2}(s)\circ s_{f_2}^{-1}\circ
(\iota_{e_2}\times\id_S)\\
&=& (p|_{V(\tilde\tau_1)\times S})\circ s_{g}^{-1}\circ
\Phi_{f_1 f_2}(s)\circ s_{f_2}^{-1}\circ
(\iota_{e_2}\times\id_S)\\
&=&(p|_{V(\tilde\tau_1)\times S})\circ\Phi_{f_1 f_2}(s)\circ
s_g^{-1}s_{f_2}^{-1}\circ (\iota_{e_2}\times\id_S)\\
&=& (p|_{V(\tilde\tau_1)\times S})\circ
\Phi_{f_1 f_2}(s)\circ 
s_{e_2}^{-1}\circ
(\iota_{e_2}\times\id_S)\\
&=&(p|_{V(\tilde\tau_2)\times S})\circ s_{e_2}^{-1}\circ
(\iota_{e_2}\times\id_S).\\
\end{eqnarray*}
This shows well-definedness.

For non-minimal $\tau\in \P$ choose a vertex $w\in \tau$ and
$g\in \Hom(\{w\}, \tau)$. Then put
\[
q_{\tau}= q_w\circ F_{S,\bar s}(g).
\]
Since $F_{S,\bar s}$ is  a functor, compatibility of these
definitions follows once we convince ourselves of independence of the
choices of $w$ and $g$. Let $w_i\in\tau$ be vertices and $g_i\in
\Hom(\{w_i\},\tau)$. Write $q_{\tau}^i$ for $q_{\tau}$ defined with
$g_i$. Let $\sigma\in \P_{\max}$ and $f\in\Hom(\tau, \sigma)$.
Putting $e_i= f\circ g_i$ we obtain the following commutative
diagram.
\[
\xymatrix@C=30pt
{w_1\ar@/^/[rrd]^{e_1}\ar[rd]_{g_1}&&\\
&\tau\ar[r]^(.3){f}&\sigma\\
w_2\ar@/_/[rru]_{e_2}\ar[ru]^{g_2}&&\\
}
\]
Each $e_i$ gives an identification of an open $V_{e_i}\subset
X_{w_i}$ with an irreducible component of $V(\sigma)$. Then
$F(g_i)^{-1} (V_{e_i})$ is the open set in $X_{\tau}$ given by the
maximal cone $K_{f}\in \Sigma_\tau$ corresponding to $f: \tau\to
\sigma$, hence does not depend on $i$. This open set is $V_f$. Since
the $V_f$ with $f\in \coprod_{\sigma \in \P_{\max}}
\Hom(\tau,\sigma)$ cover $X_{\tau}$, it is enough to show that the
$q_{\tau}^i$ agree on $V_f\times S$.

Let $K_{e_i},K_{g_i}\in \Sigma_{w_i}$ be the cones corresponding to
$e_i$ and $g_i$. Recall from Definition~\ref{functor} that
$F_{S,s}(g_i)= (F(g_i)\times \id_S) \circ \bar s_{g_i}$ where
$F(g_i)|_{V_f}$ is induced by the description of cones
\[
K_f=(K_{e_i}+\RR K_{g_i})/\RR K_{g_i}
\]
On the other hand, $q_{w_i}$ involved the closed embedding
$\iota_{e_i}: V_{e_i}\to V(\sigma)$ defined by identifying
$\dual{K_{e_i}}$ with a maximal face of $\dual{C(\sigma)}$. Therefore,
the composition
\[
V_f\stackrel{F(g_i)}{\lra} V_{e_i}\stackrel{\iota_{e_i}}{\lra}
V(\sigma)
\]
is independent of $i$. Together with the cocycle condition we see
that also
\[
s_{e_i}^{-1}\circ (\iota_{e_i}\times \id_S) \circ (F(g_i)\times \id_S)
\circ \bar s_{g_i} =
s_f^{-1}\circ \big( (\iota_{e_i}\circ F(g_i))\times\id_S\big)
\]
does not depend on $i$. The composition with $p: V(\sigma)\times S\to
X_0(B,\P,s)$ equals $q_{\tau}^i|_{V_f}$, which is hence independent
of $i$ too.
\qed
\medskip

\begin{proposition}
\label{is direct limit}
$\coprod_\tau q_\tau: \coprod_\tau X_\sigma\times S \lra
X_0(B,\P,s)$ is a direct limit for the functor $F_{S,\bar s}$.
\end{proposition}
\proof
We have to verify the universal property. Let $\psi_\tau:
X_\tau\times S \to Z$ be a system of morphisms to an algebraic space
that is compatible with $F_{S,\bar s}$. First we observe that
$V(\sigma)\times S$, $\sigma\in \P_{\max}$, is a similar direct
limit, but over the spaces $V_f\times S$ with $f\in \coprod_\tau \Hom
(\tau, \sigma)$ and morphisms $\iota_g:=F(g)|_{V_{f_2}}:V_{f_2} \to
V_{f_1}$ for $f_i:\tau_i\to\sigma$, $g:\tau_1\to \tau_2$,
$f_1=f_2\circ g$. Each $V_f$
is canonically an open subset of $X_\tau$. Consider the system of
morphisms $\chi_f=\psi_{\tau}\circ \bar s_f|_{V_f\times S}:V_f\times
S\to Z$, $f\in \coprod_{\tau} \Hom(\tau,\sigma)$. This is compatible
with the system $(V_f\times S,\iota_g)$. In fact, for $f_i, g$ as
above, 
\[
F_{S,\bar s}(g)|_{V_{f_2}\times S}=(F(g)\times \id_S)|_{V_{f_2}\times S}
\circ \bar s_g=\bar s_{f_1}\circ (\iota_g\times\id_S)\circ
\bar s_{f_2}^{-1}.
\]
Thus
\[
\chi_{f_2}=\psi_{\tau_2}\circ\bar s_{f_2}|_{V_{f_2}\times S}
=\psi_{\tau_1}\circ F_{S,\bar s}(g)\circ \bar s_{f_2}|_{V_{f_2}\times S}
=\chi_{f_1}\circ (\iota_g\times\id_S).
\]
Hence we obtain a morphism
\[
\chi(\sigma): V(\sigma)\times S\lra Z.
\]
To descend to $X_0(B,\P,s)$ it remains to show that the
$\chi(\sigma)$ commute with the local isomorphisms $\Phi_{e_1
e_2}(s)$ whose graphs define $\foR$ locally. As in the proof of
transitivity of $\foR$ we may restrict to one irreducible component
at a time. If $e_i: \tau\to\sigma_i$ this amounts to choosing a
vertex $w\in\tau$ and homomorphisms $f_i:\{w\} \to \sigma_i$ with
$f_i=e_i\circ g$ for some $g:\{w\}\to\tau$. A chart at $w$ reduces to
the case $\sigma_i\subset M\otimes_\ZZ\RR$, $w=0$ and
$\tau=\sigma_1\cap \sigma_2$. This situation gave open embeddings
$V(\tau)\subset V(\sigma_i)$, and restricted to the irreducible
component $V_g\subset V(\tau)$ labelled by $w=0\in\tau$ they defined
$\Phi_{e_1 e_2}$. We obtain the following diagram.
\[
\xymatrix@C=50pt
{
V_g\times S\ar[r]^{\hbox{open}}\ar[rd]_{s_{f_2}}\ar[dd]_{\Phi_{e_1e_2}(s)}
&V_{f_2}\times S\ar[d]_{s_{f_2}}\ar[rd]^{\chi_{f_2}}&\\
&X_w\times S\ar[r]^{\psi_w}&Z\\
V_g\times S\ar[ru]^{s_{f_1}}\ar[r]_{\hbox{open}}&
V_{f_1}\times S\ar[u]^{s_{f_1}}\ar[ru]_{\chi_{f_1}}&
}
\]
The arrows labelled ``open'' are open embeddings canonically defined
by toric geometry, while the labels $s_{f_i}$ denote similar open
embeddings but twisted by the action of $s_{f_i}$ on $X_w$. In the
case without self-intersections, $V(\tau)=V(\sigma_1)\cap
V(\sigma_2)$, $V_g=V(\tau)\cap X_w$ and $V_{f_i}=V(\sigma)\cap X_w$.
In general $V_g$ and $V_{f_i}$ are open subschemes of $X_w$.
Therefore the upper and lower triangles commute trivially. The two
triangles with common side $\psi_w$ commute by the definition of
$\chi_{f_i}$, and the triangle on the left commutes by the definition
of $\Phi_{e_1 e_2}(s)$. Therefore the whole diagram commutes, and
following the two paths from $V_g\times S$ to $Z$ via $X_w\times S$
shows the desired compatibility of $\chi_f$ with our equivalence
relation.
\qed

\begin{proposition}
$X_0(B,\P,s)$ is a separated algebraic space. If $B$ is compact
and $S$ is a separated scheme, then $X_0(B,\P,s)$ is a proper
algebraic space over $S$.
\end{proposition}

\proof An algebraic space is separated if the equivalence relation
defining it is a closed embedding (\cite{Knutson}, Def. II.1.6)
and this is the case by Lemma~\ref{equivrel}. Thus in addition if
$S$ is separated, so is the morphism $f:X_0(B,\P,s)\rightarrow S$
(\cite{Knutson}, Prop. II.3.10). If $B$ is compact, $X_0(B,\P,s)$
is clearly of finite type over $S$ as $\P$ must be a finite
set, so that $X_0(B,\P,s)$ has an \'etale cover of finite type over $S$.
Finally, to show $f$ is proper, note that
\[
q:\coprod_{v\in\P} X_v\times S\lra X_0(B,\P,s)
\]
given by $q_v:X_v\times S\rightarrow X_0(B,\P,s)$ on $X_v\times S$
is a surjective map, and $\coprod X_v\times S\rightarrow S$
is a proper map. Then $f$ is also proper (see \cite{EGA}, II 5.4.3 (ii)
in the scheme case, but the same proof works for algebraic spaces).
\qed
\bigskip

Next, we consider when two sets of open gluing data give isomorphic
spaces. Suppose $s,s'$ are open gluing data for $(B,\P)$ over $S$. We
are interested in understanding when there are isomorphisms
\[
\varphi:X_0(B,\P,s)\to X_0(B,\P,s').
\]
We will only be interested in isomorphisms $\varphi$ with
$\varphi(q_{\tau}(X_{\tau}\times S))=q_{\tau}'(X_{\tau}\times S)$; any
other isomorphism will reflect some global symmetry of $B$.
In this case, we obtain a commutative diagram for each $\tau\in\P$
\[
\begin{matrix}
X_{\tau}\times S&\mapright{\varphi_{\tau}}&X_{\tau}\times S\\
\mapdown{q_{\tau}}&&\mapdown{q_{\tau}'}\\
X_0(B,\P,s)&\mapright{\varphi}&X_0(B,\P,s')
\end{matrix}
\]
If $\varphi_{\tau}$ also preserves toric strata of $X_{\tau}\times S$,
which implies it is induced by an element $\bar t_{\tau}\in
\shQ_{\tau}\otimes\Gm(S)$, then we say the isomorphism $\varphi$
\emph{preserves $B$}. Clearly the component of the identity of
the automorphism group of $X_0(B,\P,s)$ is the set of automorphisms
preserving $B$.

A simple example of an automorphism not preserving $B$ is when
$B=\RR/\ZZ$ and $X_0(B,\P,s)$ is isomorphic to a nodal cubic, say
$y^2z=x^3+x^2z$ in $\PP^2$. Then the automorphism $y\mapsto -y$
induces a map on the normalization interchanging toric strata.

\begin{proposition}
\label{isomorphism}
\item{(1)}
Let $s,s'$ be two open gluing data for $(B,\P)$ over $S$. 
Then there exists an isomorphism
\[
\varphi:X_0(B,\P,s)\mapright{\cong} X_0(B,\P,s')
\]
preserving $B$ if and only if $s/s'\in B^1(\P,\shQ_{\P}\otimes
\Gm(S))$.
\item{(2)} The map 
\[
H^1(\P,\shQ_{\P}\otimes\Gm(S))\to H^1(\W,\shQ_{\P}\otimes\Gm(S))
\]
induced by taking open gluing data to closed gluing data is an inclusion.
\end{proposition}

\proof (1) If $s/s'\in B^1(\P,\shQ_{\P}\otimes\Gm(S))$, then there
exists $t_{\sigma}\in \check{PM}(\sigma)$ for all $\sigma\in\P$ such that for
$e:\tau\to\sigma$, $s_e'=t_{\sigma}^{-1}|_{\tau}s_et_{\tau}$. Then
the isomorphism 
\[
\coprod_{\sigma\in\P_{\max}}V(\sigma)\times
S\to\coprod_{\sigma\in\P_{\max}} V(\sigma)\times S
\]
which is given by $t_{\sigma}$ on $V(\sigma)\times S$ is compatible
with the equivalence relation defined by $s$ on the domain and the
equivalence relation defined by $s'$ on the range. Thus this
isomorphism descends to give an isomorphism
\[
X(B,\P,s)\to X(B,\P,s').
\]

Conversely, if there is an isomorphism $\varphi$ preserving $B$, then
for each $\sigma\in\P_{\max}$, $e:\tau\to\sigma$, we have
$V_e\subseteq X_{\tau}$ the open subset of $X_{\tau}$ given by
Definition~\ref{Ve}. The isomorphism $\varphi_{\tau}:X_{\tau}\times
S\to X_{\tau}\times S$ induces an isomorphism $\varphi_e:V_e\times
S\to V_e\times S$, and these patch to give an isomorphism
$\tilde\varphi_{\sigma}: V(\sigma)\times S\to V(\sigma)\times S$,
which is induced by some  $t_{\sigma}\in \check{PM}(\sigma)$. Then $\varphi$
is induced by
\[
\tilde\varphi=\coprod\tilde\varphi_{\sigma}:\coprod_{\sigma\in\P_{\max}} 
V(\sigma)\times S\to \coprod_{\sigma\in\P_{\max}} V(\sigma)\times S.
\]
Now set $s_e''=s_et_{\sigma}^{-1}$ whenever $e:\tau\to\sigma\in
\P_{\max}$, $s_e''=s_e$ otherwise.  Then $s''$ is congruent to $s
\mod B^1(\P,\shQ_{\P}\otimes \Gm(S))$, so as above we obtain an
isomorphism
\[
\psi:X_0(B,\P,s'')\to X_0(B,\P,s)
\]
induced by $t_{\sigma}^{-1}$ on $V(\sigma)$. Furthermore, the composition
$\varphi\circ\psi$ is then induced by the identity maps on
$V(\sigma)\times S$ for all $\sigma\in\P_{\max}$. Thus, by replacing $s$
by $s''$, we can assume $\tilde\varphi_{\sigma}$ is the identity for
each $\sigma\in\P_{\max}$, and thus the equivalence relations
$\foR$ and $\foR'$ defined by $s$ and $s'$ must be the same.

Thus in particular, whenever $(e_1,e_2)$ is a maximal pair,
\[
(s'_{e_1})^{-1}\circ\Phi_{e_1e_2}\circ
s'_{e_2}=s_{e_1}^{-1}\circ\Phi_{e_1e_2}\circ s_{e_2}.
\]
In particular, 
\begin{equation}
\label{samething}
s_{e_2}/s_{e_2}'=s_{e_1}/s_{e_1}'
\end{equation}

Now for each $\tau\in\P\setminus\P_{\max}$, choose some morphism
$e:\tau\to\sigma\in\P_{\max}$, and define $t_{\tau} =s_e/s_e'$. We
first need to check this is independent of $e$. Given a diagram 
\[
\xymatrix@C=30pt
{&&\sigma_1\\
\tau\ar@/^/[rru]^{e_1}\ar[r]^(.7){f}\ar@/_/[drr]_{e_2}&\omega
\ar[ru]_{g_1}\ar[rd]^{g_2}&\\
&&\sigma_2
}
\]
with $(g_1,g_2)$ a maximal pair, then using (\ref{samething}) for 
$(g_1,g_2)$, we have
\[
{s_{e_1}\over s'_{e_1}}=
{s_{g_1}|_{\tau} s_f\over s'_{g_1}|_{\tau} s'_f}
={s_{g_2}|_{\tau} s_f\over s'_{g_2}|_{\tau} s'_f}
={s_{e_2}\over s'_{e_2}}.
\]
Finally, with $t_{\sigma}=1$ for $\sigma\in\P_{\max}$, we note that
for any $e:\tau_1\to\tau_2$,
\[
t_{\tau_1}^{-1} s_e t_{\tau_2}|_{\tau_1}=s_e'.
\]
Indeed, given a map $f:\tau_2\to\sigma\in\P_{\max}$, $g=f\circ e$,
\[
t_{\tau_1}^{-1}s_et_{\tau_2}|_{\tau_1}
=t_{\tau_1}^{-1}s_g s_f^{-1}|_{\tau_1}t_{\tau_2}|_{\tau_1}=
{s_g'\over s_g}s_g s^{-1}_f|_{\tau_1} {s_f|_{\tau_1}\over
s'_f|_{\tau_1}}
=s_g'(s'_f)^{-1}|_{\tau_1}=s'_e.
\]
Thus $s$ and $s'$ are equivalent open gluing data.

(2) Clearly if $s\in B^1(\P,\shQ_{\P}\otimes\Gm(S))$ is induced by
$t=(t_{\tau})_{\tau\in\P}$, then $\bar s$ is the \v Cech 1-coboundary
of $\bar t=(\bar t_{\tau})_{\tau\in\P}$. Thus the map is
well-defined. In addition, if $s\in Z^1(\P,\shQ_{\P}\otimes\Gm(S))$
is mapped to a \v Cech 1-coboundary, say the coboundary of $(\bar
t_{\tau})_{\tau\in\P}$ with $\bar
t_{\tau}\in\Gamma(W_{\tau},\shQ_{\P}\otimes\Gm(S))
=\shQ_{\tau}\otimes\Gm(S)$, then $\bar t_{\tau}$ induces a map
$\varphi_{\tau}:X_{\tau}\to X_{\tau}$ and then an isomorphism of
$X_0(B,\P,s)=\displaystyle\lim_{\longrightarrow} F_{S,\bar s}$ with
$X_0(B,\P,1) =\displaystyle\lim_{\longrightarrow} F_{S,1}$, where $1$
denotes trivial open gluing data. But then by (1), $s\in
B^1(\P,\shQ_{\P}\otimes\Gm(S))$. Thus the map is injective.
\qed

\begin{remark}
\label{openclosedgluing}
Thus we see that $H^1(\W,\shQ_{\P}\otimes\Gm(S))$ parametrizes a
set of gluings of the components $X_{\tau}$, while
$H^1(\P,\shQ_{\P}\otimes\Gm(S))$ is a subgroup parametrizing gluings
which are ``locally trivial,'' i.e. \'etale locally isomorphic to
things of the form $V(\sigma)$.  It is only these locally trivial
gluings which will play a role in the theory. The difference between
open and closed gluing data can be exhibited locally. Given an
integral lattice polytope $\sigma\subseteq M_{\RR}$, $\shQ_{\P}$
still makes sense as a sheaf just on $\sigma$, as does the open cover
$\W$, so $H^1(\W,\shQ_{\P}\otimes \Gm(S))$ still makes sense, and
gives regluings of irreducible components of
$V(\sigma)$. If $H^1(\W,\shQ_{\P}\otimes\Gm(S)) \not=0$, there are
regluings of the irreducible components of $V(\sigma)$ which yield
spaces not isomorphic to $V(\sigma)$. One example of such a $\sigma$
is an octahedron in $\RR^3$; one can show
$H^1(\W,\shQ_{\P}\otimes\Gm(S))\not=0$ by a dimension counting
argument. Somewhat mysteriously, this group seems to have something
to do with the Brauer group of the toric variety
$X(\check\Sigma_{\sigma})$: see \cite{Miranda}.

We also note that one can easily show (using similar methods as in Lemma
\ref{cohom}) that $\W$ is an acyclic covering for $\shQ_{\P}\otimes
\Gm(S)$, and thus $H^1(\W,\shQ_{\P}\otimes\Gm(S))\cong
H^1(B,\shQ_{\P}\otimes\Gm(S))$.
\qed
\end{remark}
\bigskip

Next, we make the connection between the cone and fan pictures.
Given open gluing data $s$, $X_0(B,\P,s)$ is not necessarily projective
even if we have specified ample line bundles on each irreducible component
compatible with identifications. It turns out giving such data is equivalent
to giving a strictly convex piecewise linear function $(B,\P)$. However,
there is still an obstruction to gluing these line bundles to obtain
an ample line bundle on $X_0(B,\P,s)$. This obstruction, depending on
$s$, is the map described in the following theorem. If this obstruction
vanishes, $X_0(B,\P,s)$ will be projective and arise from the
cone picture of a discrete Legendre transform of $(B,\P)$.

\begin{theorem}
Let $B$ be an integral affine manifold with singularities with a
toric polyhedral decomposition $\P$ and a strictly convex integral
piecewise linear multi-valued function $\varphi$. Let $S=\Spec A$ be
an affine scheme. Let $(\check B,\check\P,\check\varphi)$ be the
discrete Legendre transform of $(B,\P,\varphi)$. Then there is a
group homomorphism
\[
o:H^1(\P,\shQ_{\P}\otimes\Gm(S))\to H^2(B,\Gm(S))
\]
such that if $o(s)=1$, there is a choice of gluing data $\check s$
for the cone picture for $\check\P$ over $S$
with an isomorphism
\[
\check X_0(\check B,\check\P,\check s)\cong X_0(B,\P,s).
\]
In particular, $X_0(B,\P,s)$ is a projective scheme.
\end{theorem}

\proof If $\sigma\in\P$ is a cell, then $\varphi$ induces a strictly
convex piecewise linear function $\varphi_{\sigma}$ on
$\Sigma_{\sigma}$, well-defined up to linear functions, and hence an
ample line bundle $\shL_{\sigma}$ on $X_{\sigma}$. Then it is a
standard fact of toric geometry that if $\widetilde{\check\sigma}
\subseteq\dualvs{\shQ_{\sigma,\RR}}$ is the corresponding Newton
polytope of $\shL_{\sigma}$ (as defined in the definition of discrete
Legendre transform in \S1.4), $C(\widetilde{\check\sigma})
\subseteq  \dualvs{\shQ_{\sigma,\RR}} \oplus\RR$, $\check
P_{\check\sigma}= C(\widetilde{\check\sigma}) \cap
(\dualvs\shQ_{\sigma}\oplus\ZZ)$, $R_{\check\sigma}=A[\check P_{\check
\sigma}]$, then 
\[
\Proj R_{\check\sigma}\cong X_{\sigma},
\]
with 
\[
\O_{\Proj R_{\check\sigma}}(1)\cong \shL_{\sigma}.
\]
On the other hand, $R_{\check\sigma}$ is the ring given by Definition
\ref{cones} for $\check\sigma\in\check\P$, so $\check X_{\check\sigma}
\cong X_{\sigma}$. Furthermore, it is easy to see that if $e:\tau\to
\sigma$ (or equivalently $e:\check\sigma\to\check\tau$) then
\[
F(e):X_{\sigma}\to X_{\tau}
\]
coincides with
\[
\Proj \check F(e):\check X_{\check\sigma}\to \check X_{\check\tau}.
\]

Now consider the open gluing data for the fan picture $s\in
Z^1(\P,\shQ_{\P} \otimes\Gm(S))$, with associated closed gluing data
$\bar s$, with $\bar s_e \in\Gamma(W_e,\shQ_{\P}\otimes\Gm(S))\cong
\shQ_{\sigma}\otimes\Gm(S)$ for $e:\tau\to\sigma$.  There is a natural
isomorphism between $\shQ_{\sigma}$ and  $\dualvs{(\Lambda^{\check
B}_{\check\sigma})}$.  On the other hand, choose some embedding
$\widetilde{\check\sigma}\subseteq \Lambda^{\check
B}_{\check\sigma,\RR}$. Then we can identify $\Hom(\check
P_{\check\sigma},\Gm(S))$ with $(\dualvs{\Lambda^{\check B}_{\check\sigma}
\oplus\ZZ)}\otimes\Gm(S)$. However, this is non-canonical as it
depends on the embedding of $\widetilde{\check\sigma}$; translating 
$\widetilde{\check\sigma}$ will change this identification. We can
however use the inclusion of $\dualvs{(\Lambda_{\check\sigma}^{\check
B})}\otimes\Gm(S)$ in $\dualvs{(\Lambda_{\check\sigma}^{\check
B}\oplus\ZZ)}\otimes\Gm(S)$  to identify  $\bar s_e$ with an element
$\check s_e$ of $\Hom(\check P_{\check\sigma}, \Gm(S))$. Now $(\check
s_e)$ may not satisfy the 1-cocycle condition for gluing data in the
cone picture because of the arbitrary nature of this identification.
In general, because $\bar s$ satisfied the  1-cocycle condition,
given $e_1:\check\sigma_1\to\check\sigma_2$, $e_2: \check\sigma_2
\to\check\sigma_3$, $e_3=e_2\circ e_1$,  we have
\[
\check s_{e_1}\cdot \check s_{e_2}|_{\check P_{\check \sigma_1}}
\cdot \check s_{e_3}^{-1}:\check P_{\check \sigma_1}\to \Gm(S)
\]
factors through the well-defined projection $\check
P_{\check\sigma_1} \to\ZZ$. Thus we can consider $(\check
s_{e_1}\cdot \check s_{e_2}|_{\check P_{\check \sigma_1}} \cdot
\check s_{e_3}^{-1})(1)\in \Gm(S)$. This assigns to each
two-dimensional simplex in $\Bar(\check\P)$ an element of $\Gm(S)$,
and it is easy to see from  the construction that this defines a
2-cocycle for simplicial cohomology of $B$, i.e. an element of
$H^2(B,\Gm(S))$. This defines the map $o$, and it is easy to see
$o(s)$ only depends on the class of $s$ in
$H^1(\P,\shQ_{\P}\otimes\Gm(S))$.

Now suppose $o(s)=1$. Then there exists a 1-cochain $(h_e)$,
$h_e\in\Gm(S)$, $e$ running over all edges of $\Bar(\P)$ with 
\[
(\check s_{e_1}\cdot \check s_{e_2}|_{\check P_{\check \sigma_1}}
\cdot \check s_{e_3}^{-1})(1)=h_{e_1}\cdot h_{e_2}\cdot h_{e_3}^{-1}.
\]
Interpreting $h_e$ for $e:\check\tau\to\check\sigma$ as a map
$\check P_{\check\tau}\to\Gm(S)$ given as the composition
of the projection $\check P_{\check\tau}\to\ZZ$ and the map
$n\mapsto h_e^n$, we can define
\[
\check s_e'=\check s_e\cdot h_e^{-1}.
\]
It is then clear that $\check s'=(\check s_e')$ are gluing data in the
cone picture, as the 1-cocycle condition is satisfied. Replacing
$\check s$ by $\check s'$, 
it is easy
to see that for $e:\tau\to\sigma$, we obtain a commutative diagram
\[
\begin{matrix}
\check X_{\check\sigma}\times S&\mapright{\Proj \check F_{S,\check
s}(e)} 
&\check X_{\check \tau}\times S\\
\mapdown{\cong}&&\mapdown{\cong}\\
X_{\sigma}\times S&\mapright{F_{S,\bar s}(e)}&X_{\tau}\times S
\end{matrix}
\]
Taking limits yields an isomorphism 
\[
X_0(B,\P,s)\cong \check X_0(\check B,\check\P,\check s).
\]
\qed\medskip

Thus we see that a strictly convex integral piecewise linear
multi-valued function on $B$ in particular gives rise to a choice of
ample line bundle on $X_0(B,\P,s)$ for some $s$ (though not
necessarily for all $s$ unless $H^2(B,\Gm(S))=0$). However there is
also a choice involved. There might be different choices of $(h_e)$
in the proof above, giving different choices of gluing data in the
fan picture. This choice is parametrized by $H^1(B,\Gm(S))$. While
this choice doesn't affect the underlying scheme, it does change the
line bundle.

\bigskip
\begin{example}
\label{polytope3}
If $B$ is as in Example~\ref{polytope}, then $X_0(B,\P,1)$ coincides
with $\check X_0(\check B,\check\P,1)$, where $\check B$ is obtained
from the reflexive dual $\check\Xi$ of $\Xi$. (See
Examples~\ref{Batyrevdual} and~\ref{polytope2}).
\end{example}

\begin{example}
\label{not a scheme}
According to Remark~\ref{shyalgebraic} $X_0(B,\P,s)$ is a scheme if the
cells of the polyhedral decomposition do not self-intersect. Even with
self-intersections it might be a scheme, for example in dimension one.
Here is an example where the property of being a scheme rather than an
algebraic space depends on algebraic relations between the gluing
parameters.

As in Examples~\ref{squarecone} and \ref{squarecone contd} consider
$B=\RR^2/\ZZ^2$ with only one maximal cell coming from $\tilde\sigma
=[0,1]\times[0,1] \subset \RR^2$, but now viewed in the fan picture,
with $S=\Spec k$. Then $X_v=\PP^1\times\PP^1$ for the unique vertex
$v\in B$, and closed gluing data $s=(s_e)_{e:\tau \to\sigma}$ provide
toric identifications $\{0\}\times\PP^1 \stackrel{\cong}{\to}
\{\infty\}\times\PP^1$ and $\PP^1\times\{ 0\} \stackrel{\cong}{\to}
\PP^1\times \{\infty\}$ given by the actions of some $\alpha,\beta\in
k^\times$ respectively. Here $0,\infty\in\PP^1$ denote the two
zero-dimensional torus orbits.

Now if $X=X_0(B,\P,s)$ is a scheme there exist effective divisors
$D\subset X$ through the distinguished point $o\in X$ of multiplicity
four and not containing a component of $X_\sing$. Such $D$ can be
obtained by taking the closure in $X$ of a Cartier divisor on an
affine neighbourhood of $o$. A priori this leaves us with a Weil
divisor. To obtain a Cartier divisor consider the $\ZZ/2\times
\ZZ/2$-action on $X$ defined on the normalization $\PP^1\times\PP^1$
by the toric automorphisms on each of the two factors exchanging the
two zero-dimensional torus orbits. Summing over this action we may
assume that $D$ is $\ZZ/2\times \ZZ/2$-invariant. We claim that $D$ is
now Cartier. In fact, this is only questionable at points $p\in |D\cap
X_\sing|\setminus \{o\}$. There is an \'etale neighbourhood of $p$ of
the form $\Spec k[x,y,z]/(xy)$ with one of the $\ZZ/2$-factors
exchanging the two irreducible components. (This could be made
explicit by lifting $s$ to open gluing data, which is always possible
in this case, and using the \'etale cover $V(\sigma)\to X$.) Let
$D_1+D_2$ be the decomposition of $D$ according to the two irreducible
components of this \'etale neighbourhood. Then since $D$ is
$\ZZ/2\times \ZZ/2$-invariant, $D_1\cap X_\sing =D_2\cap X_\sing$.
Therefore there exist $f_1(x,z)$, $f_2(y,z)$ defining $D_1$, $D_2$ on
their components respectively, with $f_1(0,z)=f_2(0,z)$. Then $D$ is
defined by $f_1(x,z)+f_2(y,z)-f_1(0,z)$. So $D$ is in fact Cartier.

The pull-back of $D$ to the normalization $\PP^1\times \PP^1\to X$
is now given by a bihomogeneous polynomial
\[
F=\sum_{\mu=0,\ldots,d \atop \nu=0,\ldots,e}
a_{\mu\nu}x^\mu y^{d-\mu}u^\nu v^{e-\nu}
\]
with $a_{00}\neq0$, $d,e>0$, and such that
\[
F(x,y,0,1)=c\cdot F(\alpha x,y,1,0),\quad
F(0,1,u,v)=c'\cdot F(1,0,\beta u,v),
\]
for some $c,c'\in k^\times$. Comparing coefficients gives the relations
\[
a_{\mu 0}=c\cdot\alpha^\mu a_{\mu e},\quad
a_{0\nu}=c'\cdot\beta^\nu a_{d\nu},
\]
for all $\mu,\nu$. This gives two ways to compute $a_{00}$ from $a_{de}$:
\begin{eqnarray*}
a_{00}=c a_{0e}&=& cc'\beta^e a_{de},\\
a_{00}=c' a_{d0}&=& c'c\alpha^d a_{de}.
\end{eqnarray*}
Thus a necessary condition for such a divisor $D$ to exist is an
algebraic relation of the type
\[
\alpha^d=\beta^e,
\]
with $d,e>0$, between the gluing parameters $\alpha$ and $\beta$.

Conversely, if such a relation exists the same considerations give
sufficiently many Cartier divisors on $X$ separating points and
tangents. Hence in this case $X$ is projective. In fact, such
$X$ arise from the cone picture for a rectangle with sides of
lengths $d$ and $e$.
\qed
\end{example}
\bigskip

We compute some invariants of $X_0(B,\P,s)$. 

\begin{proposition}
\label{cohoB}
Let $B$ be an integral affine manifold with singularities and toric
polyhedral decomposition $\P$. Let $s$ be open gluing data for
$(B,\P)$ over $S=\Spec(A)$, and $X=X_0(B,\P,s)$. Then
\[
H^i(X,\O_X)\cong H^i(B,A).
\]
\end{proposition}

\proof
Let 
\[
\C_{\bct}^k=\bigoplus_{\sigma_0\rightarrow\cdots\rightarrow\sigma_k
\atop \sigma_i\not=\sigma_{i+1}}
q_{\sigma_k*}\O_{X_{\sigma_k}\times S},
\]
where for $\tau\in\P$, $q_{\tau}:X_{\tau}\times S\rightarrow X$ is the
map constructed in Lemma~\ref{q_tau}. Define a differential
$d_\bct$ as in \S A.3. Then we show in 
Example \ref{Oresolution} that $\C_{\bct}^{\bullet}$ is exact. There
is also a natural inclusion $\O_X\rightarrow \C_{\bct}^0$ given by
pull-back of functions, and it is easy to check as in Example
\ref{Oresolution} that this makes $\C_{\bct}^{\bullet}$ into a resolution
of $\O_X$. 
Now 
\[
H^i(X_{\sigma_k}\times S,\O_{X_{\sigma_k}\times S})
=\begin{cases}
A & i=0\\
0 & i>0,
\end{cases}
\]
as $X_{\sigma_k}\times S$ is proper over $S$, since $\Sigma_{\sigma_k}$
is a complete fan. Thus
$H^i(X,\O_X)=H^i(\Gamma(X,\C^{\bullet}))$. However,
$\Gamma(X,\C^{\bullet})$ coincides with the complex of simplicial
cochains with coefficients in the ring $A$ with respect to the
barycentric triangulation $\Bar(\P)$. Thus
$H^i(\Gamma (X,\C^{\bullet})) =H^i(B,A)$ as desired.
\qed
\medskip

For later reference let us also record here an interpretation of
$H^1(B,\Gm(S))$. We say a line bundle on a scheme $X$ over
$S$ is \emph{numerically trivial} if its pull-back to any complete
regular curve defined over a geometric point of $S$ has degree zero.
The set of isomorphism classes of numerically trivial line bundles of
$X$ form a subgroup $\Pic^\tau(X) \subset \Pic(X)$.

\begin{proposition}\label{Pic0}
Let $B$ be an integral affine manifold with singularities and toric
polyhedral decomposition $\P$. Let $s$ be open gluing data for
$(B,\P)$ over $S=\Spec(A)$, and $X=X_0(B,\P,s)$. Then
\[
\Pic^\tau(X)=H^1(B,\Gm(S)).
\]
\end{proposition}
\proof
Given $L\in \Pic^\tau(X)$ any restriction $L_\tau=q_\tau^*(L)$ is a
numerically trivial line bundle on a toric variety defined over an
affine scheme, hence trivial. Choose isomorphisms $h_\tau: L_\tau\to
\O_{X_\tau\times S}$ and, for $e\in \Hom(\tau,\sigma)$, consider the
composition
\[
\xi_e:\O_{X_\sigma\times S}
\stackrel{h_\tau^{-1}|_{X_\sigma\times S}}{\lra}
L_\tau|_{X_\sigma\times S}= L_\sigma\stackrel{h_\sigma}{\lra}
\O_{X_\sigma\times S}.
\]
This is a locally invertible element of $\Gamma(X_\sigma\times S,
\O_{X_\sigma}\times S) = A$, that is, $\xi_e\in \Gm(S)$. For a
composition $\omega \stackrel{f}{\to} \tau \stackrel{e}{\to} \sigma$
it holds $\xi_{f\circ e}=\xi_f|_{X_\sigma\times S}\cdot \xi_e$. Hence
$(\xi_e)$ is a simplicial one-cocycle for the triangulation
$\Bar(\P)$ of $B$. Different choices of $h_\tau$ yield cohomologous
$(\xi_e)$.

Conversely, given $\xi=(\xi_e)$ we can enrich the gluing functor
$F_{S,\bar s}$ from Definition~\ref{functor} for the closed gluing
data $\bar s$ associated to $s$ to also take care of a line bundle.
The corresponding functor $G_{S,\bar s,\xi}$ maps $\tau\in\P$ to the
pair $(X_\tau,\O_{X_\tau\times S})$, and $e\in\Hom(\tau,\sigma)$ maps to
the morphism $F_{S,\bar s}(e)$ of spaces and to the morphism
\[
\xi_e\cdot F_{S,\bar s}(e)^*: F_{S,\bar s}(e)^{-1}(\O_{X_\tau\times S})
\lra \O_{X_\sigma}
\]
of line bundles. In view of the methods presented in this section it
is then clear that $G_{S,\bar s,\xi}$ has a limit $(X_0(B,\P,s),L)$ in the
category consisting of pairs of an algebraic space and a line bundle,
and that $L$ maps to $\xi\in H^1(B,\Gm(S))$ under the map described
above. Details are left to the reader.
\qed
\medskip

We end this section by discussing dualizing sheaves on 
$X=X_0(B,\P,s)$. Now we do not wish to go into the question of
duality for algebraic spaces, which is a question which does not
appear to be dealt with in the literature. Therefore, for the
discussion which follows, assume that either duality does hold, or
that $X$ is a scheme. This occurs as in Remark~\ref{shyalgebraic}, or
if $X$ arises through gluing in the cone picture, so that $X$ is in
fact a projective scheme. In any event, we take $S=\Spec k$, $k$ an
algebraically closed field. Note that $X$ is Gorenstein: it is covered by
\'etale sets of the form $V(\sigma)\subseteq U(\sigma)$, and
$V(\sigma)$ is Gorenstein by \cite{Oda}, pg. 126, Ishida's criteria,
(3). Thus there is a dualizing sheaf $\omega_X$ which is in fact a
line bundle. In addition, consider the map $\nu:\tilde
X=\coprod_{v\in\P} X_v\rightarrow X$ (induced by $q_v$ on $X_v$).
Then $\nu$ is in fact the normalization of $X$. Let $C\subseteq\tilde
X$ be the conductor locus. Then by \cite{Reid}, Proposition 2.3,
$\nu^*\omega_X \cong\omega_{\tilde X}(C)$, where the latter denotes
the ``$S_2$-isation'' of $\omega_{\tilde X}\otimes\O_{\tilde X}(C)$;
on $\tilde X\setminus \Sing(\tilde X)$, this is just $\omega_{\tilde
X\setminus\Sing(\tilde X)}(C)$. Now by construction $C\cap X_v$ is
just the reduced union of toric divisors of $X_v$, and this coincides
on $X_v\setminus\Sing(X_v)$ with the anti-canonical class of $X_v$
(see e.g. \cite{Oda}, Corollary 3.3). Thus $\omega_{\tilde
X\setminus\Sing(\tilde X)}(C)=\O_{\tilde X\setminus \Sing(\tilde
X)}$, and $\omega_{\tilde X}(C)=\O_{\tilde X}$. We shall use this to
calculate $\omega_X$:

\begin{theorem}
\label{dualizingsheaf}
Let $B$ be a {\it compact} integral affine manifold with
singularities with a toric polyhedral decomposition $\P$, and $s$
open gluing data over $S=\Spec k$, $k$ an algebraically closed
field. Let $X=X_0(B,\P,s)$.
Suppose $X$ is a scheme, or duality holds as above. Then
$\omega_{X}=\O_X$ if and only if $B$ is orientable, and
$\omega_{X}^{\otimes 2}=\O_X$ if $B$ is not orientable.
\end{theorem}

\proof First suppose $B$ is orientable. Then if $n=\dim B$,
$H^n(X,\O_X)=H^n(B,k)=k$ by orientability, compactness and 
Proposition~\ref{cohoB}. But by
Serre duality, (using the hypothesis that $B$ is compact so that $X$
is proper over $k$) we then have $H^0(X,\omega_{X})=k$. Thus we have
a non-zero global section of $\omega_{X}$. To see that it is nowhere
vanishing, use $\nu^*\omega_X=\O_{\tilde X}$. From this we conclude
that a non-zero section of $\omega_{X}$ does not vanish on any $X_v$,
and hence vanishes nowhere. Thus $\omega_X\cong\O_X$. Clearly if $B$
is not orientable, then $H^0(X,\omega_{X})=0$, so $\omega_{X}$ cannot
be trivial. However, $B$ has an oriented unramified double cover
$B'\rightarrow B$. We can pull the affine structure on $B$ back to
$B'$, and also pull back $\P$ and $s$ to $B'$ giving a polyhedral
decomposition $\P'$ of $B'$ and open gluing data $s'$, and then we
obtain an \'etale double cover $\pi:X'=X_0(B',\P',s')\rightarrow X$.
Then $\pi^*\omega_{X}=\omega_{X'}\cong \O_{X'}$ since $B'$ is
orientable. From this it follows that $\omega_{X}$ must be
two-torsion.
\qed

\section{Logarithmic structures}
\label{section3}

\subsection{Introduction to log structures}

For a normal crossing divisor $X$ in a smooth variety $V$ there is
the classical notion of differential form with logarithmic
singularities along $X$. They generate an $\O_V$-module
$\Omega^1_V(\log X)$. An important application is for a semistable
degeneration $\pi:V\to \Delta$, which is a proper map from a smooth
variety $V$ to a smooth curve $\Delta$ with $X=\pi^{-1}(0)$ a divisor
with normal crossings and $\pi|_{V\setminus X}$ smooth. One can then 
define the sheaf of log differentials $\Omega^1_X(\log X)$ by restricting
$\Omega^1_V(\log X)$ to $V$ and dividing out by $\dlog t$, where $t$
is the pull-back of a coordinate on $\Delta$.
Steenbrink
has shown that in this case $h^q(X,\Omega^p_X(\log X))$ equals
$h^{p,q}$ of a general fibre of $\pi$ \cite{steenbrink}. Logarithmic
(log-) geometry provides a means to define $\Omega^p_X(\log X)$
\emph{without knowledge of $V$}. The fundamental notions of log
geometry are due to Fontaine, Illusie and K.~Kato. Standard
references are \cite{K.Kato},\cite{F.Kato}.

Before giving the definition let us list our conventions. All sheaves
are understood in the \'etale topology. A \emph{monoid} is a set with
an associative product with a unit. We tacitly assume all monoids to
be commutative. If $M$ is a monoid then $M^\gp$ denotes the
Grothendieck group associated to $M$; this is an abelian group
together with a homomorphism $M\to M^\gp$ that has the universal
property with respect to such maps. The group of invertible elements
of $M$ is denoted $M^\times$. A monoid is \emph{integral} if $M\to
M^\gp$ is injective; and an integral monoid is \emph{saturated} if
$x\in M^\gp$, $x^n\in M$ implies $x\in M$. A finitely generated,
integral monoid is called \emph{fine}. All the finitely generated
monoids we are concerned about are \emph{toric}, which by definition
means finitely generated, saturated and $M^\gp$ free. These are
precisely the monoids isomorphic to $\ZZ^n\cap \dual{\sigma}$ with
$\sigma\subset\Hom(\ZZ^n,\RR)$ a rational, strictly convex polyhedral
cone.

\begin{definition}
A \emph{log structure} on an algebraic space $X$ is a homomorphism
\[
\alpha_X: \M_X\lra \O_X
\]
of sheaves of (multiplicative) monoids such that $\alpha_X:\alpha_X^{-1}
(\O_X^{\times})\rightarrow\O_X^{\times}$ is an isomorphism.
The triple $(X,\M_X, \alpha_X)$ is
then called a \emph{log scheme}. The quotient of $\M_X$ by the
subsheaf $\O_X^\times\subset\M_X$ is the
\emph{ghost sheaf} $\overline\M_X$ of the log structure, with monoid
structure written additively.

A morphism of log spaces $F: (X,\M_X)\to (Y,\M_Y)$ consists of a
morphism $\underline F: X\to Y$ of the underlying spaces together
with a homomorphism $F^\#: \underline{F}^{-1} (\M_Y)\to \M_X$
commuting with the structure homomorphisms: $\alpha_X\circ F^\#=
\underline F^*\circ\alpha_Y$.
\qed
\end{definition}

\noindent
By abuse of notation we often refer to $\M_X$ alone as a log
structure.

Here is the prime example. All of our log structures are locally
derived from this example.

\begin{example}\label{prime example}
Let $X$ be a scheme and $D\subset X$ a closed subset of pure
codimension $1$ and $j:X\setminus D\to X$ the embedding of the
complement. Then the natural inclusion
\[
\alpha_X: \M_{(X,D)}=j_*(\O_{X\setminus D}^\times)\cap \O_X\lra \O_X
\]
of the sheaf of regular functions with zeros contained in $D$
is a log structure on $X$. The monoid structure is multiplication. 

If $X$ is locally Noetherian and normal the stalk of
$\M_X:=\M_{(X,D)}$ at a geometric point $\bar x\to X$ has the
following form. Let $r$ be the number of irreducible components of
$D$ at $\bar x$. There is a map
\[
q:\M_{X,\bar x}\lra \NN^r
\]
given by associating to a regular function $f\in\M_{X,\bar x}$ the
vanishing orders $\ord_\mu(f)$ of $f$ along the components of $D$ at
$\bar x$. This map factors over $\M_{X,\bar x}/\O^\times_{X,\bar x}$.
On the other hand, if $f_1, f_2\in\M_{X,\bar x}\subset \O_{X,\bar x}$
have the same vanishing orders along $D$ then by normality
$f_2=h\cdot f_1$ for some $h\in \O_{X,\bar x}^\times$. This shows
that $\overline{\M}_{X,\bar x}= \M_{X,\bar x}/\O^\times_{X,\bar x}$
is a submonoid of $(\NN^r,+)$. To see that it is finitely generated
introduce the componentwise partial order on $\NN^r$, that is
$(a_1,\ldots,a_r)\ge (b_1,\ldots,b_r)$ iff $a_i\ge b_i$ for all $i$.
Then $f,g\in\M_{X,\bar x}$ and $q(f)\ge q(g)$ implies $f/g\in
\M_{X,\bar x}$, again by normality. Therefore $\im(q)$ is generated
by minimal elements. Using the fact that for any $a\in\NN^r$ the set
$\{b\in\NN^r| b\le a\}$ is finite, one can show  any infinite
sequence in $\NN^r$ contains an increasing subsequence. Therefore the
set of minimal elements in $\im(q)$ is finite.

In particular, $\overline{\M}_{X,\bar x}^{\gp}$ is a finitely
generated subgroup of $\ZZ^r$, and then 
$\Ext^1(\overline{\M}^\gp_{X,\bar x}, \O^\times_{X,\bar
x})$ is trivial. Hence the sequence
\[
1\lra \O^\times_{X,\bar x}\lra\M^\gp_{X,\bar x}\lra
\overline \M^\gp_{X,\bar x}\lra 0
\]
splits (non-canonically) and
\[
\M_{X,\bar x}= \overline \M_{X,\bar x}\times_{\overline \M_{X,\bar
x}^\gp}\M_{X,\bar x}^\gp =\overline \M_{X,\bar x}\oplus
\O^\times_{X,\bar x}.
\]
We emphasize that this splitting does not, however, hold on the sheaf
level. For example, let $k$ be a field and consider $X=\A^1_k$ and
$D=\{0\}$. Then $\overline\M_X$ is the skyscraper sheaf
$\NN_{\{0\}}$. Thus $\overline\M_X\oplus \O^\times_X$ has
non-trivial sections with support in $\{0\}$, while this is not true
for $\M_X\subset\O_X$. The reason is that unlike for coherent sheaves
the canonical map $\shHom(\overline \M_X, \O^\times_X)_{\bar x}\to 
\Hom(\overline \M_{X,\bar x}, \O^\times_{X,\bar x})$
is \emph{not} bijective at $x=0$.

If $X$ is not normal the stalks of $\overline\M_X$ might not be
finitely generated nor even countable. Here is an example that we
learned from A.~Ogus. Let $X=\Spec\CC[x,y]/(x^2-y^3)$ and $D=V(x,y)$.
Then for any $a\in \CC$ the function $x-ay$ lies in $\M_{(X,D)}$, but
if $a,b\in \CC$ the existence of $h\in\CC\lfor x,y\rfor/(x^2-y^3)$
with $x-ay=h\cdot(x-by)$ formally, implies $a=b$. The same works
for a node $X=\Spec\CC[x,y]/(xy)$.
\qed
\end{example}

For later use, we include
\begin{lemma}
\label{MXDgp}
If $X$ is normal and $D\subset X$ is the support of an effective
Cartier divisor then $\shM_{(X,D)}^\gp= j_*\O_{X\setminus D}^\times$
where $j:X\setminus D\to X$ is the inclusion.
\end{lemma}

\proof
The universal property of the morphism to the associated Grothendieck
group gives a homomorphism
\[
\shM_{(X,D)}^\gp\lra j_*\O_{X\setminus D}^\times,\quad
(f,g)\longmapsto f/g,\quad
f,g\in j_*\O_{X\setminus D}^\times\cap \O_X,
\]
which is clearly injective. Conversely, let $\bar x\to X$ be a
geometric point and $h\in (j_*\O_{X\setminus D}^\times)_{\bar x}$.
Let $t\in \O_{X,\bar x}$ be a function defining at $\bar x$ a Cartier
divisor with support $D$. If $n\ge 0$ is the maximum of the pole
orders of $h$ along the prime components of $D$ at $\bar x$ then $t^n
h$ extends in codimension one to $X$. By normality of $X$ it follows
that $f=t^n h\in \O_{X,\bar x}$. Then $f,t^n \in \shM_{X,\bar x}$ and
$(f,t^n)\in \shM_{(X,D)}^\gp$ maps to $h$. 
\qed
\medskip

The log structures obtained in this way can still be quite
pathological. For example, $\M_X$ might not be locally generated by
global sections, see Example~\ref{non-fine example}. This leads to a
certain coherence condition that requires some
preparation to explain.

\begin{definition}
\label{associated}
An arbitrary homomorphism of sheaves of monoids
\[
\varphi: \shP\lra\O_X
\]
defines an \emph{associated log structure} $\M_X$ by
\[
\M_X= \big(\shP\oplus \O^\times_X\big)\big/
\big\{(p,\varphi(p)^{-1})\,\big|\, p\in \varphi^{-1}(\O^\times_X)\big\},
\]
and $\alpha_X(p,h)= h\cdot \varphi(p)$. 

If $f:X\to Y$ is a morphism of algebraic spaces, the \emph{pull-back} 
of a log structure $\M_Y$ on $Y$ to $X$ is the log structure associated to
\[
f^*\circ \alpha_Y: f^{-1}\M_Y\lra \O_X.
\]
The notation is $f^*\M_Y$. If $f$ is an embedding we also
speak of \emph{restricting the log structure to $X$}.
\qed
\end{definition}

\begin{example}\label{log points}
For any field $k$ and monoid $Q$ (written additively) with
$Q^\times=\{0\}$ and field $k$ let $\M_{\Spec k}$ be the monoid sheaf
with sections $Q\times K^\times$ over a separable extension $k\subset
K$. Then
\[
Q\times K^\times\lra K,\quad (q,a)\longmapsto
\left\{ \begin{array}{ll}
a,&q=0\\
0,&\text{otherwise},
\end{array} \right.
\]
defines a homomorphism $\M_{\Spec k}\to \O_{\Spec k}$, which is a log
structure. Special cases are
$Q=0$ and $Q=\NN$, which give the \emph{trivial} and \emph{standard
log points} respectively. The notation for the standard log point is
$\Spec k^\ls$. This is the restriction of the log structure
associated to the origin $o\subset \AA_k^1$ to $o$.
\qed
\end{example}

Another case is when $\shP=P$ is a constant sheaf with stalks a fine
monoid, also denoted $P$. Then $\varphi$ is a \emph{chart} for its
associated log structure. More generally, a chart for a log structure
$\alpha_X:\M_X\to \O_X$ is a homomorphism $\varphi:P\rightarrow\O_X$
and a factorization of $\varphi$ via $\alpha_X$, inducing an
isomorphism of the log structure associated to $\varphi$ with $\M_X$.
A log structure possessing charts \'etale locally is called
\emph{fine}. In general it is impossible to find a global chart for a
fine log structure. This is the concept of coherence for $\M_X$. At
least in the fine saturated case there exists then for every geometric
point $\bar x\to X$ a chart in an \'etale neighbourhood with
$P=\overline\M_{X,\bar x}$, see \cite{ogus} Lemma~2.12. Another way to
think of charts is by noting that $P\to \O_X$ defines, \'etale
locally, a morphism to an affine toric variety:
\[
\Phi: X\lra \Spec \ZZ[P].
\]
Now $\Spec \ZZ[P]$ has a distinguished log structure defined either
as associated to the canonical map $P\to \O_{\Spec
\ZZ[P]}$, or by the log structure coming from the divisor defined by
the complement of the big torus orbit as in Example~\ref{prime
example}. Then $\M_X$ is isomorphic to $\Phi^* \M_{\Spec
\ZZ[P]}$. In particular, toric varieties provide an ample source of
fine log structures.

\begin{example}\label{non-fine example}
Here is an example of a naturally occuring log structure that is
\emph{not} fine. Let $k$ be an algebraically closed field and $X=\Spec
k[P]$ with $P$ the monoid generated by $e_1,\ldots, e_4$ with single
relation $e_1+e_3=e_2+e_4$. This is the $3$-dimensional
$A_1$-singularity $xy-tw=0$ in $\A^4_k$ written as a toric variety.
Let $D$ be the divisor given by $t=0$, a union of two copies of
$\A^2_k$. Then the log structure on $X$ associated to $D$ is not fine
at the zero-dimensional torus orbit $o\in X$: If $f\in \O_{X,\bar o}$
has zero locus contained in $D$ then $f=h\cdot t^n$ with $h\in
\O_{X,\bar o}^\times$. Therefore
\[
\overline\M_{X,\bar o} = \NN.
\]
On the other hand, at the generic point $\eta$ of $x=y=0$ the divisor
$D$ is normal crossing with two components. Hence
\[
\overline \M_{X,\bar \eta} =\NN^{\oplus 2}.
\]
Since $o\in \cl(\eta)$ this shows that $\overline\M_X$ is not
locally generated by global sections, which holds for fine log
structures. This is not a perverse example, but a fundamental
issue in our approach. In the example of degeneration of quartics
given in the introduction, the total space $\X$ has 24 ordinary
double points, and locally the structure of $\X_0\subseteq \X$
is exactly as described here.
\qed
\end{example}

The example also illustrates that the obstruction to being fine might
already be captured by the ghost sheaf of the log structure. This is
indeed the case as we will see instantly. Let us say a sheaf of
monoids $\overline\M$ is \emph{fine} if its stalks are fine monoids
and if for every geometric point $\bar x\to X$ there exists an
\'etale neighbourhood $f:U\to X$ of $\bar x$ and a surjection from a
constant sheaf
\[
P\lra f^{-1}\overline\M,
\]
which is an isomorphism on stalks at $\bar x$. We have the following
criterion.

\begin{proposition}\label{fine.crit}
A log structure is fine iff its ghost sheaf is fine.
\end{proposition}
\proof
Clearly if the log structure is fine, so is the ghost sheaf.
Conversely, suppose the ghost sheaf is fine.
By the same argument as in Example~\ref{prime example} the exact
sequence
\[
1\lra \O^\times_{X,\bar x}\lra\M^\gp_{X,\bar x}\lra
\overline \M^\gp_{X,\bar x}\lra 0
\]
splits, inducing a right inverse $\kappa$ of $\M_{X,\bar x}\to\overline
\M_{X,\bar x}$. Put
$P=\overline\M_{X,\bar x}$. Since $\overline\M_X$ is a fine monoid
sheaf, after going over to an \'etale neighbourhood of $\bar x$ there
exists a surjection $\psi:P\to \overline\M_X$. Composing with
$\kappa$ gives a sheaf homomorphism (in a smaller neighbourhood)
\[
\varphi: P\lra \M_X,
\]
which is a chart because it induces a surjection to the ghost sheaf.
In fact, for a geometric point $\bar y\to X$ and $m\in\M_{X,\bar y}$
let $\overline{m}$ be the image in $\overline \M_{X,\bar y}$. There
exists $p\in P$ with $\psi(p)=\overline m$. Then $\varphi(p)= h\cdot
m$ for $h\in \O_{X,\bar y}^\times$. This shows that the map $P\oplus
\O_X^\times \to\M_X$ induced by $\varphi$ is surjective. The kernel
is the submonoid of pairs $(p,h)$ with $\alpha_X (\varphi(p))\cdot
h=1$, which coincides with the submonoid divided out in the
definition of the associated log structure.
\qed
\medskip

A central concept in log geometry is (log) smoothness. The definition
runs analogously to formal smoothness for schemes. For our purposes
it is more instructive to use the characterization of log smooth
morphisms due to K.~Kato \cite{K.Kato} Theorem~3.5, see also
\cite{F.Kato} Theorem~4.1. A model log smooth map is given by the
map of log schemes induced by a morphism of monoids $Q\rightarrow P$
with finite kernel, with some subtlety in non-zero characteristic.

\begin{definition}
A morphism $f:(X,\M_X)\to (Y,\M_Y)$ of fine log schemes is
\emph{log smooth} if \'etale locally it fits into a commutative
diagram of the form
\[
\begin{CD}
X@>>> \Spec\ZZ[P]\\
@VVV @VVV\\
Y@>>> \Spec\ZZ[Q]
\end{CD}\]
such that
\begin{enumerate}
\item The horizontal arrows induce charts for the log structures on
$X$ and $Y$ that are compatible with $f^\#:f^{-1}\M_Y\to \M_X$.
\item The induced morphism $X\to Y\times_{\Spec\ZZ[Q]} \Spec\ZZ[P]$ is
a smooth morphism of schemes.
\item The right-hand vertical arrow comes from a monoid homomorphism
$Q\to P$ with $\ker (Q^\gp\to P^\gp)$ and the torsion part of
$\coker(Q^\gp\to P^\gp)$ finite groups of orders invertible on
$X$.
\end{enumerate}
\qed
\end{definition}

If the characteristic is $0$ and all monoids are toric, then the last
condition just requires the homomorphism $Q\to P$ to be injective. 

Obviously any toric morphism between affine toric varieties induced
by a monoid homomorphism $Q\to P$ as in (3) is log smooth. In
particular, log smooth morphisms need not be flat, a toric blowup for
example. This will only be the case if we also ask $f$ to be
\emph{integral}, which by definition means that for every geometric
point $\bar x\to X$ the ring homomorphism $\ZZ[\overline \M_{Y,f(\bar
x)}] \to \ZZ[\overline \M_{X,\bar x}]$ is flat.

Given a log smooth morphism it is sometimes useful to consider
certain adapted charts. An obstruction for finding such
charts in positive characteristic is torsion in $\M_X^\gp/f^*\M_Y^\gp
= \overline\M_X^\gp/f^{-1}\overline\M_Y^\gp$. In our cases torsion
never occurs and the following proposition suffices.

\begin{proposition}\label{adapted charts}
(\cite{ogus} Proposition~2.25.)
Let $f:(X,\M_X)\to (Y,\M_Y)$ be a log smooth morphism of fine log
schemes. Then for every geometric point $\bar x\to X$ with
$\M_{Y,f(\bar x)}^\gp\to \M_{X,\bar x}^\gp$ injective with
torsion-free cokernel, any chart $\varphi:Y\to \Spec
\ZZ[\overline\M_{Y,f(\bar x)}]$ inducing an isomorphism of ghost
sheaves at $f(\bar x)$ fits into a diagram
\[
\begin{CD}
U@>\psi>> \Spec\ZZ[\overline\M_{X,\bar x}\oplus\NN^r]\\
@VVV @VVV\\
Y@>\varphi>> \Spec\ZZ[\overline\M_{Y,f(\bar x)}]
\end{CD}
\]
with right-hand vertical arrow defined by the composition of
$f^\#_{\bar x}: \overline\M_{Y,f(\bar x)}\to \overline\M_{X,\bar x}$
with the inclusion $(\id,0):\overline\M_{X,\bar x}\to
\overline\M_{X,\bar x}\oplus \NN^r$, and $\psi$ \'etale 
with the induced map $\overline{\M}_{X,\bar x}\rightarrow\O_U$
defining a chart
for $(X,\M_X)$ on an \'etale neighbourhood $U\to X$ of $\bar x$. 
\end{proposition}

The factor $\NN^r$ is necessary if one insists on $\psi$ to be \'etale
rather than smooth. (Consider the example of a smooth morphism with
trivial log structures.) However, this factor has no effect on the log
structure.

Log smoothness extends the class of smooth objects considerably:

\begin{example}
1)\ \ Any toric variety over a scheme $S$ with its canonical log
structure is log smooth over $S$, viewed as log scheme with trivial
log structure.

\noindent
2)\ \ Let $\pi:V\to \Delta$ be a semistable family of algebraic
varieties as in the first paragraph of this section and $x\in\Delta$
a closed point with separably closed residue field. $\pi$ induces a
log morphism for the log structures associated to the divisors
$\{x\}\subset \Delta$ and $X=\pi^{-1}(x)\subset V$ introduced in
Example~\ref{prime example}. \'Etale locally and up to smooth factors
it has the form $\Spec k[\NN^r]\to \Spec k[\NN]$ induced by the
diagonal morphism $\NN\to\NN^r$. Hence it is log smooth and integral.

\noindent
3)\ \ Restricting the previous example to the central fibre gives
another interesting log-smooth morphism $(X,\M_X)\to \Spec k(x)^\ls$.
More generally one can show that a variety $X$ with normal crossings
supports a log smooth morphism to the standard log point iff $X$ is
\emph{d-semistable}, which by definition means
$\shExt^1(\Omega^1_X,\O_X)\cong \O_{X_\sing}$ (\cite{Friedman},
\cite{Kawamata; Namikawa 1994}, \cite{F.Kato}). So the existence of a
log smooth structure restricts not only the type of singularities of
$X$, but poses also some more subtle global analytical conditions
(here the triviality of a locally free sheaf over $X_\sing$).
\qed
\end{example}

\subsection{Sheaves of log structures}

We next study moduli of log structures on a space with given fine
ghost sheaf $\overline\M_X$. For the following discussion we make the
overall assumption that $X$ \emph{is reduced}. We
first show that the structure homomorphism $\alpha: \M_X\to \O_X$ is
then already determined by the extension
\begin{eqnarray}\label{extension}
1\lra \O_X^\times \lra \M^\gp_X\stackrel{q}{\lra}
\overline\M^\gp_X\lra 0.
\end{eqnarray}

\begin{proposition}
Let $X$ be a reduced algebraic space. Let $\alpha,\alpha': \M_X\to
\O_X$ be two fine log structures on $X$ with the same monoid sheaf
and such that $\alpha'\circ \alpha^{-1}:
\O_X^\times\to \O_X^\times$ is the identity. Then $\alpha= \alpha'$.
\end{proposition}

\proof
Let $\bar\eta$ be a generic geometric point of $X$. Let
$m\in\M_{X,\bar\eta}$. Since $\O_{X,\bar\eta}$ is a field $\alpha(m)$
is either $0$ or invertible. The first case occurs iff $q(m)\not=0$,
and then also $\alpha'(m)=0$. In the invertible case the assumption
implies that $\alpha(m)=\alpha'(m)$. Thus in any case
$\alpha(m)=\alpha'(m)$. Therefore $\alpha$ and $\alpha'$ agree
generically, and $X$ being reduced this implies $\alpha=\alpha'$.
\qed

\begin{corollary}\label{extension gives log structure}
For a reduced algebraic space $X$ the set of isomorphism classes of
log structures with ghost sheaf $\overline \M_X$ is a subset of
$\Ext^1(\overline \M^\gp_X, \O_X^\times)$.
\end{corollary}

By abuse of notation we will therefore confuse an isomorphism class
of log structures with its associated extension class. However, it is
not true that all of $\Ext^1( \overline \M^\gp_X, \O_X^\times)$
arises in this way. Let us study this in one example.

\begin{example}
\label{quadruplepoint}
Let $k$ be a separably closed field. Consider the quadruple point
$X=\Spec k[x_1,\ldots, x_4]/ (x_1 x_3, x_2 x_4)$ with $\overline
\M_X$ coming from the log structure induced by embedding into $\Spec
k[x_1,\ldots, x_4]/ (x_1 x_3 -x_2 x_4)=\Spec k[P]$ with
\[
P= \langle e_1,\dots,e_4\,|\, e_1+e_3=e_2+e_4\rangle.
\]
Denote by $o\in X$ the zero-dimensional torus orbit. This is the
unique point with $\overline \M_{X,\bar o}= P$. Note that any element of
$P$ extends uniquely to a global section of $\overline \M_X$, and
these global sections generate $\overline \M_X$ at every point.

Any extension $q:\M_X \to \overline \M_X$ by $\O_X^\times$
gives four $\O_X^\times$-torsors
\[
\shL_i =q^{-1}(e_i).
\]
There is an identification of $\O_X^\times$-torsors
\[
\big(\shL_1\times \shL_3\big)/ \O_X^\times = q^{-1}(e_1+e_3)
= q^{-1}(e_2+e_4) = \big(\shL_2\times \shL_4\big)/ \O_X^\times.
\]
By abuse of notation we are going to use tensor notation for sections
of the quotients on the left and right-hand sides. The components of
$X$ are $X_i= V(x_{i+2}, x_{i+3})$, the indices taken modulo $4$. Let
$\eta_i$ be the generic point of $X_i$. Since
$q(\shL_i)|_{X\setminus\supp(e_i)}=0$, on $X\setminus \supp(e_i)=
X\setminus (X_{i-2}\cup X_{i-3})\subseteq X_{i-1}\cup X_i=V(x_{i+2})$ the
unit in $\M_X$ lies in $\shL_i$. This gives a distinguished section
$s_i$ of $\shL_i$ over $X\setminus \supp(e_i)$. Obviously $s_i$
extends as the unit section of $\O_X^\times \subset \M_X$, but
generally not as a section of $\shL_i$. However, possibly after going
over to an \'etale neighbourhood $U\to X$ of the origin, the $\shL_i$
become trivial. Therefore $x_i^{d_i}\cdot s_i$ extends uniquely to a
section $\sigma_i$ of $\shL_i|_U$ for some $d_i\in \ZZ$ depending
only on the isomorphism class of the extension. Note that
$x_i^{d_i}s_i$ restricts to $x_i^{d_i}$ on $X\setminus \supp(e_i)$.
By the above identification there exists $h\in \O_{X,\bar o}^\times$
with
\[
\sigma_2\otimes\sigma_4=h\cdot \sigma_1\otimes \sigma_3.
\]
Any other $\sigma'_i\in \shL_{i,\bar o}$ restricting to $x_i^{d_i}$ away
from $\supp(e_i)$ has the form $\sigma'_i=h_i\cdot \sigma_i$ with
$h_i\in\O_{X,\bar o}^\times$ and $h_i|_{V(x_{i+2})}=1$. Writing
$h_i=1+x_{i+2}f_i$ this changes $h$ by $(1+x_3 f_1)(1+x_1f_3)
(1+x_4f_2)^{-1}(1+x_2f_4)^{-1}$. Terms of this form generate
$1+\mathfrak{m}_{\bar o}$. In particular, the residue of $h$ in $k=
\O_{X,\bar o}/ {\mathfrak m}_{\bar o}$ is well-defined. We will see
below that this residue together with the $d_i$ are the only
invariants of the extension. Conversely, for any $\lambda\in
k^{\times}$ the pull-back of the standard log structure by
\[
(x_1,x_2,x_3,x_4)\longmapsto (\lambda \cdot x_1,x_2,x_3,x_4)
\]
defines an extension with $d_1,\ldots,d_4=1$, $h=\lambda$. 
So any residue is possible. Using results below, one can see that in fact
\[
\Ext^1(\overline \M_X^\gp, \O_X^\times) = \ZZ^4\oplus k^\times.
\]
For such an extension $(d_i,\lambda)$ there exists a homomorphism
$\M_X\to \O_X$ inducing the identity on $\O_X^\times$ iff $d_i>0$
for every $i$.
\qed
\end{example}

This example is quite typical: To define $\alpha_X$ on $m\in
\M_{X,\bar x}$ let $\overline\eta$ be a geometric generic point of
$X$ with $\bar x\in \cl(\overline \eta)$. Then either
$q(m)_{\overline \eta}\neq 0$ and we must have
$\alpha_X(m)_{\overline \eta}=0$; or $q(m)_{\overline \eta}=0$ and
$m\in\O_{X,\overline\eta}^\times$ prescribes
$\alpha_X(m)_{\overline\eta}$. The question is then if the regular
function defined on the generic points with closure containing $\bar
x$ extends as regular function to $\bar x$. In other words, the
subspace of $\Ext^1(\overline\M^\gp_X,\O_X^\times)$ parametrizing
log structures with ghost sheaf $\overline \M_X$ is characterized
by a pointwise positivity condition.

The next aim is an explicit description of isomorphism classes of log
structures with given fine ghost sheaf, building on
Corollary~\ref{extension gives log structure}. Let $\M_X\to \O_X$ be
a fine log structure on $X$, $\bar x\to X$ a geometric point and
$P=\overline \M_{X,\bar x}$. After replacing $X$ by an \'etale open
set, we may assume that there is a chart $P\to \M_X$. In
particular, $P=\Gamma(X,\overline\M_X)$ and for every geometric point
$\bar y\in X$ the restriction map $P\to \overline\M_{X,\bar y}$ is
surjective. Define the relation sheaf $\shR$ of $\overline \M^\gp_X$
for the given chart by
\[
0\lra \shR\lra P^\gp\lra \overline \M^\gp_X\lra 0.
\]
Because $P^\gp$ is a constant sheaf with stalks a free finitely
generated abelian group, $\shExt^1(P^\gp, \O_X^\times)$ vanishes as
sheaf in the \'etale topology. We obtain an exact sequence
\[
0\lra \shHom( \overline \M^\gp_X, \O_X^\times) \lra \shHom( P^\gp,
\O_X^\times) \lra \shHom( \shR,\O_X^\times)\lra \shExt^1(
\overline \M^\gp_X, \O_X^\times)\lra 0.
\]
Therefore
\[
\shExt^1( \overline \M^\gp_X, \O_X^\times)
= \shHom(\shR, \O_X^\times)/
\shHom( P^\gp, \O_X^{\times}).
\]
Now any $p\in P^\gp$ induces a section of $\overline{\M}_X^{\gp}$
that is zero on $X\setminus\supp(p)$. (When writing $\supp(p)$ we
always mean the support of $p$ viewed as section of
$\overline\M_X^\gp$.) Thus $p|_{X\setminus \supp(p)}$ is a section of
$\shR$, and these sections generate $\shR$. Denoting
$j_p:X\setminus\supp(p)\to X$ the inclusion, a more explicit
description of the numerator is then
\[
\shHom(\shR, \O_X^\times)_{\bar x} =
\Big\{ (h_p)_{p\in P^\gp} \,\Big|\, h_p\in ({j_p}_*\O^\times_{X\setminus
\supp(p)})_{\bar x}, h_p\cdot h_q= h_{p+q} 
\text{ on } X\setminus \supp(p+q)
\Big\}.
\]
As $\shHom( P^\gp, \O_X^{\times})_{\bar x}=\Hom(\overline\M_{X,\bar
x}^{\gp}, \O_{X,\bar x}^\times)$ we obtain the following result.

\begin{proposition}\label{(h_p)}
The stalk of $\shExt^1( \overline \M^\gp_X, \O_X^\times)$ at a
geometric point $\bar x$ of $X$ is canonically isomorphic to the
quotient of
\[
\Big\{ (h_p)_{p\in \overline \M_{X,\bar x}^\gp} \,\Big|\, h_p\in
(j_{p*}\O_{X\setminus \supp(p)}^\times)_{\bar x}, h_p\cdot h_q=
h_{p+q} \text{ on } X\setminus \supp(p+q) \Big\}
\]
by $\Hom(\overline \M_{X,\bar x}^\gp, \O_{X,\bar x}^\times)$. The
extension class of the log structure associated to a homomorphism
$\varphi: \overline \M_{X,\bar x}\to \O_{X,\bar x}$ is represented by
\[
h_p=\varphi(p)|_{X\setminus \supp (p)}.
\]
Conversely, $(h_p)_p$ comes from a log structure iff $h_p$ extends
by $0$ to $X$ for all $p\in \overline \M_{X,\bar x}$.
\end{proposition}
\proof
It remains to prove the statement about log structures. Let $\tilde
\varphi: P\to \M_U$ be the chart induced by $\varphi$ on an \'etale
neighbourhood $U\to X$ of $\bar x$. Let $\shR=\ker( P^\gp\to
\overline\M_U^\gp)$ be the associated relation sheaf.
The composition
$\shR\lra P^\gp_U\lra \M^\gp_U$ factors through $\O_U^\times$. The
factorization sends the local section of $\shR$ induced by $p\in
P^\gp$ to $h_p=\alpha_X(\tilde \varphi(p))|_{U\setminus \supp(p)}
=\varphi(p)|_{U\setminus \supp(p)} \in\O^\times_{U\setminus
\supp(p)}$. This defines a map $\psi:\shR\to \O_U^\times$. By
definition $\psi$ makes the following diagram with exact rows
commutative.
\[\begin{CD}
0@>>> \shR@>>>P^\gp@>>>\overline\M_U^\gp@>>>0\\
@.\psi@VVV\tilde\varphi^\gp@VVV@|\\
0@>>> \O_U^\times@>>>\M_U^\gp@>>>\overline\M_U^\gp@>>>0.
\end{CD}\]
This shows that the class of the extension in the lower row defining
the log structure is the image of $\psi$ under the connecting
homomorphism $\Hom(\shR,\O_U^\times)\to \Ext^1
(\overline\M_U^\gp,\O_U^\times)$. Conversely, assume an extension
$\tilde h_p$ of $h_p$ by zero exists. Then
\[
p\longmapsto \tilde h_p
\]
defines a chart $P\rightarrow\O_U$
for a log-structure on $X$ with the given extension
class $(h_p)_p$.
\qed
\medskip

The proposition gives a complete description of germs of log
structures at any geometric point $\bar x\to X$ by systems of
invertible functions $(h_p)_p$, where $h_p$ is defined only on
$X\setminus \supp(p)$. Comparing the $h_p$ on $X(p):=\cl(X\setminus
\supp(p))$ (where $\cl$ denotes closure)
leads to a finer classification of log structures.

\begin{definition}
\label{logtype}
Two germs of log structures $\xi, \xi'\in \shExt^1(\overline\M^\gp_X,
\O_X^\times)_{\bar x}$ are said to be of the \emph{same type} if any
representatives $h_p$, $h'_p$ in the description of
Proposition~\ref{(h_p)} differ only by multiplication by
$e_p\in\O_{X(p),\bar x}^{\times}$, for any $p\in P$. Two log
structures on an \'etale open set $U\to X$ have the same type if
their germs at every geometric point $\bar x\to U$ have the same
type.
\end{definition}

Geometrically speaking the type of a germ of log structures is
determined by the vanishing orders of extensions of $h_p$
along the irreducible components of $X(p)\cap \supp(p)$.
\medskip

For a scheme $S$ let $S^\ls$ denote $S$ equipped with log structure
$\M_S=\NN_S\oplus \O^\times_S$. Then if $X$ is an algebraic space
over $S$ a \emph{log smooth structure} on $X$ over $S^\ls$ is the
choice of a log structure $\M_X$ on $X$ together with a log-smooth
morphism $(X,\M_X)\to S^\ls$. Isomorphisms are isomorphisms of the
log structures on $X$ commuting with the morphisms to $S^\ls$. Our
principal interest is in the case $S=\Spec k$ with $k$ separably
closed; then $S^\ls$ is the standard log point. The general case is
useful for studying locally trivial moduli of log smooth structures.
In fact, a log smooth morphism $(X,\M_X)\to S^\ls$ induces a locally
trivial log-smooth deformation of the fibre over any geometric point
of $S$.

Given that $X$ is defined over $S$, a morphism from $(X,\M_X)$ to
$S^\ls$ is equivalent to the choice of a section $\rho$ of
$\M_X$ with $\alpha(\rho)=0$. Here $\rho$ is the image of the generator
of $\NN$. Let $\bar\rho$ be the image of $\rho$ in $\overline\M_X$.
Then we obtain the commutative diagram
\[\begin{CD}
1@>>>\O_X^\times@>>>\M^\gp_X @>>> \overline\M^\gp_X@>>> 0\\
@.@|@VVV@VVV\\
1@>>>\O_X^\times@>>>\M^\gp_X/\rho @>>>
\overline\M^\gp_X/\overline{\rho}@>>> 0.
\end{CD}\]
Provided that $\rho$ is nowhere nilpotent the kernel of the middle
vertical arrow is isomorphic to $\ZZ$, with generator mapping to
$\rho$. Hence we can tell $\rho$ from the diagram as the generator of
the kernel in the middle arrow intersected with $\M_X$. Thus given a
nowhere nilpotent $\bar\rho\in\Gamma (\overline{\M}_X^{\gp})$, the
set of isomorphism classes of log structures on $X$ with a morphism
to $S^{\dagger}$ inducing $\overline{\rho}$ is a subset of
$\Ext^1(\overline\M_X^\gp/ \overline\rho, \O_X^\times)$. Indeed, an
element of $\Ext^1(\overline{\M}_X^{\gp}/\overline\rho,
\O_X^{\times})$ induces an element of $\Ext^1(\overline{\M}_X^{\gp},
\O_X^{\times})$ and a diagram as above by pull-back, hence a log
structure and a section $\rho$ of $\M_X^{\gp}$. There is then a
local description of the moduli of log structures with a morphism to
$S^{\dagger}$ analogous to Proposition~\ref{(h_p)}. The only
difference is that we take the quotient by $\Hom
(\overline\M_{X,\overline x}^\gp/ \overline \rho_{\bar x}, \O_{X,\bar
x}^\times)$ rather than $\Hom (\overline\M_{X, \overline x}^\gp,
\O_{X,\bar x}^\times)$. To impose log-smoothness we need an atlas of
log-smooth local models prescribing the type of log structures
consistently.

\begin{definition}\label{ghost structures}
Let $X$ be an algebraic space over a scheme $S$. A \emph{ghost
structure} on $X$ is a choice of a fine sheaf $\overline\M$ of toric
monoids on $X$, a nowhere zero section
$\bar\rho\in\Gamma(X,\overline\M)$, and an \'etale cover of $X$,
$\{\pi_i:U_i\to X\}$, with smooth morphisms
\[
U_i\to S\times_{\Spec \ZZ}
\Spec \ZZ[\overline\M_{\bar x_i}]/(z^{\bar\rho_{\bar x_i}})
\]
for some geometric point $\bar x_i\in U_i$. Moreover, the log 
structures $U_i^{\ls}$ induced by $\overline\M_{\bar
x_i}\to\Gamma(U_i,\O_{U_i})$ have ghost sheaves isomorphic to
$\pi_i^{-1}\overline\M$ under an isomorphism respecting
$\bar\rho_{\bar x_i}$, and are of the same type on overlaps of $U_i$ and
$U_j$
(Definition \ref{logtype}). We usually denote a ghost structure on a
scheme $X$ by $X^g$. If $(X,\M_X)\to S^\ls$ is a log smooth
structure on $X$ and $X^g$ is a ghost structure, we say $(X,\M_X)$ is
of \emph{ghost type} $X^g$ if there exists an isomorphism
$\overline\M_X\to \overline\M$ identifying the given sections, and if
on the open sets $U_i$, the log structures on $U_i$ induced by
$(X,\M_X)$ and by $X^g$ are of the same type.
\qed
\end{definition}

\begin{example}
\label{fundamental example}
Here is the key example of a ghost structure. Let $B$ be an integral
affine manifold with singularities and toric polyhedral decomposition
$\P$, and let $s$ be open gluing data for $\P$ over a scheme $S$.
Then the scheme $X=X_0(B,\P,s)$ has a canonical ghost structure as
follows. By the construction of \S2, there exists an \'etale open
cover of $X_0(B,\P,s)$ by sets $V(\sigma)\times S$ along with
embeddings $V(\sigma)\subseteq U(\sigma)$, where $\sigma$ runs over
all maximal cells. More explicitly, if $\sigma$ is a maximal cell,
$P:=P_{\sigma}$ as in Definition~\ref{fanopen},  $\bar\rho:=
\rho_{\sigma}\in P$ the element corresponding to the projection
$M\oplus\ZZ\to\ZZ$, then 
\begin{eqnarray*}
U(\sigma)&=&\Spec \ZZ[P]\\
V(\sigma)&=&\Spec \ZZ[P]/(z^{\bar\rho})=\Spec \ZZ[\partial P].
\end{eqnarray*}
Thus we obtain a canonical log structure on $V(\sigma)\times S$ from
the chart $P\to \ZZ[\partial P]$, with a ghost sheaf
$\overline\M_{V(\sigma)\times S}$. It is easy to see from the
discussion of \S\ref{section2} that on overlaps,
$\overline\M_{V(\sigma)\times S}$ match, gluing to give a sheaf
$\overline\M_X$. Furthermore the sections $\bar\rho$ of
$\overline\M_{V(\sigma)\times S}$ also glue, giving $\bar\rho\in
\Gamma(X,\overline\M_X)$. Finally the identity maps $V(\sigma)\times
S\to V(\sigma)\times S$ yield the smooth morphisms defining the ghost
structure.
\end{example}

One of the fundamental problems we need to solve in this paper
now arises. We wish to obtain log structures on $X_0(B,\P,s)$. The
naive hope would be that the log structures induced by the
above inclusions $V(\sigma)\subseteq U(\sigma)$ glue to give
a global log structure on $X_0(B,\P,s)$. Unfortunately, this
only happens when $\Delta=\emptyset$, as we shall see.

\begin{example}
In Example \ref{twotriangles}, we glue together two normal
crossings varieties. The log structures in fact glue, because the
gluing of $V(\sigma_1)$ and $V(\sigma_2)$ along $V(\tau_i)$ extends
to gluings of $U(\sigma_1)$ and $U(\sigma_2)$ along $U(\tau_i)$.
However, the log smooth structures don't. In the coordinates
used in that example, $\rho$ is given by the function $yz$ on
$U(\tau_1)$ and $U(\tau_2)$, but $\phi(yz)=x^{-1}yz$, so the 
$\rho$ do not glue. (Note that $x^{-1}$ is invertible on $U(\tau_i)$,
so the $\bar\rho$'s do glue). 

In Example \ref{twopyramids}, even the log structures don't
glue. It makes a good exercise to check this using charts.
Note in this case the gluing of $V(\sigma_1)$ and $V(\sigma_2)$
doesn't extend to a gluing of $U(\sigma_1)$ with $U(\sigma_2)$.
\qed
\end{example}

As a result of this fundamental problem, in order
to obtain a log structure on $X_0(B,\P,s)$, we
have to understand the set of all log structures with a given ghost
type. We will see that given a ghost structure $X^g$ there is a
subsheaf $\shLS_{X^g}\subset \shExt^1(\overline\M_X^{\gp}/
\bar\rho,\O_X^{\times})$ whose sections over an \'etale $U\to X$ are
in one-to-one correspondence with the log-smooth structures on $U$ of
given ghost type.

\begin{definition}
Let $X$ be a reduced algebraic space, defined over a scheme $S$, and
$X^g$ a ghost structure on $X$. Denote by  $\shLS_{X^g}$ the subsheaf
of $\shExt^1(\overline \M_X^{\gp}/ \bar\rho,\O_X^{\times})$ of germs
of extension classes of the same type as prescribed by the ghost
structure. We will usually drop the superscript $g$ as the ghost
structure will always be clear.
\end{definition}

Note that this definition makes sense because for sections of
$\shExt^1( \overline \M^\gp_X, \O_X^\times)$ to be of the same type
(Definition~\ref{logtype}) is an open condition.

\begin{proposition}\label{moduli of log smooth structures}
Let $X$ be a reduced algebraic space, of finite type over a scheme
$S$. Assume that $X^g$ is a ghost structure on $X$ with
$\overline\M_X^\gp$ generated by $\bar\rho$ at the generic points.
Then the set of isomorphism classes of log smooth structures on an
\'etale set $U\to X$ of the given ghost type is canonically in
bijection with $\Gamma(U, \shLS_{X^g})$.
\end{proposition}
\proof
On a reduced space any homomorphism $\varphi: \overline\M^\gp_X\to
\O_X^\times$ with $\varphi(\overline \rho)=1$ is trivial, because
$\overline \rho$ generates $\overline\M_X^\gp$ at the generic points.
Thus $\shHom(\overline\M_X^\gp/\rho,\O_X^\times)=0$ and the
local-global $\Ext$ spectral sequence gives
\[
\Ext^1(\overline\M_X^\gp /\overline \rho, \O_X^\times)
= \Gamma(X,\shExt^1(\overline\M_X^\gp / \overline \rho,
\O_X^\times)).
\]
In view of the discussion before Definition~\ref{ghost structures}
this shows that there is a one-to-one correspondence between
isomorphism classes of log structures on an \'etale open set $U\to X$
together with a log morphism to $S^\ls$ of the type prescribed by the
ghost structure, and $\Gamma(U,\shLS_{X^g})$. It remains to show that
the log morphisms arising in this way are log smooth. It suffices to
check this in a neighbourhood of any geometric point $\bar x\to U$
and for $S=\Spec A$ affine. Let $\xi\in (\shLS_{X^g})_{\bar x}$.

The ghost structure on $U$ and
Proposition~\ref{adapted charts} provide an \'etale map
\[
\phi: U\to \Spec A[\overline\M_{X,\bar x} \oplus\NN^r]/
(z^{\bar\rho_{\bar x}})=:V
\]
with log structure associated to $\overline\M_{X,\bar x}\to
\Gamma(U,\O_X)$ having the type prescribed by the ghost structure.
For $p\in\overline\M_{X,\bar x}$ denote by $U(p)\subset U$ the
subscheme defined by $\Ann(\phi^*(z^p))$. This is the closure of
$U\setminus \supp(p)$ when we view $p$ as a section of
$\overline\M_X$ over $U$. We may assume that the functions from
Proposition~\ref{(h_p)} representing $\xi$ are defined on $U(p)$ 
and have the form $h_p\phi^*(z^p)$ with $h_p$
invertible on $U(p)$. Define $\psi: U\to V$ by sending $(p,n)\in
\overline\M_{X,\bar x} \oplus\NN^r$ to $h_p\cdot \phi^*(z^{(p,n)})$.
The germ at $\bar x$ of the associated log structure equals $\xi$. It
remains to show that $\psi$ is \'etale.

To do this recall that \'etale maps between schemes are precisely
those locally of the form
\[
B\lra B[T_1,\ldots, T_n]/(F_1,\ldots F_n)
\]
with $\det (\partial F_i/\partial T_j)$ invertible, see for example
\cite{Milne}, Corollary~I.3.16. Clearly, to show \'etaleness in a
neighbourhood of a geometric point $\bar x$, it suffices to check
invertibility of the Jacobian determinant in the local ring at $\bar
x$. Let us now write $U=\Spec B[\mathbf{T}]/(\mathbf{F})$ with
$B=A[\overline\M_{X,\bar x}\oplus\NN^r]/ (z^{\bar\rho_{\bar x}})$ and
boldface symbols denoting indexed entities. Note that the shape of
the ghost sheaf on $V$ implies that $\bar x$ maps to the subscheme in
$U$ defined by the ideal generated by $\overline\M_{X,\bar
x}\setminus \{0\}$. We can assume without loss of generality that
$\bar x$ maps to the subscheme in $U$ defined by the ideal generated
by $(\overline\M_{X,\bar x}\oplus\NN^r) \setminus\{0\}$. Let
$p_1,\ldots,p_m$ be generators of $\partial\overline{\M}_{X,\bar
x}\oplus\NN^r$  as a monoid and write $X_k=z^{p_k}$. Their relation
ideal is generated by finitely many polynomials $G_\lambda$ of the
form $\prod X_k^{a_k}-\prod X_k^{b_k}$ with $\sum a_kp_k=\sum
b_kp_k$ or $\prod X_k^{a_k}$ with $\sum a_kp_k=\infty$. 
For each $k$ choose an extension $h_k\in
B[\mathbf{T}]/(\mathbf{F})$ of $h_{p_k}$ to $U$. Possibly after
localizing at the $h_k$ we may assume $h_k$ to be invertible, even in
$A[\mathbf{X},\mathbf{T}]/(\mathbf{F})$. Note that localization of a
ring $C$ at $a$ is isomorphic to $C[T]/(aT-1)$ and hence preserves
the standard description of \'etale maps just given.

The homomorphism defining $\psi$ is then
\[
A[\mathbf{X}]/(\mathbf{G})
\lra A[\mathbf{X},\mathbf{T}]/(\mathbf{G},\mathbf{F}),\quad
\mathbf{X}\longmapsto \mathbf{h}\cdot\mathbf{X}.
\]
To bring this into standard form introduce new variables $\mathbf{S}$
with $\mathbf{S}= \mathbf{h}\cdot\mathbf{X}$
and swap $\mathbf{X}$ and $\mathbf{S}$. Now the $h_k$ are invertible
in $A[\mathbf{X},\mathbf{T}]/(\mathbf{F})$, and $G_\lambda(\mathbf{h}
\cdot\mathbf{X})=\prod h_k^{a_k} G_\lambda(\mathbf{X})$. Hence the
ideals in $A[\mathbf{X},\mathbf{T}]/(\mathbf{F})$ generated by
$\mathbf{G}(\mathbf{X})$ and by $\mathbf{G}(\mathbf{h}\cdot
\mathbf{X})$ coincide. This shows that the introduction of
$\mathbf{S}$ turns the above homomorphism into the canonical
inclusion
\[
A[\mathbf{X}]/(\mathbf{G})
\lra \big(A[\mathbf{X}]/(\mathbf{G}(\mathbf{X})\big)
[\mathbf{S}, \mathbf{T}]/(\mathbf{X}-
\mathbf{h'}\cdot\mathbf{S},\mathbf{F'}),
\]
where $h'=h(\mathbf{S},\mathbf{T})$, $F'=F(\mathbf{S},\mathbf{T})$.
The Jacobian matrix of the relations is
\[
\left(\begin{matrix}
\big(-\mathbf{h'}-\partial\mathbf{h'}/
\partial\mathbf{S}\cdot\mathbf{S}\big)&
-\big(\partial\mathbf{h'}/ \partial\mathbf{T}\cdot\mathbf{S}\big) \\
\big(\partial \mathbf{F'}/\partial \mathbf{S}\big)&
\big(\partial \mathbf{F'}/\partial \mathbf{T}\big)
\end{matrix}\right).
\]
To show that the determinant of this matrix is invertible at $\bar x$ it
suffices to show invertibility modulo $\mathbf{S}$, that is in
\[
A[\mathbf{X},\mathbf{S}, \mathbf{T}]/(\mathbf{G},
\mathbf{X}-\mathbf{h'}\cdot\mathbf{S}, \mathbf{F'},\mathbf{S}) \cong
A[\mathbf{T}]/(\mathbf{F}(\mathbf{0},\mathbf{T})).
\]
The result is
\[
\det\left(\begin{matrix}
-\operatorname{diag}\big(\mathbf{h}(\mathbf{0},\mathbf{T})\big)&0 \\
\big(\partial \mathbf{F}/\partial \mathbf{S}
(\mathbf{0},\mathbf{T})\big)& \big(\partial \mathbf{F}/\partial
\mathbf{T}(\mathbf{0},\mathbf{T})\big) \end{matrix}\right)=
(-1)^m\left(\prod_{k=1}^m h_k(\mathbf{0},\mathbf{T})\right)
\det(\partial \mathbf{F}/\partial \mathbf{T})(\mathbf{0},\mathbf{T}).
\]
The right-hand side is invertible because all $h_k$ and
$\det(\partial \mathbf{F}/\partial \mathbf{T})$ are invertible in
$B[\mathbf{T}]/(\mathbf{F})$ $=A[\mathbf{X},\mathbf{T}]/(\mathbf{G},
\mathbf{F})$. Because the $X_k\in\ker(\O_{U,\bar x}\to k(\bar x))$
this shows that the Jacobian determinant is invertible at $\bar x$.
Hence $\psi$ is \'etale in a neighbourhood of $\bar x$.
\qed

\subsection{Log structures for the fan picture}
Our main goal now is to compute the sheaf $\shLS_{X^g}$ for 
$X=X_0(B,\P,s)$ with ghost structure given by
Example~\ref{fundamental example}, and to understand global sections
of $\shLS_{X_0(B,\P,s)}$. We will in fact find that in general, there
are no global sections if $\Delta\not=\emptyset$. However, we can
find sections over nice open subsets of $X$, which will then define a
log structure on these open subsets, which shall prove to be enough
for our purposes. We will always work over $S=\Spec k$, where $k$ is
an algebraically closed field, so we write $V(\sigma)$ instead of
$V(\sigma)\times_{\ZZ} \Spec k$, etc.

We begin by calculating the sheaf $\shLS_X$ on our standard \'etale
covering of $X$, as follows. Let $\sigma\subseteq M_{\RR}$ be a
polytope with $\dim\sigma=\dim M_{\RR}$, and let
$P:=P_{\sigma}=\dual{C(\sigma)}\cap (N\oplus\ZZ)$ as usual, with
$U(\sigma)=\Spec k[P]$ and $V(\sigma)=\Spec k[\partial P]$, also as
usual. Let $\check\Sigma$ be the normal fan of $\sigma$ in $N_{\RR}$.
Choose some non-negative integer $r$, and set 
\[
V=V(\sigma)\times \Gm^r=\Spec k[\partial P\oplus\ZZ^r].
\]
$V$ carries a canonical log structure via the smooth projection map
$V\to V(\sigma)$, or equivalently via the obvious chart
\[
P\to k[\partial P\oplus\ZZ^r].
\]
This induces a ghost structure on $V$, and we will calculate
$\shLS_V$, the sheaf of log-smooth structures on $V$ with this ghost
structure.

As in Definition~\ref{logtype} write $V(p)=\cl\big(V
\setminus\supp(p)\big)$. Let $\shF$ be the sheaf on $V$ given by
\begin{eqnarray*}
U\longmapsto&
\{(h_p)_{p\in P}| \hbox{
$h_p$ is an invertible function on $U\cap V(p)$}\\
&\hbox{ and $h_p\cdot h_q=h_{p+q}$
on $U\cap V(p+q)$}\}.
\end{eqnarray*}
Then $\shF$ can be viewed as a subsheaf of
$\shHom(\shR,\O_V^{\times})$ via the map 
\[
\Gamma(U,\shF)\ni (h_p)_{p\in
P}\longmapsto (z^ph_p)_{p\in P}\in
\Gamma(U,\shHom(\shR,\O_V^{\times})).
\]
Then it is clear from the
definition that $\shLS_{V}$ is the quotient
\[
\shF/\shHom(P^{\gp}/\rho,\O_V^{\times})
=\shF/\shHom(N,\O_V^{\times}).
\]
However, it is also clear that
$\shF$ is the same as the sheaf given by
\begin{eqnarray*}
U\longmapsto&
\{(h_p)_{p\in \partial P\setminus\{\infty\}}| \hbox{
$h_p$ is an invertible function on $U\cap V(p)$}\\
&\hbox{ and $h_p\cdot h_q=h_{p+q}$
on $U\cap V(p+q)$}\}.
\end{eqnarray*}
Note that each face $\tau$ of $\sigma$ gives a cone
$\check\tau\in\check \Sigma$ and thence an irreducible closed stratum
$V_{\tau}=\Spec k[(\check\tau\cap N)\oplus\ZZ^r]$. The scheme $V$ has
irreducible components $\{V_v|\hbox{$v$ is a vertex of $\sigma$}\}$,
and these are glued together along lower dimensional strata. As
usual, $\tau\mapsto V_{\tau}$ is an order reversing correspondence.

For this local case we introduce shorthand notation: for a face
$\tau$ of $\sigma$, we denote by $\tau^{\|}$ the intersection of the
lattice $M$ with the tangent space of $\tau$ at any interior point of
$\tau$; this is a sublattice of $M$. As in \S1.5, we can choose
primitive generators $d_{\omega}$ of  $\omega^{\|}$ for each edge
$\omega$ of $\sigma$, and this induces a choice of vertices
$v^{+}_{\omega}$ and $v^-_{\omega}$. This induces orientations on
each edge of $\sigma$.

The following definition is as in \cite{Altmann}. (This is no accident,
see Example~\ref{obstructed}.)

\begin{definition}
\label{signvector}
For each 2-face $\tau\subseteq\sigma$,  we define its \emph{sign
vector}
\[
\epsilon_{\tau}:\{\hbox{dimension one faces
of $\sigma$}\}\to \{-1,0,+1\}
\]
by
\[
\epsilon_{\tau}(\omega)=\begin{cases}
0&\hbox{if $\omega\not\subseteq\tau$}\\
\pm 1&\hbox{if $\omega\subseteq\tau$}
\end{cases}
\]
where the sign is chosen so that if $\omega_1,\ldots,\omega_N$ are the
edges of $\tau$, then the singular $1$-cycle $\sum
\epsilon_{\tau_i}(\omega_i) \omega_i$ is the oriented boundary of
$\tau$ with respect to some orientation on $\tau$ and the given chosen
orientations on $\omega_i$. This is well-defined up to sign, and we
make an arbitrary choice. 
\qed
\end{definition}

\begin{theorem}
\label{finalLS}
With the choices in Definition~\ref{signvector}, $\shLS_V$ is
isomorphic to the subsheaf of $\bigoplus_{\dim\omega=1}
\O_{V_{\omega}}^{\times}$ defined as follows. If $U\subseteq V$ is any
open set, then  $\Gamma(U,\shLS_V)$ consists of $(f_{\omega})\in 
\Gamma(U,\bigoplus_{\dim\omega=1} \O_{V_{\omega}}^{\times})$ such that
for every two-dimensional face $\tau$ of $\sigma$, we have
\begin{eqnarray}
\label{MC}
\prod_{\dim\omega=1} 
d_{\omega}\otimes f_{\omega}^{\epsilon_{\tau}(\omega)}|_{V_{\tau}}&=&1
\in M\otimes_{\ZZ}
\Gamma(U,\O_{V_{\tau}}^{\times}).
\end{eqnarray}
Here the product is taken in the group $M\otimes_{\ZZ}\Gamma(U,
\O_{V_{\tau}}^{\times})$, where $M$ is written additively and
$\Gamma(U,\O_{V_{\tau}}^{\times})$ multiplicatively.
\end{theorem}

\proof
$\shLS_V$ can be viewed as a sheaf in either the Zariski topology
or the \'etale topology. Here we will show we get an isomorphism in the
Zariski topology, though the same proof works in the \'etale topology.
For the purpose of this proof we therefore work in the Zariski
topology.

We begin by defining a homomorphism
\[
\xi:\shF\to\bigoplus_{\dim\omega=1} \O^{\times}_{V_{\omega}},
\]
and then show the kernel of $\xi$ is $\shHom(N,\O^{\times}_V)$,
and the image is the subsheaf defined by the condition (\ref{MC}).

To define $\xi$, let $h=(h_p)_{p\in\partial P\setminus\{\infty\}}$
be a section of $\shF$ over an open set $U$ of $V$. For each vertex
$v$ of $\sigma$, consider $p\in \check v\cap N$. Because $z^p
\in k[\partial P]$ does not vanish generically on the irreducible
component $V_v$ of $V$, the generic point of $V_v$ is not contained in
$\supp(p)$, so $V_v\subseteq V(p)$. Thus in particular
\[
p\longmapsto h_p|_{V_v}\in \Gamma(U,\O_{V_v}^{\times})
\]
gives a multiplicative map
\[
g_v:\check v\cap N\to\Gamma(U,\O^{\times}_{V_v})
\]
which extends multiplicatively to a map
\[
g_v:N\to \Gamma(U,\O^{\times}_{V_v}),
\]
i.e.\ an element $g_v\in M\otimes\Gamma(U,\O_{V_v}^{\times})$. Now
if $\omega$ is an edge of $\sigma$, then $\check\omega=\check v_{\omega}^+
\cap \check v_{\omega}^-$, and for $p\in\check\omega$,
$g_{v^+_{\omega}}(p)|_{V_{\omega}}=g_{v^-_{\omega}}(p)|_{V_{\omega}}$.
Thus
\[
{g_{v^-_{\omega}}|_{V_{\omega}}\over g_{v^+_{\omega}}|_{V_{\omega}}}
\in\omega^{\|}\otimes\Gamma(U,\O^{\times}_{V_{\omega}}),
\]
and we can write
\[
{g_{v^-_{\omega}}|_{V_{\omega}}\over g_{v^+_{\omega}}|_{V_{\omega}}}
=d_{\omega}\otimes\xi_{\omega}(h)
\]
for some $\xi_{\omega}(h)$. Then we define
\[
\xi(h)=(\xi_{\omega}(h))_{\dim\omega=1}
\in\Gamma\left(U,\bigoplus_{\dim\omega=1}\O_{V_{\omega}}^{\times}\right).
\]
Now by the construction of $\xi$, the kernel of $\xi$ clearly
contains $\shHom(N,\O_V^{\times})$. On the other hand, suppose $h\in
\ker \xi$ on an open set $U$. We wish to show $h$ is induced by an
element of $\shHom(N,\O_U^{\times})$, i.e.\ there is an $h'\in
\shHom(N,\O_{U}^{\times})$ such that $h'(p)|_{V(p)}=h_p$. To
construct such an $h'$, let $x\in U$ and let $\tau$ be the largest
face of $\sigma$ such that $x\in V_{\tau}$. Let $v_1,\ldots,v_r$ be
the vertices of $\tau$, so that $x\in \bigcap_{i=1}^r V_{v_i}$.
Then for $p\in N$, $g_{v_i}(p)\in \Gamma(U,\O^{\times}_{V_{v_i}})$,
where $g_{v_i}\in M\otimes \Gamma(U,\O^{\times}_{V_{v_i}})$ is as
above. If $\omega$ is an edge of $\tau$, say with vertices $v_i$ and
$v_j$, then $g_{v_i}(p)|_{V_{\omega}}=g_{v_j}(p)|_{V_{\omega}}$, since
$h\in\ker\xi$. Thus the functions $g_{v_i}(p)|_{X_{v_i}}$ glue to
give a function $h'(p)$ in a neighbourhood of $x$. Doing this for an
open covering of $U$ gives a function $h'(p)$ defined everywhere on
$U$. Since each $g_v$ is multiplicative, so is $h'$, so that
$h'\in\shHom(N, \O_V^{\times})$. Clearly by construction,
$h'(p)|_{V(p)}=h_p$. Thus $\ker\xi=\shHom(N,\O_V^{\times})$.

We next consider the image of $\xi$. If $\tau$ is any 2-face of
$\sigma$, with vertices $v_1,\ldots,v_N$ in cyclic order,
$\omega_i$ the edge connecting $v_i$ with $v_{i+1}$ (indices modulo $N$)
and $h\in\Gamma(U,\shF)$, then
\[
\prod_{i=1}^N d_{\omega_i}\otimes \xi_{\omega_i}(h)^{\epsilon_{\tau}
(\omega_i)}|_{V_\tau}=
{g_{v_1}|_{V_{\tau}}\over g_{v_2}|_{V_{\tau}}}
\cdot
{g_{v_2}|_{V_{\tau}}\over g_{v_3}|_{V_{\tau}}}
\cdots
{g_{v_N}|_{V_{\tau}}\over g_{v_1}|_{V_{\tau}}}
= 1.
\]
Thus $\xi(h)$ satisfies condition (\ref{MC}).

Conversely, we need to show that $\xi$ surjects onto the subsheaf
defined by (\ref{MC}). We will show surjectivity on stalks, so
we fix a point $x\in V$. Suppose we
are given an element
\[
f=(f_{\omega})\in\bigoplus_{\dim\omega=1}\O^{\times}_{V_{\omega},x}
\]
satisfying (\ref{MC}). For each $\omega\subseteq\sigma$, let
$M_{\omega}=M\otimes \O_{V_{\omega},x}^{\times}$, and whenever
$\omega_1\subseteq\omega_2$, let $\varphi_{\omega_2\omega_1}:M_{\omega_1}
\rightarrow M_{\omega_2}$ be given by restriction of functions. Then
$\{M_{\omega}\}_{\omega\subseteq\sigma}$ forms a system as in \S A.1.
Thus we obtain a complex $C^{\bullet}_{\phd}$ as in \S A.2.
Furthermore, $k:=(k_{\omega})=(d_{\omega}\otimes f_{\omega})_{\omega
\subseteq\sigma\atop \dim\omega=1}$ is an element of $C^1_{\phd}$ with
$d_{\phd}(\alpha)=0$, because of condition (\ref{MC}). However, the
system $\{M_{\omega}\}$  satisfies Condition~$(*)$ of \S A.1 in exactly the same
way as argued in Example \ref{Oresolution}: giving a compatible
collection on a subset of the set of faces of $\sigma$ means giving
elements of various $M\otimes\O_{V_{\omega},\bar x}^{\times}$ which
agree on intersections of strata. These elements therefore glue, and
can be lifted to an element of $M\otimes\O_{V(\sigma),\bar
x}^{\times}$. Since Condition~$(*)$ holds, there exists an element
$g:=(g_v)_{v\in\sigma}\in C^0_{\phd}$ such that $d_{\phd}(g)=k$. 

We can now construct $h=(h_p)$ as follows. If $p\in
\partial P\setminus\{\infty\}$, then $V(p)=\bigcup_{\check v\ni p}
V_v$. We want to define $h_p$ so that $h_p|_{V_v}=g_v(p)$ whenever
$p\in \check v$. Thus we need to be able to glue the $g_v(p)$ for all
$v$ with $p\in\check v$. So if $p\in\check\omega\subset\check v$, we
need to know $g_v(p)|_{V_{\omega}}$ is independent of $v\in\omega$. To
check this, it is enough to check that this is true when $\omega$ is
dimension one, with two vertices $v_{\omega}^{\pm}$.
But
\[
{g_{v_{\omega}^-}(p)|_{V_{\omega}}\over g_{v_{\omega}^+}(p)|_{V_{\omega}}}
=k_{\omega}(p)=d_{\omega}(p)\otimes f_{\omega} 
\]
But as $p\in\check\omega$, $d_{\omega}(p)=0$, so this is $1$ and we have
the desired independence.
Thus $g_{v^+_{\omega}}(p)$
and $g_{v^-_{\omega}}(p)$ glue. Thus we obtain $h_p$.

Finally, we need to check that $\xi(h)=f$. But if $\omega\subset\sigma$
is an edge, then 
\[
d_{\omega}\otimes\xi_{\omega}(h)
={g_{v^-_{\omega}}|_{V_{\omega}}
\over g_{v^+_{\omega}}|_{V_{\omega}}}
=k_{\omega}
=d_{\omega}\otimes f_{\omega}
\]
by definition of $k_{\omega}$ and $\xi$. This gives the desired result, and
surjectivity is proved.
\qed

\begin{example} 
\label{LSexamples}
(1) Let $\sigma\subseteq M_{\RR}=\RR^n$ be the standard simplex, i.e.
the  convex hull of $e_0=(0,\ldots,0)$, $e_1=(1,0,\ldots,0),\ldots,
e_n=(0,\ldots,0,1)$. Then $V:=V(\sigma)$ can be identified with the
scheme defined by $\prod z_i=0$ in $\Spec k[z_0,\ldots,z_n]$. Let
$(f_{\omega})$ be as in the above theorem, so that $f_{\omega}$ is an
invertible function on $V_{\omega}$ for $\omega\subseteq\sigma$ of
dimension one. We need to impose the condition that it be in
$\shLS_V$. Without loss of generality, we can take the two-face $\tau$
spanned by $e_0,e_1$ and $e_2$. We get the condition that
\[
(f_{e_0e_1}f_{e_1e_2}^{-1},f_{e_1e_2}
f_{e_2e_0}^{-1},1,\ldots,1)=(1,\ldots,1).
\]
on $V_{\tau}$. But this implies that the three functions
$f_{e_0e_1}$, $f_{e_0e_2}$ and $f_{e_1e_2}$ agree on $V_{\tau}$.
One can then easily see, in fact, $\shLS_V\cong \O_D^{\times}$, where
$D$ is the singular locus of $V$.

(2) We next look at the quadruple point, i.e. $\sigma\subseteq \RR^2$
with vertices $(0,0), (1,0), (1,1)$ and $(0,1)$. Starting at $(0,0)$
and proceeding counterclockwise, we have edges
$\omega_1,\ldots,\omega_4$, which we orient positively, with generators
$d_i$ of $\omega_i^{\|}$ being $d_1=(1,0)$, $d_2=(0,1)$, $d_3=(-1,0)$,
and $d_4=(0,-1)$. Then a choice of $f_i\in\O_{V_{\omega_i}}^{\times}$
must satisfy $(f_1f_3^{-1}, f_2f_4^{-1}) =(1,1)$ on $V_{\sigma}$,
i.e. $f_1=f_3$ and $f_2=f_4$ at the point $V_{\sigma}$. Note this
gives the description of $\shLS_V$ as $\O_{D_1}^{\times}\times
\O_{D_2}^{\times}$, where $D_1= V_{\omega_1}\cup V_{\omega_3}$ and
$D_2=V_{\omega_2}\cup V_{\omega_4}$.
\qed
\end{example}

We now need to understand the global nature of the sheaf of log
smooth structures. Given an isomorphism $f:X^g\to Y^g$ of
schemes with ghost structures, there is of course a natural isomorphism
\[
f^{-1}:\shLS_Y\to f_*\shLS_X
\]
defined as follows. For $U\subseteq X$, a section $\alpha\in\Gamma(
f(U),\shLS_Y)$ defines a log smooth structure $f(U)^{\ls}\to\Spec
k^{\ls}$. Pulling back this log structure to $U$, we get a composed
map $U^{\ls}\to f(U)^{\ls}\to \Spec k^{\ls}$ giving a log smooth
structure on $U$, hence a section $f^{-1}(\alpha)\in\Gamma(U,
\shLS_X)$.
More generally, if $f:X^g\rightarrow Y^g$ is any \'etale morphism such that
the pull-back of any log structure on $Y$ of ghost type $Y^g$
to $X$ is of ghost type $X^g$, we similarly get a map $f^{-1}:
\shLS_Y\rightarrow f_*\shLS_X$.

We wish to describe $f^{-1}$ in terms of the representation of
Theorem~\ref{finalLS}. Let $B$ be an integral affine manifold with
singularities, $\P$ a toric polyhedral decomposition, $k$ an
algebraically closed field, $S=\Spec k$. Let $s\in
Z^1(\P,\shQ_{\P}\otimes\Gm(S))$ be open gluing data. We first return
to the setup of Construction~\ref{basicgluing}, with two maximal
cells $\sigma_1,\sigma_2\in\P$ and $\tau=\sigma_1\cap\sigma_2$, and
use all the notation of that construction. Thus we have open sets
$V(\tau_i)\subseteq V(\sigma_i)$ and an isomorphism
\[
\Phi_{\sigma_1\sigma_2}:V(\tau_2)\to V(\tau_1).
\]
In addition, we twist this isomorphism with automorphisms
\[
s_i:V(\tau_i)\to V(\tau_i),
\]
where $s_i$ is induced by piecewise multiplicative maps
$s_i:\partial Q_i\setminus\{\infty\}\to\Gm(k)$, $Q_i$ as in
Construction~\ref{basicgluing}, and we have
\[
\Phi_{\sigma_1\sigma_2}(s)=s_1^{-1}\circ\Phi_{\sigma_1\sigma_2}\circ
s_2.
\]
We wish to describe $\Phi_{\sigma_1\sigma_2}(s)^{-1}:
\Phi_{\sigma_1\sigma_2}(s)^*\shLS_{V(\tau_1)}\rightarrow\shLS_{V(\tau_2)}$.

To state our description of this map, we introduce some notation,
which will be generalised in Definition~\ref{stwiddle} below. Assume
we have chosen generators $d_{\omega}$ of $\Lambda_{\omega}$ for each
$1$-dimensional cell $\omega$ of $\P$, as in \S1.5.

For every vertex $v$ of $\tau$, $\check v$ is
a top-dimensional cone in $\check\tau_i^{-1}\check\Sigma_i$,
and $s_i|_{\check v\cap N_i}:\check v\cap N_i\to\Gm(k)$ extends
multiplicatively to a function
\[
s_i^v:N_i\to\Gm(k).
\]
If $\omega\subseteq\tau$ is a one-dimensional cell, then
$s_i^{v^+_{\omega}}$ and $s_i^{v^-_{\omega}}$ agree on
$\check\omega$, so that $s_i^{\omega}:= s_i^{v^-_{\omega}}/
s_i^{v^+_{\omega}}$ can be viewed as an element of
$\omega^{\|}\otimes\Gm(k)$,  so in particular we can write
\[
s_i^{\omega}=d_{\omega}\otimes\D(s_i,\omega)
\]
for some $\D(s_i,\omega)\in\Gm(k)$.  Finally, a cell
$\omega\subseteq\tau$ induces closed strata $V^i_{\omega}$ of
$V(\sigma_i)$, with $V^i_{\omega} =\Spec k[\check\omega\cap N_i]$.

Then we have the following key calculation:

\begin{theorem}
\label{gluingtheorem}
Let
\[
f^1=(f^1_{\omega})_{\omega\subseteq\tau
\atop\dim\omega=1}\in\Gamma(U,\shLS_{V(\sigma_1)})
\subseteq\Gamma\left(U,\bigoplus_{\omega\subseteq
\tau\atop\dim\omega=1} \O^{\times}_{V^1_{\omega}}\right)
\]
for $U\subseteq V(\tau_1)$. Let $f^2=(f^2_{\omega})=
\Phi_{\sigma_1\sigma_2}(s)^{-1}(f^1)$. Then for each
$\omega\subseteq\tau$ of dimension one,
\begin{eqnarray}
\label{gluingformula}
\Phi_{\sigma_1\sigma_2}(s)^*f^1_{\omega}
=z^{n_{\omega}}{\D(s_1,\omega)\over\D(s_2,\omega)}s_2
(n_{\omega})f^2_{\omega},
\end{eqnarray}
where $\Phi_{\sigma_1\sigma_2}(s)^*$ is just the ordinary pull-back
of functions and $n_{\omega}:=n_{\omega}^{\sigma_1\sigma_2}$ is
as defined in \S1.5.
\end{theorem}

\proof
As in the proof of Theorem~\ref{finalLS} we work in the Zariski
topology. Let $M'\subseteq M_i$ be the parallel transport of
$\Lambda_{\tau}$ into $M_i\cong \Lambda_y$ for $y\in \Int(\sigma_i)$.
Because this space is left invariant under monodromy, it is a
canonically determined subspace of $M_i$. Because this subspace is
canonically identified under any of the parallel transports $\psi_v$
considered in Construction~\ref{basicgluing}, we don't distinguish
between $M'\subseteq M_1$ or $M'\subseteq M_2$. We then think of the
cone $C'(\tau)=\RR_{\ge 0} (\{1\}\times\tau)\subseteq \RR\oplus
M'_{\RR}$ and $\dual{C'(\tau)} \subseteq \RR\oplus N'_{\RR}$, where
$N'\cong N_i/(M')^{\perp}$, again independent of $i$.

Put $Q'=\dual{C'(\tau)}\cap (\ZZ\oplus N')$, $\partial Q'$ as usual,
which is determined by the fan $\check\Sigma_i(\check\tau_i)$
on $N'\otimes {\RR}$ (again independent of $i$). Finally, we can
choose splittings of $N_i$ given by $\beta_i$ or $\delta_i$
\[
0\mapright{} (M')^{\perp} {\mapright{\gamma_i}\atop
\mapleft{\beta_i}} N_i {\mapright{\epsilon_i}\atop \mapleft{\delta_i}}
N'\mapright{}0.
\]
This splitting is arbitrary, but the choice will prove irrelevant.
Note that $(M')^\perp$ is canonically a subspace of $N_i$ isomorphic
to $\ZZ^d$. In the following we therefore sometimes omit $\gamma_i$
from the notation. The splittings induce isomorphisms
\begin{eqnarray*}
Q_i&=&\ZZ^d\oplus Q'\\
\partial Q_i&=&\ZZ^d\oplus \partial Q'
\end{eqnarray*}
We will use $\beta_i,\gamma_i,\delta_i,\epsilon_i$ for the maps induced
by these splittings also. We also write $\epsilon_i$ for the composition
\[
P_i\hookrightarrow Q_i\to Q'.
\]
We  use  $\bar\delta_i:Q'\to \partial Q_i$ for $\delta_i$
composed with the map $Q_i\rightarrow\partial Q_i$ given
by $p\mapsto\begin{cases} p& p\in\partial Q_i\\ \infty&\hbox{otherwise.}
\end{cases}$

We will first understand how log charts behave under pull-backs. 
Note $V(\tau_i) =\Spec k[\partial Q_i]$. Suppose we are given a chart
\[
\alpha_1:Q'\to k[\partial Q_1]
\]
of the form
\begin{eqnarray}
\label{alpha1}
\alpha_1(q)&=&h_q z^{\bar\delta_1(q)}
\end{eqnarray}
where $h_q$ is an invertible function on $\cl{\{
z^{\bar\delta_1(q)}\not=0\}}$ and $h_{q_1}\cdot h_{q_2}=h_{q_1+q_2}$
where defined. Recall the map $\Phi_{\sigma_1\sigma_2}:V(\tau_2) \to
V(\tau_1)$ is induced by $\phi:\partial Q_1\to \partial Q_2$ defined
in Construction~\ref{basicgluing}. We also denote by $\phi:k[\partial
Q_1]\to k[\partial Q_2]$ the induced map.

We then pull back via $\Phi_{\sigma_1\sigma_2}(s)$ the log structure
by composing this chart with $\Phi_{\sigma_1\sigma_2}(s)^*$, i.e. set
$\alpha_2$ to be the composition
\[
\alpha_2:Q'\mapright{\alpha_1} k[\partial Q_1] 
\mapright{\Phi_{\sigma_1\sigma_2}(s)^*} k[\partial Q_2].
\]
Thus
\[
\alpha_2(q)=\Phi_{\sigma_1\sigma_2}(s)^*(z^{\bar\delta_1(q)}h_q)=
s_1(\bar\delta_1(q))^{-1}s_2(\phi(\bar\delta_1(q)))
z^{\phi\bar\delta_1(q)}\Phi_{\sigma_1\sigma_2}(s)^*(h_q).
\]
We can write
\[
\phi(\bar\delta_1(q))=\gamma_2\beta_2\phi\bar\delta_1(q)
+\delta_2\epsilon_2\phi\bar\delta_1(q).
\]
By monodromy invariance of $N'$, it follows that
$\epsilon_2\circ\phi\circ\bar\delta_1:Q'\to \partial Q'$ is just the
projection. Let
\[
\chi:=\beta_2\circ\phi\circ\bar\delta_1:Q'\to (M')^{\perp}
\subseteq \partial Q_2.
\]
Then
\begin{eqnarray}
\label{alpha2}
\alpha_2(q)&=&
\left[z^{\chi(q)}\Phi_{\sigma_1\sigma_2}(s)^*(h_q)s_1(\bar\delta_1(q))^{-1}
s_2(\phi(\bar\delta_1(q)))\right]
\cdot z^{\bar\delta_2(q)}
\end{eqnarray}

Now we need to understand how these charts correspond to sections
of $\shLS_{V(\sigma_i)}$ over $V(\tau_i)$ using the representation
given in Theorem~\ref{finalLS}. If we are given a chart
\[
\alpha_i:Q'\to k[\partial Q_i]
\]
with $\alpha_i(q)=h^i_q z^{\bar\delta_i(q)}$, 
let
\[
\alpha'_i:P_i\to k[\partial Q_i]
\]
be given by 
\[
\alpha'_i(p)=
\begin{cases}
h^i_{\epsilon_i(p)} z^p,&\hbox{if $p$ is in a face of $P_i$ containing
$\check\tau_i$}\\ 0&\hbox{otherwise.}
\end{cases}
\]
It is then easy to see that the charts $\alpha_i'$ and $\alpha_i$
induce the same log structure on $V(\tau_i)$. The chart $\alpha_i'$
then coincides with a section $h^i=(h^i_{\epsilon_i(p)})_{p\in
\partial Q_i\setminus\{\infty\}}$ of the sheaf $\shF_i$
(corresponding to the sheaf $\shF$ in the language of
Theorem~\ref{finalLS}). For $v\in\tau$ a vertex, we let $g^i_v\in
M_i\otimes\O_{V^i_v}^{\times}$ be the multiplicative extension of the
map $\check v\cap N\ni p \mapsto h^i_{\epsilon_i(p)}
|_{V^i_v}\in\O^{\times}_{V^i_v}$. Then the map
\[
\xi^i:\shF\to\shLS_{V(\sigma_i)}\subseteq
\bigoplus_{\omega\subseteq\sigma_i\atop \dim\omega=1}
\O_{V^i_{\omega}}^{\times}
\]
is given on $V^i_{\omega}$ by
\[
d_{\omega}\otimes\xi^i_{\omega}(h)=
{g^i_{v^-_{\omega}}|_{V^i_{\omega}}
\over 
g^i_{v^+_{\omega}}|_{V^i_{\omega}}}.
\]
Now we are ready to compare the answers for $\alpha_1$ and $\alpha_2$
given by formulae (\ref{alpha1}) and (\ref{alpha2}), by applying
this, for
\[
h^1_q=h_{q}
\]
and
\[
h^2_q=\left[\Phi_{\sigma_1\sigma_2}(s)^*(h_{q})\right] z^{\chi(q)}
s_1(\bar\delta_1(q))^{-1}s_2(\phi(\bar\delta_1(q))),
\]
for $q\in Q'$. If $v$ is any vertex of $\tau$ we obtain, for $p\in
\check v\cap N_1$, 
\[
g^1_v(p)=h^1_{\epsilon_1(p)}
\]
and for $p\in\check v\cap N_2$,
\[
g^2_v(p)=\left[\Phi_{\sigma_1\sigma_2}(s)^*(h^1_{\epsilon_2(p)})\right]
z^{\chi(\epsilon_2(p))}\frac{s_2(\phi(\delta_1(\epsilon_2(p))))}{
s_1(\delta_1(\epsilon_2(p)))}.
\]
Here we are using the canonical identification of the quotients  $N'$
of $N_1$ and $N_2$ to make sense of $h^1_{\epsilon_2(p)}$.  We wish
to compare $g^1_v$ and $g^2_v$, so we will for the moment  choose any
vertex $w$ of $\tau$, and use the map
\[
\psi_w^{-1}\times\Phi_{\sigma_1\sigma_2}(s)^*:M_1\times\O^{\times}_{V^1_v\cap V(\tau_1)}
\to M_2\times\O^{\times}_{V^2_v\cap V(\tau_2)}.
\]
Then for $p\in \check v\cap N_2$,
\begin{eqnarray*}
\lefteqn{[\psi_w^{-1}\times\Phi_{\sigma_1\sigma_2}(s)^*(g^1_v)](p)\big/
g^2_v(p)}\hspace{1cm}\\
&=&(g^2_v(p))^{-1}\left[\Phi_{\sigma_1\sigma_2}(s)^*
(g^1_v((\psi_w^{-1})^t(p)))\right]\\
&=&\Big\{\left[\Phi_{\sigma_1\sigma_2}(s)^*
(h^1_{\epsilon_2(p)})\right]^{-1} z^{-\chi(\epsilon_2(p))}
\frac{s_1\big(\delta_1(\epsilon_2(p))\big)}{
s_2\big(\phi(\delta_1(\epsilon_2(p)))\big)}\Big\}\cdot
\Phi_{\sigma_1\sigma_2}(s)^*\big(h^1_{\epsilon_1((\psi_w^{-1})^t(p))}
\big)\\
&=&
z^{-\chi(\epsilon_2(p))}
\frac{s_1\big(\delta_1(\epsilon_2(p))\big)}{
s_2\big(\phi(\delta_1(\epsilon_2(p)))\big)},
\end{eqnarray*}
where we used $\epsilon_1\big((\psi_w^{-1})^t(p)\big)= \epsilon_2(p)$.
This must be extended linearly to all of $N_2$, which we do as
follows. Recall that $\chi=\beta_2\circ\phi\circ\bar\delta_1$, and on
$\check v\cap N_1$, $\phi$ is linear, given by parallel transport
$\psi^t_v:N_1\to N_2$.  Furthermore, on $\check v\cap N_i$, $s_i$
extends to give the group homomorphism
\[
s_i^v:N_i\to\Gm(k).
\]
Thus for $p\in N_2$, we have
\[
\frac{[\psi_w^{-1}\times\Phi_{\sigma_1\sigma_2}(s)^*(g^1_v)](p)}
{g^2_v(p)}\ =\ 
z^{-\beta_2\circ\psi_v^t\circ\delta_1\circ\epsilon_2(p)}
s_1^v(\delta_1(\epsilon_2(p)))s_2^v(\psi^t_v(\delta_1(\epsilon_2(p))))^{-1}.
\]
We now compare $d_{\omega}\otimes\xi^1_{\omega}(h)$ and 
$d_{\omega}\otimes \xi_{\omega}^2(h)$ whenever $\omega$ is an
edge of $\tau$. Note that $\omega^{\|}\subseteq M'$, so
we can view $\omega^{\|}$ equally as a subspace of $M_1$ or
$M_2$, but to be precise, we use parallel transport 
$\psi_{v^-_{\omega}}^{-1}:M_1\to M_2$ to make the identification;
this is a convenient arbitrary choice. We then make
the comparison using
\[
\psi_{v^-_{\omega}}^{-1}\times\Phi_{\sigma_1\sigma_2}(s)^*:
\omega^{\|}\times \O^{\times}_{V^1_{\omega}}
\to
\omega^{\|}\times \O^{\times}_{V^2_{\omega}}.
\]
Then for $p\in N_2$,
\begin{eqnarray*}
&&{\left[(\psi_{v^-_{\omega}}^{-1}\times
\Phi_{\sigma_1\sigma_2}(s)^*)(
d_{\omega}\otimes\xi^1_{\omega}(h))\right](p)\over (d_{\omega}
\otimes\xi^2_{\omega}(h))(p)}\\
&=&{\left[(\psi_{v^-_{\omega}}^{-1}\times\Phi_{\sigma_1\sigma_2}(s)^*)
(g^1_{v^-_{\omega}}\cdot (g^1_{v^+_{\omega}})^{-1})\right](p)
\over
g^2_{v^-_{\omega}}(p)\cdot (g^2_{v^+_{\omega}})^{-1}(p)}\\
&=&
\Bigg(\frac
{z^{-\beta_2\circ\psi_{v^-_{\omega}}^t \circ\delta_1\circ\epsilon_2(p)}}
{z^{-\beta_2\circ\psi_{v^+_{\omega}}^t \circ\delta_1\circ\epsilon_2(p)}}
\Bigg)\cdot\Bigg(\frac
{s_1^{v^-_{\omega}}(\delta_1(\epsilon_2(p)))}
{s_1^{v^+_{\omega}}(\delta_1(\epsilon_2(p)))}
\Bigg)\cdot\Bigg(\frac
{s_2^{v^-_{\omega}}(\psi^t_{v^-_{\omega}}(\delta_1(\epsilon_2(p))))^{-1}}
{s_2^{v^+_{\omega}}(\psi^t_{v^+_{\omega}}(\delta_1(\epsilon_2(p))))^{-1}}
\Bigg).
\end{eqnarray*}
The first factor is
\[
z^{\beta_2\circ(\psi_{v^+_{\omega}}^t-\psi_{v^-_\omega}^t)\circ\delta_1\circ
\epsilon_2(p)},
\]
and using $\delta_1\circ\epsilon_1(p)=p-\gamma_1\circ\beta_1(p)$ for
$p\in N_1$ and $\epsilon_2(p)= \epsilon_1
((\psi_{v^-_{\omega}}^{-1})^t(p))$, we see that
\[
\beta_2\circ(\psi_{v^+_{\omega}}^t-\psi_{v^-_\omega}^t)\circ\delta_1\circ
\epsilon_2(p)
=\beta_2\circ(\psi_{v^+_{\omega}}^t-\psi_{v^-_\omega}^t)
((\psi_{v^-_{\omega}}^{-1})^t(p))=\beta_2\circ (\psi^t_{v^+_{\omega}}\circ
(\psi^{-1}_{v^-_{\omega}})^t-I)(p).
\]
In this computation a term $\gamma_1\circ\beta_1(p)\in M'$ cancels
because $\psi_{v_\omega^+}^t|_{M'}= \psi_{v_\omega^-}^t|_{M'}$. In
addition, the factor $\beta_2$ is unnecessary, as 
$(\psi^t_{v^+_{\omega}}\circ(\psi^{-1}_{v^-_{\omega}})^t-I)(p)$  is
already contained in $(M')^{\perp}$. But in the notation of \S1.5,
$(T_{\omega}^{\sigma_1\sigma_2})^t= \psi_{v^+_{\omega}}^t \circ 
(\psi_{v^-_{\omega}}^t)^{-1}$. Now
\[
(T_{\omega}^{\sigma_1\sigma_2})^t(n)=
n+\langle n,d_{\omega}\rangle n_{\omega}
\]
by the definition of $n_{\omega}$. Thus the first factor is
\[
z^{((T_{\omega}^{\sigma_1\sigma_2})^t-I)(p)}
=z^{\langle p,d_{\omega}\rangle n_{\omega}}.
\]
The second factor is
\[
s_1^{\omega}(\delta_1(\epsilon_1((\psi^{-1}_{v^-_{\omega}})^t(p)))
\ =\ s_1^{\omega}((\psi^{-1}_{v^-_{\omega}})^t(p)-\gamma_1\beta_1
(\psi^{-1}_{v^-_{\omega}})^t(p))
\ =\ s_1^{\omega}((\psi^{-1}_{v^-_{\omega}})^t(p)).
\]
The last factor is similarly
\[
{
s_2^{v^+_{\omega}}(\psi^t_{v^+_{\omega}}
(\psi^{-1}_{v^-_{\omega}})^t(p))
s_2^{v^+_{\omega}} (\psi^t_{v^+_{\omega}} \gamma_1\beta_1
(\psi^{-1}_{v^-_{\omega}})^t (p))^{-1}
\over
s_2^{v^-_{\omega}}(p) s_2^{v^-_{\omega}} (\psi^t_{v^-_{\omega}}
\gamma_1\beta_1 (\psi^{-1}_{v^-_{\omega}})^t (p))^{-1}
}.
\]
Using
\begin{eqnarray}
\label{agree}
\hbox{$s_2^{v^+_{\omega}}=s_2^{v^-_{\omega}}$ and
$\psi_{v^+_{\omega}}^t= \psi_{v^-_{\omega}}^t$ on
$(\omega^{\|})^{\perp}\supseteq (M')^{\perp}$,}
\end{eqnarray}
the last factors cancel and this is
\[
{s_2^{v^+_{\omega}}((T^{\sigma_1\sigma_2}_{\omega})^t(p))
\over s_2^{v^-_{\omega}}(p)}\ =\ 
{s_2^{v^+_{\omega}}(p+\langle p,d_{\omega}\rangle n_{\omega})
\over s_2^{v^-_{\omega}}(p)}\ =\ 
s_2^{\omega}(p)^{-1}(s_2^{v^+_{\omega}}
(n_{\omega}))^{\langle p,d_{\omega}\rangle}
\]
Putting this together,
\begin{eqnarray*}
{\left[(\psi_{v^-_{\omega}}^{-1}\times \Phi_{\sigma_1\sigma_2}(s)^*)(
d_{\omega}\otimes\xi^1_{\omega}(h)) \right](p) \over
(d_{\omega}\otimes\xi^2_{\omega}(h))(p)}&
=&z^{\langle p,d_{\omega}\rangle n_{\omega}}s_1^{\omega} 
((\psi_{v^-_{\omega}}^{-1})^t(p))s_2^{\omega}(p)^{-1}
s_2^{v^+_{\omega}}(n_{\omega})^{\langle p,d_{\omega}\rangle}.
\end{eqnarray*}
Thus if $f^i_{\omega}=\xi^i_{\omega}(h)$ and in view of $s_i^\omega=
d_\omega\otimes\D(s_i,\omega)$ we see that
\begin{eqnarray*}
{\Phi_{\sigma_1\sigma_2}(s)^*f^1_{\omega}\over f^2_{\omega}}
&=&z^{n_{\omega}} {\D(s_1,\omega)\over \D(s_2,\omega)}
s_2^{v^+_{\omega}}(n_{\omega}).
\end{eqnarray*}
Note that since $n_{\omega}\in (M')^{\perp}$, $s_2^{v^+_{\omega}}
(n_{\omega})= s_2(n_{\omega})$. This formula is now completely
independent of the choices of splitting.
\qed
\medskip

We are now ready to summarize the discussion of this section.
First, we define the global versions of the data appearing in
Theorem~\ref{gluingtheorem}.

\begin{definition} 
\label{stwiddle}
1)\ \ Let $B$ be an integral affine manifold with singularities and $\P$ a
toric polyhedral decomposition, and suppose $d_{\omega}$ has been
chosen as in \S1.5 for each $\omega\in\P$, $\dim\omega=1$. If
$\tau\in\P$, $s\in \check{PM}(\tau)$ (see Definition~\ref{PMdef}),
$h:v\to\tau$ with $v$ a vertex, we define $s^h\in
\Lambda_v\otimes\Gm(k)$ to be the germ of $s\in
\Gamma(\tilde\tau,\pi^* (\shQ_{\P}\otimes\Gm(k)))$ at the vertex $v_h$
of $\tilde\tau$ corresponding to $h$. 

\noindent
2)\ \ If $\dim\tau\ge 1$, and we are given morphisms $e:\omega\to\tau$,
$f:\tau\to\sigma$ with $\dim\omega=1$, $\sigma\in\P_{\max}$,
giving a diagram
\[
\xymatrix@C=30pt
{v^+_{\omega}\ar@/^/[rrd]^{h^+}\ar[rd]&&&\\
&\omega\ar[r]^(.3){e}&\tau\ar[r]^{f}&\sigma\,,\\
v^-_{\omega}\ar@/_/[rru]_{h^-}\ar[ru]&&&\\
}
\]
we can parallel transport
$s^{h^+}\in\Lambda_{v^+_{\omega}}\otimes\Gm$ and
$s^{h^-}\in\Lambda_{v^-_{\omega}}\otimes\Gm$ to $\Lambda_y\otimes
\Gm$ with $y\in \Int(\sigma)$ via $f\circ h^+$ and $f\circ h^-$. Then
we define $\D(s,e,f)\in \GG_m$ (or $\D(s,\omega,\sigma)$ when $e$ and
$f$ are understood) by
\[
d_{\omega}\otimes \D(s,e,f)={s^{h^-}\over s^{h^+}}.
\]
\qed
\end{definition}

One should think of $\D(s,e,f)$ as monodromy information for the
gluing data $s$.

\begin{remark}
\label{twiddlechange2}
We emphasize the dependence of $\D(s,e,f)$ on $f$ (or $\sigma$).
Given $f_i:\tau\to\sigma_i \in\P_{\max}$, $i=1,2$, $e:\omega\to\tau$,
$s\in \check{PM}(\tau)$, we have an identification $\psi:\Lambda_{y_1}
\otimes\Gm\to \Lambda_{y_2}\otimes\Gm$ for $y_i\in \Int(\sigma_i)$ by
parallel transport from $y_1$ to $v^-_{\omega}$ along $f_1\circ h^-$
and $v_{\omega}^-$ to $y_2$ along $f_2\circ h^-$. Then under this
identification, the parallel transport of $s^{h^-}
\in\Lambda_{v_{\omega}^-}\otimes\Gm$ to $\Lambda_{y_1}\otimes\Gm$ and
$\Lambda_{y_2}\otimes\Gm$ are identified by $\psi$. However, if
$s_i^{h^+}\in\Lambda_{y_i} \otimes\Gm$ is the parallel transport of
$s^{h^+}\in \Lambda_{v^+_{\omega}} \otimes\Gm$ to
$\Lambda_{y_i}\otimes \Gm$ along $f_i\circ h^+$, then 
\[
\psi(s_1^{h^+})=T_{\omega}^{f_1\circ e,f_2\circ e}(s_2^{h^+}).
\]
Thus 
\[
\psi(s_1^{h^-}/s_1^{h^+})=s_2^{h^-}/T_{\omega}^{f_1\circ e,f_2\circ e}(s_2^{h^+}).
\]
If $n\in\check\Lambda_{y_2}$ is such that $\langle
n,d_{\omega}\rangle=1$, then $\langle \psi^t(n),d_{\omega}\rangle=1$
also, and we see that
\begin{eqnarray*}
\D(s,e,f_1)&=& s_1^{h^-}(\psi^t(n))/ s_1^{h^+}(\psi^t(n)) =
s_2^{h^-}(n)/ (T_{\omega}^{f_1\circ e,f_2\circ e}(s_2^{h^+}))(n)\\
&=&s_2^{h^-}(n)/s_2^{h^+}(n+\langle n,d_{\omega}\rangle
n_{\omega}^{f_1\circ e,f_2\circ e})
\ =\ \D(s,e,f_2)s_2^{h^+}(n_{\omega}^{f_1\circ e,f_2\circ e})^{-1}\\
&=&\D(s,e,f_2)s(n_{\omega}^{f_1\circ e,f_2\circ e})^{-1}.
\end{eqnarray*}
The last equality follows from the fact that $n_{\omega}^{f_1\circ e,
f_2\circ e}\in\Lambda_{\tau}^{\perp}$, and hence $s_i^h
(n_{\omega}^{f_1\circ e,f_2\circ e})$ is independent of $h$ and $i$.

This explains the asymmetry in the role of $s_1$ and $s_2$ in the gluing
formulae (\ref{gluingformula}) above or (\ref{gluingform2}) below.
\qed
\end{remark}

Now suppose we are given $(B,\P)$ and $s$. How do we specify a
section of $\shLS_{X_0(B,\P,s)}$ over an open set $U\subseteq
X_0(B,\P,s)$? We have for each $\sigma\in\P_{\max}$ the projection
$p_{\sigma}: V(\sigma)\to X_0(B,\P,s)$. So a section
$f\in\Gamma(U,\shLS_{X_0(B,\P,s)})$ pulls back to sections
\[
f_{\sigma}=p_{\sigma}^{-1}(f)\in\Gamma(p_{\sigma}^{-1}(U),
\shLS_{V(\sigma)}).
\]
Furthermore, for any component $\foR_{g_1g_2}$ of $\foR$ (defined after
Definition~\ref{fanopen}),
$g_i:\tau\to \sigma_i$, we have
\[
f_{\sigma_2}=\Phi_{g_1g_2}(s)^{-1}(f_{\sigma_1})
\]
on $q_2^{-1}(U)\cap V(\tau_2)\subseteq V(\sigma_2)$.

Of course, by the sheaf gluing axiom, the converse is true:

\begin{theorem}
\label{globalgluing}
Let $B$ be an integral affine manifold with singularities with a
toric polyhedral decomposition $\P$ and open gluing data $s$ over
$\Spec k$. Choose $d_{\omega}$'s as in \S1.5. Let $U\subseteq
X_0(B,\P,s)$ be an open subset. Then to give an element
$f\in\Gamma(U,\shLS_{X_0(B,\P,s)})$, it is enough to give sections
for each $\sigma\in\P_{\max}$
\begin{eqnarray*}
f_{\sigma}=(f_{\sigma,e})_{e:\omega\to\sigma
\atop \dim\omega=1}&\in&\Gamma(p_{\sigma}^{-1}(U),\shLS_{V(\sigma)})\\
&\subseteq&\bigoplus_{e:\omega\to\sigma\atop\dim\omega=1}
\Gamma(p_{\sigma}^{-1}(U),\O_{V_e}^{\times})
\end{eqnarray*}
satisfying the following property: For every connected component
$\foR_{g_1g_2}$ of $\foR$, $g_i:\tau\to\sigma_i\in \P_{max}$, and for
every $h:\omega\to\tau$ with $\dim\omega=1$ and  $e_i=g_i\circ h$, we
have
\begin{eqnarray}
\label{gluingform2}
\Phi_{g_1g_2}(s)^*(f_{\sigma_1,e_1})&=&z^{n_{\omega}^{e_1e_2}}
{\D(s_{g_1},h,g_1)\over \D(s_{g_2},h,g_2)} 
s_{g_2}(n_{\omega}^{e_1e_2})f_{\sigma_2,e_2}
\end{eqnarray}
on $p_{\sigma_2}^{-1}(U)\cap V(\tau_2)$.
\end{theorem}

\proof
This is just a summary of the discussions, especially Theorem
\ref{gluingtheorem}, above.

Note that there is no need to define the sign vectors $\epsilon_\tau$
from Definition~\ref{signvector} globally. The reason is that these
are only used to define $\shLS_{V(\sigma)}$ as a subsheaf of
$\bigoplus_{e:\omega\to\sigma} \O_{V_e}^{\times}$, and this is
independent of choices for the sign vectors. In fact, any
$\epsilon_\tau$ is well-defined up to an overall change of sign and
this does not affect the defining equation~(\ref{MC}) for
$\shLS_{V(\sigma)}\subset
\bigoplus_{e:\omega\to\sigma}\O_{V_e}^{\times}$.
\qed\medskip

It is also helpful to give a partial description of $\shLS_{X_0(B,\P,s)}$
in terms of closed gluings.

\begin{theorem}
\label{q_tau2}
Let $B$ be an integral affine manifold with singularities, $\P$ a
toric polyhedral decomposition, and $s$ open gluing data. Choose
$d_{\omega}$'s as in \S1.5. Then for each $\omega\in\P$ of dimension
one, there is an $\O^{\times}_{X_{\omega}}$-torsor
$\shN^{\times}_{\omega}$ such that $\shLS_{X_0(B,\P,s)}$ is a
subsheaf of 
\[
\bigoplus_{\omega\in\P\atop \dim\omega=1}
q_{\omega*} \shN^{\times}_{\omega},
\]
where $q_{\omega}: X_{\omega}\to X_0(B,\P,s)$ are the maps of
Lemma~\ref{q_tau}. Furthermore, if $\shN_{\omega}$ is the line
bundle on $X_{\omega}$ associated to $\shN^{\times}_{\omega}$ (in
particular $\shN^{\times}_{\omega}$ is the sheaf of nowhere
vanishing sections of $\shN_{\omega}$) then $\shN_{\omega}$
corresponds to the piecewise linear function $\psi_{\omega}$ on the
fan $\Sigma_{\omega}$ in $\shQ_{\omega,\RR}$, well-defined up to
linear functions, defined in Remark~\ref{psiomega}.
\end{theorem}
\proof 
It is easy to see that $\shLS_{X_0(B,\P,s)}$ is a subsheaf of
$\bigoplus_{\omega\in\P\atop \dim\omega=1} q_{\omega*}
\shN^{\times}_{\omega}$, where $\shN_{\omega}^{\times}$ is defined by
the following gluing data for the open cover $\{V_e|e\in
\coprod_{\sigma\in\P_{\max}}\Hom(\omega,\sigma)\}$, with
$V_e\subseteq X_{\omega}$ the open subset corresponding to the
maximal cone $K_e$ of $\Sigma_{\omega}$ as usual. We glue
$\O^{\times}_{V_{e_1}}$ and $\O^{\times}_{V_{e_2}}$ via the
identification
\[
\O^{\times}_{V_{e_2}}|_{V_{e_1}\cap V_{e_2}}\to 
\O^{\times}_{V_{e_1}}|_{V_{e_1}\cap V_{e_2}}
\]
given by
\[
f\longmapsto z^{n^{e_1e_2}_{\omega}}\cdot s'\cdot f,
\]
where $s'\in\Gm(k)$ depends on $s$. Since $X_{\omega}$ is toric, we
can just as well assume $s'=1$, i.e.\ take the trivial open gluing
data and we will get an isomorphic $\O_{X_{\omega}
}^{\times}$-torsor.  Then the result that $\psi_{\omega}$ determines
$\shN_{\omega}$ follows from the standard recipe for obtaining a line
bundle from a piecewise linear function (see \cite{Oda}, pg. 69,
keeping in mind the sign convention of Remark \ref{convention}).
\qed

\begin{corollary}
\label{nomonodromy}
Let $B$ be an integral affine manifold with singularities,
with a toric polyhedral decomposition $\P$. If there exists a log
smooth structure on $X_0(B,\P,s)$ for some open gluing data $s$,
then the integral affine structure on $B_0\subseteq B$ can
be extended to $B$, so $B$ in fact has no singularities.
\end{corollary}

\proof If there is a log smooth structure, then each of the torsors
$\shN^{\times}_{\omega}$ has sections, hence is a trivial torsor.
Thus $\psi_{\omega}$ is a linear function. By the definition of $\psi_{\omega}$,
$n_{\omega}^{e_1e_2}$ is the difference of linear extensions of $\psi_{\omega}$
on different cones of $\Sigma_{\omega}$, hence
$n^{e_1e_2}_{\omega}=0$ for all $e_1,e_2,\omega$, and hence
by Proposition~\ref{extaff}, the affine structure extends to $B$.
\qed\medskip

In general, we now see that the existence of singularities on $B$ implies
there does not exist a global section of $\shLS_{X_0(B,\P,s)}$. Therefore,
to proceed, we have to allow the log smooth structure to break
down. We shall study this in the next section.

\begin{example}
\label{normal crossings}
If $X_0(B,\P,s)$ is normal crossings (i.e. every $\sigma\in\P_{max}$
is a standard simplex), then it follows from
Example~\ref{LSexamples}, (1), that $\shLS_{X_0(B,\P,s)}$ is an
$\O_D^{\times}$-torsor, where $D=\Sing(X_0(B,\P,s))$. One can show
\cite{Friedman}, \cite{Kawamata; Namikawa 1994} that
$\shLS_{X_0(B,\P,s)}$ is the $\O_D^{\times}$-torsor associated to the
line bundle $\shExt^1(\Omega^1_{X_0(B,\P,s)},\O_{X_0(B,\P,s)})$, the
local $T^1$-sheaf. There exists a log smooth structure on 
$X_0(B,\P,s)$ if and only if this is the trivial line bundle. In any
normal crossings case, $q_{\omega}^*\shLS_{X_0(B,\P,s)}$ is the 
$\O_{X_{\omega}}^{\times}$-torsor determined by  Theorem~\ref{q_tau2}
using monodromy data near $\omega$. Explicitly, in dimension two, if
$\omega\in\P$ is an edge containing a singular point of the affine
structure with monodromy $\begin{pmatrix}1&n\\0&1 \end{pmatrix}$ in a
suitable basis (chosen so that $n>0$ if and only if the singularity
is positive in the sense of Definition~\ref{positive def}) then
$\shN_{\omega}=\O_{\PP^1}(n)$. Even if $\shN_{\omega}=\O_{\PP^1}$ for
all $\omega\in \P$ of dimension one, $\shLS_{X_0(B,\P,s)}$ may be
non-trivial thanks to non-trivial open gluing data.

In the example given in the introduction, $n=4$ for each edge.
Let $\X_0$ carry the pull-back of the divisorial log structure
$\X_0\subseteq\X$. Then the map $\X_0^{\dagger}\rightarrow
\Spec k^{\dagger}$ fails to be log smooth precisely at the $24=4\times 6$
points where the total space $\X$ is singular, as the log structure is
not even fine at these points (Example~\ref{non-fine example}). These
24 points are the zeros of a section of the local $T^1$ sheaf of $\X_0$.
\end{example}

\section{Toric degenerations}
\label{section4}

We are now ready to study the fundamental objects of interest in our
program. These will be certain degenerations of Calabi-Yau manifolds
whose central fibres are of the form $X_0(B,\P,s)$. In
addition, these fibres come along with natural log structures of the
form studied in the previous section.

We gave a definition of toric degeneration in \cite{Announce},
Section~1; the one given here is slightly more general. From now on,
we will always work over an algebraically closed field $k$. The
algebraically closed hypothesis could be removed from this section at
the expense of a bit of additional complexity, but it's vital in the
proof of Theorem~\ref{lifted}.

\begin{definition}
\label{toric degen}
Let $R$ be a discrete valuation ring with 
residue class field $k$. A toric degeneration of Calabi-Yau
varieties over $R$ is a proper normal algebraic space $\X$ flat over
$\shS:=\Spec R$ satisfying the following properties:
\item{(1)} The generic fibre $\X_{\eta}$ is an irreducible normal
variety over $\eta$.
\item{(2)} If $\nu:\tilde\X_0\to\X_0$ is the normalization,
then $\tilde\X_0$ is a disjoint union of toric varieties,
the conductor scheme $C\subseteq\tilde\X_0$ is reduced
and the map $C\to\nu(C)$ is unramified and generically
two-to-one. The square
\[\begin{CD}
C@>>> \tilde\X_0\\
@VVV @VV{\nu}V\\
\nu(C)@>>> \X_0
\end{CD}\]
is cartesian and cocartesian.
\item{(3)} $\X_0$ is Gorenstein, and the conductor locus $C$
restricted to each irreducible component of $\tilde\X_0$ is the union
of all toric Weil divisors. (By the discussion before
Theorem~\ref{dualizingsheaf},  this a substitute for the statement
$\nu^*\omega_{\X_0}\cong\O_{\tilde\X_0}$.)
\item{(4)} There exists a closed subset $Z\subseteq\X$ of relative
codimension $\ge 2$ such that $Z$ satisfies the following properties:
$Z$ does not contain the image under $\nu$ of any toric stratum of
$\tilde\X_0$, and for any geometric point $\bar x\to \X\setminus Z$,
there is an \' etale neighbourhood $U_{\bar x}\to \X\setminus Z$ of
$\bar x$, an affine toric variety $Y_{\bar x}$, a regular function
$f_{\bar x}$ on $Y_{\bar x}$ given by a monomial, a choice of
uniformizing parameter of $R$ giving a map $k[\NN]\rightarrow R$, and
a commutative diagram
\[
\begin{matrix}
U_{\bar x}&\mapright{}&Y_{\bar x}\cr
\mapdown{f|_{U_{\bar x}}}&&\mapdown{f_{\bar x}}\cr
\Spec R&\mapright{}&\Spec k[\NN]\cr
\end{matrix}
\]
such that the induced map $U_{\bar x}\to \Spec R\times_{\Spec k[\NN]}
Y_{\bar x}$ is smooth. Furthermore, $f_{\bar x}$ vanishes on each
toric divisor of $Y_{\bar x}$.
\qed
\end{definition}

Geometrically $(2)$ and $(3)$ say that the central fibre is obtained from a
disjoint union of toric varieties by identifying pairs of irreducible
toric Weil divisors. Glued Weil divisors may lie on the same toric
variety, but it is not permitted to glue an irreducible toric Weil
divisor to itself (no pinch points).  The restriction of $\nu$ to a
toric stratum in $\tilde\X_0$ is finite and generically injective,
hence the normalization of its image. We refer to these images in
$\X_0$ also as \emph{toric strata}.

\begin{example}
\label{toric degenerations examples}
1)\ \ A straightforward source of toric degenerations are
degenerations of hypersurfaces in toric varieties. As in
Example~\ref{polytope}, let $\Xi$ be a reflexive polytope with $0\in
\Xi$ the unique interior point, and let $\Xi^*\subseteq N_{\RR}$ be
the Batyrev dual of $\Xi$ (Example~\ref{Batyrevdual}). Let
$(\PP_{\Xi^*}, \O_{\PP_{\Xi^*}}(1))$ be the projective toric variety
with Newton polytope $\Xi^*$. Define a family $\X\subseteq
\PP_{\Xi^*}\times \shS$ with $\shS=\Spec k[t]_{(t)}$ by the equation
\[
z^0+\sum_{n\in\Xi^*\cap N} a_nt z^n=0;
\]
here $z^n$ run over a basis of sections of $\O_{\PP_{\Xi^*}}(1)$. The
coefficients $a_n\in k$ are chosen to be general. The section $z^0$
vanishes precisely once on each toric divisor, so $\X_0
=\partial\PP_{\Xi^*}$. It is then easy to check $\X\to \shS$ is a
toric degeneration.

We will discuss in Example~\ref{degenexamples} below how exactly this
toric degeneration relates to the toric polyhedral decomposition of
$\partial \Xi$ introduced in Example~\ref{decompexamp},~(2).

\noindent
2)\ \ Here is an example due to Aspinwall and Morrison \cite{AM}. Take
a family $\Y\subseteq\PP^4\times\shS$, $\shS=\Spec k[t]_{(t)}$, given
by the equation
\[
t(z_0^5+\cdots+z_4^5)+z_0z_1z_2z_3z_4=0.
\]
We then divide $\Y$ by the group action of $\ZZ_5\times\ZZ_5$, with
generators acting by $(z_0,\ldots,z_4)\mapsto (z_0,\xi
z_1,\ldots,\xi^4 z_4)$ and $(z_0,\ldots,z_4)\mapsto
(z_1,\ldots,z_4,z_0)$. Here $\xi$ is a primitive fifth root of unity.
Set $\X=\Y/\ZZ_5\times\ZZ_5$; then $\X\to\shS$ is a toric
degeneration.

This example will turn out to be related to the affine manifold from
Example~\ref{aspmorr}.
\qed
\end{example}

For us, the key structure of a toric degeneration is the central fibre,
$\X_0$, along with the log structure on $\X_0$ obtained by pulling
back the log structure on $\X$ induced by the inclusion $\X_0\subseteq\X$
as in Example~\ref{prime example}. Thus we can define the logarithmic
analogue of a Calabi-Yau variety as follows.

\begin{definition}
\label{logCY}
A toric log Calabi-Yau space is a proper reduced logarithmic space $X^{\ls}$
along with a morphism of log spaces $X^{\ls}\to \Spec k^{\ls}$
to the standard log point with the following properties:
\item{(1)} If $\nu:\tilde X\to X$ is the normalization, then
$\tilde X$ is a disjoint union of toric varieties,
the conductor scheme $C\subseteq\tilde X$ is reduced
and the map $C\to\nu(C)$ is generically
two-to-one and unramified. The square formed by $C$, $\nu(C)$, $\tilde
X$ and $X$ is cartesian and cocartesian.
\item{(2)} $X$ is Gorenstein,
and the conductor locus $C$ restricted to each irreducible component
of $\tilde X$ is the union of all toric Weil divisors.
(As before, this
a substitute for the statement $\nu^*\omega_{X}\cong\O_{\tilde X}$.)
\item{(3)} There exists a closed subset $Z\subseteq \Sing(X)$ of
relative codimension $\ge 2$ and not containing any toric stratum
(called the \emph{log-singular set}), such that $X^{\ls}\setminus Z$ is
a fine log scheme and $X^{\ls}\setminus Z\to\Spec k^{\ls}$ is log
smooth. Furthermore $\overline{\M}_{X\setminus Z}$ is a sheaf of
toric monoids, and if $\bar\rho\in\Gamma(X,\overline{\M}_X)$ is the
image of $1\in\NN$ under the structure map of $X^{\ls}\to\Spec
k^{\ls}$, then the function $z^{\bar\rho}$ vanishes precisely once
along each toric Weil divisor of $\Spec k[\overline{\M}_{X,\bar x}]$ for
every geometric point $\bar x\to X\setminus Z$. 

We say $X_1^{\ls}\to\Spec k^{\ls}$ and $X_2^{\ls}\to \Spec k^{\ls}$
are isomorphic toric log Calabi-Yau spaces if there is an isomorphism
$\varphi:X_1\to X_2$ of algebraic spaces over $k$ lifting to an
isomorphism of log spaces $X_1^{\ls}\setminus Z \to
X_2^{\ls}\setminus \varphi(Z)$ over $\Spec k^{\ls}$, for some
log-singular set $Z$ for $X_1^\ls$.
\qed
\end{definition}

\begin{remark}
The spaces $X_0(B,\P,s)$ look like underlying spaces of a toric log
Calabi-Yau space. The main difficulty is to place a log structure on these
spaces. This will be done in \S 5.1 under the additional hypothesis of
$(B,\P)$ being simple.
\end{remark}

\begin{remark}\label{log structure over Z}
We never make use of the log structure over $Z$, and in fact all
relevant information is contained completely in the complement of any
$Z$ with the stated properties. This is due to the fact that there is
a canonical extension of the log structure on $X\setminus Z$ to $X$
that retrieves the original log structure at all log-smooth points:
Let $j:X\setminus Z\to X$ be the inclusion. Since $X$ is Gorenstein
and $Z$ is of codimension $\ge 2$, $j_*\O_{X\setminus Z}=\O_X$. Thus
the log structure on $X\setminus Z$ extends uniquely to a log
structure on $X$ with monoid sheaf $j_*\M_{X\setminus Z}$. If $Z$
does not contain any toric stratum the canonical map $\M_X\to
j_*\M_{X\setminus Z}$ induces an isomorphism of ghost sheaves
wherever $\M_X$ is fine, hence an isomorphism of log structures.

A canonical, minimal choice for $Z$ is the closed subspace where
$(X,j_*\M_{X\setminus Z}) \to\Spec k^\ls$ is not a log smooth
morphism of fine log spaces.
\qed
\end{remark}

\begin{proposition}
If $\X\to\Spec R$ is a toric degeneration, let $\X^{\ls}$ and $\Spec
R^{\ls}$ be the log structures induced by the inclusions
$X=\X_0\subseteq\X$ and $\Spec k\subseteq\Spec R$ as in Example 
\ref{prime example}. Let $X^{\ls}$ be the induced log structure on
$X$. Then the induced morphism $X^{\ls}\to \Spec k^{\ls}$
gives $X^{\ls}$ the structure of a toric log Calabi-Yau space.
\end{proposition}

\proof Clearly $X$ satisfies conditions (1) and (2) of Definition
\ref{logCY}, as these are (2) and (3) of Definition~\ref{toric degen}.
Note that we obtain a morphism of log spaces $\X^{\ls}
\to\Spec R^{\ls}$ induced by pull-back of functions, as
$\M_{\X}$ and $\M_{\Spec R}$ are contained in the structure sheaves
of $\X$ and $\Spec R$ respectively. This morphism
then restricts to a morphism of log spaces $X^{\ls}
\to \Spec k^{\ls}$.

Now if $P$ is a toric monoid, the log structure on $\Spec k[P]$
induced by the chart $P\to k[P]$ coincides with the log structure
induced by the inclusion $\partial\Spec k[P] \hookrightarrow \Spec
k[P]$, where by $\partial\Spec k[P]$ we mean the union of all toric
divisors of $\Spec k[P]$. On the other hand, in (4) of
Definition~\ref{toric degen}, because $f_{\bar x}$ vanishes to order
$1$ on each toric divisor of $Y_{\bar x}$, $f_{\bar
x}^{-1}(0)=\partial Y_{\bar x}$. Thus if $Y_{\bar x}=\Spec k[P_{\bar
x}]$, the induced map $P_{\bar x}\to\O_{U_{\bar x}}$ gives a chart
for the log structure on $U_{\bar x}$ induced from $\X^{\ls}$. Thus
we see that condition (4) of Definition~\ref{toric degen} implies
$\X^{\ls}\setminus Z\to\Spec R^{\ls}$ is log smooth, and hence by
restriction $X^{\ls}\setminus Z$ is also log smooth.
\qed

\begin{remark} 
In the quartic example of the introduction, the minimal singular set
$Z$ coincides with the
singular locus of the total space of $\X$, but this is not true
more generally, as $\X$ may have toric singularities. However, the
singular set $Z$ is the set of points where these singularities
are somehow ``worse.''

We will never be interested in the log
structure at the points of $Z$ itself, hence the definition of
isomorphism of log Calabi-Yau spaces is independent of the log
structure along $Z$.
\qed
\end{remark}

Toric strata also have a characterization in terms of the ghost sheaf.

\begin{lemma}\label{strata characterization}
Let $X$ be a toric log Calabi-Yau space.
A geometric point $\bar x\to X\setminus Z$ is the geometric generic
point of a toric stratum of $X$ if and only if
$\rank(\overline\M^\gp_{X,\bar x}) =\dim\O_{X,\bar x}+1$.
\end{lemma}

\proof
Log smoothness gives a smooth map $\pi:U\to V:=\Spec k[\overline
\M_{X,\bar x}]/(z^{\bar \rho_{\bar x}})$ from an \'etale neighbourhood
$U$ of $\bar x$.  Since $V$ is a subscheme inside a toric variety it
also comes with a notion of toric strata, which are the intersections
of the irreducible components of $V$.  Define toric strata of $U$ as
the components of preimages of toric strata under $U\to V$.  Now by
condition (2) of Definition~\ref{logCY}, \'etale locally toric strata
of $X$ are also precisely the intersections of irreducible components
of $X$.  Therefore, toric strata in $U$ map to toric strata in $X$ and
it suffices to prove the statement for $U$ instead of $X$.

Since $\overline\M_{X,\bar x}$ has no non-trivial invertibles there
is a distinguished closed point $0\in V$, the zero-dimensional toric
stratum. Let $\eta\in\pi^{-1}(0)$ be the generic point 
of the connected component of $\pi^{-1}(0)$ containing $\bar x$. Then
$\overline\M_{X,\bar x}=\overline \M_{X,\bar\eta}$, and the connected
component of $\pi^{-1}(0)$ containing $\bar x$ is the minimal toric
stratum of $U$ containing $\bar x$. Thus $\bar x$ is the generic
point of a toric stratum
if and only if $\bar x=\bar\eta$. By the dimension formula for
smooth morphisms this is the case if and only if $\dim\O_{X,\bar
x}=\dim(V)$. The proof is finished by noting
$\dim(V)=\rank(\overline\M^\gp_{X,\bar x})-1$.
\qed

\subsection{The dual intersection complex}

Given a toric log Calabi-Yau space $X^{\dagger}$, we now reverse the
gluing procedures described earlier in this paper in order to construct
an affine manifold with singularities. We begin by defining a category
analogous to $\Cat(\P)$, which we call $\Cat(X)$. Write 
$\tilde X=\coprod X_i$, with $\nu:\tilde X\to X$ the normalization.
The objects of $\Cat(X)$ are the \emph{strata} of $X$, i.e. the
set 
\[
\Strata(X):=\{\nu(S)|\hbox{$S$ is a toric stratum of $X_i$ for some $i$}\}.
\]
If $S_1,S_2$ are two strata of $X$, we put $\Hom(S_1,S_2)=\emptyset$ if
$S_1\not\supseteq S_2$, $\Hom(S_1,S_2)=\{id\}$ if $S_1=S_2$.
If $S_1\supset S_2$, let $\nu_1:\tilde S_1\to S_1$ be the normalization
of $S_1$, so $\tilde S_1$ is a toric variety. Then 
\[
\Hom(S_1,S_2)=\{\tilde S_2| \hbox{$\tilde S_2$ is a toric stratum of
$\tilde S_1$ with $\nu_1(\tilde S_2)=S_2$}\}.
\]
If we have a chain $S_1\supset S_2\supset S_3$, an element
of $\Hom(S_1,S_2)$ is a stratum $\tilde S_2$ of $\tilde S_1$ mapping
to $S_2$; identifying $\tilde S_2$ with the normalization of
$S_2$, an element of $\Hom(S_2,S_3)$ can then be identified also
as a substratum of $\tilde S_1$, and this defines composition of morphisms.

Next let $\LPoly$ be the category of lattice polytopes, a subcategory
of the category of topological spaces, with the objects being
polytopes with integral vertices and morphisms
$\Hom(\sigma_1,\sigma_2)$ being integral affine identifications of
$\sigma_1$ as a face of $\sigma_2$.

We can then define a functor
\[
LP:\Cat(X)\to \LPoly
\]
as follows. If $S\in \Strata(X)$, with generic point $\eta$, then the
stalk $\overline\M_{X,\bar\eta}$ is a monoid with a distinguished
element $\bar\rho$ specified by the morphism $X^{\ls}\to \Spec
k^{\ls}$: $\bar\rho$ is the germ of the image of $1\in\NN$ in 
$\overline\M_{X,\bar\eta}$. We define $LP(S)$ to be the convex
hull of
\[
\{\varphi\in\Hom(\overline\M_{X,\bar\eta},\NN)|
\varphi(\bar\rho)=1\}
\]
inside the affine space
\[
\{\varphi\in \Hom(\overline\M_{X,\bar\eta},\RR)|
\varphi(\bar\rho)=1\}.
\]
\begin{lemma}
If $S\in \Strata(X)$, then $LP(S)$ is a non-empty lattice polytope.
\end{lemma}

\proof Let $P=\overline\M_{X,\bar\eta}$ for $\eta$ the generic
point of $S$, with $\bar\rho\in P$ given by the morphism $X^{\ls}\to
\Spec k^{\ls}$. Then as $P$ is a toric monoid, the toric variety
$\Spec k[P]$ is given by a cone $K\subset \Hom(P,\RR)$, the convex
hull of $\Hom(P,\NN)$. The condition that $\bar\rho$ should be $1$ on
each irreducible component of $X$ implies that $z^{\bar\rho}\in k[P]$
vanishes to order $1$ on each toric divisor of $\Spec k[P]$. But this
is saying $\bar\rho$ evaluates to $1$ on the primitive integral
generators of the extremal rays of $K$. Thus $LP(S)$ is just the
convex hull of these integral generators, and hence is a (non-empty)
lattice polytope.
\qed
\medskip

Now suppose we are given strata $S_1\supset S_2$, and an element of
$\Hom(S_1,S_2)$, i.e. strata $\tilde S_1\supset \tilde S_2$, where
$\nu_1:\tilde S_1\to S_1$ is the normalization. We need to
define a morphism $LP(S_1)\to LP(S_2)$. This is done as follows.
If $\eta_1,\eta_2$ are the generic points of $\tilde S_1$, $\tilde S_2$
respectively, there is a well-defined cospecialization map
\[
\overline\M_{X,\bar\eta_2}=(\nu_1^*\overline{\M}_X)_{\bar\eta_2}
\to (\nu_1^*\overline{\M}_X)_{\bar\eta_1}
=\overline\M_{X,\bar\eta_1}.
\]
Note the identifications on the left and right are canonical, but the
map depends on which stratum $\tilde S_2\subset \tilde S_1$ dominating
$S_2$ we have chosen. It is easy to see that such a map is a
surjection and dually we obtain a map
\[
\Hom(\overline\M_{X,\bar\eta_1},\NN)
\to \Hom(\overline\M_{X,\bar\eta_2},\NN)
\]
which identifies the first monoid as a face of the second. In particular,
we then obtain a map $LP(S_1)\to LP(S_2)$ identifying the first
polytope as a face of the second. It is clear this construction is
compatible with composition of morphisms in $\Strata(X)$.

We can now construct the dual intersection complex $B$ of $X$ as a
union of polyhedra, explicitly $\displaystyle\lim_{\longrightarrow}
LP$ in the category of topological spaces.

\begin{proposition}
\label{intgraphman}
If $\dim X=n$, then $B$ is an $n$-dimensional manifold.
\end{proposition}

We will first need the following Lemma.

\begin{lemma}
\label{samecone}
Let $X_1$ and $X_2$ be two affine toric varieties over 
$k$, defined by strictly convex cones $\sigma_1,\sigma_2\subseteq
M_{\RR}$, with $\dim\sigma_1= \dim\sigma_2=\dim M_{\RR}$. Let $x_i\in
X_i$ be the zero-dimensional torus orbit. Suppose $X_1$ and $X_2$ are
\'etale locally isomorphic at $x_i$, i.e. there is an
isomorphism of strictly Henselian local $k$-algebras
$\psi:\O_{X_1,\bar x_1}\to \O_{X_2,\bar x_2}$. Suppose further that
this isomorphism preserves toric strata, i.e. it takes the ideal of a
toric stratum in $\O_{X_1,\bar x_1}$ to the ideal of a toric stratum
in $\O_{X_2,\bar x_2}$. Then there exists an element $\varphi:M\to M$
of $GL_n(\ZZ)$ such that $\varphi(\sigma_1)=\sigma_2$, and which is
combinatorially compatible with $\psi$, i.e. if $\psi$ identifies the
ideals of two toric strata, then $\varphi$ identifies the
corresponding cones.
\end{lemma}

\proof Let $P_i=\dual{\sigma_i}\cap N$ be the monoid defining $X_i$.
Then if $X_i^{\ls}$ is the log structure induced by the inclusion
$\partial X_i\subseteq X_i$, then $\overline{\M}_{X_i,\bar x_i} \cong
P_i$. It is then clear that the isomorphism $\psi$ induces an
isomorphism $\overline{\M}_{X_2,\bar x_2}\mapright{\cong}
\overline{\M}_{X_1,\bar x_1}$, 
and then an isomorphism $M\cong
\Hom(P_1,\ZZ)\to M\cong \Hom(P_2,\ZZ)$ taking
$\sigma_1\subseteq M_{\RR}$ to $\sigma_2\subseteq M_{\RR}$ as
desired. \qed
\medskip

\noindent
\textit{Proof of Proposition~\ref{intgraphman}.}
The vertices of $B$ correspond to irreducible components of $X$. It
is clearly enough to show $B$ is an $n$-dimensional manifold in a
neighbourhood of each vertex, since $B$ behaves uniformly along each
face.

Focusing on one vertex $v$ of $B$, there is a corresponding toric
variety $X_v$. Now for any zero-dimensional stratum $x\in X_v$, there
is a monoid $\overline\M_{X,\nu(\bar x)}=:P_{\bar x}$.
Lemma~\ref{strata characterization} shows $\dim k[P_{\bar
x}]/(z^{\bar\rho_{\bar x}})=n$. Let $K_{\bar x}$ be the convex hull
of $\Hom(P_{\bar x},\NN)$ in $\Hom(P_{\bar x},\RR)$; this is the cone
defining the toric variety $\Spec k[P_{\bar x}]$. It is then clear
that $\Hom(P_{\bar x},\RR)\cong\RR^{n+1}$ and $K_{\bar x}$ is an
$n+1$-dimensional cone, by the previous observations.
Thus $LP(\nu(\bar x))$ is in particular an $n$-dimensional lattice polytope.

If $\eta$ is the generic point of $\nu(X_v)$, then the
cospecialization map induced by $\bar x\to X_v$,
$\overline\M_{X,\nu(\bar x)} \to \overline\M_{X,\bar\eta}$, is dual
to the inclusion of some extremal ray $R_{\bar x}\subseteq K_{\bar
x}$. The tangent wedge of $LP(\nu(\bar x))$ at the vertex $p$ is just
the cone $K_{\bar x}(R_{\bar x})=(K_{\bar x}+\RR R_{\bar x})/ \RR
R_{\bar x}$, and by Lemma~\ref{samecone}, this cone is the same as
the cone corresponding to the stratum $x\in X_v$ in the fan
$\Sigma_{X_v}$ defining $X_v$, up to integral linear 
transformations. It is then easy to see that all the cones $K_{\bar
x}(R_{\bar x})$ fit together to give, at least topologically,
the fan $\Sigma_{X_v}$, as $\{x\}$ runs over all zero-dimensional
strata of $X_v$. But this is the local picture of $B$ at $v$, so $B$
is a manifold in a neighbourhood of $v$.
\qed\medskip

We can now describe an integral affine structure with singularities
on $B$ with a polyhedral decomposition $\P$. Here
\[
\P=\{\hbox{image of $LP(S)$ in $B=\displaystyle
\lim_{\longrightarrow} LP$}| S\in \Strata(X)\}.
\]
By Construction~\ref{basicconstruct}, it is enough to give a fan
structure at each vertex of $\P$ compatible with the
polyhedral decomposition. But we have already done this in
the proof of Proposition~\ref{intgraphman}: There is a
combinatorial identification of $\P$ in a neighbourhood of a vertex
$v$ with the fan $\Sigma_{X_v}$. This gives the fan structure.

Thus we have obtained an integral affine manifold $B$ with
singularities along the discriminant locus $\Delta'$ given in
Construction~\ref{basicconstruct}. We now observe:

\begin{proposition}
$\P$ is a toric polyhedral decomposition of $B$.
\end{proposition}

\proof The only thing to check is that $\P$ is toric, since it is
a polyhedral decomposition by construction. We will show this by
applying Proposition~\ref{monodromy2}. 

To do this, let $\tau\in\P$ correspond to some stratum $S$, which is
the image of toric strata $S_1\subseteq X_{v_1},\ldots,S_m\subseteq
X_{v_m}$, where $v_1,\ldots,v_m$ are the (possibly non-distinct)
vertices of $\tau$. Then $\nu$ yields canonical isomorphisms between
the $S_i$ which preserve the toric strata of the $S_i$; thus these
isomorphisms are equivariant under the torus action on $\coprod
X_{v_i}$. So the fans determining $S_1,\ldots,S_m$ can be identified.
Specifically, if $\Sigma_{v_i}$ defines $X_{v_i}$, living naturally
in $\Lambda_{\RR,v_i}$, and $\tau_{v_i}\in\Sigma_{v_i}$ is the cone
corresponding to $S_i$, then the fan for $S_i$ is
$\Sigma_{v_i}(\tau_{v_i})$, living naturally in
$\Lambda_{\RR,v_i}/\Lambda_{\tau,\RR}$.

Then parallel translation from $v_i$ to $v_j$ along some path
identifies
\[
\hbox{$\Lambda_{\RR,v_i}/\Lambda_{\tau,\RR}$ and
$\Lambda_{\RR,v_j}/\Lambda_{\tau,\RR}$}
\]
in such a way that the fans $\Sigma_{v_i}(\tau_{v_i})$ and
$\Sigma_{v_j}(\tau_{v_j})$ are identified. There is a unique such
identification, and hence this is independent of the path chosen.
Thus it follows, for any loop $\gamma$ based at $y\in
\Int(\tau)\setminus\Delta'$ in a sufficiently small neighbourhood
$U_{\tau}$ of $\Int(\tau)$, that $\tilde\rho(\gamma)$ acts trivially
on $\Lambda_{\RR,y}/\Lambda_{\tau,\RR}$, i.e.
\[
(\tilde\rho(\gamma)-I)(\Lambda_{\RR,y})\subseteq \Lambda_{\tau,\RR}
\]
as desired. \qed

\begin{definition}
If $X^{\ls}$ is a toric log Calabi-Yau space, then the integral
affine manifold with singularities $B$ along with the polyhedral
decomposition $\P$ constructed above is the \emph{dual intersection
complex} of $X^{\ls}$. 
\end{definition}

\begin{theorem} 
\label{logopen}
Let $X^{\ls}$ be a toric log Calabi-Yau space, and let $(B,\P)$
be the dual intersection complex. Then there exist open gluing
data $s$ for $\P$ over $k$ such that
\[
X\cong X_0(B,\P,s).
\]
\end{theorem}

\proof Let $X_{\tau}$ be as usual for $\tau\in\P$. Then by
construction, there is an isomorphism $\psi_{\tau}:X_{\tau}\to
S_{\tau}$, where $S_{\tau}$ is the normalization of the stratum of
$X$ corresponding to $\tau$. Furthermore, this isomorphism can be
chosen to identify the toric strata of $X_{\tau}$ with the strata of
$S_{\tau}$. Then $\psi_{\tau}$ is well-defined up to the
action of $\shQ_{\tau} \otimes\Gm(k)$ on $X_{\tau}$. Fix such
choices. In particular, for every $e:\tau\to\sigma$, we then obtain 
$\bar s_e\in \shQ_{\sigma}\otimes \Gm(k)$ such that the diagram
\[
\begin{matrix}X_{\sigma}&\mapright{\psi_{\sigma}}&S_{\sigma}\\
\mapdown{F(e)\circ \bar s_e}&&\downarrow\\
X_{\tau}&\mapright{\psi_{\tau}}&S_{\tau}
\end{matrix}
\]
commutes. Then $\bar s=(\bar s_e)$ forms \emph{closed gluing data,}
and $\displaystyle\lim_{\longrightarrow} F_{\Spec k,\bar s}=X$.
Indeed, there is a canonical map $\displaystyle\lim_{\longrightarrow}
F_{\Spec k,\bar s} \to X$ induced by the system of composed morphisms
$X_{\tau} \to S_{\tau}\to X$. Since $\displaystyle
\lim_{\longrightarrow} F_{\Spec k,\bar s}$ glues together the
components $X_v$ for $v$ a vertex in $\P$ in a normal crossings way
in codimension one, we must have $\displaystyle\lim_{\longrightarrow}
F_{\Spec k,\bar s}\to X$ being an isomorphism in codimension one. On
the other hand, as $X$ is assumed to be Gorenstein, and in particular
$S_2$, in fact this is an isomorphism everywhere by \cite{Reid},
Proposition 2.2. We then have to show this gluing arises via open
gluing data.

First, we show $X$ has an \'etale open cover by sets of the form
$V(\sigma)$ for $\sigma\in\P_{\max}$. Fix $\sigma\in\P_{\max}$,
and let $\bar x\to X$ be the corresponding zero-dimensional stratum. Then
we have a chart
\[
P_{\sigma}:=\overline{\M}_{X,\bar x}\to\O_{X,\bar x}
\]
giving the log structure at $x$. Now for each $e\in\coprod_{\tau\in\P
\atop \tau\not=\sigma} \Hom(\tau,\sigma)$, we have an affine open set
$V_e$ of $X_{\tau}$ corresponding to the cone $K_{\tau}$ in
$\Sigma_{\tau}$, and $V(\sigma)=\displaystyle\lim_{\longrightarrow}
V_e$ with respect to the canonical directed system
$F_{\sigma}(f):V_{e_2} \to V_{e_1}$ over all commutative
diagrams
\begin{eqnarray}
\label{diagram51}
\begin{matrix}
\tau_1&\mapright{f}&\tau_2\\
&\smash{\mathop{\searrow}\limits^{e_1}}&\mapdown{e_2}\\
&&\sigma
\end{matrix}
\end{eqnarray}
On the other hand, let $V'(\sigma)$ be the limit with respect to the
system
\[
F_{\sigma}(\bar s)(f):V_{e_2}\to V_{e_1}
\]
given by $F_{\sigma}(\bar s)(f)=F_{\sigma}(f)\circ \bar s_f$. Using
the universal property of colimits it is straightforward to show that
the canonical map $V'(\sigma)\to X$ is \'etale. It thus remains to
produce an isomorphism of $V(\sigma)$ and $V'(\sigma)$. 

To do so, note that $V_e\cong \Spec k[P_e]$, where $P_e$ is the face
of $P_{\sigma}$ corresponding to $e$. Thus by composing the chart
$P_{\sigma}\to\O_{X,\bar x}$ with the map $\O_{X,\bar x}
\to\O_{V_e,\bar x}$ induced by the composition $V_e\to V'(\sigma)\to
X$, we obtain a map of monoids $\phi_e:P_e\to \O_{V_e,\bar x}$ of the
form $\phi_e(p)=h_p z^p$ for $h_p\in\O_{V_e,\bar x}^{\times}$; in
particular, $p\mapsto h_p$ is a map of monoids $P_e\to \O_{V_e,\bar
x}^{\times}$. Composing this map with the residue map
\[
\O_{V_e,\bar x}^{\times}\to \O_{V_e,\bar x}^{\times}/(1+{\bf m}_{\bar x})
\cong k^{\times},
\]
we get maps
\[
\bar\phi_e:P_e\to k[P_e]
\]
given by
\[
p\longmapsto (h_p\mod 1+{\bf m}_{\bar x})\cdot z^p.
\]
This is still a map of monoids; what we are doing is removing the
``non-constant'' part of the log structure. Furthermore, given $f$ as in
(\ref{diagram51}), we necessarily have a commutative diagram
\[
\begin{matrix}
k[P_{e_1}]&\mapright{\phi_{e_1}}& \O_{V_{e_1},\bar x}\\
\mapdown{F_{\sigma}(f)^*}&&\mapdown{F_{\sigma}(\bar s)(f)^*}\\
k[P_{e_2}]&\mapright{\phi_{e_2}}& \O_{V_{e_2},\bar x}
\end{matrix}
\]
hence a commutative diagram
\[
\begin{matrix}
k[P_{e_1}]&\mapright{\bar\phi_{e_1}}& k[P_{e_1}]\\
\mapdown{F_{\sigma}(f)^*}&&\mapdown{F_{\sigma}(\bar s)(f)^*}\\
k[P_{e_2}]&\mapright{\bar\phi_{e_2}}& k[P_{e_2}]
\end{matrix}
\]
and hence an isomorphism of directed systems inducing an isomorphism
$\phi:V'(\sigma)\to V(\sigma)$, as desired.

Now let $U=\coprod_{\sigma\in\P_{\max}} V(\sigma)$. Then we obtain an
\'etale open cover $\psi:U\to X$ which is the composition of
$\phi^{-1}$ and the natural map $V'(\sigma)\to X$ on each
$V(\sigma)$. We have an \'etale equivalence relation $\foR=U\times_X U
\subseteq U\times U$. 

To obtain the open gluing data from this, choose for each
$\omega\in\P$ a morphism $e_{\omega}:\omega\to\sigma$, for some
$\sigma\in\P_{\max}$. Set $s_{e_{\omega}}=1$. If $e:\omega\to\sigma'$
is any other morphism with $\sigma'\in\P_{\max}$, we obtain
canonically $V(\omega)\subseteq V(\sigma)$ and $V(\omega)\subseteq
V(\sigma')$, with a canonical identification given by
$\Phi_{e_{\omega}e}$ as in Section~\ref{section2}. Then $\foR\cap
((V(\omega)\subseteq V(\sigma')) \times (V(\omega)\subseteq
V(\sigma))\subseteq V(\sigma')\times V(\sigma)$ can be viewed as the
graph of a morphism, and by restricting to irreducible components it
is easy to see this morphism is necessarily of the form
$\Phi_{e_{\omega}e}\circ s_e$ for some $s_e\in \check{PM}(\omega)$. This
defines $s_e$ whenever $e$ is a morphism to a maximal facet. Then
given a diagram (\ref{diagram51}), we define 
\[
s_f=s^{-1}_{e_2}|_{\tau_1}\cdot s_{e_1}.
\]
We need to check this is well-defined. Given
\[
\xymatrix@C=30pt
{&&\sigma_1\\
\tau_1\ar@/^/[rru]^{e_1}\ar[r]^(.7){f}\ar@/_/[drr]_{g_1}&\tau_2
\ar[ru]_{e_2}\ar[rd]^{g_2}&\\
&&\sigma_2
}
\]
it follows by transitivity of $\foR$ that 
\[
\foR\cap ((V(\tau_1)\subseteq V(\sigma_1))\times (V(\tau_1)\subseteq
V(\sigma_2)))
\]
is the graph of $s_{g_1}^{-1}\circ\Phi_{g_1e_1}\circ s_{e_1}$
while
\[
\foR\cap ((V(\tau_2)\subseteq V(\sigma_1))\times (V(\tau_2)\subseteq
V(\sigma_2)))
\]
is the graph of $s_{g_2}^{-1}\circ\Phi_{g_2e_2}\circ s_{e_2}$. These graphs
agree on $(V(\tau_1)\subseteq V(\sigma_1))\times (V(\tau_1)\subseteq 
V(\sigma_2))$, and hence
\[
(s_{e_2}^{-1}s_{g_2})|_{\tau_1}=s_{e_1}^{-1}s_{g_1},
\]
proving well-definedness. Then $(s_e)$ forms open gluing data, and
by construction it is clear that $X\cong X_0(B,\P,s)$.
\qed

\begin{corollary}
Let $X^{\ls}$ be a toric log Calabi-Yau space with dual
intersection complex $(B,\P)$, $X=X_0(B,\P,s)$ for some
open gluing data $s$, and let $Z\subseteq X$ be the singular set. 
If $X_0(B,\P,s)^g$ denotes the ghost structure on $X_0(B,\P,s)$ of
Example~\ref{fundamental example}, then $X^{\ls}$ is of ghost
type $X_0(B,\P,s)^g$ on $X\setminus Z$, and thus the log smooth
structure on $X\setminus Z$ is induced by a section 
$f\in\Gamma(X\setminus Z,\shLS_{X_0(B,\P,s)})$. Conversely,
given an integral affine manifold $B$ with toric polyhedral
decomposition $\P$, $s$ open gluing data, and
if $Z\subseteq X_0(B,\P,s)$ is a closed subset of codimension
$\ge 2$ not containing any toric stratum, then any section
$f\in\Gamma(X_0(B,\P,s)\setminus Z,\shLS_{X_0(B,\P,s)})$ induces
a toric log Calabi-Yau space structure on $X_0(B,\P,s)$.
\end{corollary}

\proof By log smoothness the type of log structure of $X^\ls$ along a
toric stratum minus lower-dimensional strata is constant and may
hence be read off at the generic point of the stratum. Thus the
statement follows from the construction of the dual intersection
complex, Example~\ref{fundamental example}, and
Proposition~\ref{moduli of log smooth structures}. We note for the
second part that $f\in\Gamma(X_0(B,\P,s)\setminus
Z,\shLS_{X_0(B,\P,s)})$ defines a log smooth structure on
$X_0(B,\P,s)\setminus Z$, which is all we need by Remark~\ref{log
structure over Z}.
\qed

\subsection{Polarized log Calabi-Yau spaces and the intersection complex}

Let $X^{\ls}$ be a toric log Calabi-Yau space with a line bundle
$\shL$. Then the dual intersection complex $(B,\P)$ carries a 
multi-valued piecewise linear function constructed as follows. For every 
$\sigma\in\P$, we have the map $q_{\sigma}:X_{\sigma}\to X$. Pulling
back $\shL$ to $X_{\sigma}$ gives a line bundle on the toric variety
$X_{\sigma}$, and hence an integral piecewise linear function
on the fan $\Sigma_{\sigma}$, well-defined up to linear
functions. Recall that if $m_1,\ldots,m_p$ are primitive integral
generators of the rays of $\Sigma_{\sigma}$, with corresponding toric
divisors $D_1,\ldots,D_p$, then the toric divisor $\sum a_i D_i$
corresponds to a piecewise linear function $\varphi_{\sigma}$ on
$\Sigma_{\sigma}$ with $\varphi_{\sigma}(m_i)=a_i$. By pulling this
function back to $U_{\sigma} \subseteq B$ via $s_{\sigma}$ of
Definition~\ref{PD2}, (5), one obtains piecewise linear functions on
the open covering $\{U_{\sigma}|\sigma\in\P\}$ of $B$. It is easy to
see these functions differ only by linear functions on intersections,
and thus we obtain a multi-valued piecewise linear function
$\varphi_{\shL}$ with integral slopes. If $\shL$ is ample, then it
follows from the standard characterization of ample line bundles on
toric varieties (\cite{Oda}, Corollary~2.15) that $\varphi_{\shL}$
is strictly convex.

\begin{definition}
If $\shL$ is ample, we call $(X^{\ls},\shL)$ a polarized toric log Calabi-Yau
space, and we call the triple
\[
(B,\P,\varphi_{\shL})
\]
\emph{degeneration data} associated to $(X^{\ls},\shL)$. If
$(\check B,\check\P, \check\varphi_{\shL})$ is the discrete Legendre
transform of $(B,\P,\varphi_{\shL})$, then we call $\check B,\check\P$
the \emph{intersection complex} of the polarized toric log Calabi-Yau space
$(X^{\ls},\shL)$.
\qed
\end{definition}

See \cite{Announce},~Section~4, for additional discussion of the
intersection complex.

\subsection{Positive log structures}

As we know from Corollary~\ref{nomonodromy}, it is rare that
there is no singular set $Z\subseteq X^{\ls}$. However, we need
to provide some control over these singularities. The following is a useful
restriction.

\begin{definition} 
Let $X^{\ls}$ be a toric log Calabi-Yau space with singular set $Z$,
with  $X\cong X_0(B,\P,s)$. Then the log smooth structure
$X^{\ls}\setminus Z\to \Spec k^{\ls}$ is induced by a section $f\in
\Gamma(X_0(B,\P,s) \setminus Z, \shLS_{X_0(B,\P,s)})$. We say that
$X^{\ls}$ is \emph{positive} if the corresponding sections
$f_{\omega}\in \Gamma(X_{\omega}\setminus q_{\omega}^{-1}(Z),
\shN^{\times}_{\omega})$  (see Theorem~\ref{q_tau2}) for $\dim\omega=1$
extend to sections of $\shN_{\omega}$ on $X_{\omega}$, i.e. if they have
zeros but no poles.
\end{definition}

\begin{remark}
It is easy to see that this notion of positivity is independent
of the choice of orientation of edges of $\P$.
\qed
\end{remark}

\begin{proposition}
Let $(B,\P)$ be positive, and let $\shN_{\omega}$ be the line bundle
on $X_{\omega}$ for $\dim\omega=1$ given by Theorem~\ref{q_tau2}.
Then $\shN_{\omega}$ is generated by global sections. In particular,
it has a section whose zero-locus does not contain a toric stratum of
$X_{\omega}$. Conversely, if $X^{\ls}$ is positive, then its dual
intersection complex $(B,\P)$ is positive.
\end{proposition}

\proof By Theorem~\ref{q_tau2}, $\shN_{X_{\omega}}$ corresponds to
the piecewise linear function of Remark~\ref{psiomega}. If $(B,\P)$
is positive, then $\psi_{\omega}$ is convex and hence by \cite{Oda},
Corollary~2.15, $\shN_{X_{\omega}}$ is generated by global sections.
Conversely, if $X^{\ls}$ is positive, there exist sections
$f_{\omega}\in\Gamma(X_{\omega},\shN_\omega)$ vanishing only on
$q_{\omega}^{-1}(Z)$, hence they do not vanish on any toric stratum by the
assumptions on $Z$. Since the base locus of a complete toric linear
system is a union of toric strata, $\shN_\omega$ is globally
generated, hence $\psi_{\omega}$ is convex, and $(B,\P)$ is positive.
\qed

\begin{proposition}
Let $f:\X\to \shS$ be a toric degeneration. Then $\X_0^{\ls}$ is
positive.
\end{proposition}

\proof Let $(B,\P)$ be the dual intersection complex, and
$\omega\in\P$ a one-dimensional face, corresponding to a stratum
$S_{\omega}$ of dimension $n-1$. It suffices to check the absence of
poles of $f_\omega\in \Gamma({q_\omega}_*\shN_\omega)$ in codimension
one on $S_\omega$. Without loss of generality, we can view
$\Lambda_{\omega}\cong\ZZ =M'\subseteq M$ with $\omega\subseteq\RR$
an interval with endpoints $0$ and $l>0$. Choose $d_{\omega}=1$. 
Then $Q':=\dual{C'(\omega)}\cap (N'\oplus\ZZ)$ is generated by
$p_1=(1,0)$, $p_2=(-1,l)$ and $\rho=(0,1)$, with relation
$p_1+p_2=l\rho$. Abstractly, $V(\omega)\cong (\Spec
k[u,v]/(uv))\times \Gm^{n-1}$, and the log structure on
$V(\omega)\setminus Z$ is of the same type as that given by the chart
\[
p_1\mapsto u, p_2\mapsto v, \rho\mapsto 0.
\]
If on some open subset $U$ of $V(\omega)$, we have a chart given by
\[
p_1\mapsto h_u u, p_2\mapsto h_v v, \rho\mapsto 0
\]
with $h_u, h_v$ invertible functions as usual, then by the
construction of the map $\xi$ in the proof of Theorem \ref{finalLS},
the corresponding section of $\shLS_{V(\omega)}
=\O_{V(\omega)}^{\times}$ is $(h_uh_v)^{-1}$. For the log structure
induced by the embedding $\X_0\subset \X$ it can be computed as
follows. Let $\tilde u,\tilde v\in \O_{\X,\bar x}$ be extensions of
$u,v$ at some geometric point $\bar x$ in the codimension one stratum
of $V(\omega)\setminus Z$. Then $\tilde u\tilde v$ restricts to zero
on the central fibre and hence $\tilde u\tilde v = t^l\cdot h$ for
$t$ a generator of the maximal ideal of $\Gamma(\O_{\shS})$ and $h\in
\O_{\X, \bar x}^\times$. The exponent of $t$ is $l$ because
$\M_{\X,\bar x}\cong Q'$. We may then put $h_u=1$,
$h_v=h^{-1}|_{\X_0}$. Thus $f_\omega$ is represented by $h|_{\X_0}$.

Now let $\bar y$ be a geometric point of $V(\omega)\cap Z$ and 
$\tilde u,\tilde v\in \O_{\X,\bar y}$ as before. By what we just said
$t^l$ divides $\tilde u\tilde v$ in codimension one and hence
everywhere since $\X$ is normal. The quotient $h:= \tilde u\tilde v /t^l
\in \O_{\X,\bar y}$ is the desired regular extension of
$f_\omega$ over $\bar y$.
\qed

\begin{definition}\label{LS^+_pre}
For $\sigma\subseteq M_{\RR}$ a lattice polytope,
we denote by $\shLS^+_{\pre,V(\sigma)}$ the sheaf 
\[
\bigoplus_{\omega\subseteq\sigma\atop \dim\omega=1}\O_{V_{\omega}},
\]
where $V_{\omega}:=V_{\omega \to\sigma}$ denotes the codimension one
stratum of $V(\sigma)$ corresponding to $\omega$. For $(B,\P)$ and
open gluing data $s$, $X=X_0(B,\P,s)$, we denote by $\shLS^+_{\pre,X}$
the sheaf on $X$ obtained by gluing the sheaves
$\shLS^+_{\pre,V(\sigma)}$ for $\sigma\in\P_{\max}$ using the gluing
formula of Theorem~\ref{globalgluing}, so that in particular
$\shLS_X$ is a subsheaf (of sets) of $\shLS^+_{\pre,X}$.
\qed
\end{definition}

Note that to specify a positive log Calabi-Yau structure on
$X=X_0(B,\P,s)$,  we need to give a section $f\in\Gamma
(X,\shLS^+_{\pre,X})$ satisfying: (1) $f$ satisfies the
multiplicative conditions given in Theorem~\ref{finalLS} and (2) the
components of $f$ on codimension one strata have no zero sets
containing any toric stratum. Now $\shLS^+_{\pre,X}$ is a coherent
sheaf on $X$, so if $B$ is compact, $\shLS^+_{\pre,X}$ has a  finite
dimensional space of sections. Then condition (2) determines some
open subset of this space of sections, while the multiplicative
condition (1) is a more subtle, closed condition.  This allows us to
identify the moduli space of positive log Calabi-Yau structures on
$X$ as a subvariety of an affine space. We shall examine this more
closely in the next section.

\begin{remark}
We can now determine the discriminant locus of a dual intersection
complex $B$ more precisely in the positive case. Initially, we take
$\Delta'\subseteq B$ as in Construction~\ref{basicconstruct} to be
the union of all simplices of $\Bar(\P)$ not containing a vertex of
$\P$ or the barycenter of a maximal cell. Using Proposition
\ref{extaff}, we can identify a smaller discriminant locus
$\Delta\subseteq \Delta'$ as follows if $X^{\ls}$ is positive and the
singular set $Z\subseteq X$ is taken to be minimal.

For any $\sigma\in \Bar(\P)$, $\sigma\subseteq\Delta'$, $\dim\sigma
=n-2$, there exist unique $\sigma_1,\sigma_{n-1}\in\P$ with
$\dim\sigma_i=i$ and $\sigma$ containing the barycenters of $\sigma_1$
and $\sigma_{n-1}$. In addition, there exists a unique morphism
$e_{\sigma}:\sigma_1\to\sigma_{n-1}$ corresponding to the
edge of $\sigma$ joining these two barycenters. Furthermore, there 
exist only two distinct morphisms $f_i:\sigma_{n-1}\to\tau_i$,
$i=1,2$, with $\tau_i\in\P_{\max}$. Let $e_i=f_i\circ e_{\sigma}$.

The monodromy around $\sigma$ is now given by $T^{e_1e_2}_{\sigma_1}$
of \S1.5, and hence is trivial if and only if
$n_{\sigma_1}^{e_1e_2}=0$. 

On the other hand, consider the embedding $F(e_{\sigma}):\PP^1 \cong
X_{\sigma_{n-1}}\to X_{\sigma_1}$. Then by Theorem \ref{q_tau2} and
\cite{Oda}, Corollary~2.15, the line bundle $F(e_{\sigma})^*
(\shN_{X_{\sigma_1}})$ is trivial if and only if
$n_{\sigma_1}^{e_1e_2}=0$, in which case the affine structure on $B$
will extend across $\sigma$ by Proposition~\ref{extaff}. Now the
positive log Calabi-Yau structure on $X$ determines a section
$f_{\sigma_1} \in\Gamma(X_{\sigma_1}, \shN_{X_{\sigma_1}})$ which
does not vanish on an entire toric stratum of $X_{\sigma_1}$. Thus
$F(e_{\sigma})^*(\shN_{X_{\sigma_1}}) \cong \O_{X_{\sigma_{n-1}}}$ if
and only if $f_{\sigma_1}$ does not vanish at any point of
$F(e_{\sigma})(X_{\sigma_{n-1}})$. However, the zero locus of
$f_{\sigma_1}$ is
\[
\cl(q_{\sigma_1}^{-1}(Z)\setminus \partial X_{\sigma_1})\subseteq
X_{\sigma_1}.
\]
Thus we see that $\Delta$ is the union of all $\sigma\in \Bar(\P)$, 
$\sigma\subseteq\Delta'$, $\codim\sigma=2$ satisfying
\[
\cl(q_{\sigma_1}^{-1}(Z)\setminus \partial X_{\sigma_1})
\cap F(e_{\sigma})(X_{\sigma_{n-1}})\not=\emptyset.
\]
This gives the precise relationship between $\Delta$ and $Z$.
\qed
\end{remark}

\subsection{Normalized gluing data and examples}

We saw in Proposition~\ref{isomorphism} that if $s$ and $s'$ differ by
an element of $B^1(\P,\shQ_{\P}\otimes\Gm(k))$, then $X_0(B,\P,s)
\cong X_0(B,\P,s')$. We would now like to use open gluing data
to parametrize log structures as well as the underlying schemes,
so that we can determine explicitly the moduli space of toric log
Calabi-Yau spaces with a given dual intersection complex $(B,\P)$.

We begin with the following definition:

\begin{definition}
Let $\sigma\subseteq M_{\RR}$ be an $n$-dimensional integral
polytope. Let $x\in V(\sigma)$ be the unique zero-dimensional torus
orbit. A section $f$ of $\shLS_{V(\sigma)}$ (or
$\shLS^+_{\pre,V(\sigma)}$) defined in a neighbourhood of $x$ is said
to be \emph{normalized} if, viewing
\[
f=(f_{\omega})\in\bigoplus_{\omega\subseteq\sigma\atop
\dim\omega=1}\O_{V_{\omega}}
\]
as in Theorem~\ref{finalLS}, 
$f_{\omega}$ takes the value
$1$ at $x$ for each $\omega\subseteq\sigma$, $\dim\omega=1$.
\qed
\end{definition}

Note that this condition can be rephrased in terms of a chart for
the log structure near the point $x$. Given $X=V(\sigma)$ a
scheme with the \'etale topology, and a chart
\[
\alpha:P_{\sigma}\to\O_{X,\bar x}
\]
of the form $p\mapsto h_pz^p$ as usual, with $h_p$ a germ of an invertible
function at $x\in X(p)=\cl{\{z^p\not=0\}}$, we obtain
a map 
\[
\phi:\partial P_{\sigma}\to k^{\times}
\]
given by $p\mapsto h_p/\maxid_{\bar x}$.

This map is piecewise multiplicative, and it is not difficult to see
that the corresponding section $f$ of $\shLS_{V(\sigma)}$ is
normalized if and only if this map is in fact multiplicative on $N$.

In particular, even if $f$ is not normalized, $\phi$ induces
an automorphism of $V(\sigma)$ via $z^p\mapsto\phi(p)z^p$, and it is
then clear that the pull-back of the chart $\alpha$ by $\phi^{-1}$
yields a section of $\shLS_{V(\sigma)}$ in a neighbourhood of 
$x\in V(\sigma)$ which is normalized.

\begin{definition}
Let $s$ be open gluing data for $(B,\P)$, with canonical projection maps
$p_{\sigma}:V(\sigma)\to X_0(B,\P,s)$ for each $\sigma\in\P_{\max}$.
Then a section $f\in\Gamma(X,\shLS^+_{\pre,X})$ is \emph{normalized}
if $p_{\sigma}^{-1}(f)$ is a normalized section of
$\shLS_{\pre,X}^+$ for all $\sigma\in\P_{\max}$. If in fact $f$
determines a log Calabi-Yau structure $X^{\ls}$ on $X$, i.e. 
$f\in\Gamma(X\setminus Z,\shLS_X)$ for some singular set $Z\subseteq
X$ not containing any toric stratum of $X$, and in addition
$f$ is normalized, then we say $X^{\ls}$
is normalized by the open gluing data $s$.
\end{definition}

By the above discussion, there is always some open gluing data with
respect to which a toric log Calabi-Yau space $X^{\ls}$ is normalized.
Conversely we have the following criterion.

\begin{proposition}
\label{gluingcond}
Let $(B,\P)$ be positive, and let $s\in Z^1(\P,\shQ_{\P}\otimes\Gm)$ be 
open gluing data, $X=X_0(B,\P,s)$. Then there exists a
normalized section $f\in\Gamma(X,\shLS^+_{\pre,X})$ if and only if
the following condition holds:

\emph{Condition~(LC):} 
Let $\omega\in\P$ be a one-dimensional cell,
$\sigma_1,\sigma_2\in\P_{\max}$, $\tau\in\P$ with a diagram
\begin{eqnarray}\label{transitiondiag}
\begin{minipage}{4cm}
\xymatrix@C=30pt
{&&\sigma_1\\
\omega\ar@/^/[rru]^{e_1}\ar[r]\ar@/_/[drr]_{e_2}&\tau
\ar[ru]_{g_1}\ar[rd]^{g_2}&\\
&&\sigma_2
}
\end{minipage}
\end{eqnarray}
Then if $T_{\omega}^{e_1e_2}=0$ we must have
$\D(s_{g_1},\omega,\sigma_1) =\D(s_{g_2},\omega,\sigma_2)$.
\end{proposition}

\proof Suppose there exists a normalized $f\in\Gamma
(X,\shLS^+_{\pre,X})$. This is induced by sections
$f_{\sigma}=p_{\sigma}^{-1}(f) \in\Gamma(V(\sigma),
\shLS^+_{\pre,V(\sigma)})$, $\sigma\in\P_{\max}$, which glue as given
in Theorem~\ref{globalgluing}. Now
\[
\shLS^+_{\pre,V(\sigma)}=
\bigoplus_{e\in\coprod\Hom(\omega,\sigma)\atop\dim\omega=1}
\O_{V_e},
\]
so $f_\sigma$ decomposes into components $f_{\sigma,e}\in
\Gamma(\O_{V_e})$. We can write
\[
f_{\sigma,e}=\sum_{p\in P_e} f_{\sigma,e,p}z^p,
\]
where $P_e$ is the face of $P_{\sigma}$ corresponding to $e$, with
$V_e=\Spec k[P_e]$. The normalization condition says that
$f_{\sigma,e,0}=1$. Given any diagram (\ref{transitiondiag}) such
that $(g_1,g_2)$ is a maximal pair,  Theorem~\ref{globalgluing} tells
us that if $n_{\omega}^{e_1e_2}=0$, then
\[
f_{\sigma_2,e_2}|_{V(\tilde\tau_2)} {\D(s_{g_1},\omega,\sigma_1)
\over \D(s_{g_2},\omega,\sigma_2)}
=\Phi_{g_1g_2}(s)^*(f_{\sigma_1,e_1}|_{V(\tilde\tau_1)}).
\]
Comparing constant coefficients gives
\[
f_{\sigma_2,e_2,0}\cdot {\D(s_{g_1}, \omega,\sigma_1) \over
\D(s_{g_2},\omega,\sigma_2)} =f_{\sigma_1,e_1,0}.
\]
But $f_{\sigma,e,0}=1$ for all $\sigma,e$ by the normalization
assumption, so $\D(s_{g_1},\omega,\sigma_1)
=\D(s_{g_2},\omega,\sigma_2)$ whenever $n_{\omega}^{e_1e_2}=0$.
This is the required condition when $(g_1,g_2)$ is maximal. 

Now if $(g_1,g_2)$ is not maximal, then there is a maximal pair
$(h_1,h_2)$ with $(g_1,g_2)<(h_1,h_2)$, i.e. there exists 
$f:\tau\rightarrow\rho$, $h_i:\rho\rightarrow\sigma_i$ with
$g_i=h_i\circ f$. Then
\begin{eqnarray*}
\D(s_{g_i},\omega,\sigma_i)&=&\D(s_{h_i}|_{\tau}\cdot s_f,\omega,\sigma_i)\\
&=&\D(s_{h_i}|_{\tau},\omega,\sigma_i)\D(s_f,\omega,\sigma_i)\\
&=&\D(s_{h_i},\omega,\sigma_i)\D(s_f,\omega,\sigma_i).
\end{eqnarray*}
Now $\D(s_{h_1},\omega,\sigma_1)=\D(s_{h_2},\omega,\sigma_2)$ since
$(h_1,h_2)$ is maximal and $n_{\omega}^{e_1,e_2}=0$, while
$\D(s_f,\omega,\sigma_1)=\D(s_f,\omega,\sigma_2)$ by Remark
\ref{twiddlechange2} and $n_{\omega}^{e_1e_2}=0$. Thus Condition~(LC)
holds.

Conversely, we will show how to produce a normalized section $f$ of
$\shLS^+_{\pre,X}$ by constructing $f_{\sigma}\in
\Gamma(V(\sigma),\shLS^+_{\pre,X})$ which glue correctly. To do so,
first fix for every $\omega\in\P$ with $\dim\omega=1$ a
$\sigma_{\omega}\in\P_{\max}$ with a morphism
$e_{\omega}:\omega\to\sigma_{\omega}$. Let $y\in\omega$ be a point
near $v^+_{\omega}$, and let $\check\Delta(\omega)$ be, as in 
Definition~\ref{NDelta}, the convex hull in
$\Lambda_{\omega,\RR}^{\perp} \subseteq\check\Lambda_{\RR,y}$ of the
set
\[
\{n_{\omega}^{e_{\omega}e'}|e':\omega\to\sigma'\in\P_{\max}\}
\subseteq\Lambda_{\omega}^{\perp}\subseteq\check\Lambda_y.
\]
This is the Newton polytope of the line bundle $\shN_{\omega}$ on
$X_{\omega}$, up to translation. In particular, any 
\[
f_{\sigma_{\omega},e_{\omega}}=\sum_{p\in\check\Delta(\omega)
\cap\check\Lambda_y} f_{\sigma_{\omega},e_{\omega},p} z^p
\]
gives a well-defined regular section of $\shN_{\omega}$ on
$X_{\omega}$. Explicitly, given a diagram
\smallskip
\begin{eqnarray}\label{transitiondiag2}
\begin{minipage}{4cm}
\xymatrix@C=30pt
{&&\sigma_{\omega}\\
\omega\ar@/^/[rru]^{e_{\omega}}\ar[r]\ar@/_/[drr]_{e'}&\tau
\ar[ru]_{g_1}\ar[rd]^{g_2}&\\
&&\sigma'
}
\end{minipage}
\end{eqnarray}
with $(g_1,g_2)$ maximal,
the above choice for $f_{\sigma_{\omega},e_{\omega}}$ determines
$f_{\sigma',e'}$ via the formula (coming from Theorem
\ref{globalgluing})
\begin{eqnarray*}
f_{\sigma',e'}&=&
z^{-n_{\omega}^{e_{\omega}e'}} {\D(s_{g_2},\omega,\sigma')\over\D( 
s_{g_1},\omega,\sigma_{\omega})}
s_{g_2}(n_{\omega}^{e_{\omega}e'})^{-1}\Phi_{g_1g_2}(s)^*
(f_{\sigma_{\omega},e_{\omega}})\\
&=&{\D(s_{g_2},\omega,\sigma')\over\D(s_{g_1},\omega,\sigma_\omega)}
s_{g_2}(n_{\omega}^{e_{\omega}e'})^{-1}
\sum_{p\in\check\Delta(\omega)\cap\check\Lambda_y}
f_{\sigma_{\omega},e_{\omega},p} {\Phi_{g_1g_2}(s)^*(z^p)\over
z^p} z^{p-n_{\omega}^{e_{\omega}e'}}.
\end{eqnarray*}
Here on $X_{\omega}$, $\Phi_{g_1g_2}(s)^*(z^p)\over z^p$ is just
an element of $k^\times$. Now we require that the section be 
normalized, i.e. the constant term of each $f_{\sigma,e}$ be $1$.
This implies that for each diagram (\ref{transitiondiag2}), we must have
\begin{eqnarray*}
1&=&{\D(s_{g_2},\omega,\sigma')\over
\D(s_{g_1},\omega,\sigma_{\omega})}
s_{g_2}(n_{\omega}^{e_{\omega}e'})^{-1} f_{\sigma_{\omega},
e_{\omega},n_{\omega}^{e_{\omega}e'}} {\Phi_{g_1g_2}(s)^*
(z^{n_{\omega}^{e_{\omega}e'}})\over z^{n_{\omega}^{e_{\omega}e'}}}\\
&=&{\D(s_{g_2},\omega,\sigma')\over\D(s_{g_1},\omega,\sigma_{\omega})}
s_{g_1}(n_{\omega}^{e_{\omega}e'})^{-1}f_{\sigma_{\omega},e_{\omega},
n_{\omega}^{e_{\omega}e'}},
\end{eqnarray*}
using $\Phi_{g_1g_2}(s)^*(z^{n_{\omega}^{e_{\omega}e'}})=
s_{g_1}(n_{\omega}^{e_{\omega}e'})^{-1}z^{n_{\omega}^{e_{\omega}e'}}
s_{g_2}(n_{\omega}^{e_{\omega}e'})$.
Thus
\begin{eqnarray}
\label{coefform}
f_{\sigma_{\omega},e_{\omega},n_{\omega}^{e_{\omega}e'}}=
{\D(s_{g_1},\omega,\sigma_{\omega})\over \D(s_{g_2},\omega,\sigma')}
s_{g_1}(n_{\omega}^{e_{\omega}e'}).
\end{eqnarray}
Thus the coefficients
\[
\{f_{\sigma_{\omega},e_{\omega},n_{\omega}^{e_{\omega}e'}}|
e':\omega\to\sigma'\in\P_{\max}\}
\]
of $f_{\sigma_\omega,e_\omega}$ are completely determined by the
normalization condition, while all other coefficients
$f_{\sigma_{\omega},e_{\omega},p}$ may be chosen at will.

There is one point that we must check here, which is that if
$n_{\omega}^{e_{\omega}e'}=n_{\omega}^{e_{\omega}e''}$, then the
values of $f_{\sigma_{\omega},e_{\omega},n_{\omega}^{e_{\omega}e'}}$
and $f_{\sigma_{\omega},e_{\omega},n_{\omega}^{e_{\omega}e''}}$ given
by formula (\ref{coefform}) coincide. Since $n_{\omega}^{e'e''}
=n_{\omega}^{e_{\omega}e''}-n_{\omega}^{e_{\omega}e'}$, it is enough
to check that if $n_{\omega}^{e_{\omega}e'}=0$, then (\ref{coefform})
gives $f_{\sigma_{\omega},e_{\omega},n_{\omega}^{e_{\omega}e'}}=1$.
But this follows immediately from Condition~(LC).
\qed

\begin{theorem}
\label{dimle2}
Let $(B,\P)$ be positive, with either $\Delta=\emptyset$ or $\dim B\le 2$,
and let $s\in Z^1(\P,\shQ_{\P}\otimes\Gm)$ be open gluing data satisfying
Condition~(LC) of Proposition~\ref{gluingcond}. Then if
$X=X_0(B,\P,s)$, any normalized $f\in\Gamma(X,\shLS^+_{\pre,X})$
determines a positive normalized log Calabi-Yau structure on $X$.
\end{theorem}

\proof Write $f=(f_{\sigma,e})$ with $f_{\sigma,e}\in\Gamma(V_e,\O_{V_e})$
for any $e:\omega\to\sigma\in\P_{\max}$, $\dim\omega=1$. Because $f_{\sigma,e}$
is $1$ at the zero-dimensional stratum of $V_e$, the zero set of $f_{\sigma,e}$
cannot contain any toric stratum of $V_e$. Thus if
\[
Z=\bigcup_{\sigma,e} q_{\sigma}(\{f_{\sigma,e}=0\})\subseteq X,
\]
then $Z$ does not contain any toric stratum of $X$.

If $\Delta=\emptyset$, then in fact $Z=\emptyset$ anyway, as
$f_{\sigma,e}$ is constant, hence identically $1$. Then the
multiplicative condition of Theorem~\ref{finalLS} is satisfied
automatically. If $\Delta\not=\emptyset$, then $\dim B= 2$, and
Theorem~\ref{finalLS} imposes a condition for each  two-dimensional,
i.e. maximal, cell $\sigma$. But $X_{\sigma}$ is just a point, and
the formula of Theorem~\ref{finalLS} imposes a condition on the
functions $f_{\sigma,e}|_{X_{\sigma}}$. However since $f_{\sigma,e}$
takes the value $1$ at this point by the normalization condition, the
multiplicative condition of Theorem~\ref{finalLS} is automatically
satisfied, and $f\in\Gamma(X\setminus Z,\shLS_X)$.
\qed

\begin{remark} 
In the case $\Delta=\emptyset$, the open gluing data in fact
determine the log Calabi-Yau structure uniquely. However, in the
two-dimensional case, when $\Delta\not=\emptyset$, the open gluing
data do not necessarily determine the log Calabi-Yau structure
uniquely. Indeed, if $\omega$ contains a singularity of $B$ with
holonomy $\begin{pmatrix} 1&n\\ 0&1\end{pmatrix}$ (as in
Example~\ref{normal crossings}) then
$\shN_{\omega}\cong\O_{\PP^1}(n)$, and $\check\Delta(\omega)$ is a line
segment of length $n$. The coefficients of the monomials
corresponding to the endpoints are determined by (\ref{coefform}),
but the coefficients corresponding to  the interior vertices are free
to be chosen arbitrarily in $k$. Thus for given open gluing data $s$
satisfying Condition~(LC), the set of normalized log Calabi-Yau
structures on $X_0(B,\P,s)$ is parametrized by a vector space of some
dimension. Specifically, if $\Delta=\{p_1,\ldots,p_s\}$ with
monodromy around $p_i$ given by $\begin{pmatrix}1&n_i\\
0&1\end{pmatrix}$, then the dimension of the space of positive
normalized log Calabi-Yau  structures on $X_0(B,\P,s)$ is
$\sum_{i=1}^s (n_i-1)$. In particular, if $n_i=1$ for all $i$, then
the log Calabi-Yau structure is uniquely determined by $s$.
\qed
\end{remark}

\begin{example}
\label{obstructed}
Theorem~\ref{dimle2} is not true in dimensions three and higher, and
it is an important feature of the theory that the space of possible
log structures is not even non-singular in general. We will now give
a quite general local example, which can be fitted into a global
example, featuring two maximal cells as in
Construction~\ref{basicgluing}.  Let $M_1=M_2=M=\ZZ^n$, and choose
some primitive vector  $\check d_\rho\in N$. Let $M'=\check
d_\rho^{\perp}\subseteq M$,  and let $\rho\subseteq M'_{\RR}$ be an
$n-1$-dimensional lattice polytope. Choose $\sigma_1\subseteq
M_1\otimes\RR$, $\sigma_2\subseteq M_2\otimes\RR$ to be
$n$-dimensional lattice polytopes with $\sigma_i\cap M'_\RR=\rho$ and
such that 
\[
\sigma_1\subseteq \{m\in M_1\otimes\RR| \langle \check d_\rho,m\rangle\le 0\}
\]
and
\[
\sigma_2\subseteq \{m\in M_2\otimes\RR| \langle \check
d_\rho,m\rangle\ge 0\}.
\]
Next, we define a singular affine structure with boundary
on $\sigma_1\cup\sigma_2$ as follows. Embed $M_1$ in $M\oplus\ZZ$ by
$m\mapsto (m,1)$ and $M_2$ in $M\oplus\ZZ$ by $m\mapsto
(m,1-\langle\check d_\rho,m\rangle)$. Thus we obtain
$\sigma_1\cup\sigma_2\subseteq M_{\RR}\oplus\RR$. To define the
desired affine structure, we have the canonical embeddings
$\Int(\sigma_i)\subseteq (M_i)_\RR$, giving affine coordinates on
$\Int(\sigma_i)$. In a neighbourhood $U_v$ of each vertex $v\in M'$
of $\rho$, $(v,1)\in M\oplus\ZZ$, we obtain a chart $U_v\to
(M_{\RR}\oplus\RR)/\RR(v,1)$ simply by projection. It is easy to see
this is a homeomorphism onto its image. Take $\Delta\subseteq\rho$ to
be the union of codimension one simplices of the first barycentric
subdivision of $\rho$ not containing vertices of $\rho$. We can then
choose open sets $U_v$ so that $\{\Int(\sigma_1),\Int(\sigma_2)\}
\cup \{U_v\}$ form an open covering of $(\sigma_1 \cup\sigma_2)
\setminus\Delta$ and $U_v\cap U_w=\emptyset$ if $v\not=w$. It follows
that the affine coordinate charts defined on these sets define an
integral affine structure. Parallel transport $\psi_v:M_2\to M_1$
through (or near) a vertex $v$ of $\rho$ is given by the composition
\[
M_2\hookrightarrow M\oplus\ZZ\to (M\oplus\ZZ)/\ZZ(v,1)
\to M_1,
\]
where the second map is projection and the third is the inverse of the
composed map
\[
M_1\hookrightarrow M\oplus\ZZ\to (M\oplus\ZZ)/\ZZ(v,1).
\]
Thus $\psi_v$ is the composition
\begin{eqnarray*}
M_2\ni m &\mapsto&(m,1-\langle\check d_\rho,m\rangle)\\
&\mapsto&(m,1-\langle\check d_\rho,m\rangle)\mod (v,1)\\
&=&(m+\langle\check d_\rho,m\rangle v,1)\mod (v,1)\\
&\mapsto&m+\langle\check d_\rho,m\rangle v\in M_1
\end{eqnarray*}
Thus if $\omega$ is an edge of $\rho$, $d_{\omega}$ a generator
of $\Lambda_{\omega}$, with $v^-_{\omega}
-v_{\omega}^+=pd_{\omega}$ for some positive integer $p$, then
\begin{eqnarray*}
T_{\omega\rightarrow\rho}(m)&=&\psi_{v^-_{\omega}}\circ
\psi_{v^+_{\omega}}^{-1}(m)\\
&=&m+\langle p\check d_\rho,m\rangle d_{\omega},
\end{eqnarray*}
so $n_{\omega\rightarrow\rho}=p\check d_\rho$.

We will now restrict to the case that $p=1$ for all edges $\omega$ of
$\rho$, for ease of the discussion, and to make contact with the work
of Altmann \cite{Altmann}. We now have as usual, as in
Construction~\ref{basicgluing}, $\rho=\rho_1\subseteq\sigma_1$,
$\rho=\rho_2\subseteq\sigma_2$, $V(\rho_i)\subseteq V(\sigma_i)$, and
after a choice of $s_1$ and $s_2$, an isomorphism 
$\Phi_{\sigma_1\sigma_2}(s):V(\rho_2)\to V(\rho_1)$. Gluing along
this isomorphism gives $X$.

Note that $X$ contains a one-dimensional stratum
$X_{\rho}\cong\PP^1$, and for each edge $\omega$ of $\rho$, we have
the codimension 1 stratum $X_{\omega}$ containing $X_{\rho}$. In
addition, we have the line bundle $\shN_{\omega}$ on $X_{\omega}$,
and the value of $n_{\omega\rightarrow\rho}=\check d_\rho$ tells us 
that $\shN_{\omega}|_{X_{\rho}} =\O_{\PP^1}(1)$. Thus a normalized
section  $f\in\Gamma(X,\shLS^+_{\pre,X})$ must satisfy
\[
f_{\sigma_1,\omega\to
\sigma_1}|_{X_{\rho}}=1+c_{\omega}x,
\]
where $x=z^{\check d_\rho}$. Here $c_{\omega}$ is completely
determined by $s_1$ and $s_2$ by formula (\ref{coefform}). We note
that the choice of $s_1$ and $s_2$ may not lead to an arbitrary
choice of the $c_{\omega}$'s. In addition, if $f$ is to determine a
log Calabi-Yau structure, we need the multiplicative condition of
Theorem~\ref{finalLS} to hold. In other words, for every
two-dimensional cell $\eta\subseteq \rho$, we have
\[
\prod_{\omega\subseteq\eta\atop\dim\omega=1} d_{\omega}\otimes
f_{\sigma_1,\omega\to
\sigma_1}^{\epsilon_{\eta}(\omega)}|_{X_{\eta}}=1.
\]
Restricting further to $X_{\rho}$ gives the relation
\begin{eqnarray}
\label{specialrelation}
\prod_{\omega\subseteq\eta\atop \dim\omega=1}
d_{\omega}\otimes (1+c_{\omega}x)^{\epsilon_{\eta}(\omega)}
=1.
\end{eqnarray}
Let us carry this out explicitly for the following $\rho\subseteq
M'_{\RR}$,  with orientations as depicted:

\begin{center}
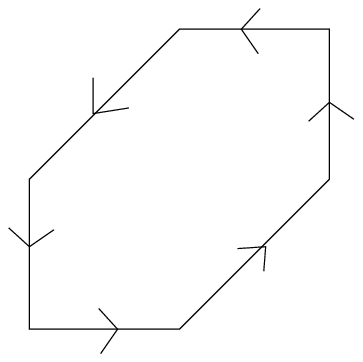
\end{center}

Then we obtain in $M'$
\[
\hbox{$d_{\omega_i}=(1,0),(1,1),(0,1),(-1,0),(-1,-1),(0,-1)$ for
$i=1,\ldots,6$.}
\]
Splitting (\ref{specialrelation}) into its two components, we get the 
equations (writing $c_i=c_{\omega_i}$)
\begin{eqnarray*}
(1+c_1x)(1+c_2x)&=&(1+c_4x)(1+c_5x)\\
(1+c_2x)(1+c_3x)&=&(1+c_5x)(1+c_6x)
\end{eqnarray*}
Comparing coefficients of powers of $x$, we get
\begin{eqnarray*}
c_1+c_2&=&c_4+c_5\\
c_2+c_3&=&c_5+c_6\\
c_1c_2&=&c_4c_5\\
c_2c_3&=&c_5c_6
\end{eqnarray*}
This defines a three-dimensional scheme $S$ in $\Gm^6$, with
irreducible components being a $\Gm^3$ given by the equations
\[
c_1=c_4,\quad c_2=c_5\quad c_3=c_6
\]
and a $\Gm^2$ given by the equations
\[
c_2=c_4=c_6,\quad c_1=c_3=c_5.
\]
As we noted earlier, because the $c_i$ depend on
$s_1$ and $s_2$, the $c_i$ cannot be chosen independently,
but it is easy to see that the only restriction on the $c_i$ is that
\[
\prod_{i=1}^6 d_{\omega_i}\otimes c_i=(1,1),
\]
but this condition has already been incorporated into the above
equations. It can also be checked that multiplying $s_1$ or $s_2$
by a multiplicative function $N_i\to\Gm$ has the effect of
replacing $c_1,\ldots,c_6$ with $c' c_1,\ldots, c' c_6$. This does
not change the isomorphism class of the gluing, and one obtains
isomorphic log Calabi-Yau structures from these gluings. Hence the
space of all open gluing data modulo this action fibres over $S/\Gm$,
where the action of $\Gm$ on $S$ is the diagonal one, and $S/\Gm$
consists of a union of a one-dimensional torus and a two-dimensional
torus meeting at a single point.

The equations determining $S/\Gm$ are precisely the same
ones defining the deformation space of a cone over a del Pezzo
surface of degree 6, as computed in \cite{Altmann}. The hexagon
we are using of course is the toric data which
defines the del Pezzo surface of degree 6.

The explanation is as follows. Let $Y$ be the four-dimensional toric
variety defined by the fan $\Sigma$ in $M_{\RR}\oplus\RR$ consisting
of the cones over $\sigma_1$ and $\sigma_2$ and their faces embedded
in $M_{\RR}\oplus\RR$ as above. Then the cone over $\rho$ is a cone
in $\Sigma$, corresponding to a one-dimensional stratum $Y_{\rho}$ of
$Y$, with $Y_{\rho}\cong\PP^1$. The singularities of $Y$ along
$Y_{\rho}$ are locally isomorphic to a cone over a del Pezzo surface
of degree 6. If we take a general anti-canonical hypersurface $H$ in
$Y$, then $H\cap Y_{\rho}$ consists of one point, and $H$ has a del
Pezzo cone singularity at that point. However, $H$ can be viewed as a
partial smoothing of the anti-canonical hypersurface $H_0$ consisting
of the union of all toric divisors of $Y$, and this is precisely
$X_0(B,\P,1)$, where $1$ denotes the trivial gluing data $s_1=s_2=1$.
The point is that we expect the moduli of log Calabi-Yau structures on
$X_0(B,\P,1)$ should be in some ways similar to (but not isomorphic to)
the moduli space of $H$, and thus it
makes sense to see the same structures appearing in $Def(H)$ and the
moduli of log Calabi-Yau structures on $H_0$.

While this example is local, one expects that there are global
examples where the space of positive log Calabi-Yau spaces coming
from general $(B,\P)$ can be singular because of the above
considerations.
\end{example}

It is not difficult to prove in three dimensions that $X_0(B,\P,1)$
(trivial gluing data) always carries a log Calabi-Yau structure.
However, it is not clear if this is true in higher dimensions, and it is
not so useful because to study mirror symmetry, we would like  to
describe the entire space of positive log Calabi-Yau spaces with dual
intersection complex $(B,\P)$. Therefore, in the next section, we
will restrict to situations where the behaviour of
Example~\ref{obstructed} does not occur.

\begin{example}
\label{degenexamples}
1)\ \ The toric degeneration of Example~\ref{toric degenerations
examples},~(1) has dual intersection complex $(B,\P)$ as defined in
Examples~\ref{polytope} and \ref{decompexamp}, (2). Details of
this case and a much more general case including the Batyrev-Borisov
construction are treated in \cite{GBB}.

\noindent
2)\ \ The dual intersection complex for the toric degeneration in
Example~\ref{toric degenerations examples},~(2), arising from the
example of Aspinwall and Morrison, has dual intersection complex as
constructed in Examples~\ref{aspmorr} and ~\ref{decompexamp},~(2).

\noindent
3)\ \ Finally, Example~\ref{enriques} gives the dual intersection
complex of a degeneration of Enriques surfaces obtained by taking a
quotient of a family of K3 surfaces in $\PP_{\Xi^*}$. We leave the
details to the reader.
\qed
\end{example}
\bigskip

We end this section with the following observation and definition.

\begin{proposition}
\label{equiv}
Let $s,s'$ be open gluing data for $(B,\P)$, and suppose there are
log Calabi-Yau structures $X=X_0(B,\P,s)^{\ls}$ and $X'
=X_0(B,\P,s')^{\ls}$ normalized by the gluing data $s$ and $s'$
respectively, and suppose there is an isomorphism of log spaces
$\varphi: X\to X'$ preserving $B$. Then there exists
$t=(t_{\sigma})_{\sigma\in\P}$, $t_{\sigma}\in \check{PM}(\sigma)$, such that
$s_e'=t_{\tau}^{-1}s_e t_{\sigma}|_\tau$ for any $e:\tau\to\sigma$. 
Furthermore, for $\sigma\in\P_{\max}$, $t_{\sigma}$ is induced by an
element of $\Lambda_y\otimes\Gm$ for $y\in \Int(\sigma)$. In this
case, we say $s$ and $s'$ are equivalent open gluing data.
\end{proposition}

\proof The only thing new over Proposition~\ref{isomorphism} is the
last statement. As in loc.\ cit.\ the isomorphism $\varphi$ is
induced by a map
\[
\tilde\varphi:\coprod_{\sigma\in\P_{\max}} V(\sigma)\to
\coprod_{\sigma\in\P_{\max}} V(\sigma)
\]
which each $\tilde \varphi_{\sigma}:V(\sigma)\to V(\sigma)$ induced
by some $t_{\sigma}\in \check{PM}(\sigma)$. However, if $p_{\sigma}:V(\sigma)
\to X$, $p_{\sigma'}:V(\sigma)\to X'$ are the projections,
$f\in\Gamma(X,\shLS^+_{\pre,X})$, $f'
\in\Gamma(X',\shLS^+_{\pre,X'})$ inducing the log Calabi-Yau
structures, then $p_{\sigma}^{-1}(f)$ and $(p_{\sigma}')^{-1}(f)$ are
both normalized, and
$\tilde\varphi_{\sigma}^{-1}p_{\sigma}^{-1}(f)=(p_{\sigma}')^{-1}(f')$.
However, as we pointed out in the discussion following
Definition~\ref{LS^+_pre}, the pull-back of a normalized section of
$\shLS^+_{\pre,V(\sigma)}$ via $t_{\sigma}$ is normalized if and only
if $t_{\sigma}$ is a multiplicative function on
$\check\Sigma_{\sigma}$ (for $\sigma\in\P_{\max}$) rather than
piecewise multiplicative, i.e. comes from an element of
$\Lambda_y\otimes\Gm$ for $y\in \Int(\sigma)$. We then complete the
proof as in Proposition~\ref{isomorphism}.
\qed

\section{Simplicity and mirror symmetry}
\label{section5}

We are now close to our first goal of realising mirror symmetry for
toric log Calabi-Yau spaces. The chief problem at this point,
however, is that for general positive $(B,\P)$, the moduli space of
all positive log Calabi-Yau spaces with dual intersection complex
$(B,\P)$ may be quite complicated. We have already seen in
Example~\ref{obstructed} that this moduli space may well be singular,
and that this reflects the fact that a partial smoothing of a log
Calabi-Yau space, if it exists, might itself have obstructed
deformation theory.  On the mirror side this should correspond to
singular moduli spaces of  complexified K\"ahler or symplectic
structures, and it is unclear what this  should mean. Also, for more
refined versions of mirror symmetry involving  Frobenius or
$A_\infty$-structures it is convenient if not indispensable to deal
with smooth moduli spaces. Therefore, we need to restrict attention
to a class of log Calabi-Yau spaces which are well-behaved. We do so
by restricting to the case that $(B,\P)$ is simple (\S1.5).

Once again, as in \S4, we always work over $S=\Spec k$, $k$ an
algebraically closed field.

\subsection{Moduli of toric log Calabi-Yau spaces in the simple case}
\label{section51}

Consider $e:\tau_1\to\tau_2$ and an element $s\in \Gamma(W_{e},
i_*\Lambda\otimes\Gm)$. This in fact defines an element $s\in
\check{PM}(\tau_1)$ as follows. Choose any $\tau_2 \to\sigma\in\P_{\max}$.
Then $\Int(\sigma)\cap W_{e}$ is non-empty, and for any point $z\in
\Int(\sigma)\cap W_{e}$, $s$ yields a germ $s_z\in\Lambda_z\otimes
\Gm$.

Now for each vertex $v$ of $\tilde\tau_1
\subseteq\tilde\sigma\subseteq \Lambda_{\RR,z}$, parallel transport
$s_z$ to $s_v\in \Lambda_{\pi(v)} \otimes\Gm=
\pi^*(\shQ_{\P}\otimes\Gm)_v$ via the image of a path in
$\tilde\sigma$ joining $z$ to $v$. It is clear that the collection
$(s_v)$ determines an element in $\check{PM}(\tau_1)$ as in
Construction~\ref{generalgluing2}. Furthermore, this element is
independent of the choice of $\sigma$: because $s_z$ is the germ of
$s$, a section of $i_*\Lambda\otimes\Gm$, it is invariant under
monodromy about any loop in $W_e$. Thus each $s_v$ is independent of
$\sigma$.

\begin{definition}
By \emph{lifted gluing data}, we mean a \v Cech 1 cocycle
$(s_e)$ for the open cover $\W=\{W_{\tau}| \tau\in\P\}$ with values
in  $i_*\Lambda\otimes\Gm$. For $e:\tau_1\to\tau_2$, 
$s_e\in\Gamma(W_e,i_*\Lambda\otimes\Gm)$ induces $s_e\in \check{PM}(\tau_1)$
by the above discussion, and hence lifted gluing data
$(s_e)$ determine open gluing data.
\end{definition}

The main point of simplicity is then

\begin{theorem} 
\label{lifted}
Let $(B,\P)$ be positive and  simple, and let $s$ be open gluing data
over an algebraically closed field $k$ satisfying Condition~(LC) of
Proposition~\ref{gluingcond}. Then
\item{(1)} $s$ is equivalent (Proposition~\ref{equiv}) 
to lifted gluing data, and
\item{(2)} there is a unique normalized section 
$f\in\Gamma(X_0(B,\P,s),\shLS^+_{\pre,X_0(B,\P,s)})$, and this
section induces a log Calabi-Yau structure on $X_0(B,\P,s)$.
\end{theorem}

\begin{remark}
Before embarking on the proof, it is worthwhile giving a quick
explanation as to why lifted gluing data are enough. Consider the
situation in Example~\ref{twotriangles}. Given
$e_i:\tau\rightarrow\sigma_i$, and lifted gluing data
$s_{e_i}\in\Gamma(W_{e_i},i_*\Lambda \otimes\Gm)$, the gluing map is
$\Phi_{e_2e_1}(s)=s_{e_1}^{-1}\circ \Phi_{e_2e_1}\circ s_{e_2}$. Now
both $s_{e_1}$ and $s_{e_2}$ can be viewed as elements of $\check{PM}(\tau)$,
and as such, we also have $\Phi_{e_2e_1}(s)=s_{e_1}^{-1}\circ
s_{e_2}\circ\Phi_{e_2e_1}$. However, to view $s_{e_2}$ as acting on
$V(\tau_1)\subseteq V(\sigma_1)$ instead of $V(\tau_2)\subseteq
V(\sigma_2)$, one parallel transports $s_{e_2}$ through the vertex
$(0,0)$ to determine the action of $s_{e_2}$ on the irreducible
component of $V(\tau_1)$ corresponding to $(0,0)$, and one parallel
transports through the vertex $(0,1)$ to describe the action on the
corresponding component. If $s_{e_2}$ is not monodromy invariant in a
neighbourhood of the singular point, the induced action on $V(\tau_1)$
does not then arise from an element of $\Gamma(W_{e_1},
i_*\Lambda\otimes\Gm)$. The choice of $s_{e_1}$ and $s_{e_2}$ then
allow all possible automorphisms of $V(\tau_1)$.

More generally, the idea will be to try to split up automorphisms
of $V(\tau)$ into a product of ones coming from lifted gluing data.
\end{remark}

\proof (1) By the definition of equivalence in Proposition
\ref{equiv}, we need to find a system $(t_{\tau})_{\tau\in\P}$ with
$t_{\tau}\in \check{PM}(\tau)$ for all $\tau\in\P\setminus\P_{\max}$, and
$t_{\tau}\in \Lambda_z\otimes\Gm$ for $z\in \Int(\tau)$, $\tau\in
\P_{\max}$, such that for each $f:\tau_1\rightarrow\tau_2$,
$s'_f=t_{\tau_1}^{-1}s_ft_{\tau_2}|_{\tau_1}$  is induced by an
element of $\Gamma(W_f,i_*\Lambda\otimes\Gm)$.
\smallskip

\noindent
\emph{Step 1}. Let $\tau\in\P$, $0<\dim\tau<\dim B$ be given, and let
$\Omega_i\subset \P_1(\tau)$, $R_i\subset \P_{n-1}(\tau)$, $\Delta_i$,
$\check\Delta_i$ be the objects associated to $\tau$ by the definition
of simplicity (Definition \ref{simplicity}). Let
$e_i:\tau\rightarrow\sigma_i\in\P_{\max}$. Then for $f\in\Omega_i$,
$n_{\omega}^{e_1\circ f,e_2\circ f}$ and 
$\D(s_{e_1},f,e_1)/\D(s_{e_2},f,e_2)$ are independent of the choice of
$f\in\Omega_i$. 

\proof The first statement is Remark~\ref{simplicityfacts}, (1). To
see that $\D(s_{e_1},f,e_1)/ \D(s_{e_2},f,e_2)$ is independent of
$f\in\Omega_i$, we note it is enough to show
$\D(s_{e_1},f,e_1)/\D(s_{e_2},f,e_2)$ is independent of $f$ in
$\Omega_i$ whenever we have a diagram
\[
\xymatrix@C=30pt
{&&&\sigma_1\\
\omega\ar[r]^{f}&\tau\ar@/^/[rru]^{e_1}\ar[r]^{g}\ar@/_/[drr]_{e_2}&\rho
\ar[ru]\ar[rd]&\\
&&&\sigma_2
}
\]
with $\dim\rho=n-1$. If $g\not\in R_i$, then $n_{g\circ
f}=n_{\omega}^{e_1\circ f, e_2\circ f}=0$ by Def.~\ref{simplicity},
(1), so by Condition~(LC), $\D(s_{e_1},f,e_1)/\D(s_{e_2},f,e_2)=1$.
Thus it is enough to show this independence for $g\in R_i$. Note that
by Remark~\ref{simplicityfacts}, (3), we can then assume
\[
\Omega_i=\{f:\omega\to\tau| n_{g\circ f}\not=0\}.
\]

Now let $z_i\in \Int(\sigma_i)$, so we can view
$\tilde\sigma_i\subseteq \Lambda_{\RR,z_i}$ and
$\tilde\tau_i\subseteq\tilde\sigma_i$, the face of $\tilde\sigma_i$
determined by $e_i$. As usual, $s_{e_i}$ can be viewed as a function
$s_{e_i}:\check\Lambda_{z_i}\to\Gm$ which is  piecewise
multiplicative on the fan $\check\tau_i^{-1}\check\Sigma_i$. It is
given by $s_{e_i}^v\in\Lambda_{z_i}\otimes\Gm$ on the cone of
$\check\tau_i^{-1}\check\Sigma_{\sigma_i}$ corresponding to a vertex
$v$ of $\tilde\tau_i$. Suppose we are given a cycle of edges
$\tilde\omega_1,\ldots,\tilde\omega_n$ of $\tilde\tau$, with vertices
of $\tilde\omega_j$ being $v_j$ and $v_{j+1}$, with $v_{n+1}=v_1$.
This cycle need not be contained in a 2-face. Specifying
$\tilde\omega_j$ is the same as specifying an $f_j:\omega_j\to\tau$.
Then we have (with the signs chosen to correspond to the orientation
of the cycle)
\[
\prod_{j=1}^{n}\D(s_{e_i},f_j,e_i)\otimes d_{\omega_j}^{\pm
1}={s^{v_1}_{e_i}\over s^{v_2}_{e_i}}{s^{v_2}_{e_i}\over s^{v_3}_{e_i}}
\cdots{s_{e_i}^{v_{n}}\over s_{e_i}^{v_1}}=1.
\]
Since $d_{\omega_j}$ is monodromy invariant in a neighbourhood of
$\Int(\tau)$, we can parallel transport this identity for $i=2$ into
$\Lambda_{z_1}\otimes\Gm$, and so obtain
\begin{eqnarray}
\label{cyclecond}
\prod_{j=1}^{n}{\D(s_{e_1},f_j,e_1)\over
\D(s_{e_2},f_j,e_2)}\otimes d_{\omega_j}^{\pm 1}&=&1.
\end{eqnarray}
Note that if $f_j\not\in\Omega_i$, then $n_{g\circ f_j}=0$,
and by Condition~(LC), $\D(s_{e_1},f_j,e_1)/\D(s_{e_2},f_j,e_2)=1$,
so only those $f_j\in\Omega_i$ contribute to the above product.

Now define a piecewise linear function $\varphi$ on the fan 
$\check\Sigma_{\tau}$ by
\[
\varphi(x)=-\inf\{\langle y,x\rangle | x\in\Delta_i\}.
\]
(See Remark~\ref{polytopefacts}.) Then $\varphi$ is given by a linear
function $q_h\in\Lambda_{\tau}$ on each maximal cone of
$\check\Sigma_{\tau}$ corresponding to $h:v\rightarrow\tau$. As in
Remark~\ref{polytopefacts}, we then have $m^{\rho}_{g\circ h_1,g\circ
h_2}=q_{h_2}-q_{h_1}\in\Lambda_{\tau}$. Now there is a fan
$\check\Sigma_{\tau}'$, consisting of (not necessarily strictly
convex) cones which are the maximal domains of linearity  of
$\varphi$. The fan $\check\Sigma_{\tau}$ is a refinement of
$\check\Sigma_{\tau}'$. There is a one-to-one order reversing
correspondence between cones of $\check\Sigma_{\tau}'$ and faces of
$\Delta_i$. Of course there is also the usual one-to-one order
reversing correspondence between cones of $\check\Sigma_{\tau}$ and
faces of $\tilde\tau$. This gives, for every face $\eta\subseteq
\tilde\tau$, a face $\eta'\subseteq\Delta_i$, as follows. The face
$\eta$ corresponds to a cone $K_{\eta}$ in $\check\Sigma_{\tau}$,
which is contained in a cone $K'_{\eta}$ of $\check\Sigma_{\tau}'$,
which corresponds to a face $\eta'\subseteq\Delta_i$. This
correspondence does not increase dimensions, i.e.\ $\Delta_i$ is
obtained by partial collapse from $\tilde\tau$. In particular, an edge
of $\tilde\tau$ corresponding to $f\in\P_1(\tau)$ is mapped to an edge
of $\Delta_i$ if and only if $f\in\Omega_i$. Otherwise it is mapped to
a vertex. If $f_j:\omega_j\rightarrow\tau$, $j=1,2$, map to the same
edge of $\Delta_i$, then we can find a cycle in $\tilde\tau$ of edges
containing $f_1$ and $f_2$ and such that all other edges of the cycle
map to vertices of $\Delta_i$. Then $d_{\omega_1} =\pm
d_{\omega_2}=:d$ is parallel to the corresponding edge of
$\Delta_i$,  and $\omega_1$ and $\omega_2$ are traversed in the cycle
in opposite directions, so that (\ref{cyclecond}) yields
\[
\left({\D(s_{e_1},f_1,e_1)\over \D(s_{e_2},f_1,
e_2)}\otimes d\right)\cdot
\left({\D(s_{e_1},f_2,e_1)\over \D(s_{e_2},f_2,
e_2)}\otimes d\right)^{-1}=1.
\]
Thus $\D(s_{e_1},f_j,e_1)/\D(s_{e_2},f_j,e_2)$ is independent
of $j$.

If $\dim\Delta_i=1$, then we are now done. Otherwise, any
two-face $\eta$ of $\Delta_i$ is an elementary, hence standard,
simplex. Thus we can find a cycle of edges of $\tilde\tau$ with only
three of its edges  $f_j:\omega_j\rightarrow\tau$ ($j=1,2,3$) in
$\Omega_i$, and such that these edges map to the three edges of
$\eta$. Also, $d_{\omega_1}$, $d_{\omega_2}$ and $d_{\omega_3}$ span
the tangent space to $\eta$, so without loss of generality we can
take (in a suitable basis) $d_{\omega_1}=(1,0)$, $d_{\omega_2}=(0,1)$
and $d_{\omega_3} =(-1,-1)$. Then (\ref{cyclecond}) is equivalent to 
\[
\left({\D(s_{e_1},f_1,e_1)\over \D(s_{e_2},f_1,
e_2)}\otimes (1,0)\right)
\left({\D(s_{e_1},f_2,e_1)\over \D(s_{e_2},f_2,
e_2)}\otimes (0,1)\right)
\left({\D(s_{e_1},f_3,e_1)\over \D(s_{e_2},f_3,
e_2)}\otimes (-1,-1)\right)=1,
\]
which implies $\D(s_{e_1},f_j,e_1)/\D(s_{e_2},f_j,e_2)$
is independent of $j$.

Since this is true of every 2-face of $\Delta_i$, one sees 
$\D(s_{e_1},f,e_1)/\D(s_{e_2},f,e_2)$ is independent of
$f\in\Omega_i$, as desired. This proves the claim. \qed
\smallskip

\noindent
\emph{Step 2}. Let $\tau\in\P$ with $0<\dim\tau<\dim B$, and fix some
$e:\tau\rightarrow\sigma\in\P_{\max}$, $z\in \Int(\sigma)$. Then
there is an element $u_{\tau}\in\Lambda_z\otimes\Gm(k)$, or
alternatively a homomorphism $u_{\tau}:\check\Lambda_z\to\Gm(k)$,
with the property that for every $e':\tau\to\sigma'\in \P_{\max}$ and
$f:\omega\to\tau$ with $\dim\omega=1$,
\begin{eqnarray}
\label{utau}
u_{\tau}(n_{\omega}^{e\circ f,e'\circ f})={\D(s_{e'},f,e')\over
\D(s_e,f,e)} s_e(n_{\omega}^{e\circ f,e'\circ f})^{-1}.
\end{eqnarray}

\proof 
Let $\Omega_i\subseteq\P_1(\tau)$, $R_i\subseteq\P_{n-1}(\tau)$,
$\Delta_i$, $\check\Delta_i$ be as in Definition~\ref{simplicity}.
First, whenever $n_{\omega}^{e\circ f,e'\circ f}=0$,
Condition~(LC) implies 
\[
\D(s_{e'},f,e')/\D(s_{e},f,e)=1,
\]
so (\ref{utau}) holds for any choice of $u_\tau$. Furthermore, if
$f:\omega\to\tau$ is not in $\bigcup_{i=1}^p\Omega_i$, then by
Remark~\ref{simplicityfacts},~(2) $n_\omega^{e\circ f,e'\circ f}=0$. Thus
we don't need to worry about such $f$. 

Next, let $T_{\check\Delta_i}\subseteq\dualvs{\shQ_{\tau,\RR}}
\subseteq\check\Lambda_{\RR,z}$ denote the tangent space to the
polytope $\check\Delta_i$. By Remark~\ref{simplicityfacts}, (4),
$T_{\check\Delta_1}+\cdots+T_{\check\Delta_p}$ form an interior direct  sum in
$\check\Lambda_{\RR,z}$. Thus for any collection of homomorphisms
$u_i:T_{\check\Delta_i}\cap \check\Lambda_z\to\Gm(k)$, $i=1,\ldots,p$, there
is a  homomorphism $u:\check\Lambda_z\to\Gm(k)$ extending each of the
$u_i$'s. (Here we use the fact that $\Gm(k)$ is divisible, so that
any homomorphism defined on a sublattice of $\check\Lambda_z$ extends
to one on $\check\Lambda_z$. It is here we use that $k$ is
algebraically closed.) Thus, since for $f:\omega\to \tau$ in
$\Omega_i$, $n_{\omega}^{e\circ f,e'\circ f}\in T_{\check\Delta_i}
\cap\check\Lambda_z$, it will be enough to show that we can construct
$u_i:T_{\check\Delta_i}\cap\check\Lambda_z\to\Gm$ such that
\begin{eqnarray}
\label{utaucond}
u_i(n_{\omega}^{e\circ f,e'\circ f})={\D(s_{e'},f,e')\over
\D(s_{e},f,e)}s_e(n_{\omega}^{e\circ f,e'\circ f})^{-1}
\end{eqnarray}
for all $f:\omega\to\tau$ in $\Omega_i$, for all $e':\tau
\to \sigma'$.

By Step 1, $n_{\omega}^{e\circ f,e'\circ f}$ and
$\D(s_{e'},f,e')/\D(s_e,f,e)$ are independent of $f\in\Omega_i$.
Because $\check\Delta_i$ is an elementary simplex which is the convex hull
of  $\{n_{\omega}^{e\circ f,e'\circ f}| e':\tau\to\sigma'\}$, the
non-zero elements of this set are linearly independent. In addition,
if $n_{\omega}^{e\circ f,e'\circ f} =n_{\omega}^{e\circ f,e''\circ
f}$, then $0=n_{\omega}^{e\circ f,e''\circ f} -n_{\omega}^{e\circ
f,e'\circ f}=n_{\omega}^{e'\circ f,e''\circ f}$, and Condition~(LC)
then shows that
\[
{\D(s_{e'},f,e')\over
\D(s_e,f,e)} s_e(n_{\omega}^{e\circ f,e'\circ f})^{-1}
={\D(s_{e''},f,e'')\over
\D(s_e,f,e)} s_e(n_{\omega}^{e\circ f,e''\circ f})^{-1}.
\]
Thus there exists a homomorphism $u_i:T_{\check\Delta_i}\cap
\check\Lambda_z\to\Gm$ satisfying (\ref{utaucond}) and hence a
homomorphism $u_{\tau}:\check\Lambda_z\to\Gm$ satisfying
(\ref{utau}).
\qed
\smallskip

\noindent
\emph{Step 3}. The open gluing data $(s_e)$ are equivalent to
lifted gluing data.

\proof
We construct $(t_{\tau})_{\tau\in\P}$ as follows. For
$\tau\in\P_{\max}$ or $\dim\tau=0$, set $t_{\tau}=1$. For any other
$\tau$, fix $e:\tau\rightarrow\sigma\in\P_{\max}$ and take
$t_{\tau}=s_eu_{\tau}$ with $u_{\tau}$ as in Step 2. We need to show
for each $f:\tau_1\rightarrow\tau_2$,
$s_f'=t_{\tau_1}^{-1}s_ft_{\tau_2}|_{\tau_1}$ is induced from
$\Gamma(W_f,i_*\Lambda\otimes\Gm(k))$.

First note that for $f:\tau_1\rightarrow\tau_2$, $s\in \check{PM}(\tau_1)$ is
induced by an element of $\Gamma(W_f,i_*\Lambda\otimes\Gm(k))$ if and
only if for all $g:\omega\rightarrow\tau_1$ with $\dim\omega=1$ and
$e:\tau_2\rightarrow \sigma\in\P_{\max}$, $\D(s,g,e\circ f)=1$.
Indeed, this says that for any vertex $v$ of $\tilde\tau_1$ and any
$\sigma\in\P_{\max}$ containing $\tau_2$, the parallel transport of
$s^v\in\Lambda_{\pi(v)}\otimes\Gm(k)$ into $z\in\Int(\sigma)$ is
independent of $v$. This gives a well-defined element of
$\Lambda_z\otimes\Gm(k)$, which is then invariant under monodromy in
$W_f$, as desired.

We first consider $e':\tau\to\sigma'\in\P_{\max}$, and wish to show
$s'_{e'}=t_{\tau}^{-1}s_{e'}$ is in
$\Gamma(W_{e'},i_*\Lambda\otimes\Gm(k))$. This is trivial if
$\dim\tau=0$. Otherwise, applying Remark~\ref{twiddlechange2} to both
$s_e$ and $u_{\tau}$ for $g:\omega\rightarrow\tau$, $\dim\omega=1$,
and the definition of $u_\tau$ (\ref{utaucond})
\begin{eqnarray*}
\D(s'_{e'},g,e')&=&\D(t_{\tau},g,e')^{-1}D(s_{e'},g,e')\\
&=&\D(s_e,g,e')^{-1}\D(u_{\tau},g,e')^{-1}\D(s_{e'},g,e')\\
&=&\D(s_e,g,e)^{-1}s_e(n_{\omega}^{e\circ g,e'\circ g})^{-1}
\D(u_{\tau},g,e)^{-1}u_{\tau}(n_{\omega}^{e\circ g,e'\circ g})^{-1}
\D(s_{e'},g,e')\\
&=&\D(s_e,g,e)^{-1}s_e(n_{\omega}^{e\circ g,e'\circ g})^{-1}
{\D(s_e,g,e)\over \D(s_{e'},g,e')} s_e(n_{\omega}^{e\circ g,e'\circ g})
\D(s_{e'},g,e')\\
&=&1
\end{eqnarray*}
as desired.

Thus $s'_{e'}\in\Gamma(W_{e'},i_*\Lambda\otimes\Gm)$ for any
${e'}:\tau\to \sigma'\in\P_{\max}$. On the other hand, for
$\tau_1\mapright{f} \tau_2\mapright{e'}\sigma'$, $s'_{{e'}\circ
f}(s'_{e'}|_{\tau_1})^{-1}=s'_f$. But then for any
$e':\tau_2\rightarrow\sigma$, $g:\omega\rightarrow\tau_1$,
\begin{eqnarray*}
\D(s'_{f},g,e'\circ f)&=&\D(s'_{e'\circ f}, g, e'\circ f) 
\D(s'_{e'}|_{\tau_1}^{-1},g,e'\circ f)\\
&=&
\D(s'_{e'\circ f}, g, e'\circ f) 
\D(s'_{e'},f\circ g,e')^{-1}\\
&=&1
\end{eqnarray*}
by the previous case. Thus $s'_f\in\Gamma(W_f,i_*\Lambda\otimes\Gm)$.
So we conclude that $s'$ is lifted gluing data.
\qed
\medskip

\noindent
(2) Given $s$, $X=X_0(B,\P,s)$, it is clear there is a unique
normalized section $f\in\Gamma(X,\shLS^+_{\pre,X})$. Indeed, using
the notation of the proof of Proposition~\ref{gluingcond}, with
$e:\omega\to \sigma$, $f_{\sigma,e}= \sum_{p\in\check\Delta
(\omega)\cap\check\Lambda} f_{\sigma,e,p} z^p,$ $f_{\sigma,e,p}$ is
determined by the normalization condition whenever $p$ is a vertex of
$\check\Delta(\omega)$. But since $\check\Delta(\omega)$ is an elementary
simplex, all integral points  of $\check\Delta(\omega)$ are vertices, so
$f_{\sigma,e}$ is completely determined, and $f$ is unique. It
remains to show $f$ defines a log Calabi-Yau structure by checking
the multiplicative condition of Theorem
\ref{finalLS}. 

Focus on one $\sigma\in\P_{\max}$ and $e:\tau\to\sigma$ with
$\dim\tau=2$. For $g:\omega\to\tau$ with $\dim\omega=1$, note that
$\D(s_{e'},g,e')=1$ for all $e':\tau\to\sigma' \in\P_{\max}$ because
the gluing data are lifted. It thus follows from (\ref{coefform})
in the proof of Proposition \ref{gluingcond} that $f_{\sigma,e\circ
g,n_{\omega}^{e\circ g,e'\circ g}}=s_e(n_{\omega}^{e\circ g, e'\circ
g})$, and in particular
\[
f_{\sigma,e\circ g}|_{V_e}=\sum_{p\in\check\Delta_i\cap \dualvs{\shQ_{\tau}}} 
s_e(p) z^p
\]
if $g\in\Omega_i$, where $\check\Delta_i$, $\Omega_i$ are associated
with $\tau$. Thus $f_{\sigma,e\circ g}|_{V_e}$ only depends on which 
$\Omega_i$ has $g$ as an element.

Now as $\tau$ is two-dimensional, we have $p\le 2$ in the definition
of simplicity for $\tau$ by Remark~\ref{simplicityfacts}, (5). If
$p=0$, then $n_{\omega}^{e\circ g,e'\circ g}=0$  and
$f_{\sigma,e\circ g}=1$ for all $g:\omega \to\tau$. Thus the
multiplicative condition is trivial. If $p=1$, then
$\dim\Delta_1=1$ or $2$, and is either a line segment of length
1 or a two-dimensional elementary simplex. Then
\[
\Omega_1=\{g_i:\omega_i\to\tau| 1\le i\le r\}
\]
with $r=2$ or $3$ in the two cases (if $r=2$, $g_1$ and $g_2$
should be thought of as opposite parallel edges of $\tau$). Note
$f_{\sigma,e\circ g}|_{V_e}=1$ for all $g\not\in\bigcup_{i=1}^p\Omega_i$.
If $\dim\Delta_1=1$, then $d_{\omega_1}=\pm d_{\omega_2}$,
and the multiplicative condition becomes
\[
(f_{\sigma,e\circ g_1}|_{V_e}\otimes d_{\omega_1})\cdot
(f_{\sigma,e\circ g_2} |_{V_e} \otimes d_{\omega_1})^{-1}=1,
\]
which holds automatically as $f_{\sigma,e\circ g}|_{V_e}$ is
independent of $g\in\Omega_i$, as remarked above. If $\dim
\Delta_1=2$, we can without loss of generality assume
$d_{\omega_1}=(1,0)$, $d_{\omega_2}=(0,1)$ and $d_{\omega_3}=(-1,-1)$
(inside the two-dimensional space $\Lambda_{\tau}$). The
multiplicative condition then becomes
\[
(f_{\sigma,e\circ g_1}|_{V_e},1)\cdot (1,f_{\sigma,e\circ g_2}|_{V_e})\cdot 
(f_{\sigma,e\circ g_3}|_{V_e}, f_{\sigma,e\circ g_3}|_{V_e})^{-1}=1,
\]
which again holds.

If $p=2$, then $\Delta_1$ and $\Delta_2$ are both line
segments, which without loss of generality we can take to be the line
segments with endpoints $(0,0)$ and $(1,0)$ for $\Delta_1$ and
$(0,0)$ and $(a,b)$ for $\Delta_2$ with $a,b\in\ZZ$,
$gcd(a,b)=1$. Again $\Omega_j=\{g_i^j:\omega_i^j\to\tau| i=1,2\}$,
$j=1,2$ (with $\omega_1^j$, $\omega_2^j$ parallel and opposite edges
of  $\tau$). Then writing $f_j=f_{\sigma,e\circ g_i^j}$ for $j=1,2$,
independent of $i$, the multiplicative condition becomes, without
loss of generality,
\[
(f_1,1)\cdot (f_2^a,f_2^b)\cdot (f_1^{-1},1)\cdot (f_2^{-a},f_2^{-b})=1,
\]
which again holds. This completes the proof. 
\qed

\begin{theorem}
\label{moduli_theorem}
Given $(B,\P)$ positive and simple, the set of positive log
Calabi-Yau spaces with dual intersection complex 
$(B,\P)$, modulo isomorphism preserving $B$, is
$H^1(\W,i_*\Lambda\otimes\Gm)$. 
\end{theorem}

\proof If $X$ is a positive log Calabi-Yau space, there exist open
gluing data $s$ such that $X=X_0(B,\P,s)$ by Theorem~\ref{logopen},
and $s$ then satisfies Condition~(LC) by Proposition~\ref{gluingcond}.
Thus by Theorem~\ref{lifted}, (1), $s$ is equivalent to lifted gluing
data.

On the other hand, any lifted gluing data trivially satisfy
Condition~(LC), hence define a positive log Calabi-Yau space by
Theorem~\ref{lifted},~(2). Now if $s,s'$ are two choices of lifted
gluing data defining log Calabi-Yau spaces $X$ and $X'$ which are
isomorphic via an isomorphism preserving $B$, then by
Proposition~\ref{equiv}, $s$ and $s'$ are  equivalent, i.e. there
exists $(t_{\tau})$ such that $s_e'=t_{\tau_1}^{-1}s_et_{\tau_2}$ for
$e:\tau_1\to\tau_2$. In particular, for $e:\tau\to\sigma
\in\P_{\max}$, $t_{\tau}=(s_e/s_e') t_{\sigma}$. But $s_e,s_e'
\in\Gamma (W_e,i_*\Lambda\otimes\Gm)$, and $t_{\sigma}
\in\Gamma(\Int(\sigma),i_*\Lambda\otimes\Gm)$. Thus
$t_{\tau}\in\Gamma(W_e\cap \Int(\sigma), i_*\Lambda\otimes\Gm)$. Since
$t_{\tau}$ is independent of $\sigma$, in fact $t_{\tau}\in
\Gamma(W_{\tau},i_*\Lambda\otimes\Gm)$. Thus $s$ and $s'$ are
cohomologous as \v Cech 1-cocycles with respect to the open covering
$\W$.

Conversely, if $s$ and $s'$ are lifted gluing data which are
cohomologous, then by Theorem~\ref{isomorphism}, $X_0(B,\P,s)$ and
$X_0(B,\P,s')$ are isomorphic. This isomorphism is induced on open
sets $V(\sigma)$ by $t_{\sigma}$, where $(t_{\tau})_{\tau\in\P}$ is
the \v Cech 0-cochain making $s$ and $s'$ cohomologous. But the
pull-back of a normalized log structure under an element of 
$\Gamma(W_{\sigma},i_*\Lambda\otimes\Gm(S))$ is still normalized, so
the pull-back of the normalized log structure on $X_0(B,\P,s')$ is a
normalized log structure on $X_0(B,\P,s)$. Hence by the uniqueness
result of Theorem \ref{lifted}, (2), $X_0(B,\P,s)$ and  $X_0(B,\P,s')$
are isomorphic as log Calabi-Yau spaces.
\qed \bigskip

We finally note that this \v Cech cohomology group computes the
cohomology of $i_*\Lambda\otimes\Gm$: 

\begin{lemma}
\label{cohom}
Suppose the discriminant locus of $B$ is straightened
(Remark~\ref{straight}). Then $\W$ is an acyclic cover for both
$i_*\Lambda$ and $i_*\Lambda\otimes \Gm(S)$. In particular,
\[
H^j(\W,i_*\Lambda)=H^j(B,i_*\Lambda)
\]
and
\[
H^j(\W,i_*\Lambda\otimes\Gm(S))=H^j(B,i_*\Lambda\otimes\Gm(S))
\]
for all $j\ge 0$.
\end{lemma}

\proof We show $H^j(W_{\tau_1\cdots\tau_p},i_*\Lambda)=0$ and
$H^j(W_{\tau_1\cdots\tau_p},i_*\Lambda\otimes\Gm(S))=1$ for $j>0$
only for $p=1$ for convenience. The more general case is the same,
keeping in mind that $W_{\tau_1\cdots\tau_p}$ may have a number of
connected components.

The first thing to observe is that $W_{\tau}$ and $\Delta\cap
W_{\tau}$ are both topological cones over the point $x=\Bar(\tau)$. In
particular, there is a fundamental system of neighbourhoods $\{W_j\}$
of the barycenter $x$ of $\tau$ and homeomorphisms $\varphi_j:W_j\to
W_{\tau}$ (just dilations) such that $\varphi_j^{-1}
(W_{\tau}\cap\Delta) =W_j\cap\Delta$ and
$\varphi_j^*\Lambda|_{W_{\tau} \setminus\Delta}=
\Lambda|_{W_j\setminus \Delta}$. Thus the restriction maps
\[
H^p(W_{\tau}, i_*\Lambda)\to H^p(W_j,i_*\Lambda)
\]
are isomorphisms for all $j$. Thus $H^p(W_{\tau},i_*\Lambda)$ is
isomorphic to the stalk of $R^p\id_*(i_*\Lambda)$ at $x$. Of course,
$R^p\id_*$ of any abelian sheaf vanishes for $p>0$ because it is
computed by applying $\id_*$ to an injective resolution. In particular
$H^p(W_{\tau},i_*\Lambda)=0$ for all $p>0$. Similarly
$H^p(W_{\tau},i_*\Lambda\otimes\Gm(S))=1$ for all $p>0$.
\qed

\begin{remark}
If $\Delta=\emptyset$, then $(B,\P)$ is automatically simple, and
the above discussion can be interpreted as follows. There is in fact
an exact sequence
\begin{equation}
\label{ses}
0\to\Lambda_{\P}\to\Lambda\to\shQ_{\P}\to 0,
\end{equation}
and hence an exact sequence
\[
H^1(B,\Lambda_{\P}\otimes\Gm(S))\mapright{f_1}
H^1(B,\Lambda\otimes\Gm(S))\mapright{f_2} 
H^1(B,\shQ_{\P}\otimes\Gm(S)).
\]
The map $f_2$ can be interpreted as taking lifted gluing
data to the corresponding closed gluing data; elements of
$H^1(B,\shQ_{\P} \otimes\Gm(S))$ represent closed gluing. The kernel
of $f_2$ can then be viewed as the space of log smooth structures on
the space $X_0(B,\P,1)$, where $1$ denotes the trivial gluing data.
This coincides with $\im f_1$. This interpretation is not so simple
for the case when $\Delta\not=\emptyset$, since then (\ref{ses})
fails to be exact on $\Delta$. However, there is always a map
\[
H^1(B,i_*\Lambda\otimes\Gm(S))\to H^1(B,\shQ_{\P}\otimes\Gm(S))
\]
taking lifted gluing data to closed gluing data, and in the
simple case the kernel again represents log Calabi-Yau structures on
the trivially glued $X_0(B,\P,1)$.
\qed
\end{remark}

\begin{example}
Following up on Example \ref{torusexamples}, (3), consider $\Gamma
\subseteq G$ defined by, for $e$ a positive integer,
\[
\Gamma=\left\{A\in\Aff(M_{\RR})\bigg| 
{\hbox{$A(m_1,m_2)=(m_1+vm_2+u+v(v-1)/2,m_2+v)$ for $(u,v)$}\atop
\hbox{in the lattice generated by $(e,0)$ and $(-e(e-1)/2,e)$}}
\right\}.
\]
The generators of $\Gamma$ are $T_1$ and $T_2$ with 
\begin{eqnarray*}
T_1(m_1,m_2)&=&(m_1+e,m_2)\\
T_2(m_1,m_2)&=&(m_1+em_2,m_2+e).
\end{eqnarray*}
These are integral not only with respect to the lattice $M\subseteq
M_{\RR}$, but also the lattice $eM\subseteq M_{\RR}$. Changing basis
of $M_{\RR}$ from the standard basis $e_1,e_2$ to $ee_1,ee_2$, we can
just as well work with the generators
\begin{eqnarray*}
T_1(m_1,m_2)&=&(m_1+1,m_2)\\
T_2(m_1,m_2)&=&(m_1+em_2,m_2+1),
\end{eqnarray*}
integral with respect to the lattice $M\subseteq M_{\RR}$.

Let $B=M_{\RR}/\Gamma$. Take as a fundamental domain for the action of
$\Gamma$ the unit square with vertices $(0,0)$, $(1,0)$, $(0,1)$ and
$(1,1)$, and obtain a polyhedral decomposition $\P$ of $B$ by
triangulating this square by adding an edge joining $(0,0)$ and
$(1,1)$. Then $\P$ contains one vertex $v$, and $\Sigma_v$ is the fan
with rays generated by $(1,0)$, $(1,1)$, $(0,1)$, $(-1,0)$, $(e-1,-1)$
and $(e,-1)$. Then $X_v$ can be viewed as a Hirzebruch surface $F_e$
blown up at two points, while $X_0(B,\P,s)$ is obtained by
identifying  ``opposite'' toric divisors of this surface. The moduli
space of log Calabi-Yau spaces with dual intersection complex $(B,\P)$
is then $H^1(B,\Lambda\otimes\Gm(k))$. Now $\Lambda$ is a local system
with monodromy $\begin{pmatrix} 1&0\\0&1\end{pmatrix}$ and 
$\begin{pmatrix} 1&e\\0&1\end{pmatrix}$ about the two generators
corresponding to $T_1$ and $T_2$, and in view of Lemma~\ref{cohom} it
is straightforward to calculate that $H^1(B,\Lambda\otimes\Gm(k))
=(\Gm(k))^2\times\mu_e$, where $\mu_e$ denotes the group of $e$th
roots of unity in $k$. This recovers a calculation of \cite{Kodaira},
and demonstrates that there may be several connected components to the
moduli space.
\end{example}

For future use, we note

\begin{corollary}
\label{Zstructure}
Let $(B,\P)$ be positive and simple, and $s$ lifted gluing
data, inducing a log Calabi-Yau structure on $X_0(B,\P,s)$. Then the
log-singular set $Z\subseteq X_0(B,\P,s)$ can be taken so that for
any $\tau\in\P$, $q_{\tau}:X_{\tau}\rightarrow X_0(B,\P,s)$,
$q_{\tau}^{-1}(Z) =Z_1\cup\cdots\cup Z_p\cup Z'$, where the
codimension of $Z'$ in $X_{\tau}$ is at least two, $Z'$ is contained
in the complement of the big torus orbit of $X_{\tau}$,
and $Z_1,\ldots,Z_p$ are irreducible and reduced Cartier divisors on
$X_{\tau}$ which are linearly independent in $\Pic(X_{\tau})$.
Furthermore the Newton polytope of $Z_i$ is $\check\Delta_i$, where
$\check\Delta_1,\ldots,\check\Delta_p$ are as in Definition~\ref{simplicity}.
\end{corollary}

\proof 
$q_{\tau}^{-1}(Z)$ consists of various irreducible components arising
as follows. The open gluing data $s$ induce a section $f\in
\Gamma(X_0(B,\P,s),\shLS^+_{\pre,X_0(B,\P,s)})$ which induces
sections of the line bundle $\shN_{\omega}$ on $X_{\omega}$ for each
$\omega\in\P$ with $\dim\omega=1$. Set $Z_{\omega}$ to be the zero
locus of this section. Then
\[
Z=\bigcup_{\omega\in\P\atop \dim\omega=1} q_{\omega}(Z_{\omega}).
\]
Thus $q_{\tau}^{-1}(Z)$
contains the codimension one components
\[
\{F_{S,\bar s}(e)^{-1}(Z_{\omega})|\hbox{$e:\omega\rightarrow\tau$ an
element of $\P_1(\tau)$}\},
\]
and a number of higher codimension components, contained in the
toric boundary of $X_{\tau}$, i.e. the complement of the big torus
orbit.

Now by definition of simplicity, $F_{S,\bar
s}(e)^{-1}(Z_{\omega})=\emptyset$ unless $e\in\bigcup_{i=1}^p
\Omega_i$. Furthermore, if $e\in\Omega_i$, then the line bundle
$F_{S,\bar s}^*(\shN_{\omega})$ has Newton  polytope given by the
elementary simplex $\check\Delta_i$. Finally, as in the  beginning of the
proof of Theorem \ref{lifted},  (2), $F_{S,\bar
s}(e)^{-1}(Z_{\omega})$ depends only on the value of $i$ for which we
have $e\in\Omega_i$. Thus we can write $q_{\tau}^{-1}(Z)
=Z_1\cup\cdots\cup Z_p\cup Z'$ where $\codim(Z')\ge 2$. Also, each
$Z_i$ is irreducible as it contains no  toric stratum of $X_{\tau}$
and comes from an elementary Newton polytope. Finally
from Remark~\ref{simplicityfacts},
$T_{\check\Delta_1}+\cdots+ T_{\check\Delta_p}$ (where $T_{\check\Delta_i}$ is the
tangent space to $\check\Delta_i$) form an interior direct sum in
$\dualvs{\shQ_{\tau,\RR}}$, from which it follows that $Z_1,\ldots,Z_p$
are linearly independent in $(\Pic X_{\tau})\otimes \RR$. 
\qed

\subsection{The logarithmic Picard group in the simple case.}

L. Illusie \cite{Illu} has suggested that if $X$ is a log scheme,
then $H^1(X,\M_X^{\gp})$ is a good candidate for a logarithmic Picard
group. Let us give some motivation for this here. Let $X$ be a normal
scheme or algebraic space, $Y\subseteq X$ a Cartier divisor. Now
suppose we have an element of $\Pic(X\setminus Y)$ which we wish to
extend to $X$. This could be done by taking a specific Weil divisor
$D$ on $X\setminus Y$ which is Cartier, representing
some element of $\Pic(X\setminus Y)$, and taking its closure $\bar
D$ in $X$. However, there are two problems. First, the linear
equivalence class of $\bar D$ might depend on which member of the
linear equivalence class of $D$ we chose. It is well-defined only up
to Weil divisors supported on $Y$. Second, $\bar D$ might be Weil,
but not Cartier. Thus we do not get a map $\Pic(X\setminus
Y)\rightarrow \Pic(X)$. 

On the other hand, if we think in terms of log geometry, and take
$j:X\setminus Y\hookrightarrow X$ the inclusion, we set $\M_{(X,Y)}
=j_*\O_{X\setminus Y}^{\times}\cap\O_X$. Then $\M_{(X,Y)}^{\gp}=
j_*\O_{X\setminus Y}^{\times}$ by Lemma~\ref{MXDgp}, and so we obtain
an injective map
\begin{equation}
\label{injmap}
H^1(X,\M_{(X,Y)}^{\gp}) \lra H^1(X\setminus
Y,\O_{X\setminus Y}^{\times}),
\end{equation}
from the Leray spectral sequence for
$j$. If in fact $R^1j_*\O_{X\setminus Y}^{\times}=0$, then we would
have $H^1(X,\M_{(X,Y)}^{\gp}) =\Pic(X\setminus Y)$. If this is the
case, it becomes clear $H^1(X,\M_{(X,Y)}^{\gp})$ is a natural notion
for the logarithmic Picard group.  This seems likely to be the case
for some natural class of pairs $Y\subseteq X$. Here is a trivial
case:

\begin{lemma}
\label{MXDzero}
Let $X$ be a toric variety over an algebraically closed field $k$,
and let $D\subseteq X$ be the union of codimension $\ge 1$ toric
strata of $X$. Let $Z\subseteq X$ be any closed subset such that $Z$
does not contain any toric stratum of $X$. Then 
\[
H^1(X\setminus
Z,\M^{\gp}_{(X,D)})=\Pic(X\setminus (D\cup Z))=0.
\]
\end{lemma}

\proof
Since $X\setminus D$ is just an algebraic torus, $\Pic(X\setminus D)
=\Pic(X\setminus (D\cup Z))=0$.  Since (\ref{injmap}) give an injection
$H^1(X\setminus Z,\M_{(X,D)}^{\gp})\rightarrow H^1(X\setminus (Z\cup
D), \O_{X\setminus D}^{\times})$, the result follows.
\qed
\medskip

The situation of interest for us is that $f:\X\rightarrow S=\Spec R$
is a toric degeneration of Calabi-Yau varieties, with singular set
$Z$. Then $f:\X\setminus Z\rightarrow S$ is log smooth, and we then
conjecturally have an isomorphism  
\[
\Pic(\X_{\bar\eta}\setminus Z)\cong
H^1(\X\setminus Z,\M_{(\X,\X_0)}^{\gp}).
\]
Since $Z$ is codimension
two, if $\X_{\bar\eta}$ is non-singular, we have $\Pic(\X_{\bar\eta}\setminus
Z)\cong \Pic(\X_{\bar\eta})$. In addition, if  $\M_{\X_0}$ is the induced
log structure on $\X_0$, there is a natural map $H^1(\X\setminus Z,
\M_{(\X,\X_0)}^{gp})\rightarrow H^1(\X_0\setminus
Z,\M_{\X_0}^{\gp})$.

This motivates the following conjecture:

\begin{conjecture}
Let $f:\X\rightarrow S$ be a toric degeneration, $X=\X_0$. Then if
$H^1(X,\O_X)=H^2(X,\O_X)=0$,
\[
\Pic(\X_{\bar\eta}\setminus Z)\cong H^1(X\setminus Z,\M_{X}^{\gp}).
\]
\end{conjecture}

If $Z=\emptyset$ and $f$ is normal crossings, then this should follow
from representability of the logarithmic Picard functor as proved by
Olsson in \cite{Olsson}. This representability probably does not
depend on $f$ being normal crossings, but the question is more subtle
when $Z\not=\emptyset$ and here we may expect to need some additional
hypotheses. 

This is meant to be motivation for computing the group
$H^1(X\setminus Z,\M_X^{\gp})$, which we will now do. This in turn
serves as motivation for the definition of the log K\"ahler moduli
space in the next section.

We will now begin to set up the technical means for a computation of
the logarithmic Picard group in our case. Let $B$ be an integral
affine manifold with singularities with a toric polyhedral
decomposition $\P$, and let $s$ be open gluing data for $(B,\P)$ over
an algebraically  closed field $k$. Let $X=X_0(B,\P,s)$, and suppose
there is a log Calabi-Yau space structure on $X$ with log-singular
set $Z$. Denote by $\bar s$ the associated closed gluing data.
Then Definition~\ref{functor} defined the gluing functor $F_{S,\bar
s}$.

For $\tau\in\P$, let $q_{\tau}:X_{\tau}\rightarrow X$ be the usual
natural map, and we will also write here $q_{\tau}$ for the
restriction $q_{\tau}: X_{\tau}\setminus q_{\tau}^{-1}(Z)\rightarrow
X\setminus Z$.  We write $\M_{\tau}:=q_{\tau}^*\M_{X\setminus Z}$,
the pull-back of the log structure (Definition~\ref{associated}).

\begin{lemma}
\label{Mresolve}
Let
\[
\C^k=\bigoplus_{\sigma_0\rightarrow\cdots\rightarrow\sigma_k
\atop \sigma_i \not=\sigma_{i+1}}
q_{\sigma_k*}(\M_{\sigma_k}^{\gp}),
\]
and define a differential
$d_{\bct}:\C^{k}\rightarrow\C^{k+1}$ by
\[
(d_{\bct}(\alpha))_{\sigma_0\rightarrow\cdots\rightarrow
\sigma_{k+1}} =\sum_{i=0}^{k}
(-1)^i\alpha_{\sigma_0\rightarrow\cdots\hat\sigma_i\rightarrow
\cdots\rightarrow\sigma_{k+1}}+(-1)^{k+1}F_{S,\bar
s}(\sigma_{k}\rightarrow
\sigma_{k+1})^*(\alpha_{\sigma_0\rightarrow\cdots
\rightarrow\sigma_{k-1}}).
\]
Then $(\C^{\bullet},d_{\bct})$ is a
resolution of the sheaf $\shM^{\gp}_{X\setminus Z}$. 
\end{lemma}

\proof Define a map $\M_{X\setminus Z}^{\gp}\rightarrow\C^0$ by
\[
\alpha\mapsto (q_{\tau}^*(\alpha))_{\tau\in\P}.
\]
We need to show
$0\rightarrow\M_{X\setminus Z}^{\gp}\rightarrow\C^{\bullet}$ is
exact. 

As for any pull-back of log-structures, $\overline\M_{\sigma}^{\gp}
=q_{\sigma}^{-1}\left(\overline{\M}_{X\setminus Z}^{\gp}\right)$
for each $\sigma\in\P$. Let 
\[
\shI^k=\bigoplus_{\sigma_0\rightarrow\cdots\rightarrow
\sigma_k}q_{\sigma_k*} \O_{X_{\sigma_k}\setminus
q_{\sigma_k}^{-1}(Z)}^{\times}
\]
and
\[
\overline{\C}^k=\bigoplus_{\sigma_0\rightarrow\cdots\rightarrow\sigma_k}
q_{\sigma_k*}\overline{\M}_{\sigma_k}^{\gp}.
\]
The differentials are defined similarly to $d_{\bct}$ (though the
differential should be written multiplicatively for
$\shI^{\bullet}$).  We then get an exact sequence of complexes
\[
\begin{matrix}
&&0&&0&&0&&\\
&&\mapdown{}&&\mapdown{}&&\mapdown{}&&\\
0&\mapright{}&\O_{X\setminus Z}^{\times}&\mapright{}&\M_{X\setminus Z}^{\gp}
&\mapright{}&\overline{\M}_{X\setminus Z}^{\gp}&\mapright{}&0\\
&&\mapdown{}&&\mapdown{}&&\mapdown{}&&\\
0&\mapright{}&\shI^{\bullet}&\mapright{}&\C^{\bullet}
&\mapright{}&\overline{\C}^{\bullet}&\mapright{}&0
\end{matrix}
\]
After taking the long exact cohomology sequence, we see the middle
complex is exact if the left and right complexes are exact. 
However, the exactness of these two complexes follows immediately from the
methods of \S A.3, Example \ref{Oresolution}, as in the proof of Proposition
\ref{cohoB}.
\qed

\begin{definition}
For any $\tau\in\P$, we denote by $D_{\tau}\subseteq X_{\tau}$
the toric boundary of $X_{\tau}$.
\end{definition}

\begin{lemma} 
\label{splitting}
For any $\tau\in\P$, there is a natural exact sequence
\[
0\lra \M_{(X_{\tau},D_{\tau})}^{\gp}
\lra \M_{\tau}^{\gp}\lra \dualvs{\Lambda_{\tau}}\oplus\ZZ
\lra 0
\]
on $X_{\tau}\setminus q_{\tau}^{-1}(Z)$. This exact sequence splits,
and the splitting is canonical if  $\dim\tau=0$. In addition,
$H^1(X_{\tau}\setminus Z,\M_{\tau}^{\gp})=0$ for all $\tau\in\P$.
Thus $H^1(X\setminus Z,\M_X^{\gp})=\HH^1(X\setminus Z,\C^{\bullet})=
H^1(\Gamma(X\setminus Z,\C^{\bullet}))$ by the hypercohomology spectral 
sequence.
\end{lemma}

\proof First consider the case of a vertex $w\in\P$. We use the
construction of $q_w$ given in Lemma \ref{q_tau}. For each $e\in
\coprod_{\sigma\in \P_{\max}} \Hom(\{w\},\sigma),$ we can view
$\tilde\sigma\subseteq \Lambda_{\RR,w}$, with $0=\tilde
w\in\tilde\sigma$ the lift of $w\in\sigma$ distinguished by $e$. We
then obtain the closed embedding $\iota_e:V_e\rightarrow V(\sigma)$,
and $q_w$ is given on the open subset $V_e=\Spec k[P_e]$ of
$X_{\omega}$ as $p\circ s_e^{-1}\circ\iota_e$. Here $P_e$ is the face
of $P_{\sigma}$ corresponding to $e$. Now the log structure on $X$
pulled back to $V(\sigma)$ is given by a section
$f\in\Gamma(V(\sigma)\setminus p^{-1}(Z),\shLS_{V(\sigma)})$. In
fact this log structure is given by a collection of charts
$P_{\sigma}\rightarrow\O_{U_i}$ for an open covering of Zariski open
subsets $\{U_i\}$ of $V(\sigma)$. The reason for this is that in the
proof of Theorem \ref{finalLS}, giving a section of $\shF$ over a
Zariski-open subset gives a chart in the Zariski topology. After
further pulling back the log structure by $s_e^{-1}$, we still obtain
charts $P_{\sigma}\rightarrow\O_{U_i}$ (replacing $U_i$ by
$s_e(U_i)$). Finally, the log structure pulled back to $V_e$ via
$\iota_e$ is given by restricting these charts to $V_e\cap U_i$,
yielding $\varphi_i: P_{\sigma}\rightarrow \O_{U_i\cap V_e}$. Such a
chart is now of the form
\[
p\mapsto \begin{cases}
0& p\not\in P_e\\
h_pz^p& p\in P_e
\end{cases}
\]
where $P_e\ni p\mapsto h_p\in \O_{U_i\cap V_e}^{\times}$ is a monoid
homomorphism. However, the log structure induced by this chart is in
fact isomorphic to the one given by the chart $\varphi_i':P_{\sigma}
\rightarrow\O_{U_i\cap V_e}$  given by
\[
p\mapsto \begin{cases}
0& p\not\in P_e\\
z^p& p\in P_e
\end{cases}
\]
Explicitly, this isomorphism is given by, for example, 
\begin{eqnarray*}
&&(P_{\sigma}\oplus\O_{V_e\cap U_i}^{\times})/
\{(p,\varphi_i(p)^{-1})|p\in\varphi_i^{-1}(\O_{U_i\cap V_e}^{\times})\}
\\
&\lra&(P_{\sigma}\oplus\O_{V_e\cap U_i}^{\times})/
\{(p,\varphi_i'(p)^{-1})|p\in\varphi_i^{-1}(\O_{U_i\cap V_e}^{\times})\}
\end{eqnarray*}
by $(p,h)\mapsto (p,h\cdot h_{\pi(p)})$, where $\pi:P^{\gp}_{\sigma}
\rightarrow P_e^{\gp}$ is any linear projection and $p\mapsto h_p$
has been extended to a homomorphism  $P_e^{\gp}\rightarrow
\O_{U_i\cap V_e}^{\times}$. 

Let $j:\eta\rightarrow X_w$ be the inclusion of the generic point
$\eta$ of $X_w$. Then there is a natural map
\[
\M_w\mapright{\beta} j_*j^*\overline{\M}_w=\NN.
\]
Let $\M'_{w}=\beta^{-1}(0)$. From the explicit description of charts
for $\M_w$ above, one sees in fact that if $\alpha:\M_w\rightarrow
\O_{X_w}$ is the structure map, then $\alpha|_{\M_w'}$ is an
inclusion of $\M_w'$ in $\O_{X_w}$, and $\alpha(\M_w')=
\M_{(X_w,D_w)}\subseteq \O_{X_w}$. Thus we obtain an exact sequence
\[
0\rightarrow \M_{(X_w,D_w)}^{\gp}\rightarrow\M_w^{\gp}\rightarrow\ZZ
\rightarrow 0.
\]
This can be split canonically via $\ZZ\ni 1\mapsto \rho\in \M_{w}^{\gp}$,
where $\rho$ is given by the morphism to $\Spec k^{\dagger}$.

If $\tau$ is arbitrary, we can choose any $e:w\rightarrow \tau$, and
then use $q_{\tau}=q_w\circ F_{S,\bar s}(e)$ to obtain charts for
$\M_{\tau}$ on open subsets of $X_{\tau}$. Now if $j:\eta\rightarrow
X_{\tau}$ is the inclusion of the generic point, then by the
construction of the dual intersection complex in \S4,
$j_*j^*\overline{\M_{\tau}}$ can be identified as the constant monoid
sheaf $\Hom(\check P_{\tau},\NN)\subseteq \dualvs{\Lambda_{\tau}}
\oplus\ZZ$. (See Definition \ref{cones}.) Thus $(j_*j^*
\overline{\M_{\tau}})^{\gp}$ can be identified with
$\dualvs{\Lambda_{\tau}}\oplus\ZZ$, though this identification is not
completely canonical, as it depends on the lifting $\tilde\tau$. The
same argument as before now shows we have an exact sequence
\[
0\mapright{} \M_{(X_{\tau},D_{\tau})}^{\gp}
\mapright{}\M_{\tau}^{\gp}\mapright{\beta} \dualvs{\Lambda_{\tau}}\oplus\ZZ
\mapright{} 0
\]
on $X_{\tau}\setminus q_{\tau}^{-1}(Z)$. 

By Lemma \ref{MXDzero}, $H^1(X_{\tau}\setminus
q_{\tau}^{-1}(Z),\M^{\gp}_{(X_{\tau}, D_{\tau})})=0$. In addition,
since $H^1(Y,\ZZ)=0$ in the \'etale topology for any normal variety
$Y$, (\cite{SGA}, IX.3.6), we see that $H^1(X_{\tau}\setminus
q_{\tau}^{-1}(Z),\M^{\gp}_{\tau})=0$. The last statement of the Lemma
then follows from Lemma \ref{Mresolve} and the hypercohomology
spectral sequence.

The only remaining statement to show is that the exact sequence of
the Lemma splits. But
$\Ext^1(\ZZ,\M^{\gp}_{(X_{\tau},D_{\tau})})=H^1(X_{\tau}\setminus
q_{\tau}^{-1}(Z),\M^{\gp}_{(X_{\tau},D_{\tau})})=0$, so the extension
class must be trivial. There is no canonical splitting, however.
\qed
\medskip

The moral of this calculation is that while there is moduli of log
structures on $V(\sigma)$, these all become a standard log structure
once restricted to irreducible components.

We now want to compute $H^1(\Gamma(X\setminus Z,\C^{\bullet}))$.

\begin{definition}
Letting
\[
Q^k=\Gamma(X\setminus Z,\C^k)/\Gamma(X\setminus Z,\shI^k),
\]
where $\shI^k=\bigoplus_{\sigma_0\rightarrow\cdots\rightarrow
\sigma_k}q_{\sigma_k*} \O_{X_{\sigma_k}\setminus
q_{\sigma_k}^{-1}(Z)}^{\times}$ is as in Lemma~\ref{Mresolve}, we
obtain an exact sequence of complexes
\[
0\rightarrow \Gamma(X\setminus Z,\shI^\bullet)
\rightarrow \Gamma(X\setminus Z,\C^{\bullet})
\rightarrow Q^{\bullet}\rightarrow 0.
\]
\qed
\end{definition}

\noindent
To think about elements of $Q^k$, we first use the sequence
\[
0\rightarrow \O_{X_{\tau}\setminus q_{\tau}^{-1}(Z)}^{\times}
\rightarrow\M_{\tau}^{\gp}\rightarrow\overline{\M}_{\tau}^{\gp}\rightarrow
0
\]
to see that
\[
\Gamma(X_{\tau}\setminus q_{\tau}^{-1}(Z),\M_{\tau}^{\gp})/
\Gamma(X_{\tau}\setminus q_{\tau}^{-1}(Z),\O_{X_{\tau}}^{\times})
\cong \ker \big(H^0(X_{\tau}\setminus q_{\tau}^{-1}(Z),
\overline{\M}_{\tau}^{\gp})
\rightarrow \Pic (X_{\tau}\setminus q_{\tau}^{-1}(Z))\big).
\]
Thus we need a way to think of elements of $H^0(X_{\tau}\setminus
q_{\tau}^{-1}(Z),\overline\M_{\tau}^{\gp})$. 

\begin{lemma}
\label{piecewise}
Let $\tau\in\P$, and choose $e:v\rightarrow\tau$. Let $\tau_e$ be
the cone of $\Sigma_v$ corresponding to $e$, with
\[
\tau_e^{-1}\Sigma_v=\{K+\RR\tau_e|\hbox{$K\in\Sigma_v$ with
$\tau_e\subseteq K$}\}.
\]
Then there is a canonical isomorphism between 
\[
H^0(X_{\tau}\setminus q_{\tau}^{-1}(Z),\overline{\M}^{\gp}_{\tau})
\]
and the group 
\begin{eqnarray*}
\PA(e)&=&\left\{\lambda:\Lambda_{\RR,v}\rightarrow\RR\bigg|
{\hbox{$\lambda$ is piecewise affine with respect to
$\tau_e^{-1}\Sigma_v$} \atop\hbox{ and takes integer values on
$\Lambda_v$}}\right\}\\
&=&\shPL_{\P}(B,\ZZ)_y
\end{eqnarray*}
for $y\in\tau$ in the interior of the edge $e$ of $\Bar(\P)$.
\end{lemma}

\proof Let $\lambda\in H^0(X_{\tau}\setminus q_{\tau}^{-1}(Z),
\overline{\M}_{\tau}^{\gp})$. Since $\overline{\M}_{\tau}^{\gp}
=q_{\tau}^{-1}(\overline{\M}_X^{\gp})$, we have stalks 
$(\overline{\M}^{\gp}_{\tau})_{\bar x}=
\overline\M^{\gp}_{X,q_{\tau}(\bar x)}$.  Now any
$f:\tau\rightarrow\sigma$ determines a  toric stratum $F_{S,\bar
s}(X_{\sigma})\subseteq X_{\tau}$ with generic point $\eta_{\sigma}$, 
and thus $\lambda$ determines a stalk 
$\lambda_{\sigma}\in\overline{\M}_{X,q_{\tau}(\bar
\eta_{\sigma})}^{\gp}$. On the other hand, the construction of the
dual intersection complex identifies $\tilde\sigma$ canonically with a
polytope in an affine hyperplane in  $\Hom( \overline{\M}_{X,
q_{\tau}(\bar \eta_{\sigma})}^{\gp},\RR)$,  and $\lambda_{\sigma}$ can
be viewed as a linear functional on this latter vector space. By
restricting $\lambda_{\sigma}$  to $\tilde\sigma$, we obtain a
canonically defined linear function $\lambda_{\sigma}:
\tilde\sigma\rightarrow\RR$. The morphism $e$ also determines a lift
$\tilde v\in\tilde\sigma$. By identifying the tangent wedge to
$\tilde\sigma$ at $\tilde v$, based at $\tilde v$, with the
corresponding cone $K$ in the fan $\Sigma_v$, we obtain just as well a
canonically defined affine linear function $\lambda_{\sigma}':
K\rightarrow \RR$ (so that $\lambda_{\sigma}'(0)=
\lambda_{\sigma}(\tilde v)$). Note that $\tau_e$ is a face of $K$.
Then $\lambda_{\sigma}'$ extends to a map $\lambda_{\sigma}':
K+\RR\tau_e\rightarrow\RR$ which is affine linear and takes integer
values on $\Lambda_v\cap (K+\RR\tau_e)$. The collection of such 
$\lambda_{\sigma}'$ defines a map $\lambda:\Lambda_{\RR,v}\rightarrow
\RR$ which is piecewise affine as desired. The fact that $\lambda$ is
continuous comes from the fact that if we have
$e_i:\omega\rightarrow\sigma_i$, $i=1,2$, the images of
$\lambda_{\sigma_1}$ and $\lambda_{\sigma_2}$ under cospecialization
maps
\[
\overline{\M}_{\tau,\bar\eta_{\sigma_1}}^{\gp}\rightarrow
\overline{\M}_{\tau,\bar\eta_{\omega}}^{\gp}\leftarrow
\overline{\M}_{\tau,\bar\eta_{\sigma_2}}^{\gp}
\]
coincide with $\lambda_{\omega}$.

Reversing this procedure, given a piecewise affine map
$\lambda:\Lambda_{\RR,v} \rightarrow\RR$, one obtains stalks
$\lambda_{\sigma}$ of $\overline{\M}_{\tau}^{\gp}$ which agree under
cospecialization, and hence give rise to a section of
$\overline{\M}_{\tau}^{\gp}$. The last equality follows from the
definition of $\shPL_{\P}(B,\ZZ)$.
\qed

\begin{remark} 
\label{PAiso}
Of course $H^0(X_{\tau}\setminus
q_{\tau}^{-1}(Z),\overline{\M}_{\tau}^{\gp})$ is independent of the
choice of $e$, so if we have $e_i:v_i \rightarrow\tau$, then the
resulting canonical isomorphism $\phi:\PA(e_1)\rightarrow \PA(e_2)$
is easily described. For any maximal cone $K_1$ of
$\tau_{e_1}^{-1}\Sigma_{v_1}$ corresponding to
$f:\tau\rightarrow\sigma\in\P_{\max}$, an element
$\lambda\in\PA(e_1)$ restricted to $K_1$ is given by an element 
$\lambda_{K_1}\in \shAff(B,\ZZ)_{v_1}$. Then $\phi(\lambda)$ restricted
to the corresponding cone $K_2$ of $\tau_{e_2}^{-1}\Sigma_{v_2}$ is
given by parallel transport of  $\lambda_{K_1}$ to
$\lambda_{K_2}\in\shAff(B,\ZZ)_{v_2}$ by first following the edge of
$\Bar(\P)$ given by $f\circ e_1$, and then following the edge $f\circ
e_2$ to arrive at $v_2$. This construction can be viewed as dual to
the construction of the map $\phi$ of Construction \ref{basicgluing}.
\qed
\end{remark}

\begin{remark}
\label{linearpart}
By Lemma \ref{splitting}, we always have a subsheaf
$\overline{\M}^{\gp}_{(X_{\tau},D_{\tau})}\subseteq
\overline{\M}_{\tau}$. We can interpret an element
$\lambda\in\Gamma(X_{\tau}\setminus q_{\tau}^{-1}(Z),
\overline{\M}_{(X_{\tau},D_{\tau})}^{\gp})$ as a piecewise linear
function on the fan $\Sigma_{\tau}$ in a slightly different way. Such
a section can be interpreted as a Cartier divisor on $X_{\tau}$
supported on $D_{\tau}$, and as such gives a well-defined piecewise
linear function $\lambda$ on the fan $\Sigma_{\tau}$ taking integer
values on $\shQ_{\tau}\subseteq \shQ_{\tau,\RR}$, by the standard
correspondence between toric Cartier divisors and piecewise linear
functions (see Remark \ref{convention}). Given
$e:v\rightarrow\sigma$, we can compose $\lambda$ with the projection
$\Lambda_{\RR,v}\rightarrow\shQ_{\tau,\RR}$, to obtain a piecewise
linear function on $\tau_e^{-1}\Sigma_v$. It is easy to check this
coincides with the piecewise affine function constructed in
Lemma~\ref{piecewise}. 

Now in general, the boundary map
\[
H^0(X_{\tau}\setminus q_{\tau}^{-1}(Z),
\overline{\M}_{(X_{\tau},D_{\tau})}^{\gp}) \rightarrow \Pic
(X_{\tau}\setminus q_{\tau}^{-1}(Z))\] induced by the standard short
exact sequence takes a Cartier divisor with support on $D_{\tau}$ to
the corresponding line bundle on $X_{\tau}\setminus
q_{\tau}^{-1}(Z)$. Thus, in particular, if $\lambda\in
H^0(X_{\tau}\setminus q_{\tau}^{-1}(Z),\overline{\M}_{\tau}^{\gp})$
is induced by a linear function on $\Sigma_{\tau}$, then the
corresponding Cartier divisor is principal and $\lambda$ lifts to an
element of $H^0(X_{\tau}\setminus q_{\tau}^{-1}(Z),
\M_{\tau}^{\gp})$.
\qed
\end{remark}

\begin{lemma} 
\label{CechAff}
Let $W_{\sigma_0\rightarrow\cdots\rightarrow\sigma_p}$ denote the 
connected component of $W_{\sigma_0\cdots\sigma_p}$ corresponding to
the $p$-dimensional simplex of $\Bar(\P)$ with edges $\sigma_i\rightarrow
\sigma_{i+1}$ (see Lemma \ref{opencovfacts}). Then there is a canonical map
\[
\Gamma(W_{\sigma_0\rightarrow\cdots\rightarrow\sigma_p},\shAff(B,\ZZ))
\rightarrow\Gamma(X_{\sigma_p}\setminus q_{\sigma_p}^{-1}(Z),
\M_{\sigma_p}^{\gp})/\Gamma(X_{\sigma_p}\setminus q_{\sigma_p}^{-1}(Z),
\O_{X_{\sigma_p}}^{\times}).
\]
This induces an injective morphism of complexes
\[
\check C^{\bullet}(\W,\shAff(B,\ZZ))\rightarrow Q^{\bullet}
\]
where the former complex denotes the \v Cech complex for $\shAff(B,\ZZ)$
with respect to the open cover $\W$.
\end{lemma}

\proof Consider first $w\in\P$ a vertex. 
By Lemma \ref{splitting}, there is a canonical isomorphism
\[
\Gamma(X_w\setminus q_w^{-1}(Z), \M_w^{\gp})
\cong\Gamma(X_w\setminus q_{w}^{-1}(Z),\M_{(X_w,D_w)}^{\gp})\oplus \ZZ\rho.
\]
Now $\M^{\gp}_{(X_w,D_w)}$ is the sheaf of rational functions on
$X_w$ with zeros and poles only along $D_w$.  Thus, since
$q_w^{-1}(Z)\subseteq D_w$, elements of $\Gamma(X_w\setminus
q_w^{-1}(Z),\M^{\gp}_{(X_w,D_w)})/\Gamma(X_w\setminus
q_w^{-1}(Z),\O_{X_w}^{\times})$ correspond to principal Cartier
divisors on $X_w$ with support on $D_w$, hence by Remark 
\ref{linearpart}, correspond to linear functions on $\Lambda_w$. In
addition, $\rho\in \Gamma(X_w\setminus q_w^{-1}(Z),\M^{\gp}_w)$
corresponds to the constant function taking the value $1$ on
$\Lambda_w$. Thus $\Gamma(X_w\setminus
q_w^{-1}(Z),\M_w^{\gp})/\Gamma(X_w\setminus
q_w^{-1}(Z),\O_{X_w}^{\times})$ is isomorphic to the subgroup of
$\PA(\id:w\rightarrow w)=\shPL_{\P}(B,\ZZ)_w$  consisting of affine
linear functions. This gives a canonical isomorphism
\[
\shAff(B,\ZZ)_w=\Gamma(W_w,\shAff(B,\ZZ))\rightarrow
\Gamma(X_w\setminus q_w^{-1}(Z),\M_w^{\gp})/
\Gamma(X_w\setminus q_w^{-1}(Z),\O_{X_w}^{\times}).
\]

We next consider the general case. Let $\sigma_0\rightarrow\cdots
\rightarrow\sigma_p$ be a chain of morphisms defining some element of
$\Bar(\P)$. Write $W=W_{\sigma_0 \rightarrow\cdots
\rightarrow\sigma_p}$. Pick a vertex $w\in\sigma_0$ and a morphism
$w\rightarrow\sigma_0$. Then $w\in\overline{W}$, and by parallel
transport from a basepoint inside $W$ to $w$ we can identify
$\Gamma(W,\shAff(B,\ZZ))$ canonically with a subgroup of
$\shAff(B,\ZZ)_w$ invariant under monodromy in $W$. We then define a
map as the composition
\begin{eqnarray*}
\Gamma(W,\shAff(B,\ZZ))\hookrightarrow \shAff(B,\ZZ)_w
=\Gamma(W_w,\shAff(B,\ZZ))&\mapright{}&
{\Gamma(X_w\setminus q_w^{-1}(Z),\M_w^{\gp})\over\Gamma(X_w
\setminus q_w^{-1}(Z),\O_{X_w}^{\times})}\\
&\mapright{F_{S,\bar s}(w\rightarrow\sigma_p)^*}&
{\Gamma(X_{\sigma_p}\setminus q_{\sigma_p}^{-1}(Z),\M_{\sigma_p}^{\gp})
\over\Gamma(X_{\sigma_p}
\setminus q_{\sigma_p}^{-1}(Z),\O_{X_{\sigma_p}}^{\times}) }.
\end{eqnarray*}
We need to show this map is independent of the choice of $w$, so
choose $e_1:w_1\rightarrow\sigma_0$ and $e_2:w_2\rightarrow\sigma_0$.
Now if we choose some $\sigma_p\rightarrow\sigma\in\P_{\max}$, we can
identify $\shAff(B,\ZZ)_{w_1}$ with $\shAff(B,\ZZ)_{w_2}$ via parallel
transport $\psi$ from $w_1$ into $\Int(\sigma)$ to $w_2$. On the
other hand, let $f_i$ be the composition $w_i\mapright{e_i}\sigma_0
\mapright{}\sigma_p$. Then we have the canonical isomorphism
$\phi:\PA(f_1) \rightarrow\PA(f_2)$ given by Remark \ref{PAiso}. This
gives a  (not necessarily commutative!) diagram
\[
\begin{matrix}
\shAff(B,\ZZ)_{w_1}&\mapright{\alpha_1}&\PA(f_1)\\
\mapdown{\psi}&&\mapdown{\phi}\\
\shAff(B,\ZZ)_{w_2}&\mapright{\alpha_2}&\PA(f_2)
\end{matrix}
\]
The second vertical map is defined on each maximal cone of
$\sigma^{-1}_{f_i} \Sigma_{w_i}$ corresponding to
$\sigma_p\rightarrow\sigma'\in\P_{\max}$  by parallel transport from
$\shAff(B,\ZZ)_{w_1}$ to $\shAff(B,\ZZ)_{w_2}$ through $\sigma'$. Thus if
the parallel transport of $\lambda\in\shAff(B,\ZZ)_{w_1}$ to
$\shAff(B,\ZZ)_{w_2}$ is independent of the choice of
$\sigma_p\rightarrow \sigma'$, then $\alpha_2(\psi(\lambda))=
\phi(\alpha_1(\lambda))$. This holds in particular if $\lambda\in
\Gamma(W,\shAff(B,\ZZ)) \subseteq\shAff(B,\ZZ)_{w_i}$. This shows
well-definedness.

Finally, it is easy to check this is sufficiently canonically defined
so as to give a map of complexes.
\qed
\medskip

\noindent
We are now in a position to carry out the calculation of
$H^1(X\setminus Z, \M_X^{\gp})$ in a way which is formally the same
as the calculation of the moduli of log Calabi-Yau spaces of \S5.1.
However, this construction will be completely dual.

\begin{theorem}
\label{logPic}
Let $(B,\P)$ be positive and simple, and $s$ lifted gluing data
determining a log Calabi-Yau space $X^{\ls}=X_0(B,\P,s)^{\ls}$ by 
Theorem~\ref{lifted}, (2). Then
there is an exact sequence
\begin{eqnarray*}
0&\rightarrow& H^0(B,k^{\times})\rightarrow H^0(X\setminus Z,\M_X^{\gp})
\rightarrow
H^0(B,\shAff(B,\ZZ))\\
&\rightarrow&H^1(B,k^{\times})\rightarrow
H^1(X\setminus Z,\M_X^{\gp})\rightarrow H^1(Q^{\bullet})
\rightarrow H^2(B,k^{\times})
\end{eqnarray*}
with an injection $H^1(B,\shAff(B,\ZZ))\hookrightarrow
H^1(Q^{\bullet})$ which is an isomorphism when tensored with  $\QQ$.
Furthermore, if $\Delta(\tau)$ of Definition~\ref{simplicity}
is a standard simplex for each $\tau\in\P$ with $\dim\tau\not=0, \dim
B$, rather than just an elementary simplex, then
$H^1(Q^{\bullet})=H^1(B,\shAff(B,\ZZ))$. 
This holds in particular if $\dim B\le 3$.
\end{theorem}

\proof
{\it Step 1}. $H^0(X_{\tau}\setminus q_{\tau}^{-1}(Z),
\O_{X_{\tau}}^{\times}) =k^{\times}$. Thus, in particular,
$\Gamma(X\setminus Z,\shI^{\bullet})$ is the simplicial cochain
complex with respect to the triangulation $\Bar(\P)$ with coefficients
in $k^{\times}$, and $H^i(\Gamma(X\setminus Z,\shI^{\bullet}))
=H^i(B,k^{\times})$.\\[1ex]
\emph{Proof.}\, If $\alpha\in H^0(X_{\tau}\setminus q_{\tau}^{-1}(Z),
\O_{X_{\tau}}^{\times})$, then the divisor of zeros and poles of
$\alpha$ is a divisor supported on $q_{\tau}^{-1}(Z)$. But then by
the linear independence statement of Corollary \ref{Zstructure},
$\alpha$ must be constant. \qed\medskip

\noindent
{\it Step 2}. Let $\check C^{\bullet} (\W,\shAff(B,\ZZ))\rightarrow
Q^{\bullet}$ be the inclusion of complexes of Lemma \ref{CechAff}.
Then $H^0(Q^{\bullet}/\check C^{\bullet}(\W,\shAff(B,\ZZ)))=0$. Thus in
particular $H^0(Q^{\bullet})=H^0(\W,\shAff(B,\ZZ))$ and
$H^1(\W,\shAff(B,\ZZ))\rightarrow H^1(Q^{\bullet})$ is injective.\\[1ex]
\emph{Proof.}\, To give an element of $H^0(Q^{\bullet}/\check
C^{\bullet}(\W, \shAff(B,\ZZ)))$ we give for each $\tau\in\P$ an
element 
\[
\lambda_{\tau}\in \Gamma(X_{\tau}\setminus
q_{\tau}^{-1}(Z),\M_{\tau}^{\gp})/ k^{\times}
\]
defined modulo $\Gamma(W_{\tau},\shAff(B,\ZZ))$, and such that  for
$e:\tau_0\rightarrow\tau_1$, $F_{S,\bar s}(e)^*(\lambda_{\tau_0})=
\lambda_{\tau_1}\mod \Gamma(W_e,\shAff(B,\ZZ))$. But for a vertex $v$,
in fact 
\[
\Gamma(X_v\setminus q_v^{-1}(Z), \M_v^{\gp})/k^{\times}=
\Gamma(W_v,\shAff(B, \ZZ)),
\]
as we saw in the proof of Lemma~\ref{CechAff},
so we can take $\lambda_v=0$, and so for
fixed $\tau$ and any $e:v\rightarrow \tau$, $\lambda_{\tau}\in
\Gamma(W_e,\shAff(B,\ZZ))= \Gamma(W_v, \shAff(B,\ZZ))=
\shAff(B,\ZZ)_v$. Now the restriction map $\Gamma (W_{\tau},
\shAff(B,\ZZ)) \rightarrow\shAff(B,\ZZ)_v$ is injective, and if
$\lambda_{\tau}$ is not in the image of this map, it is because
$\lambda_{\tau}$ is not monodromy invariant with respect to some loop
determined by $e_i:v_i\rightarrow\tau$, $i=1,2$ and $f_j:\tau
\rightarrow\sigma_j\in\P_{\max}$, $j=1,2$, passing from $v_1$ into
the interior of $\sigma_1$ to $v_2$ into the interior of $\sigma_2$
to $v_1$. Thus parallel transporting $\lambda_{\tau} \in
\shAff(B,\ZZ)_{v_1}$ from $v_1$ to $v_2$ through $\sigma_1$ and
$\sigma_2$ gives different functions. Thus $\lambda_{\tau}\in
\PA(e_1)$ is an affine linear function since $\lambda_{\tau}\in
\Gamma(W_{e_1}, \shAff(B,\ZZ))$, but as an element in $\PA(e_2)$, it is
no longer affine linear but only piecewise affine. This contradicts
$\lambda_{\tau} \in \Gamma (W_{e_2},\shAff(B,\ZZ))$. Thus
$\lambda_{\tau} \in\Gamma (W_{\tau},\shAff(B,\ZZ))$, so we can take 
$\lambda_{\tau}=0$. Thus $H^0(Q^{\bullet}/\check C^{\bullet}
(\W,\shAff(B,\ZZ)))=0$.
\qed\medskip

\noindent
{\it Step 3.} Let $\tau\in\P$ and $\lambda\in \Gamma(X_{\tau}
\setminus q_{\tau}^{-1}(Z), \overline{\M}_{\tau}^{\gp})$. For any
$e:v\rightarrow \tau$, denote by $\lambda^e$ the corresponding element
of $\PA(e)$ (Lemma~\ref{piecewise}). Suppose furthermore we have made
choices of $d_{\omega}$'s and $\check d_{\rho}$'s as in \S1.5. Then if
$f:\tau\rightarrow\rho$, with $\codim\rho=1$, we obtain
$g_{\rho}^{\pm}:\rho\rightarrow\sigma_{\rho}^{\pm}$, defining maximal
cones $K^{\pm}$ of $\tau_e^{-1}\Sigma_v$, and $\lambda^{e\pm}
\in\shAff(B,\ZZ)_v$, the restriction of $\lambda^e$ to $K^{\pm}$. Then
define $\check D(\lambda,e,f)\in\ZZ$ so that 
\[
\check D(\lambda,e,f)\check d_{\rho}=\lambda^{e-}-\lambda^{e+}.
\]
Now let $(\lambda_e)_{e:\tau_0\rightarrow\tau_1}\in Q^1$ be a
1-cocycle for the complex $Q^{\bullet}$, with $\lambda_e\in
\Gamma(X_{\tau_1} \setminus q_{\tau_1}^{-1}(Z),
\M^{\gp}_{\tau_1})/$ $k^{\times} = \Gamma(X_{\tau_1}\setminus
q_{\tau_1}^{-1}(Z), \overline\shM_{\tau_1}^\gp)$ for
$e:\tau_0\rightarrow\tau_1$. Then $(\lambda_e)$ satisfies a condition
dual to the Condition~(LC), namely: Let $\rho\in\P$ be a codimension
one cell, $v_1,v_2$ vertices, $\tau\in\P$ with a diagram 
\[
\xymatrix@C=30pt
{v_1\ar@/^/[rrd]^{e_1}\ar[rd]_{g_1}&&\\
&\tau\ar[r]^(.3){f}&\rho\\
v_2\ar@/_/[rru]_{e_2}\ar[ru]^{g_2}&&\\
}
\]
Then if $T^{\rho}_{e_1e_2}=0$, we must have
\[
\check D(\lambda_{g_1},g_1,f)=\check D(\lambda_{g_2},g_2,f).
\]
\emph{Proof.}\, First consider the case that $v_1$ and $v_2$ are
endpoints of a one-dimensional $\omega\in\P$, so that we have
\[
\xymatrix@C=30pt
{v_1\ar@/^/[rrd]^{g_1}\ar[rd]_{k_1}&&\\
&\omega\ar[r]^(.3){h}&\tau\\
v_2\ar@/_/[rru]_{g_2}\ar[ru]^{k_2}&&\\
}
\]
Then because $(\lambda_e)$ is a cocycle, we have
\[
\lambda_{g_i}=\lambda_h+F_{S,\bar s}(h)^*(\lambda_{k_i}),\quad i=1,2.
\]
We will first show that if $T^{\rho}_{e_1e_2}=0$, then  $\check
D(F_{S,\bar s}(h)^*(\lambda_{k_i}),g_i,f)=0$ for $i=1,2$. First note
that $\check D(F_{S,\bar s}(h)^*(\lambda_{k_i}),g_i,f)=\check
D(\lambda_{k_i},k_i,f\circ h)$. By Lemma \ref{splitting}, we have a
splitting $\M_{\omega}^{\gp} =\M_{(X_{\omega}, D_{\omega})}^{\gp}
\oplus (\dualvs{\Lambda_{\omega}} \oplus \ZZ\rho)$. To make this
splitting more explicit, choose an element $t_i\in \Gamma(X_{v_i},
\M^{\gp}_{(X_{v_i}, D_{v_i})})$ which vanishes exactly once along the
stratum $X_{\omega} \subseteq X_{v_i}$. This restricts to an element
$t_i\in\Gamma (X_{\omega}\setminus q_{\omega}^{-1}(Z),
\M_{\omega}^{\gp})$, and we can then write $\M_{\omega}^{\gp}=
\M_{(X_{\omega}, D_{\omega})} \oplus (\ZZ t_i\oplus\ZZ\rho)$.
Identifying $t_i$ and $\rho$ with their images in
$\overline{\M}_{\omega}^{\gp}$, we can then write $\lambda_{k_i}
=u_i+a_it_i+b_i\rho$ for $u_i\in \Gamma (X_{\omega} \setminus 
q_{\omega}^{-1}(Z), \overline{\M}^{\gp}_{ (X_{\omega}, D_{\omega})})$
and $a_i,b_i\in\ZZ$. By definition of $Q^{\bullet}$,
$\lambda_{k_i}$ is in the image of the map $\Gamma (X_{\omega}
\setminus q_{\omega}^{-1}(Z), \M_{\omega}^{\gp}) \rightarrow
\Gamma(X_{\omega} \setminus q_{\omega}^{-1}(Z),
\overline{\M}^{\gp}_{\omega})$, and since $t_i$ and $\rho$ are 
already in the image of this map, so is $u_i$, and so $u_i$ maps to
zero in $\Pic(X_{\omega} \setminus q_{\omega}^{-1}(Z))$. Then by
construction, if we view $t_i,\rho\in \PA(k_i)$, $t_i$ is a linear
function and $\rho$ is constant, so 
\[
\check D(\lambda_{k_i},k_i,f\circ h)=\check D(u_i,k_i,f\circ h).
\]
Now $u_i$ determines a Cartier divisor $C_i$ on $X_{\omega}$
supported on $D_{\omega}$, and as an element of $\PA(k_i)$, by
Remark~\ref{linearpart}, $u_i$ is just the pull-back of the
corresponding piecewise linear function on $\Sigma_{\omega}$.
Furthermore, $\check D(u_i,k_i,f\circ h)$ can be interpreted as an
element of $\Pic(X_{\rho})=\ZZ$ (see \cite{Oda}, Lemma 2.11), and as
such is the pull-back of $C_i$ to $X_{\rho}$. On the other hand,
since $u_i$ maps to zero in $\Pic(X_{\omega}\setminus
q_{\omega}^{-1}(Z))$, $C_i$ must be a principal divisor, i.e. there
exists a rational function $\alpha$ on $X_{\omega}$ with divisor of
zeros and poles  $(\alpha)=C_i+a Z_{\omega}$, where $Z_{\omega}$ is
as in the proof of Corollary \ref{Zstructure}, the single irreducible
component of $q_{\omega}^{-1}(Z)$ of codimension one. Thus $C_i+a
Z_{\omega}\sim 0$ on $X_{\omega}$, and the  restriction of $C_i+a
Z_{\omega}$ to $X_{\rho}$ is zero. But if $T^{\rho}_{e_1e_2}=0$,
$m_{e_1e_2}^{\rho}=m_{f\circ h}=0$, so $n_{f\circ h}=0$, and so it
follows from Remark \ref{psiomega} and Theorem \ref{q_tau2} that
$Z_{\omega}$ is disjoint from $X_{\rho}$. Thus $C_i$ restricts to
zero on $X_{\rho}$ and $\check D(u_i, k_i, f\circ h)=0$. Thus
$0=\check D(\lambda_{k_i},k_i,f\circ h)=\check  D(F_{S,\bar
s}(h)^*(\lambda_{k_i}), g_i,f)$ as desired.

We then conclude that $\check D(\lambda_{g_i},g_i,f)=\check
D(\lambda_h,g_i,f)$. But if $T^{\rho}_{e_1e_2}=0$, then $\check
D(\lambda_h,g_i,f)$ is also independent of $i$, so $\check
D(\lambda_{g_1},g_1,f)=\check D(\lambda_{g_2},g_2,f)$.

To conclude for the general case, we can use the above special case
if there is a sequence of one-dimensional subcells
$\omega_1,\ldots,\omega_m$ of $\tau$ with endpoints of $\omega_i$
being $w_i$ and $w_{i+1}$, with $v_1=w_1$ and $v_2=w_{m+1}$, such
that $T^{\rho}_{w_iw_{i+1}}=0$. To see this, consider the fan
$\check\Sigma_{\rho}$, with the function  $\check\psi_{\rho}$ of
Remark \ref{psiomega}. By positivity $\check \psi_{\rho}$ is convex.
However $\check\psi_{\rho}$ is determined by the same element of
$\Lambda_{\rho}$ on both the cones $\check v_1$ and $\check v_2$ of
$\check\Sigma_{\rho}$ if $m_{\rho}^{e_1e_2}=0$. By convexity of
$\check\psi_{\rho}$, one can then find a sequence of maximal cones
$\check w_1,\ldots,\check w_{m+1}$ of $\check\Sigma_{\rho}$
containing $\check\tau$ for which $\check\psi_{\rho}$ is still given
by the same element of $\Lambda_{\rho}$, and such that $\check
w_i\cap \check w_{i+1}$ is a codimension one cone, i.e.
$\check\omega_i$ for some dimension one cell $\omega_i\subseteq\tau$.
This proves the claim. \qed\medskip

\noindent
{\it Step 4.} The inclusion $H^1(\W,\shAff(B,\ZZ))\hookrightarrow 
H^1(Q^{\bullet})$ is an isomorphism when tensored with $\QQ$, and is an
isomorphism over $\ZZ$ if the additional hypotheses on the $\Delta
(\tau)$'s are satisfied.\\[1ex]
\emph{Proof.}\, Given a $1$-cocycle $(\lambda_e)$ representing an
element of  $H^1(Q^{\bullet})$, we would like to find
$(t_{\tau})_{\tau\in\P}$ with  $t_{\tau}\in\Gamma(X_{\tau}
\setminus q_{\tau}^{-1}(Z), \M_{\tau}^{\gp})/k^{\times}$  such
that for $e:\tau_1\rightarrow\tau_2$, $t_{\tau_2}+\lambda_e-F_{S,\bar
s}(e)^* (t_{\tau_1})$ is in the image of $\Gamma(W_e,\shAff(B,\ZZ))$.
This is proved in  exactly the same way as the proof of
Theorem~\ref{lifted},  (1), with everything being dual (and
additive). Since the proof would be essentially word for word the
same, we do not repeat it here, but note the one problem which
appears. The analogue of the function $u_{\tau}$ constructed in the
proof of Theorem~\ref{lifted}, (1), is first defined as a
function  $u_{\tau}:(T_{\Delta_1}\cap\Lambda_z)+\cdots+(T_{\Delta_p}
\cap\Lambda_z)\rightarrow\ZZ$. However, in general, this extends to a
function $\Lambda_z\rightarrow\ZZ$ only after multiplying $u_{\tau}$
by an integer. This is the same as multiplying the representative
$(\lambda_e)$ by an integer. Thus we only get an isomorphism over
$\QQ$. If, however, $\Delta(\tau)$ is always a standard simplex, then
$\Lambda_z$ splits as $(T_{\Delta_1}\cap\Lambda_z)\oplus \cdots\oplus
(T_{\Delta_p}\cap\Lambda_z)\oplus L$ for some free $\ZZ$-module $L$,
and then $u_{\tau}$ can be defined\qed

{\it Step 5}. To complete the proof, we now take the long exact
cohomology sequence associated to the exact sequence $0\rightarrow
\Gamma(X\setminus Z,\shI^{\bullet})\rightarrow \Gamma(X\setminus
Z,\C^{\bullet})\rightarrow Q^{\bullet}\rightarrow 0$, and use steps
1,2 and 4, as well as Lemma \ref{splitting}. Finally,
$H^p(\W,\shAff(B,\ZZ)) =H^p(B,\shAff(B,\ZZ))$, as can be shown in the
same manner as in Lemma \ref{cohom}.  Any elementary simplex of
dimension $\le 2$ is standard, hence the last statement follows if
$\dim B\le 3$.
\qed
\medskip

\begin{remark}
In case $B$ is compact every global affine function is constant by
Proposition~\ref{no global affine functions}. Hence
$H^0(B,\shAff(B,\ZZ))=\ZZ$ and $\rho\in H^0(X\setminus Z, \M_X^\gp)$ maps
to a generator. Thus the sequence in the theorem splits into two short
exact sequences, the first just giving
\[
H^0(X\setminus Z, \M_X^\gp)= k^\times\times \ZZ.
\]
Taking into account Proposition~\ref{Pic0} the second sequence reads 
\[
0\lra \Pic^\tau(X)\lra H^1(X\setminus Z,\shM_X^\gp)
\lra \ker\big(H^1(Q^\bullet)\to H^2(B,k^\times) \big) \lra 0.
\]
Thus we see that the logarithmic Picard group $H^1(X\setminus
Z,\shM_X^\gp)$ contains $\Pic^\tau(X)$ as a subgroup, and the quotient
is discrete.
\qed
\end{remark}

From the exact sequence of Proposition \ref{affexact}, we get

\begin{corollary}
If $H^i(B,\ZZ)=0$ and $H^i(B,k^{\times})=0$ for $i=1,2$ and $\Delta(\tau)$
is a standard simplex for each $\tau\in\P$, then
\[
H^1(X\setminus Z,\M_X^{\gp})\cong
H^1(B,i_*\check\Lambda).
\]
Without these hypotheses, if $H^i(B,\QQ)=0$ for $i=1,2$,
\[
H^1(X\setminus Z,\M_X^{\gp})\otimes\QQ\cong
H^1(B,i_*\check\Lambda)\otimes\QQ.
\]
\end{corollary}

\begin{example}
(1) Suppose $B=\RR/n\ZZ$, $\P$ arbitrary. Then $\shAff(B,\ZZ)$ is in
fact a local system of rank $2$ on $B$ with monodromy
$\begin{pmatrix} 1&n\\ 0&1\end{pmatrix}$. Also, $X=X_0(B,\P,s)$ is a
cycle of rational curves, and it is not difficult to see that
$H^0(X,\M_X^{\gp})=k^{\times} \times \ZZ$, generated by constant
functions and $\rho$. Then the exact sequence of Theorem \ref{logPic}
gives
\[
0\rightarrow k^{\times}\rightarrow k^{\times}\times\ZZ\rightarrow
\ZZ\rightarrow k^{\times}\rightarrow H^1(X,\M_X^{\gp})\rightarrow
\ZZ\rightarrow 0,
\]
so $H^1(X,\M_X^{\gp})=k^{\times}\times\ZZ$. See \cite{Kajiwara} and
\cite{Olsson} for
discussions of $H^1(X,\M_X^{\gp})$ for $X$ a curve.

(2) Suppose $B$ is the base of an elliptic fibration on a K3
surface (Example \ref{ellipticK3}) which is integral and has  24
singular points, and a polyhedral decomposition $\P$ with $B,\P$
simple (or alternatively, construct $B$ as the dual intersection
complex of a suitable toric degeneration of K3 surfaces). Set
$X=X_0(B,\P,s)$. One can show that $H^1(B,i_*\check\Lambda)=
\ZZ^{20}$, and $H^1(B,\shAff(B,\ZZ))$ is the kernel of a boundary map
$H^1(B,i_*\check\Lambda)\rightarrow H^2(B,\ZZ)=\ZZ$, and hence is
either $19$ or $20$ dimensional. Finally we have
\[
0\rightarrow H^1(X\setminus Z,\M_X^{\gp})\rightarrow H^1(B,\shAff(B,\ZZ))
\rightarrow k^{\times}
\]
from Theorem \ref{logPic}. The last map will depend on the choice
of open gluing data, giving a range of possible values
for $H^1(X\setminus Z,\M_X^{\gp})$.
\end{example}

\subsection{Mirror symmetry and conclusions}

We first define

\begin{definition}
Let $X$ be a toric log Calabi-Yau space over an algebraically  closed
field $k$, with dual intersection complex $B$. Then the {\it log
K\"ahler moduli space} of $X$ is defined to be
$H^1(B,i_*\check\Lambda\otimes\Gm(k))$.
\end{definition}

\begin{theorem}
Let $(X,\shL)$ be a polarized log Calabi-Yau space with degeneration
data $(B,\P,\varphi)$. Suppose $B$ is positive and simple.
If $(\check B,\check\P,\check\varphi)$ is the discrete Legendre
transform of $(B,\P,\varphi)$, then the set of log Calabi-Yau
spaces with dual intersection complex $(\check B,\check \P)$ modulo isomorphism
preserving $B$ is naturally isomorphic to the log K\"ahler moduli space
of $X$.
\end{theorem}

\proof This follows from Theorem \ref{moduli_theorem} and the fact that $\check
\Lambda^{B_0}\cong\Lambda^{\check B_0}$ by Proposition \ref{dlt}. 
\qed

Of course, at this point the reader may object that we have simply
defined the log K\"ahler moduli space to make the above theorem work.
Thus a bit more justification seems to be in order. The calculations
of \S 5.2 say that $H^1(B,i_*\check\Lambda)$ has something to do with
the logarithmic Picard group, which in turn should have something to
do with the smoothing of $X$, should such exist. This is perhaps
weaker evidence than one would like. However, in \cite{sequel}, we
will define logarithmic Hodge groups $H^{p,q}(X)$ for $X$ a log
Calabi-Yau space, and compute these groups. In certain cases, these
coincide with $H^q(B,i_*\bigwedge^p\check\Lambda\otimes_{\ZZ} k)$.
Furthermore, we will show in these cases these groups have the same
dimension as the dimension of the Hodge groups of a smoothing. Thus
$H^{1,1}(X)$ is directly related to the K\"ahler moduli space. In
addition, $H^{n-1,1}(X)=H^1(B,i_*\bigwedge^{n-1}\check\Lambda
\otimes_{\ZZ} k)=H^1(B,i_*\Lambda\otimes_{\ZZ} k)$ if the holonomy of
$B$ is contained in $\ZZ^n\rtimes \SL_n(\ZZ)$ rather than in
$\ZZ^n\rtimes \GL_n(\ZZ)$. (This property does not follow from
orientability of $B$). As $H^{n-1,1}(X)$ is then the same dimension as
the moduli space of a smoothing of $X$, we see that the moduli space
of log Calabi-Yau spaces $H^1(B, i_*\Lambda\otimes\Gm(k))$ has the
correct dimension.

Another motivation comes from \cite{SlagII}. There it was
suggested that if $f:X\rightarrow B$ was a simple torus fibration, it
was natural to consider the $B$-field living in the group
$H^1(B,R^1f_*(\RR/\ZZ))$. This was a somewhat different proposal than
is usually considered in the  physics literature, in that it depends
not just on $X$ but on the fibration. However, we feel it is the
natural group, and some moral  justification for this was given in
\cite{SlagII}. (The reader should of course keep in mind we do not
yet have a satisfactory mathematical  definition of mirror symmetry,
and therefore it does not make sense to prove such assertions).

Thus this suggests that the complexified K\"ahler moduli space lies
naturally in $H^1(B,$ $R^1f_*(\CC/\ZZ))=H^1(B,(R^1f_*\ZZ)\otimes
\CC^{\times})$. On the other hand, if $f$ was a special Lagrangian
fibration, the base $B$ would carry two affine structures (with
singularities) related by a Legendre transform \cite{Hit}. In the one
coming from the complex structure, we naturally have  $R^1f_*\ZZ\cong
i_*\check\Lambda$. Thus in this context, the  definition of the log
K\"ahler moduli space as $H^1(B,i_*\check\Lambda\otimes \Gm(k))$
appears to be the correct one. The fact that this is indeed isomorphic
to the moduli of log Calabi-Yau spaces for a ``mirror'' supports
this. 

At any rate, we have now achieved mirror symmetry for polarized log
Calabi-Yau spaces:

\begin{definition}
Let $(X,\shL)$ and $(\check X,\check\shL)$ be two polarized log
Calabi-Yau spaces with degeneration data $(B,\P,\varphi)$ and
$(\check B,\check\P,\check\varphi)$ respectively. If these are
related by a discrete Legendre transform, we say $(X,\shL)$ and
$(\check X,\check\shL)$ are a mirror pair of log Calabi-Yau spaces.
In particular, if $B$ (and hence $\check B$) are positive and simple,
then there is a canonical isomorphism between the log K\"ahler moduli
space of $X$ and the  set of log Calabi-Yau spaces with dual
intersection complex $(B,\P)$ modulo isomorphism preserving $B$.
\end{definition}

\begin{example}
Let $\Xi$ be a reflexive polytope, and consider a toric degeneration
$\X\rightarrow\shS$ constructed as in Example~\ref{degenexamples},
and let $\check\X\rightarrow\shS$ be constructed similarly using the
dual polytope $\Xi^*$. Let $X=\X_0^{\dagger}$ and $\check X=\check
\X_0^{\dagger}$. Polarize $X$ and $\check X$ by the anti-canonical
divisors $\shL=\O_{\PP_{\Xi^*}}(1)|_X$ and $\check\shL=
\O_{\PP_{\Xi}}(1)|_{\check X}$ respectively. Then by Example
\ref{Batyrevdual}, it follows that $(X,\shL)$ and $(\check X,\check\shL)$
are a mirror pair. Note in this case we don't expect $(B,\P)$ to be
simple, so we don't expect the isomorphisms between K\"ahler and
complex moduli. This is because the generic fibres $\X_{\eta}$ and
$\check\X_{\eta}$ are not in general smooth, or MPCP resolutions.
Different polarizations and degenerations are required to achieve
this, see \cite{GBB}.
\end{example}

Our goal, of course, is not just a version of mirror symmetry for log
Calabi-Yau spaces. We still need to address a number of questions to
make use of the program we have begun here. Principal among these
are:

(1) We are really interested in toric degenerations of Calabi-Yau
varieties $f:\X\rightarrow\shS$. To use the mirror symmetry above to
study mirror symmetry of non-singular Calabi-Yau varieties, we need
to  compare the log invariants of $X=\X_0^{\dagger}$ with the usual
invariants of the generic fibre $\X_{\bar\eta}$. For example, we
expect in the simple case that $h^{1,n-1}(\X_{\bar\eta})$ coincides
with the dimension of the moduli of log Calabi-Yau spaces with the
same dual intersection complex as $X$. This should follow from the
right sort of base-change results. Many of the necessary  ideas for
proving such results are already present in the literature for log
smooth morphisms; we however have the singular set $Z$ to worry
about, which complicates things.

(2) To turn these ideas into a genuine mirror symmetry construction,
we need to study the log deformation theory of log Calabi-Yau spaces,
\`a la Kato \cite{F.Kato} and Kawamata-Namikawa  \cite{Kawamata;
Namikawa 1994}. Again, we hope in the simple case that a log
Calabi-Yau space is always the central fibre of a toric degeneration.
The Bogomolov-Tian-Todorov type results of \cite{Kawamata; Namikawa
1994}  give some moral support that such should be true. However, the
presence of the singular set means we cannot apply known results off
the shelf.

(3) We believe this approach will provide a very general mirror
symmetry construction. It incorporates Batyrev-Borisov duality, as
is explained in \cite{GBB}. It also includes mirror symmetry for
abelian varieties, but it would be nice to know it gives some
genuinely new constructions.

(4) This approach should be connected more explicitly with the SYZ
approach. It is clear that our approach is a discretization of the
SYZ approach, in that the Legendre transform is replaced by the
discrete Legendre transform. However, we would like to see torus
fibrations when we work over $\CC$. Some results in this direction
were stated in \cite{Announce}. The proof of these results will
appear in \cite{tori}, along with more precise results in the simple
case.

(5) Much more generally, the log philosophy says that things you want
to compute on smoothings can be computed using logarithmic analogues
on the degenerations. For example, we hope it should be possible to
define log Gromov-Witten invariants \cite{SiebertGW} which will
coincide with those of a smoothing. The results of \cite{Mikhalkin}
and \cite{Nishinou; Siebert} suggest that this should be possible by
counting certain piecewise straight graphs on $B$ (``tropical
curves''). In addition, the philosophy presented in this paper is that
invariants of toric log Calabi-Yau spaces can be computed doing
calculations on $B$. This is the general philosophy we hope to push
further in order to understand mirror symmetry.

The sequel to this paper \cite{sequel} will address questions (1) and
(2). However, our work in this direction is still too preliminary to
give specific statements here.

So there is still much work to be done to make the logarithmic
approach to mirror symmetry realise its full potential. We hope the
results presented here already provide ample evidence that
logarithmic degeneration data capture the essential features of mirror
symmetry.
\bigskip
\appendix
\section{Simplicial and Polyhedral Complexes}\noindent
The purpose of this appendix is to collect some facts on complexes
arising from a polyhedral decomposition. We first treat the case of
abelian groups and then sheafify.


\subsection{Barycentric complexes}
Let $\P$ be the face poset of a $d$-dimensional polytope $\Xi$.
Assume that $(M_\tau)_{\tau\in\P}$ is a category of abelian groups
indexed by $\P$, that is, a functor $\Cat(\P)\to
\underline{\operatorname{Ab}}$. For $\tau,\sigma\in \P$,
$\tau\subset\sigma$ denote by $\varphi_{\sigma\tau}: M_\tau\to
M_\sigma$ the homomorphism of abelian groups thus given. Then
$\varphi_{\tau\tau}=\id$ and $\varphi_{\sigma\tau}
\circ\varphi_{\tau\omega}= \varphi_{\sigma\omega}$ whenever
$\omega\subset \tau\subset \sigma$. For notational convenience we
write $\varphi_\sigma(f)$ instead of $\varphi_{\sigma\tau}(f)$
whenever $f\in M_\tau$ and $\tau\subset\sigma$. Denote by
$\P^{[k]}\subset \P$ the subset of $k$-dimensional faces. 

The \emph{barycentric cochain complex} $(C_\bct^\bullet,
d_\bct^\bullet)$ associated to $(M_\sigma)$ is the complex of
abelian groups $C^k_\bct= \bigoplus_{\sigma_0\subsetneq
\sigma_1\subsetneq \ldots\subsetneq \sigma_k} M_{\sigma_k}$ of
\emph{simplicial cochains} with differentials
\[
\big(d^k_\bct(f_{\sigma_0\sigma_1\ldots\sigma_k})\big)_
{\sigma_0\sigma_1\ldots\sigma_{k+1}}=
\sum_{i=0}^{k+1} (-1)^i \varphi_{\sigma_{k+1}}\big(f_{\sigma_0\ldots
\widehat\sigma_i\ldots\sigma_{k+1}}\big).
\]
Here and in the following entries with a hat are to be omitted. This
is indeed a complex as one easily checks. There is generally no
reason for this complex to be acyclic, but it will be once
$(M_\sigma)$ has the following extension property. For $\P'\subset
\P$ let us call a collection $f_\tau\in M_\tau$, $\tau\in\P'$,
\emph{compatible} if $\varphi_\sigma(f_\tau) =\varphi_\sigma
(f_{\tau'})$ for any $\tau,\tau'\in \P'$, $\sigma\in\P$, $\tau,\tau'
\subset \sigma$. We consider the following 
condition.\\[2ex] $(*)$
\hfill\begin{minipage}{14cm}
\it Any compatible collection $(f_\sigma)_{\sigma\in \P'}$ indexed by
any $\P'\subset\P$ extends to a compatible collection
$(g_\sigma)_{\sigma\in\P}$, that is, $g_\sigma\in M_\sigma$ and
$g_\sigma=f_\sigma$ for $\sigma\in \P'$.
\end{minipage}\hfill
\smallskip

\begin{proposition}\label{barycentric complex}
If Condition~$(*)$ holds then the barycentric cochain complex
associated to $(M_\sigma)$ is acyclic.
\end{proposition}

\proof
We wish to write a simplicial cocycle $(f_{\sigma_0
\ldots\sigma_k})_{(\sigma_0,\ldots,\sigma_k)}$ as the coboundary of a
simplicial $(k-1)$-cochain $(g_{\sigma_0
\ldots\sigma_{k-1}})_{(\sigma_0,\ldots,\sigma_{k-1})}$. We construct
$g_{\sigma_0\cdots\sigma_{k-1}}$ by descending induction on
$m=\dim\sigma_{k-1}= d+1, \ldots,0$. The induction hypothesis is that
\begin{eqnarray}\label{f=dg}
f_{\sigma_0\ldots\sigma_k}= \sum_{i=0}^k (-1)^i
\varphi_{\sigma_k}(g_{\sigma_0\ldots \widehat\sigma_i\ldots \sigma_k}),
\end{eqnarray}
whenever $\dim\sigma_{k-1}\ge m$. The induction starts at $m=d+1$.
Condition~(\ref{f=dg}) is empty at this stage because $\P$ is the
face poset of a $d$-dimensional polytope. For the induction step
consider $\sigma_0\subsetneq \ldots \subsetneq \sigma_{k-1}$ with
$\dim \sigma_{k-1}=m-1$. We want to find $g_{\sigma_0
\ldots\sigma_{k-1}}$ such that for any $\sigma_k$ containing
$\sigma_{k-1}$
\[
(-1)^k\varphi_{\sigma_k}\big(g_{\sigma_0\ldots\sigma_{k-1}}\big)
= f_{\sigma_0\ldots\sigma_k} -\sum_{i=0}^{k-1}
(-1)^i g_{\sigma_0\ldots\widehat\sigma_i \ldots\sigma_k}.
\]
This is the required equation for $f_{\sigma_0 \ldots\sigma_k}$, and
all terms on the right-hand side are known inductively. Now view the
right-hand side as a collection of elements indexed by
$\P'=\{\sigma_k \,|\, \sigma_{k-1}\subset \sigma_k\}$. By
Condition~($*$) $g_{\sigma_0\ldots \sigma_{k-1}}$ then exists if
this collection is compatible, that is, if for any
$\sigma_{k+1}\supset \sigma_{k-1}$ the expression
\begin{eqnarray}\label{claim in step2}
\varphi_{\sigma_{k+1}}\Big( f_{\sigma_0\ldots\sigma_k}-
\sum_{i=0}^{k-1} (-1)^i g_{\sigma_0\ldots\widehat\sigma_i
\ldots\sigma_k}
\Big)
\end{eqnarray}
does not depend on $\sigma_k$ for $\sigma_{k-1}\subset
\sigma_k\subset\ \sigma_{k+1}$. For $i\le k$ the induction hypothesis
implies
\[
f_{\sigma_0\ldots\widehat\sigma_i\ldots\sigma_{k+1}}
=\sum_{j=0}^{i-1} (-1)^j g_{\sigma_0\ldots \widehat\sigma_j
\ldots\widehat \sigma_i\ldots\sigma_{k+1}} -\sum_{j=i+1}^{k+1} (-1)^j
\varphi_{\sigma_{k+1}} (g_{\sigma_0\ldots \widehat\sigma_i
\ldots\widehat \sigma_j\ldots\sigma_{k+1}}).
\]
Plugging this into the cocycle condition
\[
\varphi_{\sigma_{k+1}}\big(f_{\sigma_0\ldots\sigma_k}\big)
=(-1)^k\sum_{i=0}^k (-1)^i
f_{\sigma_0\ldots\widehat\sigma_i\ldots\sigma_{k+1}},
\]
the first term of (\ref{claim in step2}) gives $f_{\sigma_0\ldots
\sigma_{k-1}\sigma_{k+1}}$ ($i=k$) plus a sum over
$g_{\sigma_0\ldots\widehat\sigma_i\ldots\widehat\sigma_j\ldots
\sigma_{k+1}}$. For $0\le i<j<k$ the coefficient of
$g_{\sigma_0\ldots \widehat\sigma_i\ldots\widehat\sigma_j\ldots
\sigma_{k+1}}$ is $(-1)^k$ times $(-1)^i(-(-1)^j)+(-1)^j(-1)^i=0$.
Contributions involving $\varphi_{\sigma_{k+1}} (g_{\sigma_0\ldots
\widehat\sigma_i\ldots\sigma_k})$ come from the second term in
(\ref{claim in step2}) and from $j=k+1$; they cancel as well. Thus
(\ref{claim in step2}) equals
\[
f_{\sigma_0\ldots\sigma_{k-1} \sigma_{k+1}}
+(-1)^k\sum_{i=0}^{k-1} (-1)^i
(-1)^k(-g_{\sigma_0\ldots\widehat\sigma_i\ldots\widehat\sigma_k
\sigma_{k+1}}).
\]
This shows the claimed independence of (\ref{claim in step2}), and
hence the existence of $g_{\sigma_0\ldots\sigma_{k-1}}$.
\qed
\medskip


\subsection{Polyhedral complexes}
The second type of complex that we deal with is non-simplicial. We
start with the same kind of system $(M_\sigma)_{\sigma\in \P}$ of
abelian groups as in the previous subsection. The following
construction depends on a choice of \emph{orientations} of the faces
of $\Xi$. Comparison with the standard orientation of the boundary of
a face gives a sign $\sgn(\tau,\omega)$ whenever $\tau\subset\omega$
and $\dim(\omega)=\dim(\tau)+1$. Similarly, if $\sigma_0\subset\ldots
\subset\sigma_k$ and $\dim\sigma_i=i$ comparison of the natural
orientation of the simplex $(\sigma_0,\ldots,\sigma_k)$ of the
barycentric subdivision of $\sigma_k$ with the chosen orientation
defines a sign $\sgn((\sigma_0,\ldots,\sigma_k),\sigma_k)$.

We can then define the \emph{polyhedral cochain complex}
$(C_\phd^\bullet, d_\phd^\bullet)$ associated to $(M_\sigma)$ as
the complex of abelian groups $C^k_\phd= \bigoplus_{\tau\in\P^{[k]}}
M_\tau$ of \emph{polyhedral cochains} with differentials
\[
d^k_\phd:(f_\tau)_{\tau\in\P^{[k]}}\longmapsto
\left(\sum_{\tau\subset\omega,\tau\in\P^{[k]}}
\sgn(\tau,\omega)\varphi_{\omega\tau}
(f_\tau)\right)_{\omega\in \P^{[k+1]}}.
\]
To check that the composition of differentials vanishes, recall that if
$\omega$ is a $(k+1)$-dimensional polytope and $\sigma\subset \omega$
has codimension $2$ then there are exactly two $k$-dimen\-sional facets
$\tau^\pm\subset \omega$ with $\sigma\subset\tau^\pm\subset\omega$.
Because the orientations on $\sigma$ induced by the orientations of
$\tau^+,\tau^-$ as facets of $\omega$ differ, it holds
\begin{eqnarray}\label{signs work out}
\sgn(\sigma,\tau^+)\sgn(\tau^+,\omega)=
-\sgn(\sigma,\tau^-)\sgn(\tau^-,\omega).
\end{eqnarray}
Therefore,
\begin{eqnarray*}
\lefteqn{\big(d^k_\phd\circ d^{k-1}_\phd(f_\sigma)\big)_\omega
\ =\ \sum_{\sigma\subset\tau
\subset\omega} \sgn(\sigma,\tau)\sgn(\tau,\omega)
(\varphi_{\omega\tau} \circ\varphi_{\tau\sigma}) (f_\sigma)}
\hspace{1cm}\\
&=& \sum_{\scriptstyle\sigma\subset\omega \atop
\scriptstyle\dim\sigma=k-1}
\big( \sgn(\sigma,\tau^+)\sgn(\tau^+,\omega)
+\sgn(\sigma,\tau^-)\sgn(\tau^-,\omega)\big)
\varphi_{\omega\sigma}(f_\sigma)\ =\ 0.
\end{eqnarray*}

\begin{proposition}\label{polyhedral complex}
If Condition~$(*)$ holds then the polyhedral cochain complex
associated to $(M_\sigma)$ is acyclic.
\end{proposition}

\proof
Let $(f_\tau)_{\tau\in\P^{[k]}}\in C_\phd^k$ fulfill the cocycle
condition that for any $\omega\in\P^{[k+1]}$
\[
\sum_{\tau\subset\omega} \sgn(\tau,\omega)\cdot
\varphi_{\omega\tau}(f_\tau)=0.
\]
For any $(k-1)$-cell $\sigma$ we have to find $g_\sigma\in M_\sigma$
with the property that for any $\tau\in\P^{[k]}$
\[
f_\tau=\sum_{\sigma\subset\tau} \sgn(\sigma,\tau)\cdot
\varphi_{\tau\sigma}(g_\sigma).
\]
The construction of $g_\sigma$ is in three steps. Steps~1 and 2
produce a simplicial cocycle out of $g_\sigma$, which by
Proposition~\ref{barycentric complex} may then be written as a
simplicial coboundary; Step~3 translates back to the polyhedral
setting.
\smallskip

\noindent
\emph{Step~1.}
The aim of the first two steps is to transform, in a sense,
$(f_\tau)$ to a cocyle in $C^k_\bct$. For each $\tau\in \P^{[k]}$
choose a maximal chain $\tau_0\subsetneq \tau_1\subsetneq\ldots
\subsetneq\tau_k=\tau$ and write $\underline{\tau}=
(\tau_0,\ldots,\tau_k)$. For $\sigma_0\subset
\sigma_1\subset\ldots \subset\sigma_k$ with $\dim\sigma_i=i$ we
then put
\begin{eqnarray}\label{f_{sigma_0,..,sigma_k}}
f_{\sigma_0\ldots\sigma_k}=\left\{\begin{array}{ll}
(-1)^k\sgn(\underline\tau,\tau)\, f_\tau&,\ 
(\sigma_0,\ldots,\sigma_k)=(\tau_0,\ldots,\tau_k)\\
0&,\ \text{otherwise}.
\end{array}\right.
\end{eqnarray}
It remains to define $f_{\sigma_0\ldots\sigma_{k-1}\omega}$
for $\dim\omega>k$.

In this step we deal with the case $\dim\omega= k+1$. Since $\omega$
is a polytope it can be built up by consecutively attaching
$(k+1)$-simplices from its barycentric subdivision, keeping it a cell
in each step. This means that there is an order among the
$(k+1)$-simplices $(\sigma_0,\ldots,\sigma_k,\sigma_{k+1}=\omega)$
with the property that any two consecutive elements intersect along a
$k$-simplex, and all but the last cell are attached along a union of
at most $k$ of the $k+1$ facets not lying in $\partial \omega$. Let
$\Xi_\kappa$, $\kappa=0,\ldots,K$ be the cell obtained after the
$\kappa$-th attachment. Assume $f_{\sigma_0\ldots\sigma_k}$ is
constructed for all $k$-simplices contained in $\Xi_\kappa$ for $0\le
\kappa\le \kappa_0<K$, and $\Xi_{\kappa_0}$ equals $\Xi_{\kappa_0-1}$
with the $(k+1)$-simplex $(\sigma_0,\ldots,\sigma_{k+1}=\omega)$
attached. Then $f_{\sigma_0\ldots \widehat\sigma_i
\ldots\sigma_{k+1}}$ is not defined for a non-empty subset
$I\subset \{0,\ldots,k\}$. Put $f_{\sigma_0\ldots \widehat\sigma_i
\ldots\sigma_{k+1}}=0$ for $i<\max(I)$, $i\in I$, and
\[
f_{\sigma_0\ldots\widehat\sigma_{\max(I)}\ldots\sigma_{k+1}}
=(-1)^{\max(I)+1} \sum_{i\in\{0,\ldots,\widehat{\max(I)},\ldots,
k+1\}} (-1)^i\varphi_{\sigma_{k+1}}
(f_{\sigma_0\ldots \widehat\sigma_i \ldots\sigma_{k+1}}),
\]
to fulfill the cocycle condition on $(\sigma_0,\ldots,\sigma_{k+1})$.
By induction from $\kappa_0=0$ to $\kappa_0=K-1$ this determines
$f_{\sigma_0\ldots\sigma_k}$ satisfying the cocycle condition on
every $(k+1)$-simplex of the barycentric subdivision of $\omega$
except possibly on the one $(\sigma_0,\ldots,\sigma_{k+1})$ attached
in the $K$-th step. Adding the already known cocycle conditions on
all other $(k+1)$-simplices contained in $\omega$ gives the
expression
\[
\sum_{\sigma_0\subsetneq\ldots\subsetneq\sigma_{k+1}}
\sgn\big((\sigma_0,\ldots,\sigma_{k+1}),\omega\big)\,
\sum_{i=0}^{k+1} (-1)^i \varphi_\omega(f_{\sigma_0\ldots
\widehat\sigma_i \ldots\sigma_{k+1}}).
\]
For $i\le k$ the coefficient of $\varphi_\omega(f_{\sigma_0\ldots
\hat \sigma_i\ldots\sigma_{k+1}})$ is zero, because this
term occurs in the sum exactly twice with opposite signs,
cf.~(\ref{signs work out}). The remaining expression is
\[
\sum_{\sigma_0\subsetneq\ldots\subsetneq\sigma_{k+1}}
\sgn\big((\sigma_0,\ldots,\sigma_{k+1}),\omega\big)\,
(-1)^{k+1} \varphi_\omega(f_{\sigma_0\ldots\sigma_k})
=\ -\sum_{\tau\in\P^{[k]},\tau\subset\omega} 
\sgn\big( (\underline{\tau},\omega),\omega
\big)\, \sgn (\underline{\tau},\tau) \varphi_\omega(f_\tau).
\]
As $\sgn\big( (\underline{\tau},\omega),\omega
\big)\, \sgn (\underline{\tau},\tau) = \sgn(\tau,\omega)$ this is
nothing but $d_\phd(f_\tau)_\tau$, which vanishes by assumption.
\smallskip

\noindent
\emph{Step~2.}
Now let us assume $f_{\sigma_0\ldots\sigma_k}$ has been constructed
inductively for $\dim\sigma_k\le m$ fulfilling the cocycle condition
on any $m$-simplex. Let $\omega$ be of dimension $m+1>k+1$. Then for
$\sigma_0\subsetneq \ldots\subsetneq \sigma_{k+1}\subset
\partial\omega$ the cocycle condition on $(\sigma_0,\ldots,
\sigma_{k+1})$ holds by induction assumption. Applying
$\varphi_\omega$ gives $\sum_{i=0}^{k+1} (-1)^i
\varphi_\omega(f_{\sigma_0\ldots\hat \sigma_i\ldots \sigma_{k+1}})
=0$. Thus $\big(\varphi_\omega(f_{\sigma_0\ldots\hat \sigma_i\ldots
\sigma_{k+1}}) \big)_{\sigma_{k+1}\subset\omega}$ is a simplicial
$k$-cocycle on $\partial\omega\cong S^m$ with values in $M_\omega$.
But $m>k$ and hence $H^k(\partial\omega, M_\omega)=0$ by the
universal coefficient theorem. Therefore there exists
$(h_{\sigma_0\ldots \sigma_{k-1}})_{\sigma_{k-1}\subset
\partial\omega}$ with $(-1)^{k+1} \varphi_{\omega}
(f_{\sigma_0\ldots\sigma_k})= -\sum_{i=0}^k (-1)^i h_{\sigma_0\ldots
\hat\sigma_i\ldots \sigma_k}$. Putting $f_{\sigma_0\ldots
\sigma_{k-1}\omega}:= h_{\sigma_0\ldots \sigma_{k-1}}$ the cocycle
condition continues to hold on any $(\sigma_0,\ldots,\sigma_{k+1})$
with $\sigma_{k+1}\subset \omega$.
\smallskip

\noindent
\emph{Step~3.}
By Proposition~\ref{barycentric complex} there exists a simplicial
$(k-1)$-cocycle $g_{\sigma_0\ldots\sigma_{k-1}}$ with
$f_{\sigma_0\ldots \sigma_k}= \sum_{i=0}^k (-1)^i
g_{\sigma_0\ldots\hat \sigma_i\ldots \sigma_k}$. Now we want to
revert to the polyhedral setting. For $\sigma\in\P^{[k-1]}$ put
\[
g_\sigma=\sum_{\sigma_0\subsetneq\ldots\subsetneq \sigma_{k-1}=
\sigma} \sgn((\sigma_0,\ldots,\sigma_{k-1}),\sigma)
g_{\sigma_0\ldots\sigma_{k-1}}.
\]
Then for any $\tau\in\P^{[k]}$
\begin{eqnarray*}
\sum_{\sigma\subset\tau,\sigma\in \P^{[k-1]}} \sgn(\sigma,\tau)\,
\varphi_{\tau}(g_\sigma)&=& \sum_{\sigma_0\subsetneq\ldots\subsetneq \sigma_k=\tau}
\sgn((\sigma_0,\ldots,\sigma_{k-1}),\sigma_{k-1})\,\sgn(\sigma_{k-1},\tau)\,
\varphi_{\sigma_k}(g_{\sigma_0\ldots\sigma_{k-1}})\\
&=&\sum_{\sigma_0\subsetneq\ldots\subsetneq \sigma_k=\tau}
\sgn((\sigma_0,\ldots,\sigma_k),\tau)
(-1)^k\sum_{i=0}^k (-1)^i
\varphi_{\sigma_k}(g_{\sigma_0\ldots\widehat\sigma_i\ldots \sigma_k})\\
&=&\sum_{\sigma_0\subsetneq\ldots\subsetneq \sigma_k=\tau}
\sgn((\sigma_0,\ldots,\sigma_k),\tau)
(-1)^k f_{\sigma_0\ldots\sigma_k}\\
&=&\sgn(\underline\tau,\tau)\, (-1)^k\,
f_{\tau_0\ldots\tau_k}=f_\tau.
\end{eqnarray*}
To verify the second equality consider the coefficient of
$g_{\sigma_0\ldots\widehat\sigma_i\ldots \sigma_k}$. For $0\le i\le
k-1$ there are exactly two $i$-cells $\sigma_i^\pm$ between
$\sigma_{i-1}$ (empty for $i=0$) and $\sigma_{i+1}$. Now
$g_{\sigma_0\ldots \widehat{\sigma_i^\pm} \ldots\sigma_k}$ contribute
with opposite signs, just as in (\ref{signs work out}). This exhibits
$(f_\tau)_\tau$ as polyhedral coboundary as desired.
\qed
\medskip


\subsection{Complexes of Sheaves}

Let $B$ be an integral affine manifold with singularities, and let
$\P$ be a toric polyhedral decomposition on $B$. For every $\sigma\in
\P$ put $\P_\sigma=\coprod_{\tau\in\P} \Hom(\tau,\sigma)$, partially
ordered by $(e:\tau\to\sigma) \le (e':\tau' \to\sigma)$ iff there
exists $f:\tau\to\tau'$ with $e= e'\circ f$. There is a natural order
preserving map from $\P_\sigma$ to the face poset of $\tilde\sigma$. 
If we orient every $\sigma\in\P$, then as in the previous subsection
comparison of orientations gives signs $\sgn(\sigma_{k-1},\sigma_k)$
and $\sgn((\sigma_0 \to\cdots\to\sigma_k), \sigma_k)$ if $\dim
\sigma_i=i$. 

Now let $s$ be open gluing data for $(B,\P)$ over $S$, and let
$X=X_0(B,\P,s)$.  Suppose that for each $\sigma\in\P$, we are given a
sheaf $\shF_{\sigma}$ on $X_{\sigma}$, (here we write $X_{\sigma}$
for $X_{\sigma}\times S$) and for each $e\in \Hom(\tau,\sigma)$ a map
\[
\varphi_e:F_{S,\bar s}(e)^{-1}(\shF_{\tau})\rightarrow\shF_{\sigma},
\]
compatible with compositions. Analogous to the constructions for
abelian groups define the \emph{barycentric} and \emph{polyhedral
cochain complexes} associated to $(\shF_\sigma)_{\sigma\in\P}$ by
$\C^k_\bct= \bigoplus_{\sigma_0\to \cdots\to \sigma_k}
q_{\sigma_k*}\shF_{\sigma_k}$ and $\C^k_\phd=
\bigoplus_{\sigma\in\P^{[k]}} q_{\sigma*}\shF_\sigma$ respectively,
with differentials
\begin{eqnarray*}
&&d^k_\bct(f_{\sigma_0\to \cdots\to\sigma_k})\ =\ 
\left(\sum_{i=0}^{k} (-1)^i  f_{\sigma_0\to\cdots \to
\widehat\sigma_i\to\cdots \sigma_{k+1}}+(-1)^{k+1}
\varphi_{\sigma_k\rightarrow\sigma_{k+1}}
(f_{\sigma_0\to\cdots\to\sigma_k})\right)_{\sigma_0\to
\cdots\to\sigma_{k+1}},\\
&&d^k_\phd(f_\sigma)_{\sigma\in\P^{[k]}}\ =\ 
\left(\sum_{\tau\to\omega} \sgn(\tau,\omega)\varphi_{\tau\to\omega}
(f_\tau)\right)_{\omega\in \P^{[k+1]}}.
\end{eqnarray*}

We will use the previous two subsections to give a criterion for
exactness of $\C^{\bullet}_{\bct}$ and $\C^{\bullet}_{\phd}$. To check
exactness, it is enough to do so on stalks. If $\bar x\rightarrow X$
is a geometric point, let $\sigma$ be the largest $\sigma\in\P$ such
that  $\bar x$ is in  the image under $q_{\sigma}$ of a geometric
point of $X_{\sigma}$. Then $q_{\tau*} (\shF_{\tau})_{\bar x}=0$
unless $\tau\subseteq\sigma$. Furthermore, the set of geometric points
of $X_{\tau}$ mapping to  $\bar x$ is in one-to-one correspondence
with the set  $\Hom(\tau,\sigma)$. If $\bar y\to X_{\tau}$ with
$q_{\tau}(\bar y)=\bar x$ corresponds to $e:\tau\rightarrow\sigma$,
let $M_e=\shF_{\tau,\bar y}$, and if $(e:\tau\rightarrow\sigma)\le
(e':\tau'\rightarrow\sigma)$ with $e=e'\circ f$, we have a map
$\varphi_f:M_e\rightarrow M_{e'}$ induced by $\varphi_f: F_{S,\bar
s}(f)^{-1}(\shF_{\tau})\rightarrow \shF_{\tau'}$. This gives a
barycentric or polyhedral cochain complex  $(C^{\bullet}_\bct,
d_{\bct})$ or $(C^{\bullet}_\phd,d_{\phd})$ associated to the system
$(M_e)_{e\in \P_{\sigma}}$. On the other hand, since
$q_{\tau}:X_{\tau}\rightarrow X$ is always a finite map, it follows
from \cite{Milne}, II Cor. 3.5, that the stalk of
$q_{\tau*}\shF_{\tau}$ at $\bar x$ is $\bigoplus_{\bar
y}\shF_{\tau,\bar y}$, where the sum is over all geometric points
$\bar y$ of $X_{\tau}$ mapping to $\bar x$. From this one easily sees
that the stalks of the complexes $(\C^{\bullet}_\bct,d_{\bct})$ and
$(\C^{\bullet}_\phd,d_{\phd})$ at $\bar x$ coincide with the
barycentric or polyhedral cochain complexes associated to $(M_e)$.
Thus if for every point $\bar x$ of $X$ the system
$(M_e)_{e\in\sigma}$ satisfies Condition~$(*)$  (keeping in mind
$\sigma$ depends on  $\bar x$), it follows that the complexes of
sheaves $\C^{\bullet}_{\bct}$ and $\C^{\bullet}_\phd$ are exact.

\begin{example}
\label{Oresolution}
Take $\shF_{\tau}=\O_{X_{\tau}}$, and let $\varphi_e=F_{S,\bar
s}(e)^*$, the pull-back of functions via $F_{S,\bar s}(e)$. It is
easy to check Condition~$(*)$ holds. Indeed, to each
$e\in\P_{\sigma}$, we obtain a closed subscheme of $\Spec\O_{X,\bar
x}$ given by the image of $\Spec\O_{X_{\tau},\bar
y}\mapright{q_{\tau}}\Spec \O_{X,\bar x}$ where $\bar y$ is the point
of $X_{\tau}$ mapping to $\bar x$ corresponding to $e$. Then
Condition~$(*)$ is implied by the fact that given functions on some
collection of these closed subschemes which agree on intersections,
we can glue them to obtain a function on the union of these closed
subschemes, and then extend to obtain an element of $\O_{X,\bar x}$.
\end{example}
\newpage
\section*{Index of notations}
\small
\begin{tabbing}
\centerline{\bf 1.1}\\
$M$, $N=M^*$\hspace{2.2cm}
	\=free abelian group of rank $n$ and its dual:
	$N=\Hom(M,\ZZ)$\\
$M_\RR$, $N_\RR$
	\>$M\otimes_\ZZ \RR$, $N\otimes_\ZZ \RR$\\
$\Aff(M)$, $\Aff(M_\RR)$
	\>group of (integral) affine transformations of $M$\\
$B$
	\>$n$-dimensional topological manifold\\
$\shT_B$
	\>tangent bundle of $B$\\
$\pi:\tilde B\to B$
	\>universal covering\\
$\delta:\tilde B\to M_\RR$
	\>developing map\\
$\rho:\pi_1(B)\to \Aff(M_\RR)$
	\>holonomy representation\\
$\Psi_\gamma$
	\>deck transformation given by $\gamma\in\pi_1(B)$\\
$\Lin$
	\>linear part of an affine transformation\\
$\Trans$
	\>translational part of an affine transformation\\
$c_\rho$
	\>radiance obstruction of holonomy representation $\rho$\\
$\nabla$
	\>the flat connection on $\shT_B$ induced from the affine structure\\
$\Lambda_\RR,\Lambda$
	\>local system of flat (integral) vector fields\\
$\check\Lambda_\RR, \check\Lambda$
	\>duals to $\Lambda_\RR$ and $\Lambda$\\
$\shAff(B,\RR)$, $\shAff(B,\ZZ)$
	\>sheaf of (integral) affine functions on $B$\\
$g$, $K$
	\>Hessian metric on $B$ and a local potential function\\
$\check B$, $\check K$, $\check\delta$
	\>Legendre dual data\\[2ex]
\centerline{\bf 1.2}\\
$\Delta\subset B$
	\>discriminant locus\\
$i:B_0\to B$
	\>inclusion of complement of $\Delta$\\
$\Xi$
	\>reflexive lattice polytope\\[2ex]
\centerline{\bf 1.3}\\
$\P$
	\>polyhedral decomposition\\
$v$, $w$
	\>vertices of $\P$\\
$\sigma$, $\tau$, $\omega$
	\>cells of $\P$\\
$\Int\sigma$
	\>relative interior of $\sigma\in\P$\\
$\exp_v:R_v\to B$
	\>chart at vertex $v\in\P$; $R_v\subset \Lambda_{\RR,v}$\\
$\P_v$
	\>polyhedral decomposition of $R_v$\\
$\tilde\sigma\to\sigma$
	\>restriction of $\exp_v$ to lift $\tilde\sigma\in\P_v$ of
	$\sigma$\\
$S_\sigma$
	\>local submersion of neighbourhood of $\Int(\sigma)$ contracting
	$\sigma$\\
$\P_{\max}$
	\>set of maximal cells of $\P$\\
$\Bar(\P)$
	\>(first) barycentric decomposition of $B$ wrt.\ $\P$\\
$W_\tau$
	\>open star of barycenter of $\tau\in\P$ wrt.\ $\Bar(\P)$\\
$\W$
	\>open covering $\{W_\tau\,|\,\tau\in\P\}$ of $B$\\
$\Delta'\subset B$
	\>the maximal codimension two subcomplex of $\Bar(\P)$ disjoint\\
	\>from vertices and from $\Int(\sigma)$ for all
	$\sigma\in\P_{\max}$; contains $\Delta$.\\
$\Sigma_v$
	\>fan in $\Lambda_{\RR,v}$ induced from $\P$; generalizes to
	$\Sigma_\tau$ for $\tau\in\P$.\\
$\Lambda_\P$, $\Lambda_{\P,\RR}$
	\>subsheaf of $\Lambda$ ($\Lambda_\RR$) of vectors tangent to
	cells of $\P$\\
$\Lambda_\sigma$, $\Lambda_{\sigma,\RR}$
	\>stalk of $\Lambda_\P$ ($\Lambda_{\P,\RR}$) at any
	$x\in\Int(\tau)\setminus\Delta$\\
$\shQ_\P$, $\shQ_{\P,\RR}$
	\>subsheaf of $i_*(\Lambda/\Lambda_\P)$
	of locally flat sections and its realification\\
$\shQ_\sigma$, $\shQ_{\sigma,\RR}$
	\>stalk of $\shQ_\P$ ($\shQ_{\P,\RR}$) at any
	$x\in\Int(\tau)\setminus\Delta$\\
$\Sigma_\tau$
	\>fan in $\shQ_{\tau,\RR}$ defined by tangent wedges $K_\sigma$
	to $\sigma\in\P$ containing $\tau$\\
$\Cat(\P)$
	\>category with elements $\P$ and morphisms one-cells of
	$\Bar{\P}$\\
$\Hom(\tau,\sigma)$
    \>morphisms in $\Cat(\P)$; any $e:\tau\to\sigma$ corresponds
	uniquely to a face of $\tilde\sigma$\\
$\Sigma(\tau)$
	\>quotient fan of fan $\Sigma$ by $\tau\in\Sigma$\\
$\tau^{-1}\Sigma$
	\>localization of fan $\Sigma$: $\tau^{-1}\Sigma=\{\sigma+\RR\tau\,|\,
	\sigma\in\Sigma, \tau\subset\sigma\}$\\
$\check\Sigma_\tau$
	\>normal fan of $\tilde\tau$; fan in
	$\dualvs{\Lambda_{\tau,\RR}}$\\
$\shAff(B,\RR)$, $\shAff(B,\ZZ)$
	\>sheaf of continuous functions that are (integral) affine on $B_0$\\
$\shPL_{\P}(B,\RR)$, $\shPL_{\P}(B,\ZZ)$
	\>sheaf of (integral) piecewise linear functions on $B$\\
$\varphi$, $\varphi_\sigma$
	\>a piecewise linear function and the induced piecewise linear\\
	\>function on $\shQ_{\sigma,\RR}$\\
$\shMPL_{\P,\RR}$, $\shMPL_\P$
	\>sheaf of multivalued (integral) piecewise linear functions\\[2ex]
\centerline{\bf 1.4}\\
$(\check B,\check\P,\check\varphi)$
	\>Legendre transform of $(B,\P,\varphi)$\\
$\check\sigma$
	\>Legendre dual to $\sigma\in\P$; element of $\check\P$\\
$\Xi^*$
	\>dual of reflexive polytope $\Xi$\\[2ex]
\centerline{\bf 1.5}\\
$d_\omega$
	\>generator of $\Lambda_\omega$, $\dim\omega=1$\\
$\check d_\rho$
	\>generator of $\dualvs{\shQ_\rho}$, $\dim\rho=n-1$\\
$v_\omega^-$, $v_\omega^+$
	\>vertices of $\omega$ ordered by $d_\omega$\\
$e_\omega^-$, $e_\omega^+$   
	\>morphisms $v_\omega^\pm\to \omega$\\
$\sigma_\rho^-$, $\sigma_\rho^+$
	\>$n$-dimensional cells adjacent to $\rho$; ordered by $\check
	d_\rho$\\
$g_\rho^-$, $g_\rho^+$
	\>morphisms $\rho\to\sigma_\rho^\pm$\\
$\gamma^{e_1 e_2}_{f_1 f_2}$
	\>primitive loop distinguished by $f_i:v_i\to \tau$ and
	$g_i:\tau\to\sigma_i$\\
$T^{e_1 e_2}_{f_1 f_2}$
	\>monodromy along $\gamma^{e_1 e_2}_{f_1 f_2}$\\
$T^{e_1 e_2}_\omega$
	\>$T^{e_1 e_2}_{e_\omega^+ e_\omega^-}$;
	monodromy around one-dimensional $\omega\in\P$ with\\
	\>respect to adjacent $\sigma_1$, $\sigma_2\in\P_{\max}$\\
$n^{e_1 e_2}_\omega$
	\>element of $\dualvs{\shQ_\omega}$ determining $T^{e_1
	e_2}_\omega$\\
$T^\rho_{f_1 f_2}$
	\>$T^{g^+_\rho g^-_\rho}_{f_1 f_2}$; monodromy between
	$\sigma_\rho^+$, $\sigma_\rho^-$ through $v_1$, $v_2$\\
$m^\rho_{f_1 f_2}$
	\>element of $\Lambda_\rho$ determining $T^\rho_{f_1 f_2}$\\
$T_e$, $e:\omega\to\rho$
	\>distinguished primitive monodromy for $\dim\omega=1$,
	$\dim\rho=n-1$\\
$n_e$, $m_e$
	\>multiples of $d_\omega$ and $\check d_\rho$, respectively,
	determining $T_e$\\
$\psi_\omega$
	\>PL-function on $\Sigma_\omega$ determined by $n_\omega^{e_1
	e_2}$, $e_i:\omega\to\sigma_i$\\
$\check\psi_\rho$
	\>PL-function on $\check\Sigma_\rho$ determined by $m^\rho_{f_1
	f_2}$, $f_i:v_i\to\rho$\\
$\check\Delta(\omega)\subset \shQ_{\omega,\RR}$
	\>in the positive case, Newton polytope of $\psi_\omega$\\
$\Delta(\rho)\subset \Lambda_{\rho,\RR}$
	\>in the positive case, Newton polytope of $\check\psi_\rho$\\
$\check\Delta_e(\tau)\subset \check\Delta(\omega)$
	\>face determined by $e:\omega\to\tau$\\
$\Delta_e(\tau)\subset \Delta(\rho)$
	\>face determined by $e:\tau\to\rho$\\
$\P_1(\tau)$
	\>sets of morphisms $\{e:\omega\to \tau\,|\,\dim\omega=1\}$\\
$\P_{n-1}(\tau)$
	\>sets of morphisms $\{f:\tau\to\rho\,|\,\dim\rho=n-1\}$\\
$\check\Delta(\tau)$
	\>in the simple case, simplex built from
	$\check\Delta_e(\tau)$, $e\in\P_1(\tau)$\\
$\Delta(\tau)$
	\>in the simple case, simplex built from
	$\Delta_f(\tau)$, $f\in\P_{n-1}(\tau)$\\[2ex]
\centerline{\bf 2.1}\\
$A$
	\>base ring\\
$C(\sigma)$
	\>cone over $\sigma\in\P$ in $\Lambda_{\sigma,\RR}\oplus\RR$\\
$\check P_\sigma$
	\>$C(\sigma)\cap(\Lambda_\sigma	\oplus\ZZ)$\\
$R_\sigma$
	\>$\ZZ[\check P_\sigma]$\\
$\check X_\sigma$
	\>$\Proj R_\sigma$\\
$\check F$
	\>gluing functor for cone picture $\P\ni \sigma\mapsto R_\sigma$\\
$\check s=(\check s_e)_{e:\to\sigma}$
	\>gluing data for cone picture: $\check s_e\in \Hom(\check
	P_\tau,\Gm(A))$\\
$\check F_{A,\check s}$
	\>gluing functor for cone picture, twisted by $\check s$\\
$\check X_0(B,\P,\check s)$
	\>$\Proj (\displaystyle\lim_{\leftarrow} \check F_{A,\check s})$; glued
	scheme in cone picture\\[2ex]
\centerline{\bf 2.2}\\
$X(\Sigma)$
	\>toric variety for fan $\Sigma$\\
$X_\sigma$
	\>$X(\Sigma_\sigma)$\\
$\dual{\sigma}$
	\>dual of cone $\sigma$\\
$F$
	\>gluing functor for fan picture $\P\ni\sigma\mapsto X_\sigma$\\
$W_{\tau_1\ldots\tau_p}$
	\>$W_{\tau_1}\cap\ldots\cap W_{\tau_p}$\\
$W_e$
	\>connected component of $W_{\tau\sigma}$ selected by
	$e:\tau\to\sigma$\\
$S$
	\>base scheme\\
$s=(s_e)_{e:\tau\to\sigma}$
	\>closed gluing data: \v Cech $1$-cocycle for $\shQ_\P
	\otimes\Gm(S)$ wrt.\ $\W$;\\
	\>also: open gluing data, see below\\
$F_{S,s}$
	\>gluing functor for fan picture, twisted by $s$\\
$P_\sigma$
	\>$\dual{C(\sigma)}\cap (\check\Lambda_y\oplus\ZZ)$\\
$\rho_\sigma$
	\>distinguished element of $P_\sigma$\\
$\partial P$
	\>$(\partial \sigma\cap N)\cup\{\infty\}$\\
$V(\sigma)\subset U(\sigma)$
	\>$\Spec \ZZ[\partial P_\sigma]\subset \Spec \ZZ[P_\sigma]$\\
$\mathfrak R$
	\>étale equivalence relation on $\coprod_{\sigma\in\P_{\max}}
	V(\sigma)\times S$\\
$\check\omega$
	\>for $\omega\subset\sigma\in\P_{\max}$ corresponding face of
	$\dual{C(\sigma)}$;\\   
	\>also, corresponding element of $\check \Sigma_\sigma$\\
$U(\tau)$
	\>$\Spec\ZZ[\dual{C(\tau)}\cap(N\oplus\ZZ)]$\\
$V(\tau)$
	\>toric boundary of $U(\tau)$\\
$\Phi_{\sigma_1\sigma_2}$, $\Phi_{e_1 e_2}$
	\>gluing isomorphisms\\
$K_e$
	\>for $e:\tau\to\sigma$ the corresponding cone in $\Sigma_\tau$\\
$V_e$
	\>$\Spec\ZZ[\dual{K_e}\cap \dualvs{\shQ_\tau}]$ ($e:\tau\to\sigma$)\\
$\check{PM}(\tau)$
	\>piecewise multiplicative functions on $\tilde \tau$:
	$\Gamma(\tilde \tau, \pi^*(\shQ_\P\otimes\Gm(S)))$\\
$s=(s_e)_{e:\tau\to\sigma}$
	\>open gluing data: $s_e\in \check{PM}(\tau)$ for
	$e:\tau\to\sigma$;\\
	\>also: closed gluing data, see above\\
$Z^1(\P,\shQ_{\P}\otimes\Gm(S))$
	\>group of open gluing data\\
$B^1(\P,\shQ_{\P}\otimes\Gm(S))$
	\>trivial open gluing data\\
$\bar s=(\bar s_e)_{e:\tau\to\sigma}$
	\>closed gluing data associated to open gluing data $s$\\
$\Phi_{e_1 e_2}(s)$
	\>gluing isomorphism twisted by open gluing data $s$\\
$\mathfrak R_{e_1 e_2}$
	\>graph of $\Phi_{e_1e_2}$; component of $\mathfrak R$\\
$X_0(B,\P,s)$
	\>glued algebraic space in fan picture: $\coprod_{\sigma\in \P_{\max}}
	(V(\sigma)\times S)/\foR$\\
$p_\sigma$
	\>for $\sigma\in\P_{\max}$ the morphism $V(\sigma\times S) \to
	X_0(B,\P,s)$\\
$q_\tau$
	\>$X_\tau\times S\to X_0(B,\P,s)$\\
$\Pic^\tau(X)$
	\>group of isomorphism classes of numerically trivial line
	bundles on $X$\\
$\omega_X$
	\>dualizing sheaf of $X$\\[2ex]
\centerline{\bf 3.1}\\
$\alpha_X:\shM_X\to \O_X$
	\>log structure on $X$\\
$q:\shM_X\to \overline{\shM}_X$
	\>morphism to ``ghost sheaf'' $\overline\M_X=\M_X/\O_X$\\
$\M_{(X,D)}$
	\>log structure on $X$ defined by divisor $D\subset X$\\
$\Spec k^\dagger$
	\>standard log point $(\Spec k,\NN\oplus \Gm(k))$\\[2ex]
\centerline{\bf 3.2}\\
$P$
	\>toric monoid; here: $P=\overline\shM_{X,\bar x}$ for a fine log
	structure $\shM_X$\\
$\shR$
	\>local relation sheaf of $\overline\shM_X$\\
$\rho$, $\bar\rho$
	\>section of $\shM_X$ with $\alpha(\rho)=0$ and its image in
	$\overline\shM_X$\\
$X^g$
	\>space $X$ together with ``ghost structure''\\
$\shLS_{X^g}$
	\>sheaf of germs of log-smooth structures of given ghost type;\\
	\>subsheaf of $\shExt^1(\overline\shM_X^\gp/\bar\rho,
	\O_X^\times)$\\[2ex]
\centerline{\bf 3.3}\\
$V$
	\>$V(\sigma)\times\Gm^r= \Spec k[\partial P\oplus\ZZ^r]$\\
$\shF$
	\>a subsheaf of $\shHom(\shR,\O_V^\times)$ mapping onto $\shLS_V$\\
$\epsilon_\tau$
	\>for a $2$-face $\tau$ a cyclic choice of orientations of its
	edges\\
$s^v$
	\>germ of $s\in \Gamma(\tilde\tau,\pi^* (\shQ_{\P}\otimes\Gm(k)))$
	at a vertex $v$.\\
$s^\omega$
	\>$s^{v_\omega^+}/s^{v_\omega^-}$ if $\dim\omega=1$ and
	$v_\omega^\pm\to\omega$ are the two vertices\\
$D(s,\omega)$
	\>defined by $s^\omega=d_\omega\otimes D(s,\omega)$;
	global version: $D(s,e,f)$\\
$s^h$
	\>generalization of $s^v$ for the case with self-intersections;
	$h:v\to\tau$\\
$D(s,e,f)$, $D(s,\omega,\sigma)$
	\>for $e:\omega\to \tau$, $f:\tau\to\sigma$, $\dim\omega=1$,
	$\dim\sigma=n$, the element of $k^\times$\\
	\> defined by $s^{h^-}/s^{h^+}= d_\omega\otimes D(s,e,f)$\\
$\shN_\omega$
	\>line bundle on $X_\omega$ corresponding to the PL-function
	$\psi_\omega$ on $\Sigma_\omega$;\\
	\>$\shLS_{X_0(B,\P,s)}\subset \bigoplus_{\dim\omega=1}
	{q_{\omega}}_* \shN_\omega$\\[2ex]
\centerline{\bf 4.1}\\
$\shS$
	\>base of toric degeneration: spectrum of discrete valuation
	ring $R$\\
$\X$
	\>total space of toric degeneration\\
$\X_\eta$
	\>generic fiber of toric degeneration\\
$\X_0$
	\>central fiber of toric degeneration\\
$\nu:\tilde\X_0\to\X_0$
	\>normalization map\\
$Z\subset\X$, $Z\subset X$
	\>log-singular set; closed set of relative codimension two;\\
$X^\dagger=(X,\shM_X)$
	\>log Calabi-Yau space\\
$\nu:\tilde X\to X$
	\>normalization\\
$\Cat(X)$
	\>category of toric strata of $X$\\
$\Strata(X)$
	\>set of toric strata of $X$\\
$\LPoly$
	\>category of lattice polytopes\\
$LP:\Cat(X)\to\LPoly$
	\>functor used to build dual intersection complex\\[2ex]
\centerline{\bf 4.2}\\
$\shL$
	\>(ample) line bundle on $X^\dagger$\\
$\varphi_\shL$
	\>piecewise linear function on $B$ associated to $\shL$\\
$(B,\P,\varphi_\shL)$
	\>degeneration data for $(X^\ls,\shL)$ building on
	dual intersection complex of $X^\ls$\\
$(\check B,\check \P,\check \varphi_\shL)$
	\>dual degeneration data for $(X^\ls,\shL)$
	building on intersection complex of $X^\ls$\\[2ex]
\centerline{\bf 4.3}\\
$\shLS^+_{\pre,X}$
	\>$\bigoplus_{\dim\omega=1} {q_\omega}_*\shN_\omega$ where
	$\shN_\omega\in \Pic(X_\omega)$ is as in \S3.3\\
\\[2ex]
\centerline{\bf 5.1}\\
$Z_i\subset q_\tau^{-1}(Z)$
	\>irreducible components of codimension $1$ of $Z$ on $X_\tau$\\
$Z'\subset q_\tau^{-1}(Z)$
	\>higher codimension components of $Z$ on $X_\tau$\\[2ex]
\centerline{\bf 5.2}\\
$H^1(X\setminus Z, \M_X^\gp)$
	\>logarithmic Picard group of $(X,\M_X)$\\
$\M_\tau$, $\overline\M_\tau$
	\>$q_\tau^*\M_{X\setminus Z}$, $q_\tau^{-1}\overline
	\M_{X\setminus Z}$\\
$\C^\bullet$
	\>(barycentric) resolution of $\M_{X\setminus Z}^\gp$;
	$\C^k=\bigoplus_{\sigma_0\rightarrow\cdots\rightarrow\sigma_k}
	q_{\sigma_k*}(\M_{\sigma_k}^{\gp})$\\
$\shI^\bullet$
	\>(barycentric) resolution of $\O^\times_{X\setminus Z}$;
	$\shI^k=\bigoplus_{\sigma_0\rightarrow\cdots\rightarrow
	\sigma_k}q_{\sigma_k*} \O_{X_{\sigma_k}\setminus
	q_{\sigma_k}^{-1}(Z)}^{\times}$\\
$\overline{\C}^\bullet$
	\>(barycentric) resolution of $\overline\M_{X\setminus Z}^\gp$;
	$\overline{\C}^k=\bigoplus_{\sigma_0\rightarrow\cdots
	\rightarrow\sigma_k}
	q_{\sigma_k*}\overline{\M}_{\sigma_k}^{\gp}$\\
$d_\bct$
	\>(barycentric) differential in $\C^\bullet$, $\shI^\bullet$ etc.\\
$D_\tau\subset X_\tau$
	\>toric boundary of $X_\tau$\\
$Q^\bullet$
	\>$\Gamma(X\setminus Z,\C^\bullet)/
	\Gamma(X\setminus Z,\shI^\bullet)$\\
$\tau_e$
	\>cone in $\Sigma_v$ distinguished by $e:v\to \tau$\\
$\PA(e)$
	\>for $e:v\to \tau$, functions $\Lambda_{\RR,v}
	\to \RR$ piecewise affine  wrt.\ $\tau_e^{-1}\Sigma_v$\\
$W_{\sigma_0\rightarrow\cdots\rightarrow\sigma_p}$
	\>connected component of $W_{\sigma_0\cdots\sigma_p}$
	distinguished by $\sigma_0\rightarrow\cdots\rightarrow\sigma_p$\\
$\lambda$
	\>element of $\Gamma(X_\tau\setminus q_\tau^{-1}(Z),
	\overline\M_\tau^\gp)$\\
$\lambda^e$
	\>for $e:v\to\tau$ image of $\lambda$ under isomorphism\\
	\>$\Gamma(X_\tau\setminus q_\tau^{-1}(Z),\overline\M_\tau^\gp)
	\to \PA(e)$\\
$\check D(\lambda,e,f)$
	\>for $e:v\to\tau$ integer
	defined by $\check D(\lambda,e,f)\check d_{\rho}=
	\lambda^{e-}-\lambda^{e+}$,\\
	\>where $\lambda^{e\pm}\in
	\Aff(B,\ZZ)_v$ are the restrictions of $\lambda^e$ to the two\\
	\>maximal cones distinguished by
	$f:\tau\to\rho$, $\codim\rho=1$\\[2ex]
\centerline{\bf 5.3}\\
$(\check X,\check\shL)$
	\>polarized Calabi-Yau space mirror to $(X,\shL)$\\[2ex]
\centerline{\bf A.1}\\
$\Xi$
	\>$d$-dimensional polytope\\
$\P$
	\>face poset of $\Xi$\\
$\P^{[k]}$
	\>set of $k$-dimensional faces of $\Xi$\\
$(M_\tau)$
	\>system of abelian groups indexed by $\P$\\
$\varphi_{\sigma\tau}: M_\tau\to M_\sigma$
	\>homomorphism defined whenever $\tau\subset\sigma$\\
$(C^\bullet_\bct, d^\bullet_\bct)$
	\>barycentric cochain complex for $(M_\tau)$; $C^k_\bct=
	\bigoplus_{\sigma_0\sigma_1\ldots\sigma_k} M_{\sigma_k}$\\
$(f_{\sigma_0\sigma_1\ldots\sigma_k})_{\sigma_0\sigma_1
\ldots\sigma_k}$
	\>$k$-dimensional simplicial cochain, i.e.\ element of $C^k_\bct$\\[2ex]
\centerline{\bf A.2}\\
$\sgn(\tau,\omega)$
	\>for $\tau\subset\omega$ of codimension one sign
	comparing chosen orientations\\
$\sgn((\sigma_0,\ldots,\sigma_k),\sigma_k)$
	\>sign comparing orientations of
	$(\sigma_0,\ldots,\sigma_k)\in \Bar(\P)$ and of $\sigma_k$\\
$(C^\bullet_\phd, d^\bullet_\phd)$
	\>polyhedral cochain complex for $\{M_\tau\}$; $C^k_\bct=
	\bigoplus_{\tau\in\P^{[k]}} M_\tau$\\
$(f_\tau)_{\tau\in\P^{[k]}}$
	\>$k$-dimensional polyhedral cochain, i.e.\ element of $C^k_\phd$\\[2ex]
\centerline{\bf A.3}\\
$\P_\sigma$
	\>on an integral affine manifold with polyhedral decomposition
	$\P$,\\
	\>the face poset of $\tilde\sigma$\\
$\sgn(${\footnotesize $(\sigma_0 \to\cdots\to\sigma_k)$}$,
\sigma_k)$
	\>global analogue of $\sgn((\sigma_0,\ldots,\sigma_k),\sigma_k)$
	in A.2\\
$\shF_\sigma$
	\>abelian sheaf on $X_\sigma$\\
$\varphi_e$
	\>gluing maps for $(\shF_\sigma)$\\ 
$f_{\sigma_0\to \cdots\to\sigma_k}$
	\>global analogue of $f_{\sigma_0\sigma_1\ldots\sigma_k}$ in A.1\\

\end{tabbing}


\normalsize

\newpage

\begin{thebibliography}{ccc}

\bibitem{Altmann} K.~Altmann:
\emph{The versal deformation of an isolated toric  Gorenstein
singularity},
Invent.\ Math.\ {\bf 128} (1997), 443--479.

\bibitem{artin} M.\ Artin:
\emph{Algebraization of formal moduli II:
Existence of modifications},
Ann.\ Math.\ {\bf 91} (1970), 88--135.

\bibitem{SGA} M.~Artin, A.~Grothendieck, and J.-L.~Verdier: 
\emph{Th\'eorie des topos et cohomologie \'etale des sch\'emas},
(SGA~4), Lecture Notes in Mathematics, {\bf 269, 270, 305},
Springer-Verlag, Berlin (1971).

\bibitem{AM} P.\ Aspinwall and D.\ Morrison:
\emph{Chiral rings do not suffice: $N=(2,2)$ theories with nonzero
fundamental group},
Phys.\ Lett.\ B {\bf 334} (1994), 79--86.

\bibitem{Baues} O.~Baues: \emph{Gluing affine $2$-manifolds with polygons,}
Geom.\ Dedicata {\bf 75} (1999), 33--56.

\bibitem{Bat} V.~Batyrev: \emph{Dual polyhedra and mirror symmetry 
for Calabi-Yau hypersurfaces in toric varieties.}
J. Algebraic Geom.  {\bf 3}  (1994), 493--535.

\bibitem{BB} V.~Batyrev, and L.~Borisov: \emph{On Calabi-Yau complete
intersections in toric varieties,}
in {\it Higher-dimensional complex varieties (Trento, 1994)}, 39--65,
de Gruyter, Berlin, 1996.

\bibitem{Bredon} G.~Bredon:
\emph{Sheaf Theory},
2nd Edition, Springer-Verlag, 1997.

\bibitem{ChengYau} S.-Y.~Cheng and S.-T.~Yau, \emph{The real Monge-Amp\`ere 
equation and affine flat structures},
in \emph{Proceedings of the 1980 Beijing Symposium on Differential 
Geometry and Differential Equations}, Vol. 1, 2, 3 (Beijing, 1980),
339--370, Science Press, Beijing, 1982. 

\bibitem{chyn} S.~Chynoweth, and M.~Sewell:
\emph{Mesh duality and Legendre duality},
Proc.\ R.\ Soc.\ Lond.\ A, {\bf 428} (1990), 351--377.

\bibitem{Miranda} F.~DeMeyer, T.~Ford, and R.~Miranda:
\emph{The cohomological Brauer group of a toric variety},
J.\ Algebraic Geom. {\bf 2} (1993), 137--154.

\bibitem{Friedman} R.~Friedman:
\emph{Global smoothings of varieties with normal crossings}, 
Ann.\ Math.\ {\bf 118}, 75--114 (1983).

\bibitem{Fulton} W.~Fulton:
\emph{Introduction to toric varieties,}
Annals of Mathematics Studies 131, Princeton University Press,
Princeton, NJ, 1993.

\bibitem{GH} W.~Goldman, and M.~Hirsch:
\emph{The radiance obstruction  and parallel forms on affine
manifolds},
Trans.\ Amer.\ Math.\ Soc.\ {\bf 286} (1984), 629--649.

\bibitem{SlagI} M.~Gross:
\emph{Special Lagrangian Fibrations I: Topology,}
in: {\sl Integrable Systems and Algebraic Geometry},
(M.-H.\ Saito, Y.\ Shimizu and K.\ Ueno eds.),
World Scientific 1998, 156--193.

\bibitem{SlagII} M.~Gross:
\emph{Special Lagrangian Fibrations II: Geometry,}
in: {\sl Surveys in Differential Geometry}, Somerville: MA,
International Press 1999, 341--403.

\bibitem{TMS} M.~Gross:
\emph{Topological Mirror Symmetry},
Invent.\ Math.\ {\bf 144} (2001), 75--137.

\bibitem{GBB} M.~Gross:
\emph{Toric Degenerations and Batyrev-Borisov Duality},
preprint, (2004) math.AG/0406171.

\bibitem{Announce} M.~Gross, and B.~Siebert:
\emph{Affine manifolds,  log structures, and mirror symmetry},
Turkish J.\ Math.\ {\bf 27} (2003), 33-60.

\bibitem{tori} M.~Gross, and B.~Siebert:
\emph{Torus fibrations and toric degenerations,}
in preparation.

\bibitem{sequel} M.~Gross, and B.~Siebert:
\emph{Mirror symmetry via logarithmic degeneration data II},
in preparation.

\bibitem{EGA} A.~Grothendieck and J.~Dieudonn\'e:
\emph{El\'ements de G\'eom\'etrie Alg\'ebrique},
Publ.\ Math.\ IHES {\bf 8} (1961).

\bibitem{HZ} C.~Haase, and I.~Zharkov:
\emph{Integral affine structures on spheres and torus fibrations of
Calabi-Yau toric hypersurfaces I},
preprint 2002, {\tt math.AG/0205321}.

\bibitem{Hit} N.~Hitchin: 
\emph{The Moduli Space of Special Lagrangian Submanifolds},
Ann.\ Scuola Norm.\ Sup.\ Pisa Cl.\ Sci.\ (4) {\bf 25} (1997), 503--515.

\bibitem{Illu}
L.\ Illusie: \emph{Logarithmic spaces (according to K. Kato)},
in {\sl Barsotti Symposium in Algebraic Geometry} (Abano Terme 1991),
183--203, Perspect.\ Math.\ 15, Academic Press 1994.

\bibitem{Kajiwara} T.~Kajiwara: \emph{Logarithmic compactification
of the generalized Jacobian variety,} J. Fac. Sci. Univ. Tokyo, {\bf 40} (1993),
473--502.

\bibitem{F.Kato} F.~Kato:
\emph{Log smooth deformation theory},
Tohoku Math.\ J.\ {\bf 48} (1996), 317--354.

\bibitem{K.Kato} K.~Kato:
\emph{Logarithmic structures of Fontaine--Illusie},
in: {\sl Algebraic analysis, geometry, and number theory}
(J.-I.~Igusa et.~al.\ eds.), 191--224,
Johns Hopkins Univ.~Press, Baltimore, 1989.

\bibitem{Kawamata; Namikawa 1994} Y.~Kawamata, Y.~Namikawa:
\emph{Logarithmic deformations of normal crossing varieties
and smooothing of degenerate Calabi-Yau varieties},
Invent.\ Math.\ {\bf 118} (1994), 395--409.

\bibitem{Knutson} D.~Knutson:
\emph{Algebraic Spaces},
Lecture Notes in Math.\ {\bf 203}, Springer-Verlag, Heidelberg, 1971.

\bibitem{KS} M.~Kontsevich, and Y.~Soibelman: 
\emph{Homological mirror symmetry and torus fibrations},
in: {\sl Symplectic geometry and mirror symmetry} (Seoul, 2000), 203--263, 
World Sci.\ Publishing, River Edge, NJ, 2001. 

\bibitem{Leung} N.C.~Leung:
\emph{Mirror symmetry without corrections},
preprint 2000, math.DG/0009235,
to appear in Comm.\ Anal.\ Geom.

\bibitem{McLane} S.~Mac Lane:
\emph{Categories for the working mathematician},
Graduate Texts in Mathematics~{\bf 5},
Springer-Verlag 1971.

\bibitem{Mikhalkin} G.\ Mikhalkin:
\emph{Enumerative tropical  algebraic geometry in $\RR^2$},
preprint 2003, math.AG/0312530.

\bibitem{Milne} J.S.~Milne:
\emph{\'Etale Cohomology}, Princeton University Press 1980.

\bibitem{Nishinou; Siebert} T.~Nishinou, B.~Siebert:
\emph{Toric degenerations of toric varieties and tropical curves},
preprint 2004, math.AG/0409060.

\bibitem{Oda} T.~Oda:
\emph{Convex bodies and algebraic geometry. An introduction to the
theory of toric varieties},
Ergebnisse der Mathematik und ihrer Grenzgebiete (3),
{\bf 15}, Springer-Verlag, Berlin, 1988. 

\bibitem{Olsson} M.\ Olsson:
\emph{Semi-stable degenerations and period spaces for polarized
K3 surfaces}, Duke Math.\ J.\ {\bf 125} (2004), 121--203.

\bibitem{ogus} A.\ Ogus:
\emph{Logarithmic de Rham cohomology}, preprint.

\bibitem{Reid} M.\ Reid:
\emph{Nonnormal del Pezzo surfaces},
Publ.\ Res.\ Inst.\ Math.\ Sci.\ {\bf 30} (1994), 695--727.

\bibitem{Ruan} W.-D.~Ruan:
\emph{Lagrangian torus fibration and mirror symmetry of Calabi-Yau
hypersurface in toric variety},
preprint 2000, math.DG/0007028.

\bibitem{Kodaira} S.\ Schr\"oer, B.\ Siebert:
\emph{Irreducible degenerations of primary Kodaira surfaces},
in: {\sl Complex  geometry (G\"ottingen, 2000)}, 193--222,
Springer-Verlag 2002.

\bibitem{ss} S.\ Schr\"oer, B.\ Siebert:
\emph{Toroidal crossings and logarithmic structures},
preprint 2002, {\tt math.AG/0211088}, to appear in Adv.\ Math.

\bibitem{SiebertGW} B.~Siebert:
\emph{Log Gromov-Witten invariants},
unfinished manuscript. 

\bibitem{SY} H.~Shima, and K.~Yagi:
\emph{Geometry of Hessian manifolds},
Differential Geom.\ Appl.\ {\bf 7} (1997), 277--290.

\bibitem{steenbrink} J.\ Steenbrink:
\emph{Limits of Hodge structures},
Inv.\ Math.\ {\bf 31} (1976), 229--257.

\bibitem{SYZ} A.\ Strominger, S.-T.\ Yau, and E.~Zaslow, \emph{Mirror Symmetry
is $T$-duality,} Nucl.\ Phys.\ {\bf B479}, (1996) 243--259.
\end{thebibliography}
\end{document}